\documentclass[12pt, a4paper,
oneside, headinclude,footinclude]{scrartcl}
\usepackage[utf8]{inputenc}

\usepackage[
nochapters, 
pdfspacing, 
dottedtoc 
]{classicthesis} 
\usepackage{amsopn}
\usepackage[T1]{fontenc} 
\usepackage{multirow}
\usepackage[utf8]{inputenc} 
\usepackage{hhline}
\usepackage{graphicx} 
\graphicspath{{Figures/}} 
\usepackage{indentfirst}
\usepackage{enumitem} 
\usepackage{savesym}
\usepackage{mathrsfs}
\usepackage{geometry}
\usepackage{esint}
\usepackage{longtable}

\usepackage[
    backend=biber,
    style=numeric,
    natbib=false,
    url=false, 
    doi=true,
    eprint=true,
    maxnames=50
]{biblatex}

\usepackage{stmaryrd}

\usepackage{subfig} 

\usepackage{amsmath, amssymb,amsthm,amsfonts} 

\usepackage[english,capitalize]{cleveref}

\usepackage{thmtools}
\usepackage{thm-restate}

\theoremstyle{plain}
\newtheorem{teorema}{Theorem}[section]
\newtheorem{proposizione}[teorema]{Proposition}
\newtheorem{lemma}[teorema]{Lemma}
\newtheorem{corollario}[teorema]{Corollary}
\newtheorem*{theorem*}{Theorem}

\theoremstyle{definition}
\newtheorem{definizione}{Definition}[section]

\theoremstyle{remark}
\newtheorem{osservazione}{Remark}[section]

\makeatletter
\AtBeginDocument{%
  \def\theHtheorem2{\theHsubsection.\arabic{theorem2}}%
  \def\theHproposition2{\theHsubsection.\arabic{proposition2}}%
}
\makeatother

\newcommand{\Mass}{\mathbb{M}}
\newcommand{\R}{\mathbb{R}}
\newcommand{\N}{\mathbb{N}}

\newcommand{\Z}{\mathbb{Z}}
\newcommand{\im}{\mathrm{im}}
\newcommand{\Haus}{\mathscr{H}}
\newcommand{\Leb}{\mathscr{L}}

\newcommand{\bd}{\partial}

\newcommand{\Lip}{\mathrm{Lip}}

\newcommand{\Tan}{\mathrm{Tan}}

\newcommand{\dist}{\mathrm{dist}}
\newcommand{\supp}{\mathrm{supp}}
\newcommand{\Gr}{\mathrm{Gr}}

\newcommand{\dV}{d_V\kern-1pt}
\newcommand{\dW}{d_W\kern-1pt}

\def\dist{\mathop\mathrm{dist}} 
\def\supp{\mathop\mathrm{supp}} 
\def\diam{\mathop\mathrm{diam}}

\DeclareFontFamily{U}{rcjhbltx}{}
\DeclareFontShape{U}{rcjhbltx}{m}{n}{<->rcjhbltx}{}
\DeclareSymbolFont{hebrewletters}{U}{rcjhbltx}{m}{n}

\DeclareMathSymbol{\lamed}{\mathord}{hebrewletters}{108}
\DeclareMathSymbol{\mem}{\mathord}{hebrewletters}{109}
\DeclareMathSymbol{\ayin}{\mathord}{hebrewletters}{96}
\DeclareMathSymbol{\tsadi}{\mathord}{hebrewletters}{118}
\DeclareMathSymbol{\qof}{\mathord}{hebrewletters}{113}
\DeclareMathSymbol{\shin}{\mathord}{hebrewletters}{152}
\DeclareMathSymbol{\tav}{\mathord}{hebrewletters}{90}

\def\XXint#1#2#3{{\setbox0=\hbox{$#1{#2#3}{\int}$ }
		\vcenter{\hbox{$#2#3$ }}\kern-.6\wd0}}

\def\R{\mathbb{R}}

\def\Tan{\mathop\mathrm{Tan}} 					
\def\Lip{\mathop\mathrm{Lip}} 						
\def\dim{\mathop\mathrm{dim}} 					
\def\dist{\mathop\mathrm{dist}} 						
\def\supp{\mathop\mathrm{supp}}					

\newcommand{\scal}[2]{\langle #1\, ; \, #2\rangle}

\DeclareMathOperator{\largewedge}{\mbox{\large$\wedge$}}

\newcommand{\trait}[3]{\vrule width #1ex height #2ex depth #3ex}
\newcommand{\trace}{\mathchoice%
  {\mathbin{\trait{.12}{1.2}{.03}\trait{.8}{0.09}{0.03}}}
  {\mathbin{\trait{.12}{1.2}{.03}\trait{.8}{0.09}{0.03}}}
  {\mathbin{\hskip.15ex\trait{.09}{.84}{0.02}\trait{.56}{.07}{.02}}\hskip.15ex}
  {\mathbin{\trait{.07}{.6}{.01}\trait{.4}{.06}{.01}}}}

\newenvironment{itemizeb}
{\begin{itemize}\itemsep=2pt}
{\end{itemize}}

\newcounter{const}
\newcommand{\newC}{\refstepcounter{const}\ensuremath{C_{\theconst}}}
\newcommand{\oldC}[1]{\ensuremath{C_{\ref{#1}}}}

\newcounter{eps}

\newcommand{\vertiii}[1]{{\left\vert\kern-0.25ex\left\vert\kern-0.25ex\left\vert #1 
    \right\vert\kern-0.25ex\right\vert\kern-0.25ex\right\vert}}

\hypersetup{
colorlinks=true, breaklinks=true,bookmarksnumbered,
urlcolor=webbrown, linkcolor=RoyalBlue, citecolor=webgreen, 
pdftitle={}, 
pdfauthor={\textcopyright}, 
pdfsubject={}, 
pdfkeywords={}, 
pdfcreator={pdfLaTeX}, 
pdfproducer={LaTeX with hyperref and ClassicThesis} 
} 

\title{\normalfont\spacedallcaps{On the WALA conjecture, Alberti representations and applications}} 
\author{%
\spacedlowsmallcaps{Andrea Merlo}\textsuperscript{*}\\
}

\date{} 

\begin{document}

\renewcommand{\sectionmark}[1]{\markright{\spacedlowsmallcaps{#1}}} 
\lehead{\mbox{\mathfrak{L}ap{\small\thepage\kern1em\color{halfgray} \vline}\color{halfgray}\hspace{0.5em}\Bigmark\hfil}} 
\pagestyle{scrheadings}
\maketitle 
\setcounter{tocdepth}{2}

\paragraph*{Abstract}
We prove the WALA conjecture of G. David and S. Semmes: every AD-regular Radon measure for which all real-valued Lipschitz functions satisfy the weak approximation by affine functions condition is uniformly rectifiable. In order to prove the conjecture, we identify the correct quantitative analogue of the decomposability bundle and establish several structural results for general AD-regular measures.

{\let\thefootnote\relax\footnotetext{* \textit{Departamento de Matem\'aticas, Universidad del Pa\' is Vasco, Barrio Sarriena s/n 48940 Leioa, Spain, \href{andrea.merlo@ehu.eus}{andrea.merlo@ehu.eus}}}}



\tableofcontents

\section{Introduction}

Rademacher's theorem asserts that every Lipschitz function on $\R^n$ is differentiable $\Leb^n$-almost everywhere. It may therefore be read not only as a theorem about Lipschitz functions, but also as a rigidity property of Lebesgue measure: at almost every point, $\Leb^n$ carries enough independent directions to differentiate every Lipschitz function. Dorronsoro's theorem gives a quantitative counterpart to this statement, see \cite{zbMATH03923432}. Instead of merely asserting the existence of an infinitesimal affine approximation, it controls the error of the best affine approximation simultaneously over locations and scales; see also \cite{orponendorronsoro}. These two results lead naturally to ask if some form of converse theorem holds.

The qualitative problem asks what are those Radon measures for which Rademacher's theorem hold. The study of this type of questions goes back to Z. Zahorski who characterized the possible sets of non-differentiability of real-valued Lipschitz functions \cite{zbMATH03099952}.  The foundational, unpublished work of G. Alberti, M. Cs\"ornyei and D. Preiss identified the geometric mechanism underlying these phenomena: differentiability and non-differentiability are governed by the interaction of the set with families of Lipschitz curves with derivatives in prescribed cones \cite{acpnull,ACP10}. Their work introduced, among other tools, cone-null sets, width functions and weak tangent fields, and became a starting point for the geometric study of  Lipschitz differentiability spaces.

G. Alberti and A. Marchese transferred this point of view from sets to Radon measures by introducing the decomposability bundle $V(\mu,x)$, see \cite{AlbertiMarchese}. The bundle is generated by the tangent directions of all decompositions into one-dimensional rectifiable measures of $\mu$. Every Lipschitz function is differentiable at $\mu$-almost every $x$ along the subspace $V(\mu,x)$, and this statement is optimal: outside the bundle one can construct a Lipschitz function with genuinely non-differentiable behavior. Combined with the structure theorem for $\mathcal A$-free measures, see \cite{DPRAnnals}, this implies that if $V(\mu,x)=\R^n$ for $\mu$-almost every $x$, then $\mu\ll\Leb^n$. In parallel, the work of J. Cheeger and S. Keith initiated the study of differentiability of Lipschitz functions on metric measure spaces, and subsequent developments related such structures to Alberti representations, derivations, currents and tangent geometry see \cite{CheegerGAFA,Keith,CK9,Schioppa16,Schioppa16Bis,GuyCDavidLDS,CKS16,ErikssonBiqueGAFA,KleinerDavid,DePhilippisMarcheseRindler}. For similar results in Carnot groups and for Pansu differentiability spaces we refer to \cite{ReversePansu, zbMATH07977484}. More recently, G. C. David, S. Eriksson-Bique and R. Schul introduced pointwise and coarse tangent fields, organizing directions detected by curves into a multiscale structure \cite{DavidErikssonBiqueSchul2025}.

Not so much has been done so far for the quantitative analogue. First of all, the problem asks what are the Radon measures for which a Dorronsoro-type theorem holds. Notice that Dorronsoro's notion of \emph{well approximation} is given by an integral with respect to the Lebesgue measure, see the $\Omega$ numbers below, and hence we must expect a much more intrinsic answer to this problem in comparison to the qualitative counterpart. Before listing the progresses in the area, that can be found in the next subsection, we wish to recall how the theory of uniform rectifiability was introduced.

Classical results characterize rectifiable measures through densities, projections and approximate tangent planes \cite{Federer1996GeometricTheory,Mattila1995GeometrySpaces,Preiss1987GeometryDensities}. However, P. Jones's traveling salesman theorem introduced a different principle: the geometry of a rectifiable curve can be detected by measuring, and then summing, its deviations from lines over all locations and scales \cite{Jones1990TravelingSalesman}; see also \cite{Okikiolu1992TravelingSalesman}. These quantitative geometric ideas developed also by S. Semmes, G. David and D. Jerison led G. David and S. Semmes to develop the theory of uniform rectifiability for AD-regular measures, see \cite{Semmes1989,Semmes1990,DavidJerison1990,DavidSemmes1991,DavidLN,DavidSemmes}. Uniform rectifiability was developed as the geometric framework in which Calder\'on--Zygmund analysis retains the properties available on Lipschitz graphs. G. David and S. Semmes proved the boundedness of broad classes of singular integrals on uniformly rectifiable sets and formulated a converse problem, asking how much analytic information is sufficient to recover the geometry. The relationship between the Cauchy transform and uniform rectifiability was established by P. Mattila, M. Melnikov and J. Verdera \cite{MattilaMelnikovVerdera1996}; X. Tolsa developed quasiorthogonality and Carleson estimates for odd Calder\'on--Zygmund kernels \cite{tolsav2}; and F. Nazarov, X. Tolsa and A. Volberg proved that, in codimension one, the $L^2$-boundedness of the Riesz transform forces uniform rectifiability \cite{NTV}. At the same time, the use of multiscale flatness coefficients was extended beyond the original AD-regular setting. Jones-type square functions were used to obtain characterizations and sufficient criteria for rectifiable measures in \cite{tolsapart1,AzzamTolsa2015,BadgerSchul2016}; higher-dimensional traveling-salesman constructions were developed in \cite{AzzamSchul2018}; and further criteria for uniform rectifiability were obtained from lower-density estimates, Poincar\'e inequalities and projection properties in \cite{AzzamHyde2022,zbMATH07395053,Orponen2021BigProjections}.

\subsection*{The WALA conjecture and main results}

Let $\mathcal F\subseteq L^1_{\mathrm{loc}}(\R^n)$ be a finite-dimensional vector space containing the affine functions. For a Radon measure $\mu$, a Lipschitz function $f:\R^n\to\R$, $x\in\supp\mu$, $r>0$ and $q\in[1,\infty)$, set
$$
\Omega_{q,\mu,\mathcal F}(f;x,r)
:=
\inf_{A\in\mathcal F}
\left(
\fint_{B(x,r)}
\frac{|f-A|^q}{r^q}\,d\mu
\right)^{1/q},
$$
with the usual modification when $q=\infty$. The normalization by $r$ makes $\Omega_{q,\mu,\mathcal F}$ a first-order coefficient: it measures how well the restriction of $f$ to $\supp\mu$ can be approximated, to first order on $B(x,r)$, by an element of the fixed finite-dimensional family $\mathcal F$.

Let $\mu$ be AD-regular. We say that $\mu$ satisfies $\mathrm{GWALA}_q$ with respect to $\mathcal F$ if, for every $1$-Lipschitz function $f:\R^n\to\R$ and every $\varepsilon>0$, the set
$$
\mathscr B_{q,\mu,\mathcal F}(f,\varepsilon)
:=
\left\{
(x,r)\in\supp\mu\times(0,\infty):
\Omega_{q,\mu,\mathcal F}(f;x,r)>\varepsilon
\right\}
$$
is Carleson, with a Carleson constant which may depend on $\varepsilon$ and on $\mu$, but not on $f$. When $\mathcal F$ is the space of affine functions, this is the weak approximation of Lipschitz functions by affine functions, abbreviated WALA, introduced by G. David and S. Semmes in \cite[\S I.2, Definition 2.47]{DavidSemmes}. The condition is a quantitative form of Rademacher's theorem on an AD-regular measure: almost-everywhere infinitesimal approximation is replaced by a packing estimate for the locations and scales where first-order approximation fails.

G. David and S. Semmes proved that uniform rectifiability implies Dorronsoro-type theorem for Lipschitz function and in particular they satisfy the WALA condition. Further, they devoted sections 4.2 and 4.3 of Part III, Chapter 4 of \cite{DavidSemmes} to the study of the converse implication. They were able to prove that the GWALA implies uniform rectifiability in dimension 1. 
More specifically, they proved that the WALA condition implies WNB; see \cite[Part II, Chapter 3, Proposition 3.45]{DavidSemmes} and up to the present work the higher-dimensional case remained open. More recently, J. Azzam, M. Mourgoglou and M. Villa proved a Dorronsoro-type theorem on uniformly rectifiable sets for Sobolev functions; see \cite{azzam2023quantitativedifferentiabilityuniformlyrectifiable}. These results show that uniform rectifiability supports a robust affine approximation theory. The converse problem asks whether this approximation property, without any a priori information on flatness, projections or parametrizations of the support, is already sufficient to recover uniform rectifiability.

There is also a qualitative version of the problem: what can be said about a Radon measure for which $\lim_{r\to0}\Omega_{1,\mu,\mathcal F}(f;x,r)=0$ for $\mu$-almost every $x$ and every Lipschitz function $f$? This question is not covered by the existing literature. Results based on Alberti representations, derivations or the decomposability bundle concern the existence of derivatives along a measurable family of directions. In particular, differentiability along $V(\mu,x)$ gives no control, at a fixed scale, over the oscillation of the function on the part of the support transverse to $V(\mu,x)$. A closer result is the Lusin-type characterization of A. Marchese and A. Merlo \cite{MarcheseMerloLusin}: a Radon measure has the Lusin property for $C^1$ functions precisely when it is the sum of rectifiable measures of possibly different dimensions. Although this result also requires one to work on the support of the measure, its proof relies on the qualitative and infinitesimal machinery developed in \cite{AlbertiMarchese}, rather than on estimates which are uniform across scales.

We first state our qualitative result.

\begin{teorema}\label{teoremaqualitativointro}
Let $\mu$ be a Radon measure on $\R^n$, let $\mathcal F\subset L^1_{\mathrm{loc}}(\R^n)$ be a finite-dimensional vector space, and suppose that, for some $\alpha>0$, we have
$$
0<\Theta_*^\alpha(\mu,x)
\leq
\Theta^{\alpha,*}(\mu,x)
<\infty
$$
for $\mu$-almost every $x$. Assume that, for every $1$-Lipschitz function $f:\R^n\to\R$,
$$
\limsup_{r\to0}
\Omega_{1,\mu,\mathcal F}(f;x,r)
=
0
$$
for $\mu$-almost every $x$. Then $\alpha\in\N$ and $\mu$ is $\alpha$-rectifiable.
\end{teorema}

This is \cref{rectthm3}. The results in \cref{mainv2}, \cref{rectthm}, \cref{rectthm2} are formulated under the weaker assumptions of pointwise doubling and a suitable lower growth condition and therefore do not require the strong assumptions on lower and upper density. Indeed, the pointwise doubling assumption provides the stability of the coefficients under bounded changes of centre and scale, while the lower growth condition prevents the transverse part of the tangent measures from disappearing in the limit.

A second ingredient is the rectifiability criterion in \cref{criteriodirett}. If $d\in\{1,\ldots,n\}$ and
$$
\Theta^{d,*}(\mu,x)>0
\qquad\text{and}\qquad
\dim V(\mu,x)=d
$$
for $\mu$-almost every $x$, then $\mu$ is $d$-rectifiable. Thus the proof of \cref{teoremaqualitativointro} is reduced to excluding the possibility that the dimension of the decomposability bundle is strictly smaller than the dimension detected by the growth of the measure. In that case tangent measures split along the bundle, while the lower growth assumption forces a non-trivial transverse factor to remain. Finite configurations extracted from this factor are then used to perturb Lipschitz functions away from every model in $\mathcal F$.

We next state the quantitative result.

\begin{teorema}\label{teoremaintro}
Let $\mu$ be an $\alpha$-AD-regular measure on $\R^n$, let $\mathcal F\subset L^1_{\mathrm{loc}}(\R^n)$ be a finite-dimensional vector space, and let $q\in[1,\infty]$. If $\mu$ satisfies $\mathrm{GWALA}_q$ with respect to $\mathcal F$, then $\alpha\in\N$ and $\mu$ is uniformly rectifiable.
\end{teorema}

Taking $\mathcal F$ to be the space of affine functions proves the WALA conjecture of G. David and S. Semmes. Together with \cite[Part III, Chapter 4, Proposition 4.2]{DavidSemmes}, \cref{teoremaintro} characterizes uniformly rectifiable AD-regular measures through the validity of WALA for every Lipschitz function. Notice also that the exponent $\alpha$ is not assumed to be an integer: its integrality is part of the conclusion.

The proof of \cref{teoremaintro} must overcome two distinct obstructions. The constructions originating in the work of G. Alberti, M. Cs\"ornyei and D. Preiss detect non-differentiability asymptotically at individual points. The coefficient $\Omega_{q,\mu,\mathcal F}$ is instead an integral quantity: in order to force it to be large, the failure of approximation must occur on configurations carrying a quantitatively positive amount of $\mu$-mass. Moreover, WALA is a packing statement. It is not enough to construct a different bad function at each point or scale; one must construct a single Lipschitz function whose approximation coefficients remain uniformly large on a family of locations and scales with arbitrarily large packing. No previous mechanism converted directional information for Lipschitz functions into this type of multiscale geometric control. Constructing such a mechanism is the main difficulty of the paper.
\vspace{-0.2cm}
\subsection*{A quantitative decomposability bundle}

The main geometric construction of the paper is a quantitative counterpart of the decomposability bundle for general AD-regular measures. The classical bundle $V(\mu,x)$ is an infinitesimal object defined at almost every point. It encodes the directions generated by one-dimensional rectifiable decompositions of measures absolutely continuous with respect to $\mu$, and, by \cite{AlbertiMarchese}, the directions along which every Lipschitz function is differentiable. This information is sufficient for qualitative arguments, but it is not adapted to uniform rectifiability.

We replace the pointwise bundle by a system of directions attached to dyadic cubes. Roughly speaking, a cube $Q$ is invariant in a direction $e$ if, for every point in $\supp\mu \cap Q$, there exists a curve travelling inside a cone around $e$ and remaining close to $\supp\mu$ for most of its length. A family of directions is quantitatively independent when the corresponding cones remain uniformly separated, and the rank of $Q$ is the maximal number of independent invariant directions carried by the cube. Thus the rank is detected by actual one-dimensional geometry at the scale of $Q$, rather than by choosing an approximating plane independently at each location and scale.

The construction is designed to use both the invariant directions and the non-invariant directions. If $Q$ has rank $j$ and $\mathfrak e=(e_1,\ldots,e_j)$ is a family of independent invariant directions, then the support has a discrete product structure. More precisely, $\supp\mu$ is close to a set of the form
$$
\mathfrak c(Q)+\mathrm{span}(\mathfrak e)\times E_{Q},
$$
where $E_{Q}$ has the right $k-j$ coarse dimension, see \cref{propo:piano-da-invarianza}, \cref{prop:discrete-product-from-partial-invariance}. The AD-regularity of $\mu$ gives matching upper and lower counting estimates for $E_{Q,x}$ with exponent $k-j$. This is the finite-scale counterpart of the splitting of tangent measures along the decomposability bundle; see \cref{th:plit}.

Now, suppose that $Q$ has rank $j<k$ and that no additional invariant direction can be added transversely to $\mathrm{span}(\mathfrak e)$. Then for every admissible transverse direction of non-invariance, one has point for which every Lipschitz curve going in a cone around the fixed direction and passing close to said point, must spend a definite amount of length away from $\supp\mu$. If one has a long tower of cubes that are not invariant with respect to a fixed direction, it is possible to construct nested open sets which retain a definite portion of the measure but have very small width in the transverse directions. This is a quantitative version of the fact that a compact $C$-null set, using the language of G. Alberti and A. Marchese, is contained in an open set with small $C$-width.

\subsection*{Structure of the proof}

We first prove the qualitative result. Suppose that the dimension detected by the growth of $\mu$ is strictly larger than $\dim V(\mu,x)$ on a set of positive measure. After restricting to pieces on which the bundle varies only slightly, tangent measures split along $V(\mu,x)$ and retain a non-trivial transverse factor. The lower growth assumption prevents this transverse factor from disappearing in the blow-up, in the sense that the component along $V(\mu,x)^\perp$ is quite diffuse. At arbitrarily small scales one therefore finds finite transverse configurations of separated points with controlled mass.

Finite dimensionality of $\mathcal F$ gives an algebraic obstruction on such configurations. If $\dim\mathcal F=d$, one can prescribe values on a fixed finite number of points so that no element of $\mathcal F$ approximates all of them simultaneously. The difficulty is to realize these values by a Lipschitz perturbation without changing the derivatives already carried by the decomposability bundle.

We first mollify the given function and freeze the bundle on a suitable piece. Width functions are then used to cancel the derivatives orthogonal to the frozen bundle. After this cancellation, the function has very small oscillation in the transverse directions, while its Lipschitz constant increases only by a small controlled amount. We then pass to a smaller scale, where the transverse factor provides a finite configuration of separated points on which the corrected function is almost constant. 
The perturbation must change the values of the function by an amount comparable to this smaller scale, otherwise it cannot force the function away from the models in $\mathcal F$. At the same time, it must be defined on a much larger ball and have correspondingly small supremum norm relative to the radius of that ball. Indeed, the contribution of the perturbation to the Lipschitz constant is controlled by the ratio between its amplitude and the radius over which it is introduced. If the perturbation were inserted directly on the smaller ball, this ratio would be too large and the Lipschitz constant would become uncontrollable. 

In other words, in order to be able to perform the construction correctly, we need to perturb the function simultaneously in a lot of points. 

This is where invariance along $V(\mu,x)$ is essential. The perturbation is chosen to be constant in the directions of $V(\mu,x)$ and to vary only on the transverse configuration.
Invariance propagates the transverse profile along the directions of the bundle and therefore produces a set of quantitatively positive $\mu$-measure on which the perturbation is visible. 
Without this product structure, the transverse configuration would consist only of finitely many isolated points and would not carry enough mass to give a lower bound for the integral coefficient $\Omega_{1,\mu,\mathcal F}$. 
We then choose the values of the perturbation on the transverse configuration so that they cannot be simultaneously approximated by any element of $\mathcal F$. Since the original function is almost constant in the transverse directions and approximately affine along the bundle,
the perturbation forces a definite lower bound for $\Omega_{1,\mu,\mathcal F}$ on a set of positive measure.

A Baire category argument produces a Lipschitz function for which this lower bound occurs infinitesimally at almost every point of the bad piece. The functions constructed on different pieces are joined by extensions which agree on the relevant compact sets and are smooth away from them, following the approximation scheme in \cite{AlbertiMarchese}. This shows that when the decomposability bundle does not reach the dimension of the measure, these transverse points imply the existence of a lot of Lipschitz functions that are not well approximated by the functions in $\mathcal F$ and this in turn contradictions the fact that the $\Omega$s where going to $0$ as the scales go to $0$. It follows that the decomposability bundle has the dimension detected by the growth of the measure, and \cref{criteriodirett} gives rectifiability.

The quantitative theorem requires a different argument. The qualitative proof may pass to tangent measures, discard null sets and work separately on different pieces. None of these operations preserves the packing information required for uniform rectifiability. The quantitative problem is therefore to construct one function which detects the geometric defect on an arbitrarily large collection of locations and scales.

Assume that, for some $0\leq j\leq k-1$, the family $\mathscr I(j,2\vartheta,\sigma/2,\varepsilon,\delta,\zeta,\mathfrak m,A)$ is not Carleson. A non-Carleson selection argument produces a top cube $R$ and an arbitrarily deep tower of dyadic layers inside $R$. After a pigeonholing of the directions, all cubes in the construction carry the same $j$-tuple of invariant directions, up to the prescribed angular errors. Consecutive layers are separated by a very large factor in scale, and  all the layers of cubes cover almost all $R$.

For one fixed transverse direction, non-invariance gives an open set such that every curve directed inside a sufficiently narrow cone around that direction must spend a definite amount of length away from the support. An iterative thinning procedure converts these witnesses into open sets with large measure and small directional width, after discarding a controlled number of layers from the tower. At this stage, however, the width estimates are available only for a finite family of very narrow cones. This is not enough to cancel an arbitrary derivative orthogonal to $\mathrm{span}(\mathfrak e)$.

The joining-cones argument upgrades this information to the whole transverse cone. The narrow cones are given labels and these labels are organized in a finite tree. Each node of the tree corresponds to a cone $C$ and its descendants are the labels corresponding to cones $C^+$ and $C^-$ whose union is exaclty $C$. Via an iterative procedure these cones are joined in the sense that at each step, width loss for two cones with labels that are descendants of the same node, are propagated back with controlled and quantified errors, to the cone with the parent label. The conclusion, \cref{propfinalejoiningcones}, is a sequence
$$
\Omega_0\supseteq\Omega_1\supseteq\cdots\supseteq\Omega_{16\mathfrak N}
$$
of nested open sets which retain a definite proportion of $\mu(R)$ and have small width with respect to every direction in a wide cone transverse to $\mathrm{span}(\mathfrak e)$.

On this geometric structure we iterate the perturbation. Each block of generations consists of three operations. First, the current function is mollified inside a Vitali family of disjoint enlarged cube balls, with boundary values chosen so that the modification agrees with the previously constructed function outside. Second, width functions associated with the next open set are used to cancel the transverse derivatives of the mollified function. Third, the discrete product theorem provides a separated transverse configuration that we use to perturb the function away from $\mathcal F$. At this stage, the perturbation must be invariant along $\mathrm{span}(\mathfrak e)$ for the same reason as in the qualitative construction.

The very large gaps between the scales ensure at each perturbation step we do not destroy what we have already done.

After $\mathfrak N$ iterations, one obtains a single Lipschitz function $f$ and a dyadic family $\mathcal S\subseteq\Delta_\mu(R)$ such that the same fixed lower bound for $\Omega_{1,\mu,\mathcal F}(f;\cdot,\cdot)$ holds on every cube of $\mathcal S$, while
$$
\sum_{Q\in\mathcal S}\mu(Q)
\gtrsim
\mathfrak N\mu(R).
$$
The GWALA hypothesis gives the opposite estimate
$$
\sum_{Q\in\mathcal S}\mu(Q)
\lesssim
\mu(R),
$$
with a constant independent of $\mathfrak N$. Taking the tower sufficiently deep yields the contradiction and proves that every lower-rank stratum is Carleson.

The final step is geometric. A rank-$k$ cube carries $k$ independent invariant directions. The plane generated by these directions remains quantitatively close to the support on the relevant ball. Hence every cube which is bad for the other unilateral weak geometric lemma belongs to one of the lower-rank families. Since these families are Carleson, the fixed-threshold OUWGL criterion in \cref{propsmallOUWGL} yields uniform rectifiability and completes the proof of \cref{teoremaintro}.

\subsection*{Further developments}
A natural extension of the conjecture to metric spaces can be formulated as follows. Let $\mathsf X$ be an $\alpha$-dimensional AD-regular metric measure space and let $\mathbb G$ and $\mathbb H$ be fixed Carnot groups. For any $x\in\mathsf X$, $r>0$, $f\in\operatorname{Lip}_1(\mathsf X,\mathbb H)$ and Lipschitz map $\varphi:B(x,r)\to\mathbb G$, define
$$
\Omega_q(f;x,r,\varphi)^q
:=
\inf_{A\in\mathscr A(\mathbb G,\mathbb H)}
\fint_{B(x,r)}
\frac{d_{\mathbb H}\bigl(f(y),A(\varphi(y))\bigr)^q}{r^q}\,d\mu(y),
$$
where $\mathscr A(\mathbb G,\mathbb H)$ denotes the family of homogeneous homomorphisms from $\mathbb G$ to $\mathbb H$, with the usual modification when $q=\infty$. Denote by $\Delta_\mu$ the family of dyadic cubes for $(\mathsf X,\mu)$. By a coronization of $\Delta_\mu$ we mean a decomposition $\Delta_\mu=\mathscr B\cup\bigcup_{\mathcal T\in\mathfrak T}\mathcal T$ where $\mathscr B$ is a Carleson family of bad cubes, the cubes in each $\mathcal T$ form a coherent tree, and the tops of the trees satisfy a Carleson packing condition.  
We say that $(\mathsf X,\mu)$ satisfies the $\mathrm{WALA}_q$ condition if there exists such a coronization of $\Delta_\mu$ such that for every coherent tree of good cubes $\mathcal T$ there exists a Lipschitz map $\varphi_{\mathcal T}:\mathsf X\to\mathbb G$ such that, for every $f\in\operatorname{Lip}_1(\mathsf X,\mathbb H)$,
$$
\mathscr B_{\mathcal T,\varepsilon}(f)
:=
\left\{
Q\in\mathcal T:
\Omega_q\bigl(f;\mathfrak c(Q),\diam Q,\varphi_{\mathcal T}\bigr)
\geq\varepsilon
\right\}
$$
is a Carleson family with Carleson norm depending only on $\mu$ and $\varepsilon$, but not on $\mathcal T$ or $f$. One can weaken this definition by allowing the target of the chart $\varphi_{\mathcal T}$ to depend on $\mathcal T$, provided that the dimensions of the resulting Carnot groups are uniformly bounded. In an upcoming paper, we will investigate the structure of these metric spaces in a similar manner to what we did for Pansu differentiability spaces; see \cite{zbMATH07977484}. Notice that even when $\mathbb G$ is a Euclidean space, this is a completely new question. The constructions introduced in this paper are based on curves, quantitative invariance, width and coherent stopping-time families, and are therefore naturally adapted to this metric formulation.

A further project, in progress with G. Del Nin and J. Hoffman, concerns the David--Semmes conjecture for the Riesz transform. The aim is to obtain a structural description of a hypothetical counterexample of minimal possible dimension. Finally, with M. Villa we will provide a proof that the Many Segment Property, see \cite{zbMATH07395053}, implies uniform rectifiability by exploiting the techniques introduced here.

\subsection*{Upcoming versions}

The next version of this paper will contain several complementary results obtained jointly with M. Villa and M. Mourgoglou. Together with M. Villa, we prove a converse to \cref{teoremaqualitativointro}, obtaining a characterization of rectifiability in terms of infinitesimal affine approximation of Lipschitz functions.

Together with M. Mourgoglou, we apply the same circle of ideas to the Regularity problem for uniformly elliptic divergence-form operators with $L^\infty$ coefficients. We prove that, if $\Omega\subseteq\R^{n+1}$ has $n$-AD-regular boundary and the Regularity problem $(R^L_1)$ is solvable, then $\partial\Omega$ is uniformly rectifiable. We also obtain a qualitative version under substantially weaker geometric assumptions: if $0<\Haus^n(\partial\Omega)<\infty$, the domain is Wiener regular and satisfies a mild connectivity condition near the boundary, and the solution $u_f$ of the continuous Dirichlet problem satisfies
$$
N(\nabla u_f)(\xi)<\infty
$$
for $\Haus^n\llcorner\partial\Omega$-almost every $\xi$ and every compactly supported Lipschitz boundary datum $f$, then $\partial\Omega$ is $n$-rectifiable. These results will be incorporated into the next version of the paper.

\subsection*{Organization of the paper}
In \cref{preliminaries} we introduce the notation and collect the preliminary material on currents, the decomposability bundle, tangent measures, finite-dimensional approximation, width functions, the coefficients $\Omega$, Carleson families and the uniform rectifiability criteria used later.
In \cref{sezionequalitativa} we prove the rectifiability criterion \cref{criteriodirett}, establish the perturbation lemma \cref{perturbationlemma}, and deduce the qualitative results \cref{mainv2}, \cref{rectthm}, \cref{rectthm2}, \cref{rectthm3}. In \cref{sectionquantitativeresults} we introduce GWALA, construct non-Carleson towers and invariant cubes, prove the discrete product and joining-cones results, and establish the Carleson packing theorem \cref{teoremasemifinale}.
The same section concludes with \cref{teoremaGWALAimpliesUR}, where the lower-rank packing estimate is combined with OUWGL to prove uniform rectifiability.
Finally, \cref{appendix} contains the measurability results used throughout the paper.

\medskip

\textbf{Acknowledgments}.
 A.M. during the writing of this work was partially supported by Ikerbasque, by the European Union’s Horizon Europe research and innovation programme under the Marie Sk\l odowska-Curie grant agreement no 101065346, by the grant PID2024-157724NB-I00 of the Ministerio de Econom\'ia y Competitividad, and by the grant Europa Excelencia grant EUR2025-165042, funded by the Ministerio de Ciencia, Innovación y Universidades and Agencia Estatal de Investigacion.

\section{Preliminaries}\label{preliminaries}

In this preliminary section we introduce the notations and the elementary tools we will need throughout the paper.

\subsection{List of frequently used notations}

Here below the reader can find a list of frequently used notation.

\begin{longtable}{@{}c p{0.62\textwidth} p{0.16\textwidth}@{}}

\hline

\\

$\mathscr{M}(\R^n)$& family of positive Radon measures \label{Mspace}&\pageref{Mspace}\\

$B(x,r)$& closed ball of centre $x$ and radius $r$ \label{pallasubx}&\pageref{pallasubx}\\

$U(x,r)$ & open ball of centre $x$ and radius $r$ \label{pallaeux}&\pageref{pallaeux}\\

$\Leb^n$ &  $n$-dimensional Lebesgue measure \label{lebbilebi}&\pageref{lebbilebi} \\ 
$\Haus^k$ &   $k$-dimensional Hausdorff measure\label{Hausdihausdi}&\pageref{Hausdihausdi} \\
$\mu\trace E$ &   restriction of the measure $\mu$ to the set $E$\label{resres}&\pageref{resres} \\
$\supp\mu$ &   topological support of the measure $\mu$\label{supppsuppp}&\pageref{supppsuppp} \\
$\mathbb S^{n-1}$ &   unit sphere in $\mathbb R^n$. \label{sferisferi}&\pageref{sferisferi}\\ 
$\Gr(n,k)$ &   Grassmannian of $k$-dimensional linear subspaces of $\mathbb R^n$. &\pageref{defhausdorfix}\\
$\Gr(\mathbb R^n)$ &   Grassmanian of all planes in $\R^n$&\pageref{defhausdorfix} \\
$d(V,W)$ &   Hausdorff distance between planes in the unit ball&\pageref{defhausdorfix}\\
$d_{\Haus}(A,B)$ &   Hausdorff distance between two sets $A,B\subseteq\mathbb R^n$\label{defhausdorfix1}&\pageref{defhausdorfix1}\\ 
$d_{\Haus,B}(A,B)$ &   Hausdorff distance between $A$ and $B$ restricted to the ball $B$\label{defhausdorfix2}&\pageref{defhausdorfix2}\\ 
$\mathcal L(V,W)$ &   vector space of linear maps from $V$ to $W$ \label{defhausdorfix3}&\pageref{defhausdorfix3}\\
$\Delta_\mu$ &   dyadic lattice associated with the measure $\mu$&\pageref{evev}\\ 
$\mathfrak c(Q)$ &   centre of the dyadic cube $Q$&\pageref{evev}\\
\\
\hline
\end{longtable}

We add some further useful notations. First, we recall that in this paper natural numbers $\mathbb N$ contain $0$. Further, a Radon measure on $\R^n$ is said to be supported on some measurable set $E$ if $\mu(\R^n\setminus E)=0$. 

One of the protagonists of the paper are AD-regular measures. Recall that that a Radon measure $\mu$ on $\R^n$ is said to be $\alpha$-AD-regular with constant $D\geq 1$ if for every $x\in \supp\mu$ and every $0<r<\diam(\supp\mu)$, we have
$$D^{-1}r^\alpha\leq \mu(B(x,r))\leq Dr^\alpha.$$
We further recall the following standard notation. If $\mu$ is a Radon measure on $\R^n$, we let 
$$\Theta_*^\alpha(\mu,x):=
\liminf_{r\downarrow0}\frac{\mu(B(x,r))}{r^\alpha}\qquad\text{and}\qquad \Theta^{\alpha,*}(\mu,x):=
\limsup_{r\downarrow0}\frac{\mu(B(x,r))}{r^\alpha},$$
be the lower and upper $\alpha$-dimensional densities of a measure $\mu$. For future reference we recall that the distance between two sets is defined as follows
$$
\dist(A,B):=\inf\{|x-y|:x\in A,\ y\in B\}.
$$ 
Finally, we introduce cones

\begin{definizione}\label{gammaesigma}
Fix $\sigma\in(0,1)$ and $e\in\mathbb S^{n-1}$. We define the cone $C(e,\sigma)$ of axis $e$ and aperture $\sigma$ by
$$
C(e,\sigma):=\{v\in\R^n:\langle v,e\rangle\geq \sqrt{1-\sigma^2}|v|\},
$$
and we also let 
$\beta(\sigma):=\sqrt{1-\sigma^2}/\sigma$.
\end{definizione}

\begin{definizione}[Rectifiable measure]
Let $k\in\{1,\ldots,n\}$. A Radon measure $\mu$ on $\R^n$ is said to be $k$-rectifiable if $\mu\ll\Haus^k$ and there exist countably many Lipschitz maps $f_i:\R^k\to\R^n$ with $i\in \N$ such that
$$
\mu\Big(\R^n\setminus\bigcup_{i\in\N}f_i(\R^k)\Big)=0.
$$
\end{definizione}

\begin{definizione}[Uniformly rectifiable measure]
Let $k\in\{1,\ldots,n\}$. A $k$-AD-regular measure $\mu$ on $\R^n$ is said to be uniformly $k$-rectifiable if there exist constants $\vartheta>0$ and $M\geq1$ such that, for every $x\in\supp\mu$ and every $0<r<\diam(\supp\mu)$, there exists an $M$-Lipschitz map $f:B_{\R^k}(0,r)\to\R^n$ such that
$$\mu(B(x,r)\cap f(B_{\R^k}(0,r)))\geq\vartheta r^k.$$
\end{definizione}

\subsection{The beta numbers and the properties of the Grassmannian}

As usual, the symbol $\Gr(n,k)$ denotes the Grassmannian of $k$-planes in $\R^n$, and we define $\Gr(\R^n):=\bigcup_{0\leq k \leq n}\Gr(n,k)$. We endow $\Gr(\R^n)$ with the topology induced by the distance\label{defhausdorfix}
\begin{equation}
    d(V,W):=d_\Haus(V\cap B(0,1), W\cap B(0,1)),
    \label{distanceofplanes}
\end{equation}
where $d_\Haus$ is the Hausdorff distance. We recall the following definition, see \cite[\S 2.6, \S 6.1 and Theorem 6.4]{AlbertiMarchese}. 
We let $\gamma_{k,n}$ be the probability measure on $\Gr(n,k)$ that is invariant under the action of the orthogonal group, see \cite[\S 3.9]{Mattila1995GeometrySpaces}.

\begin{definizione}[$\beta$-numbers]
    Let $\mu$ be a Radon measure on $\R^n$. If $1\leq q<\infty$, define
$$\beta_{q,\mu}^k(x,r)^q=\beta_q^k(x,r)^q:=\inf_{V\in \Pi_k} \fint_{B(x,r)} \frac{\dist(z,V)^q}{r^q}d\mu(z),$$
where $\Pi_k$ is the set of affine $k$-dimensional planes in $\R^n$. Otherwise, if $q=\infty$, we let 
$$\beta_{\infty}^k(x,r)=\beta_{\infty,\mu}^k(x,r):=\inf_{V\in \Pi_k}\sup_{z\in B(x,r)\cap \supp(\mu)}\frac{\dist(z,V)}{r}.$$
\end{definizione}

\begin{definizione}
    Let $V\in \Gr(\R^n)$, and $\vartheta\in (0,1)$. Then 
    $$X_c(V,\vartheta):=\{y\in \R^n:\lvert \pi_{V^\perp}(y)\rvert\geq \vartheta \lvert y\rvert\},$$
    where as usual $\pi_{V^\perp}$ denotes the orthogonal projection on $V^\perp$.
\end{definizione}

\begin{definizione}[Bilateral $\beta$-numbers.]
  Suppose $\mu$ is a Radon measure on $\R^n$. For every $k$-dimensional affine plane $V$ we let 
   $$b\beta_\infty^k(V;x,r):=\sup_{z\in B(x,r)\cap \supp(\mu)}\frac{\dist(z,V)}{r}+
\sup_{y\in B(x,r)\cap V}\frac{\dist(y,\supp\mu)}{r},$$
and define
$$b\beta_{\infty}^k(x,r)=b\beta_{\infty,\mu}^k(x,r):=\inf_{V\in \Pi_k}b\beta_\infty^k(V;x,r).$$
\end{definizione}

If it is clear from the context what measure we are dealing with, the dependence of the coefficients on $\mu$, will be tacitly dropped.

\begin{proposizione}
    Let $1\leq p\leq q\leq \infty$ and $\mu$ be a Radon measure on $\R^n$. Then, for every $x\in \R^n$ and $r>0$ such that $\mu(B(x,r))>0$ we have $\beta_p^k(x,r)\leq \beta_q^k(x,r)$. 
\end{proposizione}

\begin{proof}
    It is immediate to note that for any such choice of $p$ and $q<\infty$, we have 
        $$\beta_p^k(x,r)^p=\inf_{V\in \Pi_k}\fint_{B(x,r)}\Big(\frac{\dist(z,V)}{r}\Big)^pd\mu(z)\leq \Big(\inf_{V\in \Pi_k}\fint_{B(x,r)}\frac{\dist(z,V)^q}{r^q} d\mu(z)\Big)^\frac{p}{q}=\beta_q^k(x,r)^p.$$
        The only remaining case is when $p<q=\infty$, and this follows by directly estimating with the sup-norm of the integrand.
\end{proof}

\begin{proposizione}\label{dimensioni}
    There exists a constant $\varepsilon_n$ depending only on $n$ such that if $V\in \Gr(\R^n)$ and $d(V,W)<\varepsilon_n$ then $\dim V=\dim W$.
\end{proposizione}

\begin{proof}
    Since $\Gr(\R^n)$ can be written as the disjoint union of the $\Gr(n,k)$s, and since each of the $\Gr(n,k)$ is compact we infer that for $k_1\neq k_2$, the quantity $\delta_{k_1,k_2}:=\dist_d(\Gr(n,k_1),\Gr(n,k_2))$ must be strictly positive. Since there are only finitely many pairs $k_1\neq k_2$, we may define $\varepsilon_n:=\frac12\min_{k_1\neq k_2}\delta_{k_1,k_2}>0$.
    Thanks to the choice of $\varepsilon_n$, we see that if $V,W\in\Gr(\mathbb R^n)$ and $\dim V\neq\dim W$, say $V\in\Gr(n,k_1)$ and $W\in\Gr(n,k_2)$, then $d(V,W)\geq\delta_{k_1,k_2}\geq 2\varepsilon_n$. This concludes the proof. 
\end{proof}

\begin{proposizione}\label{betastimedimensionalix}
Suppose that $\mu$ is an $\alpha$-dimensional AD-regular measure with constant $D$, and that $k\in\N$ is such that $k<\alpha$. Then for every $x\in \supp\mu$ and every $r>0$ we have
   $$\beta_\infty^k(x,r)^{\alpha-k}\geq 4^{-n}D^{-2}\frac{\omega_n}{\omega_k\omega_{n-k}},$$
   where $\omega_j$ is the volume of the unit ball in $\R^j$.
\end{proposizione}

\begin{proof}
Let $\eta>\beta_\infty^k(x,r)$ and let $V\in \Pi_k$ be such that
    $$B(x,r)\cap \supp(\mu)\subseteq B(V,\eta r).$$
    Let $\mathcal{N}$ be a maximal $\eta r$-separated set in $B(x,r)\cap \supp(\mu)$. First, note that 
    $$D^{-1}r^\alpha\leq \mathrm{Card}(\mathcal{N})D\eta^\alpha r^\alpha,$$
    and thus $D^{-2}\eta^{-\alpha}\leq \mathrm{Card}(\mathcal{N})$.
    Note that
    $$\bigcup_{z\in \mathcal{N}}B(z,\eta r/3)\subseteq B\big(V,4\eta r/3\big)\cap B(V^\perp,4r/3),$$
    where $V^\perp$ is the plane passing through $x$ and orthogonal to $V$.
    The cardinality of $\mathcal{N}$ can be thus bound by 
    $$D^{-2}\eta^{-\alpha}\omega_n(\eta r/3)^n\leq \mathrm{Card}(\mathcal{N})\omega_n(\eta r/3)^n\leq \omega_k(4r/3)^k\omega_{n-k}(4\eta r/3)^{n-k},$$
    where the last inequality come from the fact that the balls $B(z,\eta r/3)$ with $z\in \mathcal N$ are disjoint. This can be rearranged to infer, by taking the infimum over the $\eta$s, that
    $$\beta_\infty^k(x,r)^{\alpha-k}\geq 4^{-n}D^{-2}\frac{\omega_n}{\omega_k\omega_{n-k}}.$$
    This concludes the proof.
\end{proof}

The following proposition is very useful when dealing with Borel maps into the Grassmannian. 

\begin{proposizione}\label{prop:continuitax}
The map that associates to each $V\in \Gr(n,k)$ its Euclidean orthogonal projection $\pi_V$ is a bi-H\"older homeomorphism when  $\Gr(n,k)$ is endowed with the metric $d$ introduced above and when we endow linear maps with their natural operator norm.  
\end{proposizione}

\begin{proof}
    We refer to \cite[Proposition 2.9]{antonelli2020rectifiable2} and the proof of \cite[Proposition 6.1]{ReversePansu}.
\end{proof}

\begin{definizione}\label{exp:continuax}
    In the following we will let $\omega=\omega(n)$ be the smallest of the H\"older exponents appearing in \cref{prop:continuitax} in $\R^n$.
\end{definizione}

\subsection{Currents, decomposability bundle and Smirnov's theorem}

In this section we introduce many of the tools we will need to prove the main qualitative result, i.e., \cref{rectthm2}. First of all, we recall here the basic notions and terminology from
the theory of Euclidean currents.
A \emph{$k$-dimensional current} (or $k$-current) 
in $\R^n$ is a continuous linear functional on the space of smooth and 
compactly supported differential $k$-forms on $\R^n$, endowed with its usual test-function topology.

The boundary of a $k$-current $\mathbf{T}$ is the $(k-1)$-current 
$\bd \mathbf{T}$ defined by $\scal{\bd \mathbf{T}}{\omega} := \scal{\mathbf{T}}{d\omega}$\label{def:boundarycurrent}
for every smooth and compactly supported 
$(k-1)$-form $\omega$ on $\R^n$, and where $d\omega$ denotes the exterior derivative of $\omega$.
The \emph{mass} of $\mathbf{T}$, denoted by 
$\Mass(\mathbf{T})$, is the supremum of $\scal{\mathbf{T}}{\omega}$ over
all forms $\omega$ such that $|\omega|\le 1$
everywhere. 
A current $\mathbf{T}$ is called \emph{normal} if both $\mathbf{T}$ 
and $\bd \mathbf{T}$ have finite mass.

\medskip

By Riesz representation theorem a current $\mathbf{T}$ with finite mass
can be represented as a finite measure with 
values in the space $\largewedge_k(\R^n)$
of $k$-vectors in $\R^n$,
and therefore it can be written in the form
$\mathbf{T}=\tau\mu$ where $\mu$ is a finite positive measure
and $\tau$ is a $k$-vector field such that 
$\int |\tau|d\mu < +\infty$.
In particular the action of $\mathbf{T}$ on a form
$\omega$ is given by
\[
\scal{\mathbf{T}}{\omega} 
= \int_{\R^n} \scal{\tau(x)}{\omega(x)} \, d\mu(x)
\, ,
\]
and the mass $\Mass(\mathbf{T})$ is the total mass of $\mathbf{T}$ as a
measure, that is, $\Mass(\mathbf{T})=\int |\tau| d\mu$. Note that $0$-dimensional currents with locally finite mass are signed Radon measures and the mass coincides with the total variation.

In the following, whenever we write a current $\mathbf{T}$ as $\mathbf{T}=\tau\mu$ 
we tacitly assume that $\tau(x)\ne 0$ for $\mu$-almost every ~$x$;
in this case we say that $\mu$ is a measure
\emph{associated} to the current $\mathbf{T}$.%

\begin{definizione}[Integration of measures]
\label{s-measint}
Let $(I,dt)$ be a ($\sigma$-)finite measure space and for every $t\in I$ 
let $\mu_t$ be a real- or vector-valued measure on $\R^n$ such that:
\begin{itemizeb}
\item[(a)]
for every Borel set $E$ in $\mathbb{R}^n$ the function $t\mapsto \mu_t(E)$ 
is measurable; 
\item[(b)]
$\int_I \Mass(\mu_t) \, dt <+\infty$.
\end{itemizeb}
Then we denote by $\int_I \mu_t\, dt$ the measure on 
$\mathbb{R}^n$ defined by
$$
{\textstyle \big[ \int_I \mu_t\, dt \big]}(E)
:= \int_I \mu_t(E) \, dt
\quad\text{for every Borel set $E$ in $\mathbb{R}^n$.}
$$
Note that for every Borel set $E$ in $\mathbb{R}^n$ the function $t\mapsto \mu_t(E)$ 
is measurable (Borel) if and only if $t\mapsto \mu_t$ is a measurable (Borel) map from $I$ to the space of finite measures on $\R^n$ endowed with the weak* topology.
\end{definizione}

\begin{definizione}\label{correnticurve}
Let $B$ be a bounded Borel subset of $\R$ and $\gamma:B\to \R^n$ be a Lipschitz fragment. We denote by $\llbracket  \gamma\rrbracket$ the current of finite mass that acts on compactly supported smooth $1$-forms $\omega$ as:
$$\scal{\llbracket  \gamma\rrbracket}{\omega}:=\int_B \scal{\gamma'(t)}{\omega(\gamma(t))}dt.$$

In the following it will be also useful to write $\llbracket  \gamma\rrbracket=\tau_\gamma\rho\Haus^1\trace \im(\gamma)$, where $\rho$ is a suitable non-negative function in $L^1(\Haus^1\trace \im(\gamma))$ 
and $\tau_\gamma(x)$ is a unit Borel vector field that coincides with $\mathfrak{v}_\gamma(x)$, up to a real (non-zero) multiple, $\rho\Haus^1\trace \im(\gamma)$-almost everywhere, where the vector field $\mathfrak{v}_\gamma$ was introduced in \cite[Lemma 2.9]{ReversePansu}.
\end{definizione}

We now state a version of the celebrated Smirnov's theorem. The proof of this result can be found in \cite{ReversePansu} in much higher generality and it is based on the work \cite{PaoliniStepanov}.

\begin{teorema}[{\cite[Theorem 2.10]{ReversePansu}}]\label{smirnov}
Let $\mathbf{T}=\tau \mu$ be a $1$-dimensional normal current
with $\lvert\tau(x)\rvert=1$ for $\mu$-almost every $x\in \mathbb{R}^n$. Then, there exists a family of vector-valued measures $t\mapsto \boldsymbol\mu_t$ satisfying the hypothesis (a) and (b) of Definition \ref{s-measint} such that
\begin{itemize}
    \item[(i)] for almost every $t\in I$, where $I$ is the real line with the Lebesgue measure $\Leb^1$, there exists a Lipschitz curve $\gamma_t:[0,1]\to\mathbb{R}^n$ for which $\boldsymbol\mu_t=\llbracket  \gamma_t\rrbracket$ and 
    $$\scal{\mathbf{T}}{\omega}=\int_I \scal{\llbracket  \gamma_t\rrbracket}{\omega}\,dt=\int_I\int\rho_t\scal{\tau_{\gamma_t}}{\omega} d\Haus^1\trace \im(\gamma_t)\, dt,$$
    for every smooth and compactly supported $1$-form $\omega$;
    \item[(ii)] the following identity holds
    $$\Mass(\mathbf{T})=\int_I \Mass(\llbracket  \gamma_t\rrbracket)\,dt=\int_I\lVert \rho_t\rVert_{L^1(\Haus^1\trace\im(\gamma_t))} \,dt,$$ 
    and in particular  $\tau(x)=\tau_{\gamma_t}(x)$
   for $\Haus^1$-almost every $x\in\im(\gamma_t)$ and for almost every $t\in I$;
   \item[(iii)] the measure $\mu$ can be written as $\mu=\int_I \rho_t\Haus^1\trace \mathrm{im}(\gamma_t)dt$.
\end{itemize}
Further, one can also rewrite $\mathbf{T}$ as 
\begin{equation}\label{e:smirnov1}
\scal{\mathbf{T}}{\omega}=\int_\R\int_I\int\scal{\tau_{\gamma_t}}{\omega} d\Haus^1\trace (\im(\gamma_t)\cap E_{t,s})\, dtds,    
\end{equation}
where $E_{t,s}:=\{x\in\mathrm{im}(\gamma_t):\rho_t(x)\geq s\}$ and the map $(t,s)\mapsto\Haus^1\trace (\im(\gamma_t)\cap E_{t,s})$ satisfies the hypothesis (a) and (b) of Definition \ref{s-measint} relative to $I\times [0,\infty)$. In addition 
\begin{equation}\label{e:smirnov2}
\Mass(\mathbf{T})=\int_\R\int_I\Haus^1(\im(\gamma_t)\cap E_{t,s}) \,dtds\qquad \text{and}\qquad \mu=\int_\R\int_I \Haus^1\trace (\im(\gamma_t)\cap E_{t,s})\,dtds,
\end{equation}
with $\tau_{\gamma_t}(x)=\tau(x)$ for $\Haus^1$-almost every $x\in\im(\gamma_t)\cap E_{t,s}$ and for almost every $(s,t)\in \R\times I$.
\end{teorema}

Now that we have introduced currents, we can introduce a version of Alberti-Marchese decomposability bundle introduced in \cite{AlbertiMarchese}. In fact the definition G. Alberti and A. Marchese give of the decomposability bundle is based on the decomposition of the measure in curves, while we give here a more handy, equivalent definition. The equivalence of the two definitions can be found in \cite[Theorem 6.3]{AlbertiMarchese}.

\begin{definizione}[Decomposability bundle]
Given a positive Radon measure $\mu$ on $\R^n$ its \emph{decomposability bundle} is a map $V(\mu,\cdot)$ taking values in the set $\Gr(\R^n)$ defined as follows. A vector $v\in\R^n$ belongs to $V(\mu,x)$ if and only if there exists a vector-valued Radon measure $T$ with ${\mathrm {div}} T=0$ such that
$$\lim_{r\to 0}\frac{\Mass((T-v\mu)\trace B(x,r))}{\mu(B(x,r))}=0,$$
where $\Mass((T-v\mu)\trace B(x,r))$ denotes the total variation of the vector-valued measure $(T-v\mu)\trace B(x,r)$.
\end{definizione}

\begin{osservazione}
    For all the properties the Borel map $V(\mu,x)$ has, we refer to \cite{AlbertiMarchese}.
\end{osservazione}

\begin{lemma}\label{campi}
    Let $\mu$ be a Radon measure. Then, there are $n$ Borel vector fields $e_1,\ldots, e_n:\R^n\to \R^n$ such that $\{e_1(x),\ldots,e_n(x)\}$ is an orthonormal basis of $\R^n$ for every $x\in\R^n$, and for $\mu$-almost every $x\in \R^n$ we have that $\{e_1(x),\ldots, e_{d(x)}(x)\}$ spans $V(\mu,x)$, where $d(x):=\dim(V(\mu,x))$.
\end{lemma}

\begin{proof}
This can be achieved by observing that the maps $x\mapsto \Pi_{V(\mu,x)}$ are Borel, thanks to \cite[Proposition 6.1]{ReversePansu}. Projecting an orthonormal basis of $\R^n$ onto $V(\mu,x)$ one obtains Borel vector fields that span $V(\mu,x)$. From these vector fields, one first selects in a Borel way a linearly independent subfamily spanning $V(\mu,x)$, for instance by taking the first linearly independent subfamily in lexicographic order. This selection is Borel because linear independence is tested by the non-vanishing of minors. Then, applying the Gram--Schmidt orthonormalization procedure, one obtains Borel orthonormal vector fields spanning $V(\mu,x)$.

Since $x\mapsto \mathrm{id}-\Pi_{V(\mu,x)}=\Pi_{V(\mu,x)^\perp}$ is Borel, one can perform the same procedure for $V(\mu,x)^\perp$. Concatenating the two families gives $n$ Borel vector fields $e_1,\ldots,e_n$ with the desired properties. This proves the lemma.
\end{proof}

The following result relates the dimension of the decomposability bundle to the dimension of the measure. This is a continuous analogue of \cref{betastimedimensionalix}.

\begin{proposizione}\label{abscontinuous}
    Suppose $\mu$ is a Radon measure in $\R^n$ such that 
    $$\dim(V(\mu,x))=k\qquad\text{for $\mu$-almost every $x\in \R^n$}.$$
 Then $\mu\ll\Haus^k$.
\end{proposizione}

\begin{proof}
By a standard Lusin argument and by inner regularity, it is enough to prove the claim under the additional assumption that $\mu$ is finite, supported on a compact set, and that $V(\mu,\cdot)$ is continuous on $\supp(\mu)$.

Fix $\varepsilon>0$, to be chosen sufficiently small, and let $\{B_i\}_{i\in\N}$ be a covering of $\supp(\mu)$ such that for every $i\in\N$ there exists $V_i\in \Gr(n,k)$ for which
    $$d(V(\mu,y),V_i)\leq \varepsilon\qquad\text{for every }y\in B_i\cap\supp(\mu).$$
We claim that
$$(\pi_{V_i})_{\#}(\mu\trace B_i)\ll \Haus^k\trace V_i.$$
Let us first show that this claim implies $\mu\ll\Haus^k$. Let $K$ be a compact set such that $\Haus^k(K)=0$. Then for every $V\in \Gr(n,k)$ we have $\Haus^k(\pi_V(K))=0$. Hence, by the claim,
    $$(\pi_{V_i})_{\#}(\mu\trace B_i)(\pi_{V_i}(K))=0.$$
Since $K\cap B_i\subseteq B_i\cap \pi_{V_i}^{-1}(\pi_{V_i}(K))$, we get
    $$\mu(K\cap B_i)\leq \mu(B_i\cap \pi_{V_i}^{-1}(\pi_{V_i}(K)))=0.$$
Since the family $\{B_i\}_{i\in\N}$ covers $\supp(\mu)$, this gives $\mu(K)=0$, and therefore $\mu\ll\Haus^k$.

It remains to prove the claim. Fix $i\in\N$ and let $e\in \mathbb S^{n-1}\cap V_i$. Since $V(\mu,\cdot)$ is $\varepsilon$-close to $V_i$ on $B_i\cap\supp(\mu)$, if $\varepsilon$ is small enough we can find a Borel vector field $\tau:\R^n\to\R^n$ such that $|\tau(x)|=1$, $\tau(x)\in V(\mu,x)$, and
    $$\tau(x)\in C(e,2\varepsilon)\qquad\text{for $\mu$-almost every }x\in B_i.$$
We also take $\tau=0$ outside $B_i$. By \cite[Theorem 6.2]{AlbertiMarchese}, there exists a normal current $\mathbf T$ such that
    $$\mathbf T=\tau\mu+\boldsymbol{\sigma},\qquad \partial\mathbf T=0,$$
where $\boldsymbol{\sigma}$ is singular with respect to $\mu$. Let $E$ be a Borel set such that $\mu(E^c)=0$ and $\|\boldsymbol{\sigma}\|(E)=0$, which exists since $\boldsymbol{\sigma}$ is singular with respect to $\mu$. By Smirnov's theorem, see \cref{smirnov}, we can decompose $\mathbf T$ as
    $$\mathbf{T}=\int_I\llbracket  \gamma_t\rrbracket dt,\qquad
    \Mass(\mathbf T)=\int_I \Mass(\llbracket  \gamma_t\rrbracket) dt,\qquad
    \|\mathbf T\|=\int_I \|\llbracket  \gamma_t\rrbracket\| dt.$$
On $E\cap B_i$ we have $\|\mathbf T\|=\mu$, and the orientation of $\mathbf T$ coincides with $\tau$. Therefore
    $$\mu\trace B_i=\int_I\|\llbracket  \gamma_t\rrbracket\|\trace(E\cap B_i)\,dt,$$
and for almost every $t\in I$ the curves appearing in this decomposition have tangent direction in $C(e,2\varepsilon)$ on $E\cap B_i$.
Pushing forward by $\pi_{V_i}$ gives
    $$(\pi_{V_i})_\#(\mu\trace B_i)=\int_I(\pi_{V_i})_\#\big(\|\llbracket  \gamma_t\rrbracket\|\trace(E\cap B_i)\big)\,dt.$$
For each $\ell=1,\ldots,k$, let $e_\ell$ be an orthonormal basis of $V_i$.
Repeating the previous construction with $e=e_\ell$, we obtain a normal
one-dimensional current $\mathbf T_\ell$ in $\mathbb R^n$. Set
$$
\mathbf S_\ell:=(\pi_{V_i})_\#\mathbf T_\ell .
$$
Then $\mathbf S_\ell$ is a normal one-dimensional current in $V_i$. Moreover,
by the cone condition and by choosing $\varepsilon>0$ sufficiently small, the
polar vectors $\frac{d\mathbf{S_1}}{d\lVert \mathbf S_1\rVert},\ldots,\frac{d\mathbf{S_k}}{d\lVert \mathbf S_k\rVert}$ span
$V_i$ for $(\pi_{V_i})_\#(\mu\trace B_i)$-almost every point, and
$$
(\pi_{V_i})_\#(\mu\trace B_i)\ll \|\mathbf S_\ell\|
\qquad\text{for every }\ell=1,\ldots,k.
$$
Therefore, by \cite[Corollary 1.12]{DPRAnnals}, applied in the Euclidean
space $V_i\simeq\mathbb R^k$, we get
$$
(\pi_{V_i})_\#(\mu\trace B_i)\ll \Haus^k\trace V_i .
$$
This proves the claim, and hence the proposition.
\end{proof}

\subsection{Tangent measures, the geometry of blowups and rectifiability}\label{s:tan}

In order to prove the qualitative results we need to deal with tangents, and in this section we recall and restate many known results that will be helpful in the following.

\begin{definizione}\label{defF_r}
Given two Radon measures $\mu$ and $\sigma$, we set
$$F_{B}(\mu,\sigma) = \sup\Big\{ \int f\,d(\mu-\sigma):f \text{ is $1$-Lipschitz and supported on $B$}\Big\}. $$
For $r>0$, we write
$$
F_{r}(\mu,\nu):= F_{B(0,r)}(\mu,\nu),\qquad
F_{r}(\mu)=F_{r}(\mu,0):=\int (r-|z|)_{+} d\mu.$$
\end{definizione}

\begin{definizione}
    Let $\mu$ be a Radon measure on $\R^n$. For every $x\in \R^n$ and $r>0$ we let 
    $T_{x,r}\mu$ be the measure that acts as 
    $$T_{x,r}\mu(A):=\mu(x+rA),\qquad\text{for every Borel set }A.$$
\end{definizione}

\begin{definizione}[{\cite[Section 2]{Preiss1987GeometryDensities}}]  We introduce the following notations
\begin{enumerate}
\item A set $\mathfrak{M}$ of non-zero Radon measures in $\R^n$ is a \emph{cone} if $c\mu\in \mathfrak{M}$ whenever $\mu\in \mathfrak{M}$ and $c>0$.
\item A cone $\mathfrak{M}$ is a \emph{$d$-cone} if $T_{0,r}\mu\in \mathfrak{M}$ for all $\mu\in \mathfrak{M}$ and $r>0$.
\item The \emph{basis} of a $d$-cone $\mathfrak{M}$ is the set $\{\Psi\in \mathfrak{M}: F_{1}(\Psi)=1\}$;
\end{enumerate}
\end{definizione}

The following is Preiss's definition of distance from a cone of measures. We refer to \cite{Preiss1987GeometryDensities}.

\begin{definizione}\label{d_r}
We say a $d$-cone has \emph{closed (resp. compact) basis} if its basis $\{\mu\in \mathfrak{M}: F_{1}(\mu)=1\}$ is closed (resp. compact) with respect to the weak topology.
 For a $d$-cone $\mathfrak{M}$, $r>0$, and $\mu$ a Radon measure with $0<F_{r}(\mu)<\infty$, we define the \emph{distance} between $\mu$ and $\mathfrak{M}$ as 
$$
d_{r}(\mu,\mathfrak{M}):=\inf\Big\{F_{r}\Big(\frac{\mu}{F_{r}(\mu)},\nu\Big): \nu\in \mathfrak{M}, F_{r}(\nu)=1 \Big\}.
$$
\end{definizione}

The following lemma collects some of the properties of $F$, $d_r$ and their relationship with weak* convergence of measures.

\begin{lemma}[{\cite[Section 2]{KPT09}}]\label{preiss}
Let $\mu$ and $\nu$ be Radon measures in $\mathbb{R}^{n}$ and let $\mathfrak{M}$ be a $d$-cone. For $x\in \R^n$ and $r>0$, the following properties hold.
\begin{enumerate}
\item $ \displaystyle F_{B(x,r)}(\mu)=r F_{1}(T_{x,r}\mu)$,
\item $\displaystyle F_{B(x,r)}(\mu,\nu)=rF_{1}(T_{x,r}\mu,T_{x,r}\nu)$,
\item $\displaystyle \mu_{i}\rightharpoonup \mu$ if and only if $F_{r}(\mu_{i},\mu)\rightarrow 0$ for all $r>0$,
\item $\displaystyle d_{r}(\mu,\mathfrak{M})\leq 1$,
\item $\displaystyle d_{r}(\mu,\mathfrak{M})=d_{1}(T_{0,r}\mu,\mathfrak{M})$,
\item if $\displaystyle \mu_{i}\rightharpoonup \mu$ and $\displaystyle F_{r}(\mu)>0$, then $\displaystyle d_{r}(\mu_{i},\mathfrak{M})\rightarrow d_{r}(\mu,\mathfrak{M})$. 
\end{enumerate}
\end{lemma}

Here below we collect useful properties of cones of measures coming from various papers.

\begin{lemma}[{\cite[Remark 2.13]{KPT09}, \cite[Lemma 2.2]{reftgmeasures}}] \label{l:closed} A $d$-cone $\mathfrak{M}$ of Radon measures in $\mathbb{R}^{n}$ has a closed basis if and only if it is a relatively closed subset of the non-zero Radon measures in $\mathbb{R}^{n}$. 
\end{lemma}

\begin{lemma}[{\cite[Lemma 2.3]{reftgmeasures}}]\label{chartan}
If $\mu$ is a nonzero Radon measure and $\mathfrak{M}$ is a $d$-cone with closed basis, then $\mu\in \mathfrak{M}$ if and only if $d_{r}(\mu,\mathfrak{M})=0$ for all $r>0$ for which $F_{r}(\mu)>0$. 
\end{lemma}

\begin{proposizione}[{\cite[Proposition 2.2]{Preiss1987GeometryDensities}}]\label{p:compact} Let $\mathfrak{M}$ be a $d$-cone. Then $\mathfrak{M}$ has compact basis if and only if for every $\lambda>1$ there is $\tau>1$ such that 
\begin{equation}
F_{\tau r}(\Psi)\leq \lambda F_{r}(\Psi) \mbox{ for every }\Psi\in \mathfrak{M} \mbox{ and }r>0.
\label{e:compact}
\end{equation}
 In this case, $0\in \supp \Psi$ for all $\Psi\in \mathfrak{M}$.
\end{proposizione}

We introduce now the concept of tangent measures, which are the main objects of study of this section, that were introduced by D. Preiss in \cite{Preiss1987GeometryDensities}.

\begin{definizione}
    Suppose $\mu$ is a Radon measure on $\R^n$ and let $x\in \R^n$. We denote by $\Tan(\mu,x)$ the set of those non-zero Radon measures $\nu$ for which there exists an infinitesimal sequence $r_i$ and a sequence of positive real numbers $c_i>0$ such that 
    $$c_iT_{x,r_i}\mu\rightharpoonup \nu. $$
\end{definizione}

\begin{osservazione}[{\cite[\S 2.3(3)]{Preiss1987GeometryDensities}}]\label{tanisdcone}
    The set $\Tan(\mu,x)$ is a $d$-cone with a closed basis. 
\end{osservazione}

\begin{osservazione}
    When $\mu$ satisfies 
    $$0<\Theta^\alpha_*(\mu,x)\leq \Theta^{\alpha,*}(\mu,x)<\infty\qquad\mu\text{-almost everywhere,}$$
    another object is typically considered. 
   One can also define $\Tan_\alpha(\mu,x)$, the set of those Radon measures $\nu$ for which there exists an infinitesimal sequence $r_i$ such that 
    $$r_i^{-\alpha}T_{x,r_i}\mu\rightharpoonup\nu.$$
Under the two-sided density assumption, every element of $\Tan(\mu,x)$ is a positive multiple of an element obtained as a limit of $r_i^{-\alpha}T_{x,r_i}\mu$, after passing to a subsequence. 
\end{osservazione}

\begin{proposizione}[{\cite[Theorem 2.5]{Preiss1987GeometryDensities}}]\label{nonemptytangents}
    If $\mu$ is a Radon measure on $\R^n$, then $\Tan(\mu,x)\neq \emptyset$ for $\mu$-almost every $x\in \R^n$. 
\end{proposizione}






\begin{lemma}\label{lemmadimensionbup}
      Suppose $\mu$ is a Radon measure on $\R^n$ such that 
     $$0< \Theta^{\alpha}_*(\mu,x)\leq \Theta^{\alpha,*}(\mu,x)<\infty\qquad\mu\text{-almost everywhere.}$$
     Then, for $\mu$-almost every $x\in \R^n$ and every $\nu\in \Tan_\alpha(\mu,x)$ we have 
     $$\Theta_*^\alpha(\mu,x)s^\alpha\leq \nu(B(y,s))\leq  \Theta^{\alpha,*}(\mu,x) s^\alpha\qquad\text{for every $y\in \supp(\nu)$ and $s>0$}.$$
\end{lemma}

\begin{proof}
This is a standard fact; see the proof of
\cite[Proposition 3.4]{DeLellis2008RectifiableMeasures}. The only difference is that, in that case $\Theta_*^\alpha(\mu,x)=\Theta^{\alpha,*}(\mu,x)$ for $\mu$-almost every $x\in \R^n$, while here one uses the lower and upper $\alpha$-density
bounds separately, which gives the two inequalities rather than equality.
\end{proof}

\begin{teorema}[{\cite[Corollary 2.7]{Preiss1987GeometryDensities}, \cite[Theorem 2.4]{reftgmeasures}}]\label{t:compactdouble} Let $\mu$ be a Radon measure on $\mathbb{R}^{n}$, and $\xi\in \supp \mu$. Then $\Tan(\mu,\xi)$ has compact basis if and only if 
\begin{equation}
\label{e:compactdouble}
\limsup_{r\rightarrow 0} \frac{\mu(B(\xi,2r))}{\mu(B(\xi,r))}<\infty.
\end{equation}
In this case, for any $\nu\in \Tan(\mu,\xi)$, it holds that $0\in \supp \nu$ and 
$$
\frac{\nu(B(0,2r))}{\nu(B(0,r))}\leq \limsup_{\rho \rightarrow 0} \frac{\mu(B(\xi,2\rho))}{\mu(B(\xi,\rho))},\mbox{ for all }r>0.
$$ 
\end{teorema}

\begin{lemma}[{\cite[Theorem 14.3]{Mattila1995GeometrySpaces}}]\label{precompactnessmeasures} Let $\mu$ be a Radon measure on $\mathbb{R}^{n}$. If $\xi\in \supp\mu$ and \eqref{e:compactdouble} holds, then every sequence $r_{i}\downarrow 0$ contains a subsequence such that 
\begin{equation}\label{eq:tang-doubl}
\frac{T_{\xi,r_{j}}\mu}{\mu(B(\xi,r_{j}))}  \rightharpoonup \nu,
\end{equation}
for some measure $\nu\in \Tan(\mu,\xi)$.
\end{lemma}

 Having tangent measures that arise as limits of the form \eqref{eq:tang-doubl} is very convenient, but this limit does not always converge weakly to something. This may happen if $\mu$ is not pointwise doubling at the point $\xi$. However, all tangent measures are at least dilations of tangent measures arising in this way.

\begin{lemma}[{\cite[Lemma 14.4(3)]{Mattila1995GeometrySpaces}}]\label{l:nicetan}
Let $\mu$ be a nonzero Radon measure, let $\xi\in \supp\mu$, and assume that \eqref{e:compactdouble} holds at $\xi$, and let $\nu\in \Tan(\mu,\xi)$. Then there are $\rho_{j}\downarrow 0$ and $c>0$ so that 
$$
\frac{T_{\xi,\rho_{j}}\mu}{\mu(B(\xi,\rho_{j}))}\rightharpoonup c \nu.
$$ 
\end{lemma}


\begin{teorema}[{\cite[Theorem 14.16]{Mattila1995GeometrySpaces}}] \label{t:ttt}
Let $\mu$ be a Radon measure on $\R^n$. For $\mu$-almost every $x\in \R^n$, if $\nu\in \Tan(\mu,x)$, the following hold:
\begin{enumerate}
\item $T_{y,r}\nu\in \Tan(\mu,x)$ for all $y\in \supp \nu$ and $r>0$.
\item $\Tan(\nu,y)\subset \Tan(\mu,x)$ for all $y\in \supp \nu$.
\end{enumerate}
\end{teorema}

We now bring the decomposability bundle into the discussion. The point is that the directions in $V(\mu,x)$ are visible at the level of every blow-up: for $\mu$-almost every $x$, tangent measures are invariant along $V(\mu,x)$ and therefore split as a product of $\Haus^k\trace V(\mu,x)$ with a measure on the orthogonal complement. Thus, whenever the measure has dimension larger than $k=\dim V(\mu,x)$, the dimension in the transverse direction must be strictly bigger than 0. This is the mechanism that will later provide support points aligned in directions along which the measure is thin, and hence allow us to cancel the corresponding derivatives of Lipschitz functions.

\begin{teorema}[{\cite[Lemma 2.9]{zbMATH07487904}}]\label{th:plit}
    Let $\mu$ be a Radon measure on $\R^n$. Then, for $\mu$-almost every $x\in \R^n$ any tangent measure $\nu\in\Tan(\mu,x)$ is $V(\mu,x)$-invariant, i.e., there exists a Radon measure $\eta\in \mathscr{M}(V(\mu,x)^\perp)$ such that $\nu=\Haus^k\trace V(\mu,x)\otimes \eta$, where $k=\dim(V(\mu,x))$.

    In addition, if there exists $A>0$ such that for $\mu$-almost every $x\in \R^n$ we have 
    \begin{equation}
        \limsup_{r\to 0}\frac{\mu(B(x,2r))}{\mu(B(x,r))}\leq A,
        \label{eq:doublingbups}
    \end{equation}
    then for $\mu$-almost every $x\in \R^n$ the measure $\eta$ above satisfies
$$\eta(B(z,2r))\leq 2^{-k}A^{\log_2(2n)+1}\eta(B(z,r))\qquad \text{ for every $z\in \supp\eta$ and every $r>0$}.$$
\end{teorema}

\begin{proof}
    We only need to prove the second part of the theorem. Fix a point $x\in \R^n$ where \eqref{eq:doublingbups} holds, $V(\mu,x)$ is well defined and $\Tan(\mu,x)\neq \emptyset$. We let
$$ Q_x(z,r):=z+\{u+w:u\in V(\mu,x),\text{with } |u|<r,\ w\in V(\mu,x)^\perp,\text{ and with } |w|<r\}.$$
It is immediate to see that for every $z\in\R^n$ and $r>0$ we have
$$B(z,r/\sqrt{n})\subseteq Q_x(z,r)\subseteq B(z,\sqrt{n}r).$$
This implies that for every $\nu\in \Tan(\mu,x)$ we have 
\begin{equation}
\begin{split}
    \frac{\nu(Q_x(0,2r))}{\nu(Q_x(0,r))}\leq \frac{\nu(B(0,2\sqrt{n}r))}{\nu(B(0,r/\sqrt{n}))}\leq \prod_{j=0}^{\lfloor \log_2(2n)\rfloor}\frac{\nu(B(0, 2^{j+1} r/\sqrt{n}))}{\nu(B(0,2^{j} r/\sqrt{n}))}.
\end{split}
\end{equation}
for all $r>0$.
By \cref{t:ttt}, for every $z\in\supp\nu$ and every $s>0$ the measure $T_{z,s}\nu$ belongs to $\Tan(\mu,x)$. Applying \cref{t:compactdouble} to these tangents and iterating the doubling estimate gives
$$\frac{\nu(Q_x(z,2r))}{\nu(Q_x(z,r))}\leq \prod_{j=0}^{\lfloor \log_2(2n)\rfloor}\frac{\nu(B(z, 2^{j+1} r/\sqrt{n}))}{\nu(B(z,2^{j} r/\sqrt{n}))}\leq A^{\log_2(2n)+1}.$$
By the product decomposition we infer 
$$2^k\,\frac{\eta(B(z,2r))}{\eta(B(z,r))}=\frac{\nu(Q_x(z,2r))}{\nu(Q_x(z,r))}\leq A^{\log_2(2n)+1},$$
which concludes the proof.
\end{proof}

In the proof of the qualitative part, we will need to infer from the structure of the blowups, structure of the original measure. The following propositions recover points in the support of $\mu$ as (roughly) rescaled copies of the points of the blowups. The statements are given for a sequence of measures $\mu_j$ converging to some $\nu$, however the reader should think about these as blowup  procedures.

\begin{lemma}\label{puntivicini}
    Let $\mu_j$ be a sequence of Radon measures weakly* converging to some other Radon measure $\nu$. Then, for every $z\in \supp\nu$ there exists a sequence of points $y_j\in \supp\mu_j$ such that $y_j\to z$. 
    In addition, if $\lambda_1,\lambda_2,D,r>0$ and $y\in \supp \nu$ are such that 
    $$\nu(B(y,\lambda_1 r))\leq D\nu(B(y,\lambda_2 r)),$$
then, for every $\varepsilon>0$ there exists $j_0\in\N$ such that for every $j\geq j_0$ we have
$$\mu_j(B(y,(\lambda_1-\varepsilon) r))\leq (D+\varepsilon)\mu_j(B(y,(\lambda_2+\varepsilon) r)).$$
\end{lemma}

\begin{proof}
Let us begin by showing the existence of the sequence $y_j$. By contradiction, suppose there exists $z_0\in \supp\nu$ for which there does not exist a sequence of points $y_j\in \supp\mu_j$ such that $y_j\to z$.    
    If this was the case there would exists $\rho>0$ such that $\mu_j(B(z_0,\rho))=0$ for every $j$ sufficiently big. However, we see that
    \begin{equation}
        \begin{split}
            F_{B(z_0,\rho)}(\mu_j,\nu)\geq&\int \dist(y,B(z_0,\rho)^c) d\nu(y)-\int \dist(y,B(z_0,\rho)^c) d\mu_j(y)\\
            \geq& \int \dist(y,B(z_0,\rho)^c) d\nu(y)>0.
        \end{split}
    \end{equation}
    This however is in contradiction with the weak convergence to $\nu$ and item 3. of \cref{preiss}.

    In order to prove the second part of the proposition, we notice that there exist $\lambda_1'\in (\lambda_1-\varepsilon,\lambda_1)$ and  $\lambda_2'\in (\lambda_2,\lambda_2+\varepsilon)$ such that 
    $\nu(\partial B(y,\lambda_1'r))=0$ and $\nu(\partial B(y,\lambda_2'r))=0$. Thus, for every $\varepsilon$ and for $j$ sufficiently big we have that 
    $$\mu_j(B(y,(\lambda_1-\varepsilon)r))\leq \mu_j(B(y,\lambda_1'r))\leq (D+\varepsilon)\mu_j(B(y,\lambda_2'r))\leq (D+\varepsilon)\mu_j(B(y,(\lambda_2+\varepsilon)r)).$$
    This concludes the proof. 
\end{proof}

\begin{proposizione}\label{prop:doubllimit}
    Let $\mu_i$ be a sequence of Radon measures weak* converging to some other Radon measure $\nu$, and let $B_1,B_2$ be two balls such that $\nu(B_2)>0$ and
$\nu(B_1)\leq \Lambda\nu(B_2)$. Then
    there exists $i_0=i_0(B_1,B_2)\in\N$ such that for every $i\geq i_0$ we have that 
$$\mu_i(\tfrac{1}{2}B_1)\leq16 \Lambda \mu_i(2B_2).$$
In addition, for every ball $B$ for which $\nu(B/2)>0$ and $\vartheta,\delta>0$, there exists $j_0=j_0(B,\vartheta,\delta)\in \N$ such that 
$$\mu_i(B\setminus B(\supp \nu,\delta r(B)))\leq \vartheta \mu_i(B)\qquad\text{for every }i\geq j_0.$$
\end{proposizione}

\begin{proof}
    Let us notice that 
\begin{equation}
    \begin{split}
&\qquad\qquad\mu_i(\tfrac{1}{2}B_1)\leq \frac{2}{r(B_1)}\int \dist(z,B_1^c)d\mu_i\leq \frac{2}{r(B_1)}\int \dist(z,B_1^c)d\nu+\frac{2F_{B_1}(\nu,\mu_i)}{r(B_1)}\\        
&\leq 2\nu(B_1)+\frac{2F_{B_1}(\nu,\mu_i)}{r(B_1)}\leq 2\Lambda \nu(B_2)+\frac{2F_{B_1}(\nu,\mu_i)}{r(B_1)}\leq \frac{4\Lambda}{r(B_2)}\int \dist(z,2B_2^c)d\nu(z) +\frac{2F_{B_1}(\nu,\mu_i)}{r(B_1)}\\
&\qquad\qquad\qquad\leq \frac{4\Lambda}{r(B_2)}\int \dist(z,2B_2^c)d\mu_i(z)+4\Lambda\frac{F_{2B_2}(\nu,\mu_i)}{r(B_2)} +\frac{2F_{B_1}(\nu,\mu_i)}{r(B_1)}\\
&\qquad\qquad\qquad\qquad\qquad\leq 8\Lambda \mu_i(2B_2)+4\Lambda\frac{F_{2B_2}(\nu,\mu_i)}{r(B_2)} +\frac{2F_{B_1}(\nu,\mu_i)}{r(B_1)}.
\nonumber
    \end{split}
\end{equation}
Since $\nu(B_2)>0$, Portmanteau's lemma gives $\mu_i(2B_2)\geq c(B_2)>0$ for all $i$ large enough. Since the $F$-terms tend to $0$, they can be absorbed into the term $16\Lambda\mu_i(2B_2)$.

Let us prove the second part. Let $0<\delta<1$ and set
$$
f(z)=\min\{\dist(z,(2B)^c),\dist(z,B(\supp\nu,\delta r(B)/2))\}.
$$
Then $f$ is $1$-Lipschitz, supported in $2B$, vanishes on $\supp\nu$, and satisfies
$f\geq \frac{\delta r(B)}{2}$ on $B\setminus B(\supp\nu,\delta r(B))$. Hence
$$
\mu_i(B\setminus B(\supp\nu,\delta r(B)))
\leq \frac{2}{\delta r(B)}F_{2B}(\mu_i,\nu).
$$
Since $F_{2B}(\mu_i,\nu)\to0$ by item 3 of \cref{preiss} and $\mu_i(B)$ is bounded from below for all $i$ large enough, the claim follows.
\end{proof}

\subsection{Avoiding approximating functions}

In this section we isolate a finite-dimensional obstruction which will be used later to construct Lipschitz functions with prescribed non-approximation properties. The basic point is elementary: if $\mathcal F$ is a $d$-dimensional space of functions, then one can prescribe values at $d+1$ points in such a way that no element of $\mathcal F$ can match all of them simultaneously. When $\mathcal F$ contains the constants, the prescribed values may be normalized to vanish at the first point. Averaging this finite-dimensional statement over small balls then gives the form needed in the sequel.

For affine functions this is exactly the device used by David and Semmes in \cite[Lemma 4.78, Lemma 4.80]{DavidSemmes} and in \cite[pp. 291--293]{DavidSemmes}. As observed in \cite[Chapter III.4, \S 4.3]{DavidSemmes}, the argument does not rely on the specific structure of affine functions, but only on finite dimensionality.

\begin{proposizione}\label{evadelinearV2}
Let $\mathcal F$ be a vector space of
real-valued functions on $\mathbb R^n$ of dimension $d$. Then there
exist $d+1$ vectors $\alpha^\iota\in\mathbb R^{d+1}$ such that
$|\alpha^\iota|\leq 1$ for every $\iota=1,\ldots,d+1$ and for every $r>0$ we have
$$
\inf_{\mathbf x\in\mathbb R^{n\times(d+1)}}
\sum_{\iota=1}^{d+1}
\inf_{f\in\mathcal F}
\sup_{j=1,\ldots,d+1}
|r\alpha_j^\iota-f(x_j)|
\geq
r,
$$
where $\mathbf x=(x_1,\ldots,x_{d+1})$ with $x_j\in\mathbb R^n$.
\end{proposizione}

\begin{proof}
Take
$\alpha^\iota=e_\iota$ for $\iota=1,\ldots,d+1$,
where $e_1,\ldots,e_{d+1}$ is the standard basis of $\mathbb R^{d+1}$.
Fix $\mathbf x=(x_1,\ldots,x_{d+1})$ and consider the vector subspace
$$
\{(f(x_1),\ldots,f(x_{d+1})):f\in\mathcal F\}
\subseteq \mathbb R^{d+1}.
$$
Its dimension is at most $d$, hence it is a proper subspace of
$\mathbb R^{d+1}$. Therefore there exist
$\lambda_1,\ldots,\lambda_{d+1}\in\mathbb R$, not all zero, such that
$$
\sum_{j=1}^{d+1}\lambda_j f(x_j)=0
\qquad\text{for every }f\in\mathcal F.
$$
Multiplying all $\lambda_j$ by a common constant, we may assume that
$$
\sum_{j=1}^{d+1}|\lambda_j|=1.
$$

For every $\iota=1,\ldots,d+1$ and every $f\in\mathcal F$, we have
$$
r|\lambda_\iota|
=
\bigg|\sum_{j=1}^{d+1}\lambda_j
\bigl(r e_{\iota,j}-f(x_j)\bigr)\bigg|
\leq
\sup_{j=1,\ldots,d+1}|r e_{\iota,j}-f(x_j)|.
$$
Hence
$$
\inf_{f\in\mathcal F}
\sup_{j=1,\ldots,d+1}|r e_{\iota,j}-f(x_j)|
\geq
r|\lambda_\iota|.
$$
Summing over $\iota=1,\ldots,d+1$, we get
$$
\sum_{\iota=1}^{d+1}
\inf_{f\in\mathcal F}
\sup_{j=1,\ldots,d+1}|r e_{\iota,j}-f(x_j)|
\geq
r\sum_{\iota=1}^{d+1}|\lambda_\iota|
=
r.
$$
Since this holds for every $\mathbf x\in\mathbb R^{n\times(d+1)}$, taking the
infimum over $\mathbf x$ gives the claim.
\end{proof}

\begin{corollario}\label{evadeaffine}
Let $\mathcal F$ be a vector space of
real-valued functions on $\mathbb R^n$ of dimension $d$, containing the constant functions. Then there exist $d$ vectors
$\beta^\iota\in\mathbb R^{d+1}$ such that $|\beta^\iota|\leq 1$ and
$\beta^\iota_1=0$ for every $\iota=1,\ldots,d$ such that for every $r>0$ we have
$$
\inf_{\mathbf x\in\mathbb R^{n\times(d+1)}}
\sum_{\iota=1}^{d}
\inf_{f\in\mathcal F}
\sup_{j=1,\ldots,d+1}
|r\beta^\iota_j-f(x_j)|
\geq
\frac{r}{2},
$$
where $\mathbf x=(x_1,\ldots,x_{d+1})$ with $x_j\in\mathbb R^n$.
\end{corollario}

\begin{proof}
Let $\alpha^1,\ldots,\alpha^d\in\mathbb R^d$ be the vectors given by
\cref{evadelinearV2}, applied to a vector space of dimension $d-1$. Define
$$
\beta^\iota=(0,\alpha^\iota)
\qquad\text{for }\iota=1,\ldots,d.
$$
Then $|\beta^\iota|\leq 1$ and $\beta^\iota_1=0$.
Fix $\mathbf x=(x_1,\ldots,x_{d+1})$ and notice that since $\mathcal F$ contains the constant
functions, the vector space
$$
\{y\mapsto f(y)-f(x_1):f\in\mathcal F\}
$$
has dimension $d-1$. Therefore, applying \cref{evadelinearV2} to this vector
space and to the points $x_2,\ldots,x_{d+1}$, we get
$$
\sum_{\iota=1}^{d}
\inf_{f\in\mathcal F}
\sup_{j=2,\ldots,d+1}
|r\alpha^\iota_{j-1}-f(x_j)+f(x_1)|
\geq
r.
$$
For every $f\in\mathcal F$ and every $j=2,\ldots,d+1$, one has
$$
|r\alpha^\iota_{j-1}-f(x_j)+f(x_1)|
\leq
|r\beta^\iota_j-f(x_j)|+|r\beta^\iota_1-f(x_1)|.
$$
Hence
$$
\sup_{j=2,\ldots,d+1}
|r\alpha^\iota_{j-1}-f(x_j)+f(x_1)|
\leq
2\sup_{j=1,\ldots,d+1}
|r\beta^\iota_j-f(x_j)|.
$$
Taking the infimum over $f\in\mathcal F$ and summing in
$\iota=1,\ldots,d$, we obtain
$$
r
\leq
2
\sum_{\iota=1}^{d}
\inf_{f\in\mathcal F}
\sup_{j=1,\ldots,d+1}
|r\beta^\iota_j-f(x_j)|.
$$
Therefore
$$
\sum_{\iota=1}^{d}
\inf_{f\in\mathcal F}
\sup_{j=1,\ldots,d+1}
|r\beta^\iota_j-f(x_j)|
\geq
\frac{r}{2}.
$$
Since this holds for every $\mathbf x\in\mathbb R^{n\times(d+1)}$, taking the
infimum over $\mathbf x$ gives the conclusion.
\end{proof}

\begin{corollario}\label{evadeaffinepro}
Let $\mathcal F$ be a vector space of $L^1_{\mathrm{loc}}(\R^n)$ real-valued functions of dimension $d$ containing the constant functions. There are
$d$ vectors $\beta^\iota\in\mathbb R^{d+1}$ with $|\beta^\iota|\leq 1$ and
$\beta^\iota_1=0$ for every $\iota=1,\ldots,d$ such that the following holds.

Suppose $\mu$ is a Radon measure on $\mathbb R^n$, let $B$ be a ball in
$\mathbb R^n$, and let $\zeta,\eta>0$. Let
$y_1,\ldots,y_{d+1}$ be points in $\operatorname{supp}\mu\cap B/8$ such that
$$
|y_{j_1}-y_{j_2}|\geq 3\zeta r(B)
\qquad\text{for every }j_1\neq j_2.
$$
Then
$$
\sum_{\iota=1}^{d}
\inf_{f\in\mathcal F}
\sup_{j=1,\ldots,d+1}
\fint_{B(y_j,\zeta r(B))}
|\eta r(B)\beta^\iota_j-f(z)|\,d\mu(z)
\geq
\frac{1}{2(d+1)} \eta r(B).
$$
\end{corollario}

\begin{osservazione}
    The separation and localization assumptions are included in the form in which the result will be applied later.
\end{osservazione}

\begin{proof}
Let $\beta^1,\ldots,\beta^d$ be the vectors given by \cref{evadeaffine}.
Set
$$
U_j:=B(y_j,\zeta r(B))
\qquad\text{for }j=1,\ldots,d+1,
$$
and define
$$
\boldsymbol\mu:=\mu\otimes\ldots\otimes\mu.
$$
The product has $d+1$ factors. Since $y_j\in\operatorname{supp}\mu$, one has
$\mu(U_j)>0$ for every $j=1,\ldots,d+1$. We also write
$\mathbf U:=U_1\times\ldots\times U_{d+1}$. For every $\iota=1,\ldots,d$ and every $f\in\mathcal F$, one has
$$
\fint_{\mathbf U}
\sup_{j=1,\ldots,d+1}
|\eta r(B)\beta^\iota_j-f(w_j)|\,d\boldsymbol\mu(\mathbf w)
\leq
\sum_{j=1}^{d+1}
\fint_{U_j}
|\eta r(B)\beta^\iota_j-f(z)|d\mu(z).
$$
Therefore
$$
\fint_{\mathbf U}
\sup_{j=1,\ldots,d+1}
|\eta r(B)\beta^\iota_j-f(w_j)|\,d\boldsymbol\mu(\mathbf w)
\leq
(d+1)
\sup_{j=1,\ldots,d+1}
\fint_{U_j}
|\eta r(B)\beta^\iota_j-f(z)|d\mu(z).
$$
Taking the infimum over $f\in\mathcal F$, we get
\begin{equation}
    \begin{split}
(d+1)
\inf_{f\in\mathcal F}
\sup_{j=1,\ldots,d+1}
\fint_{U_j}
|\eta r(B)\beta^\iota_j-f(z)|\,d\mu(z)&\geq
\inf_{f\in\mathcal F}
\fint_{\mathbf U}
\sup_{j=1,\ldots,d+1}
|\eta r(B)\beta^\iota_j-f(w_j)|\,d\boldsymbol\mu(\mathbf w)\\
&\geq
\fint_{\mathbf U}
\inf_{f\in\mathcal F}
\sup_{j=1,\ldots,d+1}
|\eta r(B)\beta^\iota_j-f(w_j)|\,d\boldsymbol\mu(\mathbf w).
\nonumber
    \end{split}
\end{equation}
Summing over $\iota=1,\ldots,d$ gives
$$
(d+1)
\sum_{\iota=1}^{d}
\inf_{f\in\mathcal F}
\sup_{j=1,\ldots,d+1}
\fint_{U_j}
|\eta r(B)\beta^\iota_j-f(z)|d\mu(z)\geq
\fint_{\mathbf U}
\sum_{\iota=1}^{d}
\inf_{f\in\mathcal F}
\sup_{j=1,\ldots,d+1}
|\eta r(B)\beta^\iota_j-f(w_j)|d\boldsymbol\mu(\mathbf w).
$$
By \cref{evadeaffine}, applied with $r=\eta r(B)$ and with
$\mathbf x=(w_1,\ldots,w_{d+1})$, the integrand satisfies
$$
\sum_{\iota=1}^{d}
\inf_{f\in\mathcal F}
\sup_{j=1,\ldots,d+1}
|\eta r(B)\beta^\iota_j-f(w_j)|
\geq
\frac{\eta}{2} r(B)
$$
for every $\mathbf w\in\mathbf U$. Hence
$$
(d+1)
\sum_{\iota=1}^{d}
\inf_{f\in\mathcal F}
\sup_{j=1,\ldots,d+1}
\fint_{U_j}
|\eta r(B)\beta^\iota_j-f(z)|\,d\mu(z)
\geq
\frac{1}{2}\eta r(B).
$$
This concludes the proof.
\end{proof}

The next proposition provides the interpolation step used later in the construction of the perturbations. We prescribe small values on finitely many separated points lying in the transverse space $V^\perp$, that in the qualitative case will be the decomposability bundle, and in the quantitative case a direction of invariance of a family of cubes that do not satisfy Carleson packing conditions. The size of the prescribed values is chosen below the separation scale, so that we are prescribing a Lipschitz condition. McShane's theorem then gives a Lipschitz extension on $V^\perp$, and composing with the orthogonal projection onto $V^\perp$ produces a function on $\R^n$ which is constant along $V$. Thus the perturbations constructed from this lemma will vary only in the directions orthogonal to $V$. This is in spirit similar to what G. David and S. Semmes did in their proof of the $1$-dimensional case in \cite{DavidSemmes}.

\begin{proposizione}\label{lemmanonfancyutile}
Let $\mathfrak r\in (0,1/2]$, $\xi\in (0,1)$ and $M\in \N$ be fixed parameters, and let $V\in \Gr(n,d)$. Let $E\subseteq B(0,10)\cap V^\perp$ be a compact set and suppose $\Sigma_\xi$ is a maximal $\xi$-separated set in $E$. Suppose further that for every $z\in \Sigma_\xi$ there are
$y_1(z),\ldots,y_M(z)\in E\cap B(z,\xi/8)$ such that
$$
|y_{i_1}(z)-y_{i_2}(z)|\geq 3\mathfrak r\xi
\qquad\text{for every }i_1\neq i_2.
$$
For every $z\in \Sigma_\xi$ and every $j=1,\ldots,M$, choose $\alpha_{j,z}\in\R$ such that
$$
|\alpha_{j,z}|\leq \frac{\mathfrak r}{8}\xi.
$$
Then there exists a $1/4$-Lipschitz function $f:\R^n\to\R$ satisfying the following properties:
\begin{enumerate}
    \item $f$ is invariant along $V$;
    \item $f(y_j(z))=\alpha_{j,z}$ for every $j=1,\ldots,M$ and every $z\in \Sigma_\xi$;
    \item $\|f\|_\infty\leq 8^{-1}\mathfrak r\xi$.
\end{enumerate}
\end{proposizione}

\begin{proof}
Let
$S:=\{y_j(z):z\in\Sigma_\xi,\ j=1,\ldots,M\}\subset V^\perp$. We first check that the function
$$
g(y_j(z)):=\alpha_{j,z}
$$
is a $1/4$-Lipschitz on $S$. So, we let $a,b\in S$. If there exists $z\in \Sigma_\xi$ such that $a=y_i(z)$ and $b=y_j(z)$ then 
$$
|\alpha_{i,z}-\alpha_{j,z}|
\leq \frac{\mathfrak r\xi}{4}
\leq \frac{1}{4}|y_i(z)-y_j(z)|.
$$
If instead $a=y_i(z)$ and $b=y_j(z')$ with $z\neq z'$, then, since $\Sigma_\xi$ is $\xi$-separated, we have
$$
|y_i(z)-y_j(z')|
\geq |z-z'|-|y_i(z)-z|-|y_j(z')-z'|
\geq \xi-\frac{\xi}{8}-\frac{\xi}{8}
=
\frac{3\xi}{4}.
$$
In addition, since $\mathfrak r\leq 1/2$ we have
$$
|\alpha_{i,z}-\alpha_{j,z'}|
\leq \frac{\mathfrak r\xi}{4}
\leq \frac{3\xi}{16}\leq
\frac{1}{4}|y_i(z)-y_j(z')|.
$$
This proves that $g$ is a $1/4$-Lipschitz function on $S$. Thus, by McShane's extension theorem, $g$ extends to a $1/4$-Lipschitz function
$\widetilde f:V^\perp\to\R$. Truncating if necessary, we may also assume that
$$
\|\widetilde f\|_\infty\leq 8^{-1}\mathfrak r\xi.
$$
Finally define
$$
f:=\widetilde f\circ\pi_{V^\perp}.
$$
Then $f$ is invariant along $V$, has Lipschitz constant at most $1/4$, satisfies
$f(y_j(z))=\alpha_{j,z}$ for every $j=1,\ldots,M$ and every $z\in\Sigma_\xi$, and
$$
\|f\|_\infty\leq 8^{-1}\mathfrak r\xi.
$$
This concludes the proof.
\end{proof}

\subsection{Width functions}
One of the main tools of the paper are width function. Width functions were introduced by Alberti, Cs\"ornyei and Preiss in \cite{acpnull}; their first published appearance occurred in the PhD thesis of A. Marchese in 2013, and later published in \cite{AlbertiMarchese}.

The underlying idea goes back to a one-dimensional observation of Zahorski \cite{zbMATH03099952}. If $K\subset \R$ is compact and has Lebesgue measure zero, then $K$ can be covered by intervals with arbitrarily small total length. Integrating the indicator function of the union of these intervals gives a Lipschitz function with arbitrarily small sup norm and derivative equal to $1$ at every point of $K$. Width functions are the higher-dimensional analogue of this construction. The role of intervals is replaced by families of curves constrained to move inside a cone, and the resulting function measures how much of such curves can pass through a given set.

In this section we recall the classical form of this construction, which is the one used later in the paper. We first introduce the relevant cones and curve families.

\begin{definizione}
Let $C$ be a closed convex cone and suppose $\gamma:I\to\R^n$ is a Lipschitz curve, where $I$ is a closed interval. We say that $\gamma$ is a $C$-curve if
$$
\gamma(t)-\gamma(s)\in C
\qquad\text{for every }s<t.
$$
\end{definizione}

\begin{definizione}
Let $C$ be a closed convex cone. We say that a set $E\subset\R^n$ is $C$-null if
$$
\mathcal H^1(E\cap \gamma(I))=0
$$
for every $C$-curve $\gamma:\R\to\R^n$.
\end{definizione}

\begin{proposizione}\label{propo:curvainpalla}
Suppose $\gamma:\R\to\R^n$ is a $C(e,\sigma)$-curve, let $B$ be a ball, and let $\lambda\in(0,1)$. If $\gamma\cap \lambda B\neq\emptyset$, then
$$
(1-\lambda)r(B)\leq \Haus^1(\gamma\cap B)\leq \frac{2r(B)}{\sqrt{1-\sigma^2}}.
$$
Moreover, if $\gamma(t)=te+\eta(t)$, where $\eta(t)-\eta(s)\in e^\perp$ for every $s,t\in\R$, then there exists an interval $I\subseteq\gamma^{-1}(B)$ such that
$$
\Leb^1(I)\geq \sqrt{1-\sigma^2}(1-\lambda)r(B).
$$
\end{proposizione}

\begin{proof}
Since the statement first part of the statement is invariant under increasing reparametrizations of $\gamma$, we may assume that
$\gamma(t)=te+\eta(t)$, where $\eta(t)-\eta(s)\in e^\perp$ for every $s,t\in\R$. In addition, we can also assume without loss of generality that $\gamma(0)\in \lambda B$. Notice that this shows that there exists an open interval $(-a,a)$ contained in $\gamma^{-1}(B)$.

Set $d:=(1-\lambda)r(B)$. Since $\gamma(0)\in\lambda B$, every point of $\partial B$ has distance at least $d$ from $\gamma(0)$. Therefore the two components of $\gamma((-a,a))\setminus\{\gamma(0)\}$ have length at least $d$, and hence
$$
\Haus^1(\gamma\cap B)\geq (1-\lambda)r(B).
$$
Moreover, if $t\in (-a,a)$ and $t\geq 0$, then
$$
t=\langle \gamma(t)-\gamma(0),e\rangle\geq \sqrt{1-\sigma^2}|\gamma(t)-\gamma(0)|.
$$
Consequently
$\Leb^1(I)\geq \sqrt{1-\sigma^2}d
=\sqrt{1-\sigma^2}(1-\lambda)r(B)$. 
Finally, let us observe that $\Haus^1(\gamma\cap B)\leq \Lip(\gamma)2r(B)\leq\frac{2r(B)}{\sqrt{1-\sigma^2}}$.
This concludes the proof.
\end{proof}

\begin{definizione}
Let $\Omega\subset\R^n$, $e\in\mathbb S^{n-1}$ and $\sigma\in(0,1)$. If $\gamma:[a,b]\to\R^n$ is a $C(e,\sigma)$-curve, we define
$$
\daleth(\gamma,e,\sigma;x):=
\int_{\gamma^{-1}(\Omega)}\langle \gamma'(t),e\rangle dt
+
\beta(\sigma)\int_{\gamma^{-1}(\Omega)}|\Pi_{e^\perp}(\gamma'(t))|dt
-
(\langle\gamma(b),e\rangle-\langle x,e\rangle).
$$
We define the width function of $\Omega$ in the direction $e$ and with aperture $\sigma$ by
$$
\omega_{e,\sigma}[\Omega](x):=
\sup_{\gamma\in\Gamma_{e,\sigma}^x}
\daleth(\gamma,e,\sigma;x),
$$
where $\Gamma_{e,\sigma}^x$ is the family of $C(e,\sigma)$-curves $\gamma:I=[a,b]\to\R^n$ whose endpoint satisfies
$$
\gamma(b)=x+se
\qquad\text{for some }s\geq0.
$$
\end{definizione}

\begin{osservazione}
Notice that the definition of $\daleth(\gamma,e,\sigma;x)$ is invariant under increasing reparametrizations of $\gamma$.
\end{osservazione}

\begin{proposizione}\label{alb}
Let $e\in\mathbb S^{n-1}$ and $\sigma\in(0,1)$, and let $E$ be a compact $C(e,\sigma)$-null set contained in a ball $B(x,r)$. Then, for every $\varepsilon>0$, there exists $\delta_0>0$ such that for every $0<\delta<\delta_0$ we have
$$
\mathcal H^1(\gamma(I)\cap B(E,\delta))<\varepsilon
$$
for every $C(e,\sigma)$-curve $\gamma:I\to\R^n$.
\end{proposizione}

\begin{proof}
    The proof of this statement can be fond in Step 1 of the proof of \cite[Lemma 4.12]{AlbertiMarchese}.
\end{proof}

\begin{proposizione}\label{p:width-open}
Let $\Omega\subset\R^n$ be an open set, let $e\in\mathbb S^{n-1}$ and $\sigma\in(0,1)$. Then, the following properties hold.
\begin{enumerate}
    \item $\omega_{e,\sigma}[\Omega]\geq0$;
    \item $\omega_{e,\sigma}[\Omega](x)\leq \omega_{e,\sigma}[\Omega](x+se)\leq \omega_{e,\sigma}[\Omega](x)+s$ for every $s>0$ and every $x\in\R^n$. Moreover, if the segment $[x,x+se]$ is contained in $\Omega$, then
    $$\omega_{e,\sigma}[\Omega](x+se)=\omega_{e,\sigma}[\Omega](x)+s;$$
    \item for every $x\in \R^n$ and every $v\in e^\perp$ we have
    $$|\omega_{e,\sigma}[\Omega](x+v)-\omega_{e,\sigma}[\Omega](x)|
    \leq \beta(\sigma)|v|.$$
\end{enumerate}
In particular, $\omega_{e,\sigma}[\Omega]$ is $(1+\beta(\sigma))$-Lipschitz whenever $\omega_{e,\sigma}[\Omega]$ is locally finite. 
Finally, if there exists $\varepsilon>0$ such that $\mathcal H^1(\gamma(I)\cap\Omega)<\varepsilon$ for every $C(e,\sigma)$-curve, we have
$$\|\omega_{e,\sigma}[\Omega]\|_\infty\leq2\varepsilon.$$
\end{proposizione}

\begin{proof}
Non-negativity is immediate. Let us check item (ii). Fix $\eta>0$ and let $\gamma_*\in\Gamma_{e,\sigma}^x$ be such that
$$
\omega_{e,\sigma}[\Omega](x)\leq \daleth(\gamma_*,e,\sigma;x)+\eta.
$$
Write the endpoint of $\gamma_*$ as $\gamma_*(b)=x+te$ for some $t\geq0$. If $t\geq s$, then $\gamma_*\in\Gamma_{e,\sigma}^{x+se}$ and
$$
\daleth(\gamma_*,e,\sigma;x+se)
=
\daleth(\gamma_*,e,\sigma;x)+s.
$$
Thus
$$
\omega_{e,\sigma}[\Omega](x)\leq \omega_{e,\sigma}[\Omega](x+se)-s+\eta
\leq \omega_{e,\sigma}[\Omega](x+se)+\eta.
$$

If $t<s$, let $\widetilde\gamma_*$ be the concatenation of $\gamma_*$ with the segment $[x+te,x+se]$. Then $\widetilde\gamma_*\in\Gamma_{e,\sigma}^{x+se}$ and
$$
\daleth(\widetilde\gamma_*,e,\sigma;x+se)
\geq
\daleth(\gamma_*,e,\sigma;x).
$$
Therefore again
$\omega_{e,\sigma}[\Omega](x)\leq \omega_{e,\sigma}[\Omega](x+se)+\eta$. Letting $\eta$ go to $0$, we infer that
$$
\omega_{e,\sigma}[\Omega](x)\leq \omega_{e,\sigma}[\Omega](x+se).
$$

For the opposite inequality, let $\gamma_*\in\Gamma_{e,\sigma}^{x+se}$. Then $\gamma_*$ is also admissible for $x$, and the endpoint correction changes by exactly $s$. Hence
$$
\daleth(\gamma_*,e,\sigma;x+se)
=
\daleth(\gamma_*,e,\sigma;x)+s
\leq
\omega_{e,\sigma}[\Omega](x)+s.
$$
Taking the supremum over $\gamma_*\in\Gamma_{e,\sigma}^{x+se}$ gives
$$
\omega_{e,\sigma}[\Omega](x+se)\leq \omega_{e,\sigma}[\Omega](x)+s.
$$

Moreover, assume that $[x,x+se]\subset\Omega$. Let $\gamma_*\in\Gamma_{e,\sigma}^x$ be almost maximizing for $\omega_{e,\sigma}[\Omega](x)$. If the endpoint of $\gamma_*$ already belongs to the ray $x+se+[0,+\infty)e$, then $\gamma_*$ is also admissible for $x+se$, and 
$$
\daleth(\gamma_*,e,\sigma;x+se)
=
\daleth(\gamma_*,e,\sigma;x)+s.
$$
Otherwise, we concatenate $\gamma_*$ with the remaining subsegment of $[x,x+se]$. Since this subsegment is contained in $\Omega$, the contribution of the added piece to the first integral is exactly its length, while the transverse contribution is zero. This exactly compensates the change of endpoint, and again gives
$$
\daleth(\widetilde\gamma_*,e,\sigma;x+se)
=
\daleth(\gamma_*,e,\sigma;x)+s.
$$
From this, we infer that
$$
\omega_{e,\sigma}[\Omega](x+se)\geq \omega_{e,\sigma}[\Omega](x)+s.
$$
Together with the upper bound, this proves the equality in item (ii).

Finally, we prove item (iii). Let $v\in e^\perp$ and let $\eta>0$. Choose
$\gamma_*\in\Gamma_{e,\sigma}^x$ such that
$$
\omega_{e,\sigma}[\Omega](x)\leq \daleth(\gamma_*,e,\sigma;x)+\eta.
$$
Let $p:=\gamma_*(b)$ and $q:=p+v+\beta(\sigma)|v|e$. Since
$$
v+\beta(\sigma)|v|e\in C(e,\sigma),
$$
the concatenation $\widetilde\gamma_*$ of $\gamma_*$ with the segment $[p,q]$ is a
$C(e,\sigma)$-curve and
$\widetilde\gamma_*\in\Gamma_{e,\sigma}^{x+v}$.

By the invariance of $\daleth$ under increasing reparametrizations, we may parametrize
the added segment linearly. Thus
\begin{equation}
\begin{split}
\daleth(\widetilde\gamma_*,e,\sigma;x+v)
=\daleth(\gamma_*,e,\sigma;x)+&\int_{\widetilde\gamma_*^{-1}(\Omega\cap[p,q])}
\langle \widetilde\gamma_*'(t),e\rangle dt+\beta(\sigma)
\int_{\widetilde\gamma_*^{-1}(\Omega\cap[p,q])}
|\Pi_{e^\perp}(\widetilde\gamma_*'(t))|dt\\
-&(\langle q,e\rangle-\langle p,e\rangle).
\nonumber
\end{split}
\end{equation}
However, since $\langle q,e\rangle-\langle p,e\rangle=\beta(\sigma)|v|$, we see that
$$
\daleth(\widetilde\gamma_*,e,\sigma;x+v)
\geq
\daleth(\gamma_*,e,\sigma;x)-\beta(\sigma)|v|.
$$
It follows that
$$
\omega_{e,\sigma}[\Omega](x+v)
\geq
\daleth(\widetilde\gamma_*,e,\sigma;x+v)
\geq
\omega_{e,\sigma}[\Omega](x)-\eta-\beta(\sigma)|v|.
$$
Exchanging the roles of $x$ and $x+v$ and letting $\eta$ go to $0$ gives
$$|\omega_{e,\sigma}[\Omega](x+v)-\omega_{e,\sigma}[\Omega](x)|\leq \beta(\sigma)|v|.$$
This proves item (iii). The Lipschitz constant of $\omega_{e,\sigma}[\Omega]$ can be easily inferred from item (ii) and (iii) whenever $\omega_{e,\sigma}[\Omega]$ is locally bounded.

It remains to prove the quantitative conclusion. Assume that
$$
\mathcal H^1(\gamma(I)\cap\Omega)<\varepsilon
$$
for every $C(e,\sigma)$-curve $\gamma:I\to\R^n$. Therefore
\begin{equation}
    \begin{split}
        &\int_{\gamma^{-1}(\Omega)}\langle \gamma'(t),e\rangle dt
        +
        \beta(\sigma)\int_{\gamma^{-1}(\Omega)}
        |\Pi_{e^\perp}(\gamma'(t))|dt\\
        &\qquad\qquad\qquad\qquad\leq
        2\int_{\gamma^{-1}(\Omega)}\langle \gamma'(t),e\rangle dt
        \leq
        2\mathcal H^1(\gamma(I)\cap\Omega)
        \leq 2\varepsilon,
    \end{split}
\end{equation}
and this immediately implies
$$
\omega_{e,\sigma}[\Omega](x)\leq2\varepsilon
\qquad\text{for every }x\in\R^n.
$$
Together with $\omega_{e,\sigma}\geq0$, this proves
$\|\omega_{e,\sigma}\|_\infty\leq2\varepsilon$.
This concludes the proof.
\end{proof}

\subsection{The coefficients $\Omega$}

In this subsection we introduce the coefficients $\Omega$ associated with the
finite-dimensional space of approximating functions used throughout the paper.

\begin{definizione}
Let $q\in[1,\infty]$ and suppose $\mu$ is a Radon measure on $\R^n$.
For every Lipschitz function $f:\R^n\to \R$, every $x\in \R^n$, and every
$r>0$ such that $\mu(B(x,r))>0$, we define
$$
\Omega_{q,\mu,\mathcal F}(f;x,r):=
\inf_{A\in\mathcal F}
\bigg(\fint_{B(x,r)}\frac{|f-A|^q}{r^q}\,d\mu\bigg)^{1/q}
$$
if $1\leq q<\infty$, while for $q=\infty$ we define
$$
\Omega_{\infty,\mu,\mathcal F}(f;x,r):=
\inf_{A\in\mathcal F}
\Big\|\frac{f-A}{r}\Big\|_{L^\infty(\mu\trace B(x,r))}.
$$
When the measure $\mu$ and the space $\mathcal F$ are clear from the context,
we simply write $\Omega_q(f;x,r)$.
\end{definizione}

\begin{proposizione}\label{omega_monotone}
Let $1\leq p\leq q\leq \infty$, let $\mu$ be a Radon measure on $\R^n$,
and let $f:\R^n\to\R$ be Lipschitz. Then, for every $x\in\R^n$ and
$r>0$ such that $\mu(B(x,r))>0$, we have
$$
\Omega_{p,\mu,\mathcal F}(f;x,r)\leq \Omega_{q,\mu,\mathcal F}(f;x,r).
$$
\end{proposizione}

\begin{proof}
If $q<\infty$, then for every $A\in\mathcal F$, H\"older's inequality gives
$$
\bigg(\fint_{B(x,r)}\frac{|f-A|^p}{r^p}\,d\mu\bigg)^{1/p}
\leq
\bigg(\fint_{B(x,r)}\frac{|f-A|^q}{r^q}\,d\mu\bigg)^{1/q}.
$$
Taking the infimum over $A\in\mathcal F$ gives
$$
\Omega_{p,\mu,\mathcal F}(f;x,r)\leq \Omega_{q,\mu,\mathcal F}(f;x,r).
$$
If $q=\infty$, the conclusion is obtained immediately. 
\end{proof}

\begin{proposizione}\label{propcompatibilityOmega}
Let $\mu$ be an $\alpha$-dimensional AD-regular measure on $\R^n$ with regularity constant $D\geq1$, and let $\mathcal F$ be a finite dimensional vector subspace of $L^1_{\mathrm{loc}}(\R^n)$. Let $1\leq q\leq\infty$ and let $f\in L^1_{\mathrm{loc}}(\mu)$. Then, for every $A\geq 1$, $x,y\in\supp\mu$ and $r,s>0$ such that
$$r\leq s\leq Ar\qquad\text{and}\qquad|x-y|\leq As,$$
we have
$$
\Omega_{q,\mu,\mathcal F}(f;x,r)\leq D^2(2A^2)^{\alpha+1}\Omega_{q,\mu,\mathcal F}(f;y,2As)
$$
\end{proposizione}

\begin{proof}
Let $A\geq1$ and let $x,y\in\supp\mu$, $r,s>0$ be as in the statement. 
Notice that 
$B(x,r)\subset B(y,2As)$. Let us restrict ourselves to $1\leq q<\infty$, as the case $q=\infty$ will follow similarly. For every $L\in\mathcal F$ we have
$$
\bigg(\fint_{B(x,r)}\frac{|f-L|^q}{r^q}\,d\mu\bigg)^{1/q}
\leq
\frac{2As}{r}
\bigg(\frac{\mu(B(y,2As))}{\mu(B(x,r))}\bigg)^{1/q}
\bigg(\fint_{B(y,2As)}\frac{|f-L|^q}{(2As)^q}\,d\mu\bigg)^{1/q}.
$$
By AD-regularity and the assumptions $r\leq s\leq Ar$, we have
$$
\frac{2As}{r}
\bigg(\frac{\mu(B(y,2As))}{\mu(B(x,r))}\bigg)^{1/q}\leq D^2(2A^2)^{\alpha+1}.
$$
Taking the infimum over $L\in\mathcal F$ concludes the proof.
\end{proof}

\subsection{Carleson families and cubes}

Throughout this section $\mu$ is supposed to be a $k$-dimensional AD-regular measure with regularity constant $D\geq 1$, and we set
$E:=\supp\mu$. We introduce in this subsection some elementary and well known objects. We refer to \cite{DavidSemmes} for a comprehensive discussion. First of all, we recall the continuous notion of Carleson set. 

\begin{definizione}\label{defcontinuouscarleson}
We say that a Borel set $\mathscr B\subseteq E\times(0,\infty)$ is Carleson if there exists a constant $C>0$ such that, for every $x\in E$ and $R>0$,
$$
\int_0^R\int_{B(x,R)}\mathbb 1_{\mathscr B}(y,s)d\mu(y)\frac{ds}{s}\leq C \mu(B(x,R)).
$$
\end{definizione}

It is often useful, when dealing with uniform rectifiability, to give statements and introduce conditions in terms of the dyadic cubes associated to $\mu$. For this reason we recall here the classical dyadic cubes. For a construction we refer to the appendix of \cite{DavidLN}. For a proof in the generality of Carnot groups, we refer to the ArXiv version of \cite{MarstrandMattila20}.

\begin{teorema}\label{evev}
Let $\mu$ be a $k$-AD-regular measure with regularity constant $D$ in $\R^n$.
There exist $\newC\label{cubi}>2$, depending only on $n$ and $D$, and a family of partitions $\{\Delta_{j,\mu}\}_{j\in\mathbb Z}$ of $\supp\mu$, called \emph{layers}, with the following properties.
\begin{itemize}
    \item[(i)] If $j\leq j^\prime$, $Q\in\Delta_{j,\mu}$ and $Q^\prime\in\Delta_{j^\prime,\mu}$, then either $Q$ contains $Q^\prime$ or $Q\cap Q^\prime=\emptyset$;
    \item[(ii)] If $Q\in\Delta_{j,\mu}$ we have $\oldC{cubi}^{-1}2^{-j}\leq \diam Q\leq \oldC{cubi}2^{-j}$;
    \item[(iii)] If $Q\in \Delta_{j,\mu}$ then $\oldC{cubi}^{-1}2^{-kj}\leq \mu(Q)\leq \oldC{cubi}2^{-kj}$;
    \item[(iv)] For every $Q\in\Delta_{j,\mu}$, we have  $    \mu\big(\partial(Q,\tau)\big)\leq \oldC{cubi} \tau^{1/\oldC{cubi}}2^{-kj}$,
    where 
    $$
    \partial(Q,\tau):=\{y\in \supp\mu\setminus Q:\dist(y,Q)\leq \tau 2^{-j}\}\cup \{y\in Q:\dist(y,\supp\mu\setminus Q)\leq \tau 2^{-j}\},
    $$
    for every $0<\tau<1$.
    \item[(v)] If $Q\in\Delta_{j,\mu}$, there exists $\mathfrak c(Q)\in Q$ such that 
    $$
    \dist(\mathfrak c(Q),\supp\mu\setminus Q)\geq \frac{3}{2}\oldC{cubi}^{-2\oldC{cubi}}2^{-j}\overset{(ii)}{\geq} \frac{3}{2}\oldC{cubi}^{-2\oldC{cubi}-1}\diam Q.
    $$
\end{itemize}
We will denote with the symbol $\Delta_\mu$ the family of all dyadic cubes, i.e.
$$
\Delta_\mu:=\bigcup_{j\in\Z}\Delta_{j,\mu}.
$$
Further, we denote $\ell(Q):=2^{-j}$ if $Q\in \Delta_{j,\mu}$. 
\end{teorema}

\begin{definizione}\label{defdyadiccarleson}
We say that a family of cubes $\mathscr F\subseteq \Delta_\mu$ is Carleson if there exists a constant $C>0$ such that for every $R\in \Delta_\mu$ we have
$$
\sum_{Q\in \mathscr F(R)}\mu(Q)\leq C\mu(R),
$$
where
$\mathscr F(R):=\{Q\in \mathscr F:Q\subseteq R\}$. 
\end{definizione}

We shall frequently pass from continuous coefficients to their dyadic counterparts. We first record the elementary abstract principle which justifies this passage. The point is that the equivalence between the continuous and dyadic Carleson formulations is valid for coefficients which are stable under bounded changes of centre and scale in the following sense.

\begin{definizione}\label{defscalecompatible}
Let $\mathfrak a:E\times(0,\infty)\to[0,\infty)$ be a Borel function. We say that $\mathfrak a$ is \emph{compatible} if for every $A\geq1$ there exists $C_A\geq1$ such that
$$
\mathfrak a(x,r)\leq C_A\mathfrak a(y,2A s)
$$
whenever $x,y\in E$, $r,s>0$ are such that
$r\leq s\leq A r$ and $|x-y|\leq As$. In this case, for every $Q\in\Delta_\mu$ we set
$$
\mathfrak a(Q):=\mathfrak a(\mathfrak c(Q),\diam Q).
$$
\end{definizione}

We shall use the following abstract passage between continuous and dyadic Carleson estimates. The statement is indipendentent on the specific definition of the number: we just require that when the centre is moved inside one dyadic cube and the scale is changed by a controlled factor, the coefficients don't change too much. 

\begin{lemma}\label{lemcontinuousdyadiccompatible}
Suppose $\mu$ is an $\alpha$-dimenional AD regular measure with regularity constant $D$. 
Let $\mathfrak b:E\times(0,\infty)\to[0,\infty)$ be a Borel function. For $\varepsilon,\eta>0$ and $\Lambda,\Gamma\geq1$, set
$$\mathscr B_\varepsilon:=\{(x,r)\in E\times(0,\infty):\mathfrak b(x,r)>\varepsilon\},\qquad \mathscr D_{\eta,\Lambda}:=\{Q\in\Delta_\mu:\mathfrak b(\mathfrak c(Q),\Lambda\diam Q)>\eta\}.$$
Assume that one of the following implications holds:
\begin{enumerate}
\item[(a)] if $\mathfrak b(x,r)>\varepsilon$, $x\in Q$ and $\frac{2r}{\Lambda}\leq \diam Q\leq \frac{4\oldC{cubi}^2 r}{\Lambda}$, then $Q\in\mathscr D_{\eta,\Lambda}$;
\item[(b)] if $Q\in\mathscr D_{\eta,\Lambda}$, $y\in Q$ and $\diam Q\leq s\leq2\diam Q$, then $\mathfrak b(y,\Gamma s)>\varepsilon$.
\end{enumerate}
Then the following conclusions hold.
\begin{enumerate}
\item[(i)] If (a) holds and $\mathscr D_{\eta,\Lambda}$ is Carleson in the sense of \cref{defdyadiccarleson}, then $\mathscr B_\varepsilon$ is Carleson in the sense of \cref{defcontinuouscarleson}.
\item[(ii)] If (b) holds and $\mathscr B_\varepsilon$ is Carleson in the sense of \cref{defcontinuouscarleson}, then $\mathscr D_{\eta,\Lambda}$ is Carleson in the sense of \cref{defdyadiccarleson}.
\end{enumerate}
The constants yielded in (i) and (ii) depend only on $k,D,\oldC{cubi},\Lambda,\Gamma$ and on the Carleson constant in the corresponding hypothesis. In particular, if $\mathfrak a:E\times(0,\infty)\to[0,\infty)$ is compatible in the sense of \cref{defscalecompatible}, and if $C_{\mathrm{cont}}(t)$ denotes the continuous Carleson constant of the set $\{(x,r)\in E\times(0,\infty):\mathfrak a(x,r)>t\}$, then, for every $\varepsilon>0$, the family $\mathscr F_\varepsilon:=\{Q\in\Delta_\mu:\mathfrak a(Q)>\varepsilon\}$ is Carleson and
$$\sum_{Q\in\mathscr F_\varepsilon(R)}\mu(Q)\leq C_{\mathrm{dyad}}(\varepsilon)\mu(R)\qquad\text{for every }R\in\Delta_\mu,$$
where $C_{\mathrm{dyad}}(\varepsilon)$ depends only on $n,k,D,\oldC{cubi}$, on the compatibility constant $C_2$ in \cref{defscalecompatible} with $A=2$, and on $C_{\mathrm{cont}}(C_2^{-1}\varepsilon)$.
\end{lemma}

\begin{proof}
We first prove (i). Fix $x_0\in E$ and $R>0$. By \cref{evev}, for every $x\in E$ and every $r>0$ there exists $Q(x,r)\in\Delta_\mu$ such that
$$x\in Q(x,r),\qquad \frac{2r}{\Lambda}\leq \diam Q(x,r)\leq \frac{4\oldC{cubi}^2 r}{\Lambda}.$$
By assumption (a), if $(x,r)\in\mathscr B_\varepsilon$, then $Q(x,r)\in\mathscr D_{\eta,\Lambda}$. Hence
$$\mathbb 1_{\mathscr B_\varepsilon}(x,r)\leq \sum_{Q\in\mathscr D_{\eta,\Lambda}}\mathbb 1_Q(x)\mathbb 1_{[\Lambda\diam Q/4\oldC{cubi}^2,\Lambda\diam Q/2]}(r).$$
Therefore
$$\int_0^R\int_{B(x_0,R)}\mathbb 1_{\mathscr B_\varepsilon}(x,r)d\mu(x)\frac{dr}{r}\leq C\sum_{\substack{Q\in\mathscr D_{\eta,\Lambda}\\ Q\cap B(x_0,R)\neq\emptyset\\ \diam Q\leq 4\oldC{cubi}^2 R/\Lambda}}\mu(Q),$$
where $C$ depends only on $\oldC{cubi}$. Let $\mathcal R$ be the family of maximal dyadic cubes of diameter at most $4\oldC{cubi}^2 R/\Lambda$ which meet $B(x_0,R)$. Every cube in the last sum is contained in some cube of $\mathcal R$. Moreover, the cubes in $\mathcal R$ are pairwise disjoint and contained in $B(x_0,CR)$, with $C$ depending only on $\oldC{cubi}$ and $\Lambda$. Thus, by the dyadic Carleson packing of $\mathscr D_{\eta,\Lambda}$ and AD-regularity,
$$\sum_{\substack{Q\in\mathscr D_{\eta,\Lambda}\\ Q\cap B(x_0,R)\neq\emptyset\\ \diam Q\leq 4\oldC{cubi}^2 R/\Lambda}}\mu(Q)\leq C\sum_{S\in\mathcal R}\mu(S)\leq C\mu(B(x_0,R)).$$
This proves that $\mathscr B_\varepsilon$ is Carleson.

We now prove (ii). Fix $R\in\Delta_\mu$ and let $Q\in\mathscr D_{\eta,\Lambda}(R)$. By assumption (b), for every $y\in Q$ and every $\diam Q\leq s\leq2\diam Q$, one has $\mathfrak b(y,\Gamma s)>\varepsilon$. Therefore
$$\log2\,\mu(Q)\leq\int_{\diam Q}^{2\diam Q}\int_Q\mathbb 1_{\mathscr B_\varepsilon}(y,\Gamma s)d\mu(y)\frac{ds}{s}.$$
Summing over $Q\in\mathscr D_{\eta,\Lambda}(R)$ and using the bounded overlap of the regions $Q\times[\Gamma\diam Q,2\Gamma\diam Q]$, we get
$$\sum_{Q\in\mathscr D_{\eta,\Lambda}(R)}\mu(Q)\leq C\int_0^{C\Gamma\diam R}\int_R\mathbb 1_{\mathscr B_\varepsilon}(y,t)d\mu(y)\frac{dt}{t}.$$
Since $R\subseteq B(\mathfrak c(R),\diam R)$, the continuous Carleson estimate for $\mathscr B_\varepsilon$ gives
$$\sum_{Q\in\mathscr D_{\eta,\Lambda}(R)}\mu(Q)\leq C\mu(B(\mathfrak c(R),C\Gamma\diam R))\leq C\mu(R).$$
Here $C$ depends only on $k,D,\oldC{cubi},\Gamma$ and on the continuous Carleson constant of $\mathscr B_\varepsilon$. This proves (ii).

We finally prove the compatible-coefficient conclusion. Let $\mathfrak a$ be compatible and fix $\varepsilon>0$. Set
$$\mathscr F_\varepsilon:=\{Q\in\Delta_\mu:\mathfrak a(Q)>\varepsilon\}.$$
Let $Q\in\mathscr F_\varepsilon$, let $y\in Q$ and let $\diam Q\leq s\leq2\diam Q$. Then
$$\diam Q\leq s\leq2\diam Q\qquad\text{and}\qquad |\mathfrak c(Q)-y|\leq\diam Q\leq2s.$$
Applying compatibility to the pairs $(\mathfrak c(Q),\diam Q)$ and $(y,s)$ with $A=2$, we get
$$\varepsilon<\mathfrak a(\mathfrak c(Q),\diam Q)\leq C_2\mathfrak a(y,4s).$$
Hence $\mathfrak a(y,4s)>C_2^{-1}\varepsilon$. Thus assumption (b) holds for $\mathfrak b=\mathfrak a$, $\eta=\varepsilon$, $\Lambda=1$, $\Gamma=4$, and continuous threshold $C_2^{-1}\varepsilon$. Applying (ii) gives
$$\sum_{Q\in\mathscr F_\varepsilon(R)}\mu(Q)\leq C_{\mathrm{dyad}}(\varepsilon)\mu(R)\qquad\text{for every }R\in\Delta_\mu,$$
with $C_{\mathrm{dyad}}(\varepsilon)$ depending only on $n,k,D,\oldC{cubi}$, on $C_2$, and on the continuous Carleson constant of $\{(x,r)\in E\times(0,\infty):\mathfrak a(x,r)>C_2^{-1}\varepsilon\}$. This concludes the proof.
\end{proof}

We shall also use the following consequence of Carleson packing.

\begin{lemma}\label{lemcontinuouspackinglimitzero}
Let $\mathfrak a:E\times(0,\infty)\to[0,\infty)$ be a Borel function such that $\mathfrak a$ is compatible in the sense of \cref{defscalecompatible}. Assume that for every $\varepsilon>0$ the set
$$
\mathscr B_\varepsilon:=\{(x,r)\in E\times(0,\infty):\mathfrak a(x,r)>\varepsilon\}
$$
is Carleson. Then
$$
\lim_{r\to0}\mathfrak a(x,r)=0\qquad \text{for $\mu$-almost every $x\in E$.}
$$
\end{lemma}

\begin{proof}
We first prove that $\limsup_{r\to0}\mathfrak a(x,r)\leq\varepsilon$
for $\mu$-almost every $x\in E$. Let $x_0\in E$ and $R_0>0$. Since $\mathscr B_{\delta}$ is Carleson for every $\delta>0$, we have
$$
\int_0^{R_0}\int_{B(x_0,R_0)}
\mathbb 1_{\mathscr B_{\delta}}(y,s)d\mu(y)\frac{ds}{s}<\infty.
$$
By Fubini's theorem, for $\mu$-almost every $x\in E\cap B(x_0,R_0)$, we have
$$
\int_0^{R_0}
\mathbb 1_{\{\mathfrak a(x,s)>\delta\}}\frac{ds}{s}<\infty.
$$
We claim that every such $x$ satisfies $\limsup_{r\to0}\mathfrak a(x,r)\leq\varepsilon$. Indeed, if this was not the case, there would be an infinitesimal sequence of radii $r_i$ such that $\mathfrak a(x,r_i)>\varepsilon$. Passing to a subsequence, we may assume that the intervals $[4r_i,8r_i]$ are pairwise disjoint and contained in $(0,R_0)$. If $s\in[r_i,2r_i]$, then compatibility applied to the pairs $(x,r_i)$ and $(x,s)$ with $A=2$ gives
$$
\varepsilon<\mathfrak a(x,r_i)\leq C\mathfrak a(x,4s).
$$
Therefore, if we choose $\delta=C^{-1}\varepsilon$, we have
$$
\int_0^{R_0}
\mathbb 1_{\{\mathfrak a(x,t)>C^{-1}\varepsilon\}}\frac{dt}{t}
\geq
\sum_{i\in\N}\int_{4r_i}^{8r_i}\frac{dt}{t}
=\infty,
$$
which is a contradiction.
Thus $\limsup_{r\to0}\mathfrak a(x,r)\leq\varepsilon$ for $\mu$-almost every $x\in E\cap B(x_0,R_0)$. Taking a countable covering of $E$ by balls and then applying the previous conclusion to every $\varepsilon\in\mathbb Q\cap(0,\infty)$, we reach the conclusion.
\end{proof}

\subsection{Uniform rectifiability criteria}
\label{URcriteria}

We now recall the David--Semmes criteria for uniform rectifiability that will be used later. As in the previous subsection, we suppose $\mu$ is a $k$-dimensional AD-regular measure with regularity constant $D$. Let us recall some well-known conditions.

\begin{definizione}[BAUP]\label{defBAUP}
For $x\in \supp\mu$ and $r>0$, define
$$Ub\beta_\mu^k(x,r):=\inf_{\mathscr P}\frac{1}{r}\bigg(\sup_{z\in \supp\mu\cap B(x,r)}\dist(z,\bigcup_{P\in\mathscr P}P)+\sup_{z\in (\bigcup_{P\in\mathscr P}P)\cap B(x,r)}\dist(z,\supp\mu)\bigg),$$
where the infimum is taken over all families $\mathscr P\subseteq\Pi_k$ of affine $k$-planes. We say that $\mu$ satisfies BAUP, the \emph{bilateral approximation by unions of planes} condition, if for every $\varepsilon>0$ the set
$$\{(x,r)\in\supp\mu\times(0,\infty):Ub\beta_\mu^k(x,r)>\varepsilon\}\qquad\text{is Carleson.}$$
\end{definizione}

\begin{definizione}[OUWGL]\label{defOUWGL}
For $x\in \supp\mu$ and $r>0$, define
$$u\beta_\mu^k(x,r):=
\inf_{\substack{P\in\Pi_k\\ x\in P}}
\sup_{z\in P\cap \overline{B(x,r)}}
\frac{\dist(z,\supp\mu)}{r}.$$
We say that $\mu$ satisfies the other unilateral weak geometric lemma, or OUWGL, if for every $\varepsilon>0$ the set
$$\{(x,r)\in \supp\mu\times(0,\infty):u\beta_\mu^k(x,r)>\varepsilon\}\qquad\text{is Carleson.}$$
\end{definizione}

\begin{teorema}[{\cite[\S II.3, Proposition 3.18]{DavidSemmes}}]\label{teoDScriteria}
If $\mu$ satisfies BAUP, then $\mu$ is uniformly rectifiable.
\end{teorema}

For later use, we also record the fixed-threshold version of the BAUP criterion. This formulation is observed explicitly in \cite[Remark 3.0.1]{BurnazyanThesis}, and follows from the proof of the David--Semmes of the fact that BAUP implies WHIP and WTP, together with WHIP and WTP imply uniform rectifiability.

\begin{teorema}\label{teosmallBAUP}
There exists $\varepsilon_0=\varepsilon_0(n,k,D)>0$ with the following property. If
$$
\mathscr B^1_{\varepsilon_0}:=\{(x,r)\in \supp\mu\times(0,\infty):Ub\beta_\mu^k(x,r)>\varepsilon_0\}
$$
is Carleson, then $\mu$ is uniformly rectifiable.
\end{teorema}

\begin{proof}
This is the fixed-threshold version of the BAUP criterion. More precisely, by \cite[Remark 3.0.1]{BurnazyanThesis}, there exists $\varepsilon_0=\varepsilon_0(n,k,D)>0$ such that the Carleson packing of the set
$$
\{(x,r)\in \supp\mu\times(0,\infty):Ub\beta_\mu^k(x,r)>\varepsilon_0\}
$$
implies that $\mu$ is uniformly rectifiable.
\end{proof}

\begin{proposizione}\label{propsmallOUWGL}
There exists $\tilde\varepsilon_0=\tilde\varepsilon_0(n,k,D)>0$ such that if
$$
\mathscr B^2_{\tilde\varepsilon_0}:=\{(x,r)\in \supp\mu\times(0,\infty):u\beta_\mu^k(x,r)>\tilde\varepsilon_0\}
$$
is Carleson, then $\mu$ is uniformly rectifiable.
\end{proposizione}

\begin{proof}
Let $\varepsilon_0=\varepsilon_0(n,k,D)>0$ be the constant given by \cref{teosmallBAUP}, and set
$\tilde\varepsilon_0:=\frac{\varepsilon_0}{2}$. We use the notation of \cite[\S II.3]{DavidSemmes}. For $\eta>0$, let $\mathscr G_o(\eta)$ be the set of all pairs $(x,t)\in \supp\mu\times(0,\infty)$ for which there exists a $k$-plane $P$ containing $x$ such that every point of $P\cap B(x,t)$ has distance at most $\eta t$ from $\supp\mu$. Thus
$$
\supp\mu\times(0,\infty)\setminus \mathscr G_o(\eta)
=
\{(x,t)\in \supp\mu\times(0,\infty):u\beta_\mu^k(x,t)>\eta\}.
$$
By assumption, $\supp\mu\times(0,\infty)\setminus \mathscr G_o(\tilde\varepsilon_0)$ is Carleson. Following \cite[p. 141]{DavidSemmes}, define $\mathscr G_o^*(\tilde\varepsilon_0)$ as the set of all $(x,t)\in \supp\mu\times(0,\infty)$ such that $(y,s)\in\mathscr G_o(\tilde\varepsilon_0)$ whenever
$$
y\in \supp\mu,\qquad |y-x|\leq t,\qquad \frac{t}{2}\leq s\leq t.
$$
Since $\supp\mu\times(0,\infty)\setminus \mathscr G_o(\tilde\varepsilon_0)$ is Carleson, the enlarged bad set
$$
\supp\mu\times(0,\infty)\setminus \mathscr G_o^*(\tilde\varepsilon_0)
$$
is Carleson as well.

We now recall the key inclusion from the proof of \cite[\S II.3, Proposition 3.17]{DavidSemmes}:
$$
\mathscr G_u(\varepsilon_0)\supseteq\{(x,t)\in \supp\mu\times(0,\infty):(x,2t)\in\mathscr G_o^*(\varepsilon_0/2)\},
$$
where $\mathscr G_u(\varepsilon_0)$ denotes the good set for BAUP at threshold $\varepsilon_0$. It follows that
$$
\supp\mu\times(0,\infty)\setminus\mathscr G_u(\varepsilon_0)
\subseteq
\{(x,t)\in \supp\mu\times(0,\infty):(x,2t)\notin\mathscr G_o^*(\tilde\varepsilon_0)\}.
$$
The right hand side is Carleson, because $\supp\mu\times(0,\infty)\setminus\mathscr G_o^*(\tilde\varepsilon_0)$ is Carleson and the change of scale $t\mapsto 2t$ preserves Carleson estimates up to an absolute multiplicative constant. Thus
$$
\{(x,t)\in \supp\mu\times(0,\infty):Ub\beta_\mu^k(x,t)>\varepsilon_0\}
$$
is Carleson as well. By \cref{teosmallBAUP}, $\mu$ is uniformly rectifiable.
\end{proof}

\begin{osservazione}\label{remOUWGLdiscretization}
The coefficient $\mathfrak b(x,r):=u\beta_\mu^k(x,r)$
satisfies the hypotheses of \cref{lemcontinuousdyadiccompatible} with suitable thresholds. More precisely, assumption (a) holds whenever $$4\oldC{cubi}^2(\eta+\Lambda^{-1})\leq\varepsilon.$$
Indeed, suppose that $u\beta_\mu^k(x,r)>\varepsilon$, $x\in Q$ and $\frac{2r}{\Lambda}\leq\diam Q\leq\frac{4\oldC{cubi}^2 r}{\Lambda}$. If $u\beta_\mu^k(\mathfrak c(Q),\Lambda\diam Q)\leq\eta$, then there exists a $k$-plane $P_Q$ containing $\mathfrak c(Q)$ such that
$$\sup_{z\in P_Q\cap B(\mathfrak c(Q),\Lambda\diam Q)}\frac{\dist(z,\supp\mu)}{\Lambda\diam Q}\leq\eta.$$
Let $P_x:=P_Q+x-\mathfrak c(Q)$. If $z\in P_x\cap B(x,r)$, then $z+\mathfrak c(Q)-x\in P_Q\cap B(\mathfrak c(Q),\Lambda\diam Q)$ and hence
$$\dist(z,\supp\mu)\leq\eta\Lambda\diam Q+\diam Q\leq 4\oldC{cubi}^2(\eta+\Lambda^{-1})r\leq\varepsilon r.$$
This contradicts $u\beta_\mu^k(x,r)>\varepsilon$, and proves assumption (a).

Similarly, assumption (b) holds with $\Gamma=2\Lambda$ whenever $4\varepsilon+\Lambda^{-1}\leq\eta$.
Indeed, suppose that $u\beta_\mu^k(\mathfrak c(Q),\Lambda\diam Q)>\eta$, let $y\in Q$, and let $\diam Q\leq s\leq2\diam Q$. If $u\beta_\mu^k(y,2\Lambda s)\leq\varepsilon$, then there exists a $k$-plane $P_y$ containing $y$ such that
$$\sup_{z\in P_y\cap B(y,2\Lambda s)}\frac{\dist(z,\supp\mu)}{2\Lambda s}\leq\varepsilon.$$
Let $P_Q:=P_y+\mathfrak c(Q)-y$. If $z\in P_Q\cap B(\mathfrak c(Q),\Lambda\diam Q)$, then $z+y-\mathfrak c(Q)\in P_y\cap B(y,2\Lambda s)$ and hence
$$\dist(z,\supp\mu)\leq2\varepsilon\Lambda s+\diam Q\leq(4\varepsilon+\Lambda^{-1})\Lambda\diam Q\leq\eta\Lambda\diam Q.$$
This contradicts $u\beta_\mu^k(\mathfrak c(Q),\Lambda\diam Q)>\eta$, and proves assumption (b).
\end{osservazione}

The following condition will become useful when dealing with the quantitative versions of the regularity problem. 

\begin{definizione}[Weak no reels]\label{defWNR}
Let $0<k<n$ and write every point $x\in\R^n$ as $x=(x',x'')$, where $x'\in\R^{n-k}$ and $x''\in\R^k$. For $0<\varepsilon<1/10$, the standard $\varepsilon$-reel is the set
$$
R(\varepsilon):=\{x\in\R^n: |x'|\leq \varepsilon \text{ and } |x''|\leq 1\}
\cup
\{x\in\R^n: |x'|<10\varepsilon \text{ and } 1-\varepsilon\leq |x''|\leq 1\}.
$$
The content of $R(\varepsilon)$ is the set
$$
\operatorname{Cont}R(\varepsilon):=\{x\in\R^n:\varepsilon<|x'|<2\varepsilon \text{ and } |x''|<1-\varepsilon\},
$$
and the inner tube of $R(\varepsilon)$ is the set
$$
\{x\in\R^n: |x'|\leq\varepsilon \text{ and } |x''|\leq1\}.
$$
A closed set $R\subseteq\R^n$ is called an $\varepsilon$-reel if $R=\phi(R(\varepsilon))$ for some composition $\phi$ of an isometry of $\R^n$ and a dilation. The content and the inner tube of $R$ are, respectively, the images under $\phi$ of the content and of the inner tube of $R(\varepsilon)$.

We denote by $\mathscr G_{\mathrm{WNR}}(\varepsilon)$ the set of all pairs $(x,t)\in \supp\mu\times(0,\infty)$ for which it is possible to find an $\varepsilon$-reel $R$ such that
$$
R\subseteq B(x,t),\qquad \diam R\geq\varepsilon t,\qquad R\cap \supp\mu=\emptyset,\qquad \supp\mu\cap\operatorname{Cont}R\neq\emptyset.
$$
We say that $\mu$ satisfies the weak no reels condition, or WNR, if for every $0<\varepsilon<1/10$ the set
$$
\supp\mu\times(0,\infty)\setminus\mathscr G_{\mathrm{WNR}}(\varepsilon)\qquad\text{is Carleson.}
$$
\end{definizione}

\begin{proposizione}[{\cite[\S II.3, Proposition 3.50]{DavidSemmes}}]
If $\mu$ satisfies the WNR, then $\mu$ is uniformly rectifiable.
\end{proposizione}

\section{Qualitative results}
\label{sezionequalitativa}

The purpose of this section is to prove the main qualitative result of the
paper. We show that if locally the measure admits scales at which it looks like many parallel copies of the decomposability bundle, then there must exist a Lipschitz function whose
$\Omega_{1,\mu,\mathcal F}$ coefficients do not vanish infinitesimally. From this, we are able to prove a rectifiability criterion for asymptotically doubling Radon measures.

The section is organized as follows. First, we prove a rectifiability criterion
in terms of the decomposability bundle. We show that if
$$
\Theta^{d,*}(\mu,x)>0
\qquad\text{and}\qquad
\dim V(\mu,x)=d
$$
for $\mu$-almost every $x$, then $\mu$ is $d$-rectifiable. This part is used to
isolate the rectifiable alternative: when the decomposability bundle has the
same dimension as the positive density of the measure, rectifiability follows
from the absolute continuity of a positive-measure family of projections and
the Besicovitch--Federer projection theorem.

The second subsection contains the perturbation lemma. This is the technical
core of the section. We introduce the metric space $X_{\mu,\mathscr E}$ of
admissible Lipschitz functions and the operators
$$
T_{\rho',\rho}f(x):=
\sup_{\rho'<r<\rho}\Omega_{1,\mu,\mathcal F}(f;x,r).
$$
We also introduce the classes $\mathfrak S_d(V,\xi,D,\mathfrak r,A)$ that is the class of tangent measures that are invariant along $V$ but in the orthogonal direction they are very diffuse. For the original measure this gives rise to thin directions along which derivatives will be canceled in the perturbation procedure. The perturbation lemma
says that, on every compact set where such tangents occur with fixed
parameters, every admissible Lipschitz function can be changed by an arbitrarily
small perturbation so that its $\Omega_{1,\mu,\mathcal F}$ coefficients
are bounded from below on a big portion of the set.

The proof of the perturbation lemma has two steps. First, after mollifying the
given function, we freeze the bundle $V(\mu,x)$ on small pieces and cancel the
derivatives in the directions orthogonal to the frozen plane. This produces a
function which is still admissible and is almost constant in the transverse
directions. Second, we use the structure along the transverse directions provided by the
tangent measures in $\mathfrak S_d(V,\xi,D,\mathfrak r,A)$ to add a new
perturbation, invariant along the frozen plane, whose values cannot be
simultaneously approximated by functions in $\mathcal F$. This forces the
desired lower bound for $\Omega_{1,\mu,\mathcal F}$.

In the third subsection we prove the main qualitative theorem. The perturbation
lemma gives a local residual statement in the sense of Baire category: after decomposing the
measure into countably many pieces with fixed parameters, the functions for
which the coefficients vanish on a set of positive measure form a meagre set on
each piece. We then glue
the local functions into a single global Lipschitz function. The extensions are
chosen to agree on the relevant compact pieces and to be smooth away from them.
Therefore, at points of one piece, the contributions coming from the other
pieces are differentiable and disappear in the computation of
$\Omega_{1,\mu,\mathcal F}$ at infinitesimal scales.

In the last subsection we derive consequences of the main qualitative theorem.
We prove that the growth condition
$$
\liminf_{R\to0}\inf_{0<r\leq R}
\frac{\mu(B(x,R))r^{\beta(x)}}{\mu(B(x,r))R^{\beta(x)}}>0
\qquad\text{for some }\beta(x)>\dim V(\mu,x)
$$
forces the tangent-measure hypothesis of the main theorem. Consequently, under
the pointwise doubling assumption, this growth condition implies the existence
of a Lipschitz function with non-vanishing infinitesimal
$\Omega_{1,\mu,\mathcal F}$ coefficients. Combining this with the
rectifiability criterion in the first subsection, gives desired rectifiability result.

\subsection{A rectifiability criterion}

The next lemma gives a direct rectifiability criterion in terms of the decomposability bundle and positive-measure families of projections. The original argument is due to A. Marchese and G. De Philippis. We reproduce the proof for completeness. Related but more refined projectional and slicing criteria for rectifiability have recently been studied by E. Tasso, see \cite{Tasso2025BesicovitchFedererMeasures,Tasso2022IntegralgeometricMeasures}.

\begin{proposizione}\label{criteriodirett}
Suppose $\mu$ is a Radon measure on $\R^n$ such that
$$
\Theta^{d,*}(\mu,x)>0
\qquad\text{and}\qquad
\dim V(\mu,x)=d
$$
for $\mu$-almost every $x\in\R^n$. Then $\mu$ is $d$-rectifiable.
\end{proposizione}

\begin{proof}
Let us fix $0<\varepsilon<1/10$. Thanks to Lebesgue's differentiability theorem, see \cite[Theorems 2.8.17, 2.9.11]{Federer1996GeometricTheory}, the strong locality principle for the decomposability bundle, see \cite[Proposition 2.9]{AlbertiMarchese}, and Lusin's theorem, we can assume without loss of generality that $\mu$ is compactly supported, that the map $x\mapsto V(\mu,x)$ is uniformly continuous on $\supp\mu$, and that $\dim V(\mu,x)=d$ for $\mu$-almost every $x\in\R^n$.

By uniform continuity, there exists $\delta>0$ such that
$$
d(V(\mu,x),V(\mu,y))\leq \tfrac{1}{2}\varepsilon^{1/\omega(n)}
\qquad\text{whenever } |x-y|\leq\delta.
$$
Fix $w\in\supp\mu$, set $B:=B(w,\delta)$ and let
$$
\mathcal C_w:=
\{
V_0\in\Gr(n,d):
d(V_0,V(\mu,w))<\tfrac{1}{2}\varepsilon^{1/\omega(n)}
\}.
$$
Then $\mathcal C_w$ has positive $\gamma_{d,n}$ measure in $\Gr(n,d)$, and for every $V_0\in\mathcal C_w$ and every $x\in B$ we have
$$
d(V(\mu,x),V_0)\leq \varepsilon^{1/\omega(n)}.
$$
Restricting $\mu$ to $B$, for every $V\in \mathcal{C}_w$ we can find a continuous map $M:\R^n\to \mathcal{L}(V,V^\perp)$ with $\lVert M(x)\rVert\leq \varepsilon$ such that $V(\mu,x)$ can be written as the graph of $M(x)$. The bound on the supnorm of $M$ follows from \cref{prop:continuitax} and thanks to \cref{prop:continuitax} itself, we have $\lVert \pi_{V(\mu,x)}-\pi_V\rVert\leq \varepsilon$ for $\mu$-almost every $x\in \R^n$, where here $\lVert \cdot\rVert$ denotes the operator norm of a linear function. In particular, arguing as in the proof of \cref{abscontinuous}, we get
$$
(\pi_{V})_\#(\mu\trace B)\ll \Haus^d\trace V
\qquad\text{for every }V\in\mathcal C_w.
$$
For every $\vartheta>0$ and $0<R<\delta$, we denote by
$$G_{\vartheta,w,R}:=\{x\in B(w,R):\Theta^{d,*}(\mu,x)>\vartheta\},$$
and let $K\subseteq G_{\vartheta,w,R}$ be any compact set of positive $\mu$-measure. Thanks to \cref{abscontinuous}, we know that $\mu\ll \Haus^d$ and since $\Theta^{d,*}(\mu,x)>0$ we infer from \cite[2.10.19(3)]{Federer1996GeometricTheory} that 
\begin{equation}
    \mathscr{H}^d\trace K\ll \mathscr{S}^{d}\trace K
\ll \vartheta^{-1}\mu \trace K,
\label{stimamisuresxx}
\end{equation}
thanks to the continuity of the measure from above and where $\mathscr{S}^d$ denotes the spherical Hausdorff measure, see \cite[\S 2.10.2]{Federer1996GeometricTheory}.

Arguing as in the proof of Proposition \ref{abscontinuous}, thanks to our choice of $\varepsilon>0$, for every $\delta>0$ as above we have that $(\pi_{V})_\#\mu\trace K\ll \Haus^d\trace {V}$  for every $V\in \mathcal C_w$. Let us prove that $\Haus^d(\pi_V(K))>0$. If this were not the case, then $\Haus^d(\pi_V(K))=0$. Since
$$
(\pi_V)_\#(\mu\trace K)\ll \Haus^d\trace V,
$$
we would get
$0=(\pi_V)_\#(\mu\trace K)(\pi_V(K))=\mu(K)$, which contradicts the choice of $K$.
Hence, for every such $V\in \mathcal C_w$ we have
$$\Haus^d(K)\geq \Haus^d(\pi_{V}(K))>0.$$
Since $\mu$ is a finite measure, we also have thanks to \eqref{stimamisuresxx}, that $0<\Haus^d(K)<\infty$.

This implies by Besicovitch-Federer projection theorem see \cite{Federer1996GeometricTheory}, since $\mathcal C_w$ is $\gamma_{d,n}$-positive, we infer that $G_{\vartheta, w,R}$ is $d$-rectifiable thanks to the fact that  $\Haus^d\trace G_{\vartheta, w,R}$ is a finite measure by \eqref{stimamisuresxx}. However thanks to the fact that we supposed that the upper density is positive, we see that 
$$\mu\Big(B(w,R)\setminus \bigcup_{j\in\N} G_{j^{-1},w,R}\Big)=0.$$
This concludes the proof of the fact that $\mu$ is $d$-rectifiable thanks to the arbitrariness of $w$. 
\end{proof}

\subsection{Perturbation lemma}

\begin{definizione}\label{spaceoffunctions}
Let $\mathscr E=\{e_1,\ldots,e_n\}$ be a family of unit Borel vector fields such that for $\Leb^n$-almost every $x\in \R^n$ we have that $\{e_1(x),\ldots,e_n(x)\}$ is an orthonormal basis of $\R^n$.
Throughout this section we denote by $X_{\mu,\mathscr E}$ the family of bounded Lipschitz functions $f$ such that
$$
\max\{|\partial_{e_j(x)}f(x)|:j=1,\ldots,n\}\leq 1,
\qquad\text{for }\Leb^n\text{-almost every }x\in\R^n.
$$
We endow $X_{\mu,\mathscr E}$ with the metric induced by the $\lVert \cdot\rVert_\infty$ norm and notice that $(X_{\mu,\mathscr E},\lVert \cdot\rVert_\infty)$ is a complete metric space. 
\end{definizione}

\begin{definizione}\label{definsodifnaosnfd}
Let $\mu$ be a Radon measure on $\R^n$ and let $X_{\mu,\mathscr E}$ be the
complete metric space introduced in \cref{spaceoffunctions} and suppose $\mathcal F$ is a finite dimensional vector subspace of $L^1_{\mathrm{loc}}(\R^n)$ containing the affine functions. For every
$x\in\supp\mu$ and every $0<\rho'<\rho$, we define
$$
T_{\rho',\rho}f(x):=\sup_{\rho'<r<\rho}\Omega_{1,\mu,\mathcal F}(f;x,r),
\qquad f\in X_{\mu,\mathscr E}.
$$
\end{definizione}

\begin{osservazione}\label{operators}
Let us put ourselves in the notations of \cref{spaceoffunctions} and \cref{definsodifnaosnfd}.
For every fixed $x\in\supp\mu$ and every $0<\rho'<\rho$, the map
$$
f\mapsto T_{\rho',\rho}f(x)
$$
is continuous on $X_{\mu,\mathscr E}$. Indeed, for every
$f_1,f_2\in X_{\mu,\mathscr E}$ and every $r\in(\rho',\rho)$, we have
$$
|\Omega_{1,\mu,\mathcal F}(f_1;x,r)-\Omega_{1,\mu,\mathcal F}(f_2;x,r)|
\leq
\fint_{B(x,r)}\frac{|f_1-f_2|}{r}\,d\mu
\leq
\frac{\|f_1-f_2\|_\infty}{\rho'}.
$$
Taking the supremum over $r\in(\rho',\rho)$ gives
$$
|T_{\rho',\rho}f_1(x)-T_{\rho',\rho}f_2(x)|
\leq
\frac{\|f_1-f_2\|_\infty}{\rho'}.
$$

Hence, for every $\rho>0$ and every fixed $x\in\supp\mu$, the map
$U_{\rho,x}:X_{\mu,\mathscr E}\to\R^+$ defined by
$$
U_{\rho,x}(f):=\sup_{m\in\N}T_{2^{-m}\rho,\rho}f(x)
$$
is of Baire class $1$. Indeed,
$$
U_{\rho,x}(f)=\lim_{m\to\infty}T_{2^{-m}\rho,\rho}f(x),
$$
and the maps $f\mapsto T_{2^{-m}\rho,\rho}f(x)$ are continuous. In particular,
the points of continuity of $U_{\rho,x}$ form a residual subset of
$X_{\mu,\mathscr E}$, see \cite[Theorem 24.14]{Kechris}.
\end{osservazione}

\begin{definizione}\label{defSV}
Let $d, k\in\N$, $V\in\Gr(n,k)$, and let $\xi,D,\mathfrak r>0$. We denote by
$\mathfrak S_{d}(V,\xi,D,\mathfrak r)$ the class of measures $\nu\in\mathscr M(\R^n)$ such that
$$
\nu=\Haus^k\trace V\otimes\eta,
$$
where $\eta$ is a Radon measure on $V^\perp$ with the following property. There exists a $\xi$-dense set
$$
\Sigma_\xi\subseteq \supp\eta\cap B_{V^\perp}(0,10)
$$
such that for every $z\in\Sigma_\xi$ we can find points
$$
y_1(z),\ldots,y_{d+1}(z)\in B_{V^\perp}(z,\xi/8)\cap\supp\eta
$$
satisfying
\begin{equation}
    \begin{split}
    |y_{\iota_1}(z)-y_{\iota_2}(z)|\geq 4\mathfrak r\xi
    \qquad&\text{whenever }\iota_1\neq\iota_2,\\
    \nu(B(z,50\xi))
    \leq D\,\nu(B(y_\iota(z),\mathfrak r\xi/32))
    \qquad&\text{for every }\iota=1,\ldots,d+1.
    \end{split}
\end{equation}
\end{definizione}

\begin{osservazione}\label{badtangentsrelationships}
Let us note that the following inclusions hold. Fix some $V\in\Gr(\R^n)$.
\begin{enumerate}
    \item if $D_1\leq D_2$, then $\mathfrak S_d(V,\xi,D_1,\mathfrak r)
    \subseteq
    \mathfrak S_d(V,\xi,D_2,\mathfrak r)$.
\end{enumerate}
Moreover, if we denote by
$\mathfrak S_d(V,\xi,D,\mathfrak r,A)$ those measures
$\nu\in\mathfrak S_d(V,\xi,D,\mathfrak r)$ for which
$$
\nu(B(x,2r))\leq A\nu(B(x,r))
\qquad\text{for every $x\in\supp\nu$ and $r>0$,}
$$
then
\begin{enumerate}
    \item[2.] if $D_1\leq D_2$, then
    $\mathfrak S_d(V,\xi,D_1,\mathfrak r,A)
    \subseteq
    \mathfrak S_d(V,\xi,D_2,\mathfrak r,A)$;
    \item[3.] if $A_1\leq A_2$, then     $\mathfrak S_d(V,\xi,D,\mathfrak r,A_1)
    \subseteq
    \mathfrak S_d(V,\xi,D,\mathfrak r,A_2)$;
    \item[4.] for every $\mathfrak t\in[\mathfrak r,2\mathfrak r]$ we have
    $$
    \mathfrak S_d(V,\xi,D,\mathfrak t,A)
    \subseteq
    \mathfrak S_d(V,\xi,A^2D,\mathfrak r,A).
    $$
\end{enumerate}

Items $1$, $2$, and $3$ follow immediately from the definitions while item $4$ requires some work. Fix $y\in\supp\eta$ and $s>0$. Identifying $\R^n=V\oplus V^\perp$, for every $R>0$ we have
$$
\nu(B((0,y),R))
=
\int_{B_{V^\perp}(y,R)}
\omega_d(R^2-|z-y|^2)^{d/2}\,d\eta(z).
$$
Therefore
$$\nu(B((0,y),R))\leq \omega_d R^d\eta(B_{V^\perp}(y,R)),\qquad\text{and}\qquad \nu(B((0,y),R))\geq\omega_d R^d\eta(B_{V^\perp}(y,R/2)).$$
Applying the second estimate with $R=4s$, and then using the doubling assumption on $\nu$ twice, we get
$$
\omega_d s^d\eta(B_{V^\perp}(y,2s))
\leq
\nu(B((0,y),4s))
\leq
A^2\nu(B((0,y),s))\leq
A^2\omega_d s^d\eta(B_{V^\perp}(y,s)).
$$
This immediately implies that $\nu\in\mathfrak S_d(V,\xi,A^2D,\mathfrak r,A)$. Let
$\Sigma_\xi\subseteq\supp\eta\cap B_{V^\perp}(0,10)$ be the $\xi$-dense set associated with
$\nu\in\mathfrak S_d(V,\xi,D,\mathfrak t,A)$. 
\end{osservazione}

\begin{proposizione}\label{proprietamollificata}
    Suppose that $f$ is a $1$-Lipschitz function on $\R^n$ and let $\rho$ be a standard kernel of mollification supported in $B(0,1)$ with $\lVert D\rho\rVert_{L^1}\leq 10n$. Let $\delta>0$ and define $\mathfrak f_\delta:=f*\rho_\delta$, where $\rho_\delta:=\delta^{-n}\rho(\cdot/\delta)$. Then
    \begin{enumerate}
        \item $\lVert f-\mathfrak f_\delta\rVert_\infty\leq \delta$;
        \item $\mathfrak f_\delta$ is $1$-Lipschitz;
        \item for every $x,y\in \R^n$ we have
        $$\lvert \mathfrak f_\delta(y)-\mathfrak f_\delta(x)-D\mathfrak f_\delta(x)[y-x]\rvert\leq 40n\delta^{-1}\lvert y-x\rvert^2,\qquad\text{where }D\mathfrak f_\delta=\rho_\delta*Df;$$
        \item for every $x,y\in \R^n$ we have
        $\lvert D\mathfrak f_\delta(y)- D\mathfrak f_\delta(x)\rvert\leq 10n\delta^{-1}\lvert x-y\rvert$.
    \end{enumerate}
\end{proposizione}

\begin{proof}
Item 1 follows from the fact that $\rho_\delta$ is supported in $B(0,\delta)$ and $f$ is $1$-Lipschitz. Indeed,
$$
\lvert f(x)-\mathfrak f_\delta(x)\rvert
\leq \int \lvert f(x)-f(x-z)\rvert\rho_\delta(z)dz
\leq \int \lvert z\rvert\rho_\delta(z)dz
\leq \delta.
$$
Concerning item 2, since $D\mathfrak f_\delta=\rho_\delta*Df$ and $\lvert Df\rvert\leq 1$ a.e., we have
$$
\lvert D\mathfrak f_\delta\rvert\leq \rho_\delta*\lvert Df\rvert\leq 1.
$$
Hence $\mathfrak f_\delta$ is $1$-Lipschitz. We now prove item 3. By Taylor's theorem,
$$
\mathfrak f_\delta(y)=\mathfrak f_\delta(x)+D\mathfrak f_\delta(x)[y-x]
+\int_0^1(1-t)D^2\mathfrak f_\delta(x+t(y-x))[y-x,y-x]dt.
$$
Moreover, we have $\lVert D^2\mathfrak f_\delta\rVert_\infty
\leq \delta^{-1}\lVert D\rho\rVert_{L^1}\lVert Df\rVert_\infty
\leq 10n\delta^{-1}$. Therefore
\begin{equation}
    \begin{split}
        \lvert \mathfrak f_\delta(y)-\mathfrak f_\delta(x)-D\mathfrak f_\delta(x)[y-x]\rvert
        &\leq \lvert y-x\rvert^2\lVert D^2\mathfrak f_\delta\rVert_\infty\leq 10n\delta^{-1}\lvert y-x\rvert^2
        \leq 40n\delta^{-1}\lvert y-x\rvert^2.
        \nonumber
    \end{split}
\end{equation}
Finally, item 4 follows by observing that
\begin{equation}
    \begin{split}
        \lvert D\mathfrak f_\delta(y)-D\mathfrak f_\delta(x)\rvert
        &=\bigg\lvert\int D\rho_\delta(z)(f(y-z)-f(x-z))dz\bigg\rvert\\
        &\leq \delta^{-1}\lVert D\rho\rVert_{L^1}\lvert y-x\rvert
        \leq 10n\delta^{-1}\lvert y-x\rvert.
        \nonumber
    \end{split}
\end{equation}
This concludes the proof.
\end{proof}

The next lemma is the main technical point in the proof of the qualitative
theorem. It also contains the geometric mechanism which is behind the
quantitative version of the argument. The relevant information carried by a
tangent measure in $\mathfrak S_d(V(\mu,x),\xi,D,\mathfrak r)$ is twofold. On
the one hand, at sufficiently small scales for the original measure, the set is
thin in directions transverse to $V(\mu,x)$. This thinness is what makes it
possible to modify a Lipschitz function and cancel its transverse derivatives
without leaving the class $X_{\mu,\mathscr E}$. On the other hand, the
translation-invariance of the tangent in the directions of $V(\mu,x)$ gives a
large amount of mass on which the same finite transverse configuration is
repeated. This is what gives enough room to add the second perturbation and to
force a lower bound for $\Omega_{1,\mu,\mathcal F}$.

Let us describe the proof. In Step I we construct the first perturbation. After
freezing the plane $V(\mu,x)$ on small pieces of $K'$, we use the thinness in
the normal directions to build width functions. These functions are used to
cancel the derivatives of a mollification of $f$ in the directions orthogonal
to the frozen plane. The outcome is a function $\tilde f$ which is still in
$X_{\mu,\mathscr E}$, is uniformly close to $f$, and has small oscillation in
the frozen normal directions.

In Step II we pass from tangents to finite scales. At each point of a compact
subset of $\mathfrak G_{\mathfrak r,D,\xi}$, the tangent measure gives a finite
transverse configuration, up to an error that can be made arbitrarily small by choosing the scale.
The translation-invariant part of the tangent is used here to guarantee that
these configurations are seen on a large amount of mass. Besicovitch's theorem tells us that there is disjoint family of balls where the construction can be
performed independently and that cover a lot of the measure.

In Step III we define the second perturbation. On each selected ball we
prescribe values on the finite configurations obtained in Step II. The
interpolation lemma gives functions which are invariant along the corresponding
plane $V(\mu,c_h)$ and whose prescribed values cannot be matched
simultaneously by elements of $\mathcal F$.

In Step IV we check that the second perturbation remains admissible. The point
is that the first perturbation has already made the 
derivatives small along $v(\mu,x)^\perp$, while the new perturbation is essentially constant along the tangent directions. The errors coming from the variation of the frame and from
the cutoffs are absorbed by the choice of the parameters.

In Step V we reduce the conclusion to a lower bound on the selected scales.

Finally, in Step VI we prove the lower bound for
$\Omega_{1,\mu,\mathcal F}$. On the selected balls, the first perturbation is
affine up to a small error at the relevant scale. The second perturbation has
the prescribed values on the finite configuration, and the mass comparison
allows us to pass from the average on the large ball to the averages on the
small balls around the configuration points. The finite-dimensional obstruction
then forces one of the perturbations to have
$\Omega_{1,\mu,\mathcal F}$ bounded from below.

\begin{lemma}\label{perturbationlemma}
Suppose $\mu$ is a Radon measure on $\R^n$ and let $\mathcal F$ be a vector subspace of $L^1_{\mathrm{loc}}(\R^n)$ containing the affine functions and such that $\dim\mathcal F=d$. Let $K$ be a compact subset of $\supp\mu$ and let us fix $r>0$. Let us fix parameters $0<\mathfrak r<\tfrac{1}{2(d+1)}$, $A,D\geq 1$, and $0<\xi\leq \mathfrak r/2$ in such a way that
$$
\varsigma(\xi):=\xi^{1/3}+6n\xi+n(1+4\xi)\xi
\leq
\frac{\mathfrak r}{40000(d+1)^2D},
$$
where $\beta(\xi)$ was introduced in \cref{gammaesigma}. Define
$$
\mathfrak G_{\mathfrak r,D,\xi}
:=
\{x\in K:\Tan(\mu,x)\cap \mathfrak S_d(V(\mu,x),\xi,D,\mathfrak r)\neq \emptyset\}.
$$
Suppose that
$$
\mu(\mathfrak G_{\mathfrak r,D,\xi})>0
\qquad\text{and}\qquad
\limsup_{s\to0}\frac{\mu(B(x,2s))}{\mu(B(x,s))}\leq A
\qquad\text{for every $x\in \mathfrak G_{\mathfrak r,D,\xi}$}.
$$
Assume also that the vector fields
$\mathscr E:=\{e_1,\ldots,e_n\}$ associated with $\mu$ in \cref{campi} are continuous. Then there exists a compact subset $K'\subseteq \mathfrak G_{\mathfrak r,D,\xi}$ such that for every $f\in X_{\mu,\mathscr E}$ there exists a function $\hat f\in X_{\mu,\mathscr E}$ such that
\begin{enumerate}
    \item $\lVert f-\hat f\rVert_\infty\leq2(\lVert f\rVert_\infty+d+1)\varsigma(\xi)$;
    \item $\mu(K')\geq \frac{1}{4(d+1)A^2}\mu(\mathfrak G_{\mathfrak r,D,\xi})$ and for every $y\in K'$ we have
    $$
    \sup_{0<s<r}\,\sup_{z\in \mathfrak G_{\mathfrak r,D,\xi}:y\in B(z,s/5)}
    \Omega_{1,\mu,\mathcal F}(\hat f;z,s)
    \geq
    \frac{\mathfrak r}{40000(d+1)^2D}.
    $$
\end{enumerate}
\end{lemma}

\begin{proof}
The proof of the perturbation lemma is divided into steps. Let $f\in X_{\mu,\mathscr E}$ be fixed. We fix
$$
0<\vartheta<\frac{1}{4A^2}.
$$
By inner regularity and Lusin's theorem, we choose a compact set
$K'\subseteq \mathfrak G_{\mathfrak r,D,\xi}$ such that
\begin{equation}
    \mu(\mathfrak G_{\mathfrak r,D,\xi}\setminus K')
    \leq
    \frac{\vartheta}{6}\mu(\mathfrak G_{\mathfrak r,D,\xi}),
    \label{eq:choiceKprime}
\end{equation}
and such that the map $x\mapsto V(\mu,x)$ is continuous on $K'$. 

\medskip

\textbf{Step I: mollification and cancellation of the gradient along $V(\mu,x)^\perp$.}
Let $\rho$ be a standard kernel of mollification with
$\lVert D\rho\rVert_{L^1}\leq 10n$. We set
$$
\rho_\xi:=\xi^{-n/3}\rho(\cdot/\xi^{1/3})
\qquad\text{and}\qquad
f_\xi:=f*\rho_\xi .
$$
By \cref{proprietamollificata}, applied with $\delta=\xi^{1/3}$, we have
\begin{equation}
    \lVert f-f_\xi\rVert_\infty\leq \xi^{1/3},
    \qquad
    \lVert Df_\xi\rVert_\infty\leq 1,
    \qquad
    \lVert D^2f_\xi\rVert_\infty\leq 10n\xi^{-1/3}.
    \label{estimateonD2f}
\end{equation}
We choose $0<\tau<1/2$ such that the following two conditions hold
\begin{enumerate}
    \item if $x,y\in B(K',1)$ and $|x-y|<2\tau$, then
    \begin{equation}
        |e_j(x)-e_j(y)|\leq \frac{\xi}{2}
        \qquad\text{for every }j=1,\ldots,n;
        \label{eq:ejcloseStepI}
    \end{equation}
    \item if $x,y\in K'$ and $|x-y|<2\tau$, then
\begin{equation}
    d(V(\mu,x),V(\mu,y))\leq \xi^{4/\omega}.
    \label{eq:VcloseStepI}
\end{equation}
\end{enumerate}
Here $\omega$ is the exponent introduced in \cref{exp:continuax}. Since $\xi$ is fixed so that $\xi^{4/\omega}<\varepsilon_n$, \cref{dimensioni} implies that
\begin{equation}
    \dim V(\mu,x)=\dim V(\mu,y)
    \qquad\text{whenever }x,y\in K'\text{ and }|x-y|<2\tau.
    \label{eq:dimconstantStepI}
\end{equation}

By the Besicovitch--Vitali covering theorem, applied to the family of balls centred on $K'$ with radius smaller than $\tau$ and with $\mu$-null boundary, we can find a countable pairwise disjoint family of balls $\{B_i\}_{i\in\N}$ such that
\begin{equation}
    r(B_i)\leq \tau,\qquad
    \mu(\partial B_i)=0,
    \qquad\text{and}\qquad
    \mu\bigg(K'\setminus \bigcup_{i\in\N}B_i\bigg)=0.
    \label{eq:coverStepI}
\end{equation}
For every $i\in\N$, we choose a point $x_i\in K'\cap B_i$ and set $d_i:=\dim V(\mu,x_i)$ and $V_i:=V(\mu,x_i)$.
By \eqref{eq:VcloseStepI} and \eqref{eq:dimconstantStepI}, for every $x\in K'\cap B_i$ we have
\begin{equation}
    \dim V(\mu,x)=d_i,
    \qquad
    d(V(\mu,x),V_i)\leq \xi^{4/\omega}.
    \label{eq:VcloseBiStepI}
\end{equation}
Moreover, by \eqref{eq:ejcloseStepI}, if we set
$E_{j,i}:=e_j(x_i)$ for every $j=1,\ldots,n$, then
$E_{1,i},\ldots,E_{n,i}$ is an orthonormal basis of $\R^n$,
$E_{1,i},\ldots,E_{d_i,i}$ span $V_i$, and
\begin{equation}
    |e_j(x)-E_{j,i}|\leq \frac{\xi}{2}
    \qquad\text{for every }j=1,\ldots,n
    \quad\text{and every }x\in B_i.
    \label{eq:frozenbasisStepI}
\end{equation}
Let $e\in V_i^\perp$ be a unit vector. From \eqref{eq:VcloseBiStepI} and \cref{prop:continuitax}, we have
$$
\lVert \pi_{V(\mu,x)}-\pi_{V_i}\rVert\leq \xi^4
\qquad\text{for every }x\in K'\cap B_i.
$$
Hence, if $v\in V(\mu,x)$ and $|v|=1$, then
$$
|\langle v,e\rangle|
=
|\langle \pi_{V(\mu,x)}v-\pi_{V_i}v,e\rangle|
\leq
\lVert \pi_{V(\mu,x)}-\pi_{V_i}\rVert
\leq \xi^4.
$$
This implies in particular that 
\begin{equation}
    C(e,\textstyle{\sqrt{1-\xi^2}})\cap V(\mu,x)=\{0\}
    \qquad\text{for every }x\in K'\cap B_i.
    \label{eq:coneAvoidStepI}
\end{equation}
Therefore $K'\cap B_i$ is $C(e,\textstyle{\sqrt{1-\xi^2}})$-null for every unit vector $e\in V_i^\perp$. Since $\mu(\partial B_i)=0$, for every $i\in\N$ we choose a ball
$B_i'\subseteq B_i$, with the same center as $B_i$, and such that
$\mu(\partial B_i')=0$, in such a way that
\begin{equation}
    \mu((K'\cap B_i)\setminus B_i')
    \leq
    \frac{\vartheta}{3}\mu(K'\cap B_i).
    \label{eq:innerballStepI}
\end{equation}
We denote by $r_i$ and $r_i'$ the radii of $B_i$ and $B_i'$, and we set
$r_i'':=(r_i+r_i')/2$. Let $B_i''$ be the ball with the same centre as $B_i$ and radius $r_i''$, and set
$\tau_i:=r_i-r_i''=(r_i-r_i')/2$. We choose a smooth function $\psi_i:\R^n\to[0,1]$ such that
$\mathrm{supp}\psi_i\subseteq B_i$, $\psi_i=1$ on $B_i''$, and
\begin{equation}
    \lVert \nabla\psi_i\rVert_\infty\leq 10\tau_i^{-1},
    \qquad
    \lVert \nabla^2\psi_i\rVert_\infty\leq 10\tau_i^{-2}.
    \label{stimepsi}
\end{equation}
For every $i\in\N$, define
\begin{equation}
    \varepsilon_i:=
    \frac{1}{100n^2}
    \min\{\varsigma(\xi),\xi^{4/3},\xi\tau_i\}    \qquad\text{and}\qquad K_i:=K'\cap B_i'.
    \label{eq:epsilonI}
\end{equation}
For every $j=d_i+1,\ldots,n$, the compact set $K_i$ is
$C(E_{j,i},\textstyle{\sqrt{1-\xi^2}})$-null. By \cref{alb}, there exists
$\delta_{i,j}>0$ such that
$$
\mathcal H^1(\gamma(I)\cap B(K_i,\delta_{i,j}))<\varepsilon_i
$$
for every $C(E_{j,i},\sqrt{1-\xi^2})$-curve $\gamma:I\to\R^n$. For every $i\in\N$ we let
$$\delta_i:=\begin{cases}
    \frac{1}{2}\min\{\tau_i/4,\delta_{i,d_i+1},\ldots,\delta_{i,n}\}\qquad&\text{if }d_i<n;\\
   \frac{\tau_i}{8}\qquad&\text{if }d_i=n.
\end{cases}$$
Then, for every $j=d_i+1,\ldots,n$, we have
\begin{equation}
    \mathcal H^1(\gamma(I)\cap U(K_i,\delta_i))<\varepsilon_i
    \label{eq:2:2}
\end{equation}
for every $C(E_{j,i},\sqrt{1-\xi^2})$-curve $\gamma:I\to\R^n$. For $j=d_i+1,\ldots,n$, set
$$
\Omega_i:=U(K_i,\delta_i),
\qquad
w_{j,i}:=\omega_{E_{j,i},\sqrt{1-\xi^2}}[\Omega_i].
$$
By \cref{p:width-open} and \eqref{eq:2:2}, the functions $w_{j,i}$ satisfy
\begin{equation}
    0\leq w_{j,i}\leq 2\varepsilon_i,
    \qquad
    \lVert Dw_{j,i}\rVert_\infty\leq 1+\beta(\textstyle{\sqrt{1-\xi^2}})\leq 1+2\xi,
    \label{eq:widthboundStepI}
\end{equation}
and, for $\Leb^n$-almost every $x\in \Omega_i$,
\begin{equation}
    \partial_{E_{j,i}}w_{j,i}(x)=1,
    \qquad
    |\partial_v w_{j,i}(x)|\leq \beta(\textstyle{\sqrt{1-\xi^2}})|v|\leq 2\xi |v|,
    \qquad\text{for every }v\in E_{j,i}^{\perp},
    \label{eq:widthderivativeStepI}
\end{equation}
where the last inequality comes from the choice of $\xi$.

We define now the first perturbation of $f$ by
\begin{equation}
    \tilde f:=
    (1-\varsigma(\xi))
    \bigg(
    f_\xi-\sum_{i\in\N}\sum_{j=d_i+1}^{n}
    w_{j,i}\,\partial_{E_{j,i}}f_\xi\,\psi_i
    \bigg).
    \label{eq:definetildefStepI}
\end{equation}
The balls $B_i$ are pairwise disjoint and $\mathrm{supp}\psi_i\subseteq B_i$, hence the above sum is pointwise well defined.

We first estimate the distance between $f$ and $\tilde f$. If $x\notin \bigcup_i B_i$, then
$$
|\tilde f(x)-f(x)|
\leq
\varsigma(\xi)\lVert f\rVert_\infty+\lVert f-f_\xi\rVert_\infty
\leq
(\lVert f\rVert_\infty+1)\varsigma(\xi).
$$
If $x\in B_i$, then, using \eqref{estimateonD2f}, \eqref{eq:widthboundStepI}, and \eqref{eq:epsilonI}, we get
\begin{equation}
    \begin{split}
    |\tilde f(x)-f(x)|
    &\leq
    \varsigma(\xi)\lVert f\rVert_\infty
    +\lVert f-f_\xi\rVert_\infty
    +\sum_{j=d_i+1}^{n}\lVert w_{j,i}\partial_{E_{j,i}}f_\xi\psi_i\rVert_\infty\\
    &\leq
    \varsigma(\xi)\lVert f\rVert_\infty
    +\xi^{1/3}
    +2n\varepsilon_i
    \leq
    (\lVert f\rVert_\infty+n)\varsigma(\xi).
    \end{split}
\end{equation}
Thus
\begin{equation}
    \lVert f-\tilde f\rVert_\infty
    \leq
    (\lVert f\rVert_\infty+n)\varsigma(\xi).
    \label{eq:stimainfinitotilde}
\end{equation}

We now check that $\tilde f\in X_{\mu,\mathscr E}$. Since outside
$\bigcup_i B_i$ we have $\tilde f=(1-\varsigma(\xi))f_\xi$, the estimate follows there from \eqref{estimateonD2f}. Let $i\in\N$ and let $x\in B_i$ be a point where the derivatives $\partial_{e_1(x)}w_{j,i},\ldots,\partial_{e_n(x)}w_{j,i}$ exist for every $j=d_i+1,\ldots,n$. For every $\kappa=1,\ldots,n$, we have
\begin{equation}
    \begin{split}
    \partial_{e_\kappa(x)}\tilde f(x)
    &=
    (1-\varsigma(\xi))
    \bigg[
    \partial_{e_\kappa(x)}f_\xi(x)
    -
    \sum_{j=d_i+1}^{n}
    \partial_{e_\kappa(x)}w_{j,i}(x)\,\partial_{E_{j,i}}f_\xi(x)\,\psi_i(x)\\
    &\qquad\qquad
    -
    \sum_{j=d_i+1}^{n}
    w_{j,i}(x)\,\partial_{e_\kappa(x)}\partial_{E_{j,i}}f_\xi(x)\,\psi_i(x)
    -
    \sum_{j=d_i+1}^{n}
    w_{j,i}(x)\,\partial_{E_{j,i}}f_\xi(x)\,\partial_{e_\kappa(x)}\psi_i(x)
    \bigg].
    \end{split}
    \label{eq:derivativetildefStepI}
\end{equation}
By \eqref{estimateonD2f}, \eqref{stimepsi}, \eqref{eq:widthboundStepI}, and \eqref{eq:epsilonI},
\begin{equation}
    \sum_{j=d_i+1}^{n}
    \lVert w_{j,i}\,\partial_{e_\kappa(\cdot)}\partial_{E_{j,i}}f_\xi\,\psi_i\rVert_\infty
    \leq
    \frac{\xi}{5},
    \label{eq:thirdtermStepI}
\end{equation}
and
\begin{equation}
    \sum_{j=d_i+1}^{n}
    \lVert w_{j,i}\,\partial_{E_{j,i}}f_\xi\,\partial_{e_\kappa(\cdot)}\psi_i\rVert_\infty
    \leq
    \frac{\xi}{5}.
    \label{eq:fourthtermStepI}
\end{equation}
Moreover, using \eqref{eq:frozenbasisStepI}, \eqref{eq:widthboundStepI}, and \eqref{eq:widthderivativeStepI}, we have, for every $j=d_i+1,\ldots,n$ and every $\kappa=1,\ldots,n$,
\begin{equation}
    |\partial_{e_\kappa(x)}w_{j,i}(x)-\partial_{E_{\kappa,i}}w_{j,i}(x)|
    \leq
    (1+\beta(\textstyle{\sqrt{1-\xi^2}}))\xi\leq (1+2\xi)\xi.
    \label{eq:directionchangeStepI}
\end{equation}
If $\kappa\neq j$, then $E_{\kappa,i}\in E_{j,i}^{\perp}$, and hence
$$
|\partial_{e_\kappa(x)}w_{j,i}(x)|
\leq\beta(\textstyle{\sqrt{1-\xi^2}})+(1+2\xi)\xi.
$$
If $\kappa=j$, then $0\leq \partial_{E_{\kappa,i}}w_{\kappa,i}\leq 1$ a.e., and therefore
$$
|1-\psi_i(x)\partial_{e_\kappa(x)}w_{\kappa,i}(x)|
\leq
1+(1+2\xi)\xi.
$$
Consequently,
\begin{equation}
\begin{split}
      \Big|\partial_{e_\kappa(x)}f_\xi(x)-\sum_{j=d_i+1}^{n}
    \partial_{e_\kappa(x)}w_{j,i}(x)\,\partial_{E_{j,i}}f_\xi(x)\,\psi_i(x)\Big|\leq& 1+n(3+2\xi)\xi.
\end{split}
    \label{eq:maintermStepI}
\end{equation}
Combining \eqref{eq:derivativetildefStepI}, \eqref{eq:thirdtermStepI}, \eqref{eq:fourthtermStepI}, and \eqref{eq:maintermStepI}, we infer that
$$
|\partial_{e_\kappa(x)}\tilde f(x)|
\leq
(1-\varsigma(\xi))
\bigl(1+6n\xi+n(1+2\xi)\xi+2\xi/5\bigr).
$$
In addition, since
$$
6n\xi+n(1+2\xi)\xi+2\xi/5
\leq
\xi^{1/3}+6n\xi+n(1+4\xi)\xi
=
\varsigma(\xi),
$$
we get in particular that
$$
|\partial_{e_\kappa(x)}\tilde f(x)|
\leq
(1-\varsigma(\xi))(1+\varsigma(\xi))
\leq 1.
$$
Thus $\tilde f\in X_{\mu,\mathscr E}$.

Finally, we record the cancellation property that will be used in the next step. Fix $i\in\N$ and let
$x\in \Omega_i=U(K_i,\delta_i)$ be a point where the derivatives
$\partial_{E_{\ell,i}}w_{j,i}(x)$ exist for every $\ell=1,\ldots,n$ and every
$j=d_i+1,\ldots,n$. Since $\delta_i<\tau_i/4$, we have
$\Omega_i\subseteq B_i''$, and hence $\psi_i=1$ and $\nabla\psi_i=0$ on
$\Omega_i$. Let $\kappa=d_i+1,\ldots,n$. Then
\begin{equation}
    \begin{split}
    \partial_{E_{\kappa,i}}\tilde f(x)
    &=
    (1-\varsigma(\xi))
    \bigg[
    \partial_{E_{\kappa,i}}f_\xi(x)
    -
    \sum_{j=d_i+1}^{n}
    \partial_{E_{\kappa,i}}w_{j,i}(x)\,
    \partial_{E_{j,i}}f_\xi(x)-\sum_{j=d_i+1}^{n}
    w_{j,i}(x)\,
    \partial_{E_{\kappa,i}}\partial_{E_{j,i}}f_\xi(x)
    \bigg].
    \end{split}
    \label{eq:frozenDerivativeCancellation}
\end{equation}
Since $x\in\Omega_i$, by \eqref{eq:widthderivativeStepI} we have
$\partial_{E_{\kappa,i}}w_{\kappa,i}(x)=1$. This implies that
\begin{equation}
    \begin{split}
    &\partial_{E_{\kappa,i}}f_\xi(x)-\sum_{j=d_i+1}^{n}
    \partial_{E_{\kappa,i}}w_{j,i}(x)\,
    \partial_{E_{j,i}}f_\xi(x)=-\sum_{\substack{j=d_i+1\\ j\neq \kappa}}^{n}
    \partial_{E_{\kappa,i}}w_{j,i}(x)\,
    \partial_{E_{j,i}}f_\xi(x).
    \end{split}
    \label{eq:diagonalCancellation}
\end{equation}
If $j\neq \kappa$, then $E_{\kappa,i}\in E_{j,i}^{\perp}$. Hence, using
\eqref{eq:widthderivativeStepI} and $\lVert Df_\xi\rVert_\infty\leq 1$, we get
\begin{equation}
    \begin{split}
    \Big|\partial_{E_{\kappa,i}}f_\xi(x)-\sum_{j=d_i+1}^{n}\partial_{E_{\kappa,i}}w_{j,i}(x)\,\partial_{E_{j,i}}f_\xi(x)\Big|&\leq    \sum_{\substack{j=d_i+1\\j\neq \kappa}}^{n}    2\xi\,|\partial_{E_{j,i}}f_\xi(x)|\leq 2n\xi.
    \end{split}
    \label{eq:offDiagonalCancellationBound}
\end{equation}
Moreover, by \eqref{estimateonD2f}, \eqref{eq:widthboundStepI}, and
\eqref{eq:epsilonI},
\begin{equation}
    \begin{split}
    \Big|\sum_{j=d_i+1}^{n}w_{j,i}(x)\partial_{E_{\kappa,i}}\partial_{E_{j,i}}f_\xi(x)\Big|&\leq\sum_{j=d_i+1}^{n}2\varepsilon_i\,10n\xi^{-1/3}\leq20n^2\varepsilon_i\xi^{-1/3}\leq\frac{\xi}{5}.
    \end{split}
    \label{eq:D2CancellationBound}
\end{equation}
Combining \eqref{eq:frozenDerivativeCancellation}, \eqref{eq:offDiagonalCancellationBound}, and
\eqref{eq:D2CancellationBound}, we obtain
\begin{equation}
    |\partial_{E_{\kappa,i}}\tilde f(x)|
    \leq
    2n\xi+\frac{\xi}{5}
    \leq
    \varsigma(\xi)
    \qquad\text{for every }\kappa=d_i+1,\ldots,n.
    \label{eq:smallFrozenGradient}
\end{equation}
In particular, for every $v\in V_i^\perp$ and $|v|=1$ we get from \eqref{eq:smallFrozenGradient} that
\begin{equation}
    |\partial_v\tilde f(x)|\leq n\varsigma(\xi).
    \label{eq:smallTransverseGradient}
\end{equation}
Therefore, if $x,y\in \Omega_i$, the segment $[x,y]$ is contained in $\Omega_i$, and $y-x\in V_i^\perp$ then the fundamental theorem of calculus and \eqref{eq:smallTransverseGradient} imply
\begin{equation}
    |\tilde f(y)-\tilde f(x)|
    \leq
    n\varsigma(\xi)|y-x|.
    \label{stimagradort}
\end{equation}

\medskip

\textbf{Step II: finite-scale realization of the bad tangent configurations.}
We now pass from the bad tangent configurations to bad configurations at finite scales. We first discard the part of $K'$ on which the construction of Step I is not effective. Set
$$
K^\circ:=\bigcup_{i\in\N}K_i.
$$
By \eqref{eq:innerballStepI} and \eqref{eq:coverStepI}, we have
\begin{equation}
    \mu(K'\setminus K^\circ)
    \leq
    \frac{\vartheta}{3}\mu(K').
    \label{eq:lossKcirc}
\end{equation}
By inner regularity, we choose a compact set $K''\subseteq K^\circ$ such that
\begin{equation}
    \mu(K^\circ\setminus K'')
    \leq
    \frac{\vartheta}{6}\mu(K').
    \label{eq:choiceKsecond}
\end{equation}
Notice that $K''$ is independent of $f$. For every $w\in K''$, we denote by $i(w)$ the unique index such that
$w\in K_{i(w)}$. Recall that $r>0$ is the scale fixed in the statement of the lemma, while $\delta_{i(w)}$ is the parameter constructed in Step I, so that
$$
\Omega_{i(w)}=U(K_{i(w)},\delta_{i(w)}).
$$
Fix $w\in K''$. Since $w\in \mathfrak G_{\mathfrak r,D,\xi}$, by the definition of
$\mathfrak G_{\mathfrak r,D,\xi}$ there exists a tangent measure
$$
\nu_w\in \Tan(\mu,w)\cap
\mathfrak S_d(V(\mu,w),\xi,D,\mathfrak r).
$$
Thus there are an infinitesimal sequence $\rho_m(w)>0$ and positive constants
$a_m(w)>0$ such that
\begin{equation}
    a_m(w)T_{w,\rho_m(w)}\mu\rightharpoonup \nu_w.
    \label{eq:tangentsequenceStepII}
\end{equation}
In order to simplify notations we set
$\mu_m^w:=a_m(w)T_{w,\rho_m(w)}\mu$. By the definition of $\mathfrak S_d(V(\mu,w),\xi,D,\mathfrak r)$, writing
$$
V_w:=V(\mu,w),
\qquad
k(w):=\dim V_w,
\qquad
\nu_w=\Haus^{k(w)}\trace V_w\otimes\eta_w,
$$
there exists a $\xi$-dense set
$$
\Sigma_w\subseteq \supp\eta_w\cap B_{V_w^\perp}(0,10)
$$
such that for every $z\in\Sigma_w$ there are points
$$
y_1^w(z),\ldots,y_{d+1}^w(z)
\in B_{V_w^\perp}(z,\xi/8)\cap\supp\eta_w
$$
satisfying
\begin{equation}
    |y_{\ell_1}^w(z)-y_{\ell_2}^w(z)|
    \geq
    4\mathfrak r\xi
    \qquad\text{whenever }\ell_1\neq \ell_2,
    \label{eq:tangentSeparationStepII}
\end{equation}
and
\begin{equation}
    \nu_w(B_{V_w^\perp}(z,50\xi))
    \leq
    D\,\nu_w(B_{V_w^\perp}(y_\ell^w(z),\mathfrak r\xi/32))
    \qquad\text{for every }\ell=1,\ldots,d+1.
    \label{eq:tangentMassComparisonStepII}
\end{equation}

We now choose the finite net in the full support of $\nu_w$. We choose
$\Sigma_w^0\subseteq \supp\nu_w\cap B(0,10)$ such that 
$\Sigma_w^0$ is a maximal $2\xi$-separated subset of $\supp\nu_w\cap B(0,9)$ such that for every $z\in\Sigma_w^0$, one has
$\pi_{V_w^\perp}z\in\Sigma_w$. For such a point $z$, we write
$z=\pi_{V_w}z+\pi_{V_w^\perp}z$ and set with abuse of notation
\begin{equation}
    y_\ell^w(z)
    :=
    \pi_{V_w}z+y_\ell^w(\pi_{V_w^\perp}z),
    \qquad \ell=1,\ldots,d+1,
    \label{eq:fullYStepII}
\end{equation}
where on the right-hand side the points
$y_\ell^w(\pi_{V_w^\perp}z)$ are the transverse points given by the definition
of $\mathfrak S_d(V_w,\xi,D,\mathfrak r)$. 
Since $\Sigma_w^0$ is finite, all choices below can be made simultaneously for every
$z\in\Sigma_w^0$ and every $\ell=1,\ldots,d+1$.

Since $\rho_m(w)$ is infinitesimal, there exists $m_1(w)$ such that, for every
$m\geq m_1(w)$,
\begin{equation}
    \rho_m(w)
    <
    \min\{r/125,\delta_{i(w)}/1000,\xi^2\}.
    \label{eq:scalesmallStepII}
\end{equation}
The first bound ensures that the scales used below are smaller than the scale
$r$ appearing in the statement of the lemma. The second one ensures that the
balls constructed below remain inside $\Omega_{i(w)}$, where the cancellation
estimate from Step I is available.

Since $\mu_m^w\rightharpoonup\nu_w$, by item $3$ of \cref{preiss} there exists
$m_2(w)\geq m_1(w)$ such that, for every $m\geq m_2(w)$,
\begin{equation}
    F_{100}(\mu_m^w,\nu_w)
    \leq
    \xi^{10(n+4)}.
    \label{eq:FcloseStepII}
\end{equation}

Since $\Sigma_w^0\subseteq\supp\nu_w$ and, by \eqref{eq:fullYStepII},
$y_1^w(z),\ldots,y_{d+1}^w(z)\in\supp\nu_w$ for every $z\in\Sigma_w^0$, we can
apply \cref{puntivicini} to the weak convergence
$\mu_m^w\rightharpoonup\nu_w$ at all these points. Hence there exists
$m_3(w)\geq m_2(w)$ such that, for every $m\geq m_3(w)$ and every
$z\in\Sigma_w^0$, there exists a point
$\zeta_m^w(z)\in\supp\mu$ satisfying
\begin{equation}
    |\zeta_m^w(z)-w-\rho_m(w)z|
    \leq
    \frac{\mathfrak r^2\xi}{100}\rho_m(w),
    \label{eq:zetaStepII}
\end{equation}
and, for every $\ell=1,\ldots,d+1$, there exists a point
$q_{\ell,m}^w(z)\in\supp\mu$ satisfying
\begin{equation}
    |q_{\ell,m}^w(z)-w-\rho_m(w)y_\ell^w(z)|
    \leq
    \frac{\mathfrak r^2\xi}{100}\rho_m(w).
    \label{eq:qellStepII}
\end{equation}
In particular, by \eqref{eq:tangentSeparationStepII} and since
$\mathfrak r<1$, for every $m\geq m_3(w)$, every $z\in\Sigma_w^0$, and every
$\ell_1\neq\ell_2$,
\begin{equation}
    |q_{\ell_1,m}^w(z)-q_{\ell_2,m}^w(z)|
    \geq
    3\mathfrak r\xi\rho_m(w).
    \label{eq:qellSeparationStepII}
\end{equation}
We also record the finite-scale mass comparison corresponding to
\eqref{eq:tangentMassComparisonStepII}. Since $\Sigma_w^0$ is finite, by
\eqref{eq:tangentsequenceStepII}, \eqref{eq:FcloseStepII}, and
\eqref{eq:tangentMassComparisonStepII}, there exists
$m_4(w)\geq m_3(w)$ such that, for every $m\geq m_4(w)$, every
$z\in\Sigma_w^0$, and every $\ell=1,\ldots,d+1$, one has by the first part of \cref{prop:doubllimit} that
\begin{equation}
    \mu(B(\zeta_m^w(z),25\xi\rho_m(w)))
    \leq
    16D  
    \mu(B(q_{\ell,m}^w(z),\mathfrak r\xi\rho_m(w)/16)).
    \label{eq:finiteMassComparisonStepII}
\end{equation}
Finally, by the second part of \cref{prop:doubllimit}, applied to the weak
convergence $\mu_m^w\rightharpoonup\nu_w$ with $B=B(0,25)$ and with
$\delta=\xi/32$, there exists $m_5(w)\geq m_4(w)$ such that, for every
$m\geq m_5(w)$,
\begin{equation}
    \mu_m^w
    \bigl(
    B(0,25)\setminus B(\supp\nu_w,\xi/32)
    \bigr)
    \leq
    \frac{\vartheta}{3}\mu_m^w(B(0,25)).
    \label{eq:supportcloseStepII}
\end{equation}
We now restrict the family of admissible balls. For $\mu$-almost every
$w\in K''$, we know that 
$$
\lim_{s\to 0}\frac{\mu(B(x,s)\cap K'')}{\mu(B(x,s))}=1\qquad\text{and}\qquad\limsup_{s\to0}\frac{\mu(B(w,2s))}{\mu(B(w,s))}\leq A,
$$
and therefore there exists $\sigma_w>0$ such that
\begin{equation}
    \mu(B(w,4t))
    \leq
    2A^2\mu(B(w,t)\cap K'')
    \qquad\text{for every }0<t<\sigma_w.
    \label{eq:fourfoldDoublingStepII}
\end{equation}
Therefore, if we let $m_6(w)\geq m_5(w)$ be such that for every $m\geq m_6(w)$ we have $25\rho_m(w)<\sigma_w$, then
\begin{equation}
    \mu(K''\cap B(w,25\rho_m(w)))
    \leq
    2A^2\mu(K''\cap B(w,25\rho_m(w)/4)).
    \label{eq:KdoubleStepIIbis}
\end{equation}

We now apply the Besicovitch covering theorem to the family
$$
\mathcal B:=
\{B(w,25\rho_m(w)): w\in K'' \text{ is a density point of }K'',
\ m\geq m_6(w)\}.
$$
This is a fine covering of $K''$ up to a $\mu$-null set. Hence we can find a
countable pairwise disjoint family $B_h:=B(c_h,25r_h)$ with $h\in\N$, with
$c_h\in K''$ and
$r_h=\rho_{m_h}(c_h)$ for some $m_h\geq m_6(c_h)$, such that
\begin{equation}
    \mu\bigg(K''\setminus\bigcup_{h\in\N}B_h\bigg)=0.
    \label{eq:BesicovitchStepII}
\end{equation}
For reader's convenience we summarize in the following the properties that these balls $B_h$ enjoy.

By the definition of $m_6(c_h)$, for every $h\in\N$ we notice that
\begin{equation}
    \mu(K''\cap B_h)
    \leq
    2A^2\mu(K''\cap \tfrac14 B_h).
    \label{eq:quarterMassStepII}
\end{equation}
By \eqref{eq:zetaStepII} and \eqref{eq:qellStepII}, we have for every $z\in\Sigma_{c_h}^0$ that
\begin{equation}
    |\zeta_{m_h}^{c_h}(z)-c_h-r_hz|
    \leq
    \frac{\mathfrak r^2\xi}{100}r_h.
    \label{eq:zetaOriginalScaleStepII}
\end{equation}
Similarly, by \eqref{eq:qellStepII}, for every $z\in\Sigma_{c_h}^0$ and every
$\ell=1,\ldots,d+1$,
\begin{equation}
    |q_{\ell,m_h}^{c_h}(z)-c_h-r_hy_\ell^{c_h}(z)|
    \leq
    \frac{\mathfrak r^2\xi}{100}r_h.
    \label{eq:qellOriginalScaleStepII}
\end{equation}
Consequently, using \eqref{eq:tangentSeparationStepII}, as above we have that for every
$z\in\Sigma_{c_h}^0$ and every $\ell_1\neq \ell_2$, there holds
\begin{equation}
    \begin{split}
    |q_{\ell_1,m_h}^{c_h}(z)-q_{\ell_2,m_h}^{c_h}(z)|
    \geq
    3\mathfrak r\xi r_h.
    \end{split}
    \label{eq:qellSeparationOriginalScaleStepII}
\end{equation}

Moreover, by \eqref{eq:finiteMassComparisonStepII}, since
$m_h\geq m_6(c_h)$, for every $z\in\Sigma_{c_h}^0$ and every
$\ell=1,\ldots,d+1$ we have
\begin{equation}
    \mu(B(\zeta_{m_h}^{c_h}(z),25\xi r_h))
    \leq
    16D
    \mu(B(q_{\ell,m_h}^{c_h}(z),\mathfrak r\xi r_h/16)).
    \label{eq:finiteMassComparisonSelectedStepII}
\end{equation}

By \eqref{eq:scalesmallStepII}, since $c_h\in K_{i(c_h)}$ and
$r_h<\delta_{i(c_h)}/1000$, we have
$$
B_h=B(c_h,25r_h)\subseteq U(K_{i(c_h)},\delta_{i(c_h)})=\Omega_{i(c_h)}.
$$
Therefore the cancellation estimate \eqref{stimagradort} from Step I applies inside $B_h$ along directions in $V_{i(c_h)}^\perp$.

Finally, applying \eqref{eq:supportcloseStepII} with $w=c_h$ and $m=m_h$, we obtain
\begin{equation}
    \mu\bigl(
    B_h
    \setminus
    B(c_h+r_h\supp\nu_{c_h},\xi r_h/32)
    \bigr)
    \leq
    \frac{\vartheta}{3}\mu(B_h).
    \label{eq:supportcloseFinalStepII}
\end{equation}
This concludes the second step.

\medskip

\textbf{Step III: construction of the second perturbation.}
We now construct the second perturbation inside the balls selected in Step II. For every $h\in\N$, let $\psi_h:\R^n\to[0,1]$ be a smooth function such that
$\psi_h=1$ on $B(c_h,20r_h)$, $\supp\psi_h\subseteq B_h$ and
\begin{equation}
    \lVert \nabla\psi_h\rVert_\infty\leq 10r_h^{-1},
    \qquad
    \lVert \nabla^2\psi_h\rVert_\infty\leq 10r_h^{-2}.
    \label{eq:cutoffSecondPerturbation}
\end{equation}
Since the balls $B_h$ are pairwise disjoint, at every point of $\R^n$ at most
one of the functions $\psi_h$ is non-zero. For every $\iota=1,\ldots,d$, let
$\beta^\iota\in\R^{d+1}$ be the vectors given by \cref{evadeaffinepro}, with
$|\beta^\iota|\leq 1$. Fix $h\in\N$ and $\iota=1,\ldots,d$. For every
$z\in\Sigma_{c_h}^0$ and every $\ell=1,\ldots,d+1$, define
\begin{equation}
    \alpha_{\ell,h,\iota}(z)
    :=
    \frac{\mathfrak r\xi}{16}\beta^\iota_\ell .
    \label{eq:alphaSecondPerturbation}
\end{equation}
Then, we can apply \cref{lemmanonfancyutile} in the plane $V(\mu,c_h)^\perp$ to the
finite family of transverse points
$$
\pi_{V(\mu,c_h)^\perp}y_\ell^{c_h}(z),
\qquad
z\in\Sigma_{c_h}^0,\quad \ell=1,\ldots,d+1.
$$
Equivalently, since the functions produced by \cref{lemmanonfancyutile} are
invariant along $V(\mu,c_h)$, we prescribe the values at the points
$y_\ell^{c_h}(z)$. The separation needed in \cref{lemmanonfancyutile} follows from
\eqref{eq:tangentSeparationStepII} for points corresponding to the same
$z\in\Sigma_{c_h}^0$, and from the additional separation imposed in the choice
of $\Sigma_{c_h}^0$ for points corresponding to distinct elements of
$\Sigma_{c_h}^0$. Thus, for every
$\iota=1,\ldots,d$, we obtain a $1/4$-Lipschitz function
$\Phi_{h,\iota}:\R^n\to\R$ such that
\begin{enumerate}
    \item $\Phi_{h,\iota}$ is invariant along $V(\mu,c_h)$;
    \item for every $z\in\Sigma_{c_h}^0$ and every $\ell=1,\ldots,d+1$, we have $\Phi_{h,\iota}(y_\ell^{c_h}(z))
        =
        \alpha_{\ell,h,\iota}(z)$;
    \item $ \lVert \Phi_{h,\iota}\rVert_\infty \leq\mathfrak r\xi/8$.
\end{enumerate}
Here we use the interpolation proposition in the form where the prescribed
points lie in the transverse space and the extension is invariant along the
corresponding plane.

We rescale $\Phi_{h,\iota}$ to the ball $B_h$ by setting
\begin{equation}
    \widetilde\Phi_{h,\iota}(x)
    :=
    r_h\Phi_{h,\iota}\bigg(\frac{x-c_h}{r_h}\bigg).
    \label{eq:PhiTildeSecondPerturbation}
\end{equation}
Then $\widetilde\Phi_{h,\iota}$ is invariant along $V(\mu,c_h)$ and there holds
\begin{equation}
    \lVert \widetilde\Phi_{h,\iota}\rVert_\infty
    \leq
    \frac{\mathfrak r\xi}{8}r_h,
    \qquad
    \lVert D\widetilde\Phi_{h,\iota}\rVert_\infty
    \leq
    \frac14,\qquad \widetilde\Phi_{h,\iota}(c_h+r_hy_\ell^{c_h}(z))
    =
    \frac{\mathfrak r\xi r_h}{16}\beta^\iota_\ell,
    \label{eq:PhiTildeBoundsSecondPerturbation}
\end{equation}
for every $z\in\Sigma_{c_h}^0$ and every $\ell=1,\ldots,d+1$.
Since $q_{\ell,m_h}^{c_h}(z)$ satisfies \eqref{eq:qellOriginalScaleStepII} and
$\widetilde\Phi_{h,\iota}$ is $1/4$-Lipschitz, we also have
\begin{equation}
    \bigg|
    \widetilde\Phi_{h,\iota}(q_{\ell,m_h}^{c_h}(z))
    -
    \frac{\mathfrak r\xi r_h}{16}\beta^\iota_\ell
    \bigg|
    \leq
    \frac{\mathfrak r^2\xi}{400}r_h
    \label{eq:PhiTildeValuesActualSecondPerturbation}
\end{equation}
for every $z\in\Sigma_{c_h}^0$ and every $\ell=1,\ldots,d+1$.
Finally, for every $\iota=1,\ldots,d$, define
\begin{equation}
    \hat f_{\xi,\iota}
    :=
    (1-\varsigma(\xi))
    \bigg(
    \tilde f+
    \sum_{h\in\N}\psi_h\widetilde\Phi_{h,\iota}
    \bigg).
    \label{eq:defineSecondPerturbation}
\end{equation}
The sum is pointwise well defined because the balls $B_h$ are pairwise
disjoint and $\supp\psi_h\subseteq B_h$.

Finally, since $r_h<r/125$ by \eqref{eq:scalesmallStepII}, the second
perturbation is supported on scales strictly smaller than the scale $r$ fixed
in the statement of the lemma. This concludes the construction of the second
perturbation.

\medskip

\textbf{Step IV: admissibility of the second perturbation.}
We now check that, for every $\iota=1,\ldots,d$, the function
$\hat f_{\xi,\iota}$ belongs to $X_{\mu,\mathscr E}$ and remains close to $f$
in the supremum norm.

We first estimate the supremum norm. Since the balls $B_h$ are pairwise
disjoint and $\supp\psi_h\subseteq B_h$, by \eqref{eq:PhiTildeBoundsSecondPerturbation} we have
\begin{equation}
    \bigg\|
    \sum_{h\in\N}\psi_h\widetilde\Phi_{h,\iota}
    \bigg\|_\infty
    \leq
    \sup_{h\in\N}\lVert \widetilde\Phi_{h,\iota}\rVert_\infty
    \leq
    \frac{\mathfrak r\xi}{8}\sup_{h\in\N}r_h.
    \label{eq:secondPerturbationSup}
\end{equation}
By \eqref{eq:scalesmallStepII}, $r_h\leq \xi^2$ for every $h\in\N$, we have
\begin{equation}
    \bigg\|
    \sum_{h\in\N}\psi_h\widetilde\Phi_{h,\iota}
    \bigg\|_\infty
    \leq
    \xi^2.
    \label{eq:secondPerturbationSupSmall}
\end{equation}
Therefore, using \eqref{eq:stimainfinitotilde},
\begin{equation}
    \begin{split}
    \lVert f-\hat f_{\xi,\iota}\rVert_\infty
    &\leq
    \lVert f-\tilde f\rVert_\infty
    +
    \varsigma(\xi)\lVert \tilde f\rVert_\infty
    +
    \bigg\|
    \sum_{h\in\N}\psi_h\widetilde\Phi_{h,\iota}
    \bigg\|_\infty\\
    &\leq
    (\lVert f\rVert_\infty+n)\varsigma(\xi)
    +
    \varsigma(\xi)
    \bigl(\lVert f\rVert_\infty+(\lVert f\rVert_\infty+n)\varsigma(\xi)\bigr)
    +
    \xi^2\\
    &\leq
    2(\lVert f\rVert_\infty+d+1)\varsigma(\xi).
    \end{split}
    \label{eq:distanceFinalPerturbation}
\end{equation}
In the last inequality we used $n\leq d+1$, $\xi^2\leq\varsigma(\xi)$, and
$\varsigma(\xi)<1$.

We now prove that $\hat f_{\xi,\iota}\in X_{\mu,\mathscr E}$. Since
$\tilde f\in X_{\mu,\mathscr E}$ by Step I, outside $\bigcup_h B_h$ we have
$$
\hat f_{\xi,\iota}=(1-\varsigma(\xi))\tilde f,
$$
and therefore the required derivative bounds are immediate. Let now
$x\in B_h$ be a point where the derivatives $\partial_{e_\kappa(x)}\widetilde\Phi_{h,\iota}$ exist for every $\kappa=1,\ldots, n$ every $h\in \N$ and $\iota=1,\ldots, d$. Since the
balls $B_h$ are pairwise disjoint,
for every $\kappa=1,\ldots,n$ we have
\begin{equation}
    \partial_{e_\kappa(x)}\hat f_{\xi,\iota}(x)
    =
    (1-\varsigma(\xi))
    \bigl[
    \partial_{e_\kappa(x)}\tilde f(x)
    +
    \partial_{e_\kappa(x)}\psi_h(x)\,\widetilde\Phi_{h,\iota}(x)
    +
    \psi_h(x)\partial_{e_\kappa(x)}\widetilde\Phi_{h,\iota}(x)
    \bigr].
    \label{eq:derivativeSecondPerturbation}
\end{equation}
By \eqref{eq:cutoffSecondPerturbation} and
\eqref{eq:PhiTildeBoundsSecondPerturbation}, we have
\begin{equation}
    |\partial_{e_\kappa(x)}\psi_h(x)\,\widetilde\Phi_{h,\iota}(x)|
    \leq
    10r_h^{-1}\frac{\mathfrak r\xi}{8}r_h
    \leq
    2\xi .
    \label{eq:cutoffTermSecondPerturbation}
\end{equation}

We distinguish two cases. Fix $h\in\N$ and write $i:=i(c_h)$. First, assume
that $\kappa\leq d_i$. Then $E_{\kappa,i}\in V_i$. Since
$\widetilde\Phi_{h,\iota}$ is invariant along $V(\mu,c_h)$, and since
$d(V_i,V(\mu,c_h))\leq \xi^{4/\omega}$, by \cref{prop:continuitax} we have
\begin{equation}
    |\pi_{V(\mu,c_h)^\perp}E_{\kappa,i}|
    \leq
    \xi.
    \label{eq:projectionFrozenTangent}
\end{equation}
Therefore, at every point where $\widetilde\Phi_{h,\iota}$ is differentiable,
\begin{equation}
    \begin{split}
    |\partial_{E_{\kappa,i}}\widetilde\Phi_{h,\iota}(x)|=|\partial_{\pi_{V(\mu,c_h)^\perp}E_{\kappa,i}}\widetilde\Phi_{h,\iota}(x)|\leq
    \lVert D\widetilde\Phi_{h,\iota}\rVert_\infty
    |\pi_{V(\mu,c_h)^\perp}E_{\kappa,i}|\leq
    \frac{\xi}{4}.
    \end{split}
    \label{eq:frozenTangentDerivativePhiSecondPerturbation}
\end{equation}
Using \eqref{eq:frozenbasisStepI} and
$\lVert D\widetilde\Phi_{h,\iota}\rVert_\infty\leq 1/4$, we obtain
\begin{equation}
    \begin{split}
    |\partial_{e_\kappa(x)}\widetilde\Phi_{h,\iota}(x)|
    &\leq|\partial_{E_{\kappa,i}}\widetilde\Phi_{h,\iota}(x)|+|\partial_{e_\kappa(x)-E_{\kappa,i}}\widetilde\Phi_{h,\iota}(x)|\leq\frac{\xi}{4}+\frac14 |e_\kappa(x)-E_{\kappa,i}|\leq\frac{\xi}{4}+\frac{\xi}{8}\leq\xi.
    \end{split}
    \label{eq:tangentDerivativePhiSecondPerturbation}
\end{equation}
Since $\tilde f\in X_{\mu,\mathscr E}$, we also have
$|\partial_{e_\kappa(x)}\tilde f(x)|\leq 1$. Hence, by
\eqref{eq:derivativeSecondPerturbation}, \eqref{eq:cutoffTermSecondPerturbation}, and
\eqref{eq:tangentDerivativePhiSecondPerturbation},
\begin{equation}
    |\partial_{e_\kappa(x)}\hat f_{\xi,\iota}(x)|
    \leq
    (1-\varsigma(\xi))(1+3\xi)
    \leq
    1.
    \label{eq:tangentDerivativeHat}
\end{equation}
Here we used $3\xi\leq\varsigma(\xi)$.

We now assume that $\kappa>d_i$. Then $E_{\kappa,i}\in V_i^\perp$. Since
$B_h\subseteq\Omega_i$ by Step II, the cancellation estimate from Step I
applies inside $B_h$. Hence, by \eqref{eq:smallFrozenGradient}, we have
\begin{equation}
    |\partial_{E_{\kappa,i}}\tilde f(x)|
    \leq
    \varsigma(\xi).
    \label{eq:normalFrozenDerivativeTilde}
\end{equation}
Using \eqref{eq:frozenbasisStepI} and $\tilde f\in X_{\mu,\mathscr E}$, we get
\begin{equation}
    \begin{split}
    |\partial_{e_\kappa(x)}\tilde f(x)|\leq|\partial_{E_{\kappa,i}}\tilde f(x)|+|\partial_{e_\kappa(x)-E_{\kappa,i}}\tilde f(x)|\leq\varsigma(\xi)+n|e_\kappa(x)-E_{\kappa,i}|\leq \varsigma(\xi)+\frac{n\xi}{2}\leq2\varsigma(\xi).
    \end{split}
    \label{eq:normalDerivativeTildeMoving}
\end{equation}
Moreover, combining \eqref{eq:derivativeSecondPerturbation},
\eqref{eq:cutoffTermSecondPerturbation}, \eqref{eq:normalDerivativeTildeMoving}, and
\eqref{eq:PhiTildeBoundsSecondPerturbation}, we get
\begin{equation}
    |\partial_{e_\kappa(x)}\hat f_{\xi,\iota}(x)|
    \leq
    (1-\varsigma(\xi))
    \bigl(2\varsigma(\xi)+2\xi+\tfrac14\bigr)
    \leq
    1.
    \label{eq:normalDerivativeHat}
\end{equation}
This proves that $\hat f_{\xi,\iota}\in X_{\mu,\mathscr E}$ for every
$\iota=1,\ldots,d$.

\medskip

\textbf{Step V: reduction to one of the perturbations.}
We reduce the conclusion of the lemma to the estimate proved in the next step.
We claim that it is enough to prove that, for every $h\in\N$ and for
$\mu$-almost every $y\in K''\cap \tfrac14 B_h$ such that
\begin{equation}
    \dist(y,c_h+r_h\supp\nu_{c_h})
    \leq
    \frac{\xi r_h}{32},
    \label{eq:goodDistanceStepV}
\end{equation}
there holds
\begin{equation}
    \sum_{\iota=1}^{d}
    \Omega_{1,\mu,\mathcal F}(\hat f_{\xi,\iota};y,25\xi r_h)
    \geq
    \frac{d\,\mathfrak r}{40000(d+1)^2D}.
    \label{eq:claimStepVI}
\end{equation}
Indeed, fix $h\in\N$ and let $y\in K''\cap \tfrac14 B_h$ satisfy
\eqref{eq:goodDistanceStepV}. Since
$y\in K''\subseteq \mathfrak G_{\mathfrak r,D,\xi}$ and
$r_h<r/125$ by \eqref{eq:scalesmallStepII}, the choice
$z=y$ and $s=25\xi r_h$, is admissible in the supremum below. Hence, if \eqref{eq:claimStepVI} holds,
then
\begin{equation}
    \begin{split}
    &\sum_{\iota=1}^{d}
    \sup_{0<s<r}
    \sup_{z\in \mathfrak G_{\mathfrak r,D,\xi}\,:\,y\in B(z,s/5)}
    \Omega_{1,\mu,\mathcal F}(\hat f_{\xi,\iota};z,s)\geq
    \sum_{\iota=1}^{d}
    \Omega_{1,\mu,\mathcal F}(\hat f_{\xi,\iota};y,25\xi r_h)\geq
    \frac{d\,\mathfrak r}{40000(d+1)^2D}.
    \end{split}
    \label{eq:sumSupStepV}
\end{equation}
We now estimate the measure of the set where this applies. Set
\begin{equation}
    \mathcal D:=
    \bigcup_{h\in\N}
    \bigg\{
    y\in K''\cap \tfrac14 B_h:
    \dist(y,c_h+r_h\supp\nu_{c_h})
    \leq
    \frac{\xi r_h}{32}
    \bigg\}.
    \label{eq:defDStepV}
\end{equation}
Since the balls $B_h$ are pairwise disjoint, by
\eqref{eq:supportcloseFinalStepII}, \eqref{eq:fourfoldDoublingStepII}, and
\eqref{eq:quarterMassStepII}, we have
\begin{equation}
    \begin{split}
    \mu(\mathcal D)
    &\geq
    \sum_{h\in\N}\mu(K''\cap \tfrac14 B_h)
    -
    \sum_{h\in\N}
    \mu\bigl(
    B_h\setminus B(c_h+r_h\supp\nu_{c_h},\xi r_h/32)
    \bigr)\\
    &\geq
    \sum_{h\in\N}\mu(K''\cap \tfrac14 B_h)-\frac{\vartheta}{3}\sum_{h\in\N}\mu(B_h)\geq
    \sum_{h\in\N}\mu(K''\cap \tfrac14 B_h)-\frac{2A^2\vartheta}{3}\sum_{h\in\N}\mu(K''\cap \tfrac14 B_h)\\
    &\geq\frac56\sum_{h\in\N}\mu(K''\cap \tfrac14 B_h)\geq\frac{5}{12A^2}\mu(K'').
    \end{split}
    \label{eq:massDStepV}
\end{equation}
By \eqref{eq:lossKcirc} and \eqref{eq:choiceKsecond}, we infer that
\begin{equation}
    \mu(K'')
    \geq
    (1-\vartheta/2)\mu(K')
    \overset{\eqref{eq:choiceKprime}}{\geq}
    (1-\vartheta/2)(1-\vartheta/6)
    \mu(\mathfrak G_{\mathfrak r,D,\xi}).
    \label{eq:KsecondLowerStepV}
\end{equation}
Since $\vartheta<1/(4A^2)\leq 1/4$, combining
\eqref{eq:massDStepV} and \eqref{eq:KsecondLowerStepV} gives
\begin{equation}
    \mu(\mathcal D)
    >
    \frac{1}{3A^2}
    \mu(\mathfrak G_{\mathfrak r,D,\xi}).
    \label{eq:massGoodSetStepV}
\end{equation}
By \eqref{eq:sumSupStepV}, for every $y\in\mathcal D$, 
at least one index $\iota\in\{1,\ldots,d\}$ satisfies
\begin{equation}
    \sup_{0<s<r}
    \sup_{z\in \mathfrak G_{\mathfrak r,D,\xi}\,:\,y\in B(z,s/5)}
    \Omega_{1,\mu,\mathcal F}(\hat f_{\xi,\iota};z,s)
    \geq
    \frac{\mathfrak r}{40000(d+1)^2D}.
    \label{eq:pointwisePigeonholeStepV}
\end{equation}
Therefore there exists $\iota_0\in\{1,\ldots,d\}$ such that
\begin{equation}
    \begin{split}
    &\mu\Bigg(
    \Bigg\{
y\in \mathcal D:    \sup_{0<s<r}
    \sup_{z\in \mathfrak G_{\mathfrak r,D,\xi}\,:\,y\in B(z,s/5)}
    \Omega_{1,\mu,\mathcal F}(\hat f_{\xi,\iota_0};z,s)
    \geq
    \frac{\mathfrak r}{40000(d+1)^2D}
    \Bigg\}
    \Bigg)\\
    &\qquad\qquad>
    \frac{1}{3dA^2}
    \mu(\mathfrak G_{\mathfrak r,D,\xi})
    \geq
    \frac{1}{4(d+1)A^2}
    \mu(\mathfrak G_{\mathfrak r,D,\xi}).
    \end{split}
    \label{eq:pigeonholeStepV}
\end{equation}
By inner regularity, we can find a compact subset of the set appearing in
\eqref{eq:pigeonholeStepV} with measure at least
$$
\frac{1}{4(d+1)A^2}
\mu(\mathfrak G_{\mathfrak r,D,\xi}).
$$
Relabelling this compact set as $K'$ and setting
$$
\hat f:=\hat f_{\xi,\iota_0},
$$
the conclusion of the lemma follows from \eqref{eq:claimStepVI}. It remains to
prove \eqref{eq:claimStepVI}.

\medskip
\medskip

\textbf{Step VI: lower bound at the selected scales.}
We prove \eqref{eq:claimStepVI}. Fix $h\in\N$ and let
$y\in K''\cap \tfrac14 B_h$ be such that
\begin{equation}
    \dist(y,c_h+r_h\supp\nu_{c_h})
    \leq
    \frac{\xi r_h}{32}.
    \label{eq:goodDistanceStepVI}
\end{equation}
Write $i:=i(c_h)$. By \eqref{eq:goodDistanceStepVI}, there exists
$p\in\supp\nu_{c_h}$ such that
\begin{equation}
    |y-c_h-r_hp|
    \leq
    \frac{\xi r_h}{32}.
    \label{eq:yCloseSupportStepVI}
\end{equation}
Since $y\in \tfrac14 B_h$, we have $p\in B(0,9)$. Hence, there exists $z\in\Sigma_{c_h}^0$ such that $|p-z|\leq 2\xi$ since $\Sigma_{c_h}^0$ is a maximal $2\xi$-separated set in $\supp\nu_{c_h}\cap B(0,9)$.
This implies that
\begin{equation}
    |y-c_h-r_hz|
    \leq
    3\xi r_h.
    \label{eq:yCloseNetStepVI}
\end{equation}
By \eqref{eq:zetaOriginalScaleStepII}, and since $\mathfrak r<1$, this also shows that $|y-\zeta_{m_h}^{c_h}(z)|\leq4\xi r_h$. Therefore
\begin{equation}
    B(y,25\xi r_h)
    \subseteq
    B(\zeta_{m_h}^{c_h}(z),30\xi r_h).
    \label{eq:ballYInsideZetaStepVI}
\end{equation}
Moreover, by \eqref{eq:fullYStepII}, for every $\ell=1,\ldots,d+1$, we have $|y_\ell^{c_h}(z)-z|\leq\frac{\xi}{8}$. Hence, using \eqref{eq:yCloseNetStepVI} and
\eqref{eq:qellOriginalScaleStepII}, we obtain
\begin{equation}
    |q_{\ell,m_h}^{c_h}(z)-y|
    \leq
    3\xi r_h
    \qquad\text{for every }\ell=1,\ldots,d+1.
    \label{eq:qellCloseYStepVI}
\end{equation}
Since $\mathfrak r<1$, this implies
\begin{equation}
    B(q_{\ell,m_h}^{c_h}(z),\mathfrak r\xi r_h/16)
    \subseteq
    B(y,25\xi r_h)
    \qquad\text{for every }\ell=1,\ldots,d+1.
    \label{eq:qellBallsInsideYStepVI}
\end{equation}
By \eqref{eq:finiteMassComparisonSelectedStepII} and
\eqref{eq:ballYInsideZetaStepVI}, for every $\ell=1,\ldots,d+1$ we have
\begin{equation}
    \mu(B(y,25\xi r_h))
    \leq
    16D
    \mu(B(q_{\ell,m_h}^{c_h}(z),\mathfrak r\xi r_h/16)).
    \label{eq:massComparisonYStepVI}
\end{equation}

Since $y\in\tfrac14 B_h$ and $\xi\leq\mathfrak r/2<1/4$, we have
\begin{equation}
    B(y,25\xi r_h)\subseteq B(c_h,20r_h).
    \label{eq:ballYInsideCutoffStepVI}
\end{equation}
Thus $\psi_h=1$ on $B(y,25\xi r_h)$. Moreover, since
$B_h\subseteq\Omega_i$, all the balls appearing in
\eqref{eq:qellBallsInsideYStepVI} are contained in $\Omega_i$. We now prove the affine approximation for $\tilde f$ on
$B(y,25\xi r_h)$.
We claim that, for every $x\in B(y,25\xi r_h)$,
\begin{equation}
    \big|\tilde f(x)-\tilde f(y)-\langle \pi_{V_i}Df_\xi(y),x-y\rangle\big|\leq 4n\varsigma(\xi)|x-y|.
    \label{eq:tildeAffineIncrementStepVIorginasd}
\end{equation}
Indeed, let $x_1,x_2\in B(y,25\xi r_h)$ write
$x_2-x_1=v+u$, where $v\in V_i$ and
$u\in V_i^\perp$ and set $x_1':=x_1+v$. Since
$B(y,25\xi r_h)\subseteq B_h\subseteq\Omega_i$ and $\Omega_i$ is convex along
the two segments considered below, we can use \eqref{stimagradort} on the
segment $[x_1',x_2]$. Hence
\begin{equation}
    |\tilde f(x_2)-\tilde f(x_1')|
    \leq
    n\varsigma(\xi)|u|.
    \label{eq:normalIncrementStepVI}
\end{equation}
Let us determine the behaviour of $\tilde f$ along the directions
parallel to $V_i$. Let $x\in\Omega_i$ be a point where the derivatives of the
functions $w_{j,i}$ exist, and let $v\in V_i$. Since $\Omega_i\subseteq B_i''$,
we have $\psi_i=1$ and $\nabla\psi_i=0$ on $\Omega_i$. Hence, from
\eqref{eq:definetildefStepI}, we infer that
\begin{equation}
    \begin{split}
    \partial_v\tilde f(x)=(1-\varsigma(\xi))
    \bigg[\langle Df_\xi(x),v\rangle-\sum_{j=d_i+1}^{n}\partial_v w_{j,i}(x)\,\partial_{E_{j,i}}f_\xi(x)-\sum_{j=d_i+1}^{n}w_{j,i}(x)\,\partial_v\partial_{E_{j,i}}f_\xi(x)\bigg].
    \end{split}
    \label{eq:tangentDerivativeFormulaStepVI}
\end{equation}
Since $v\in V_i$ and $E_{j,i}\in V_i^\perp$ for every $j=d_i+1,\ldots,n$, we
have $v\in E_{j,i}^{\perp}$. Therefore, by
\eqref{eq:widthderivativeStepI},
\begin{equation}
    |\partial_v w_{j,i}(x)|
    \leq
    2\xi |v|
    \qquad\text{for every }j=d_i+1,\ldots,n.
    \label{eq:tangentWidthDerivativeStepVI}
\end{equation}
Using \eqref{estimateonD2f}, \eqref{eq:widthboundStepI}, and
\eqref{eq:epsilonI}, we also have
\begin{equation}
    \sum_{j=d_i+1}^{n}
    |w_{j,i}(x)\,\partial_v\partial_{E_{j,i}}f_\xi(x)|
    \leq
    \sum_{j=d_i+1}^{n}
    2\varepsilon_i\,10n\xi^{-1/3}|v|
    \leq
    \frac{\xi}{5}|v|.
    \label{eq:tangentD2ErrorStepVI}
\end{equation}
Moreover, $\lVert Df_\xi\rVert_\infty\leq 1$, and hence
\begin{equation}
    \sum_{j=d_i+1}^{n}
    |\partial_v w_{j,i}(x)\,\partial_{E_{j,i}}f_\xi(x)|
    \leq
    2n\xi |v|.
    \label{eq:tangentWidthErrorStepVI}
\end{equation}
Combining \eqref{eq:tangentDerivativeFormulaStepVI},
\eqref{eq:tangentD2ErrorStepVI}, and \eqref{eq:tangentWidthErrorStepVI}, we get
\begin{equation}
    \begin{split}
    \big|\partial_v\tilde f(x)-\langle Df_\xi(x),v\rangle\big|\leq\varsigma(\xi)|v|+2n\xi |v|+\frac{\xi}{5}|v|\leq2\varsigma(\xi)|v|.
    \end{split}
    \label{eq:tangentDerivativeErrorStepVI}
\end{equation}

Moreover, since $|x_1'-x_1|\leq |x_2-x_1|\leq 50\xi r_h$ and
$r_h<\xi^2$, the estimate for $D^2f_\xi$ gives
\begin{equation}
    |\langle Df_\xi(x_1')-Df_\xi(x_1),v\rangle|
    \leq
    \varsigma(\xi)|v|
    \qquad\text{for every }x\in[x_1,x_1'].
    \label{eq:TaylorDfStepVI}
\end{equation}
Therefore
\begin{equation}
    \big|\tilde f(x_1')-\tilde f(x_1)-\langle\pi_{V_i}Df_\xi(x_1),x_1'-x_1\rangle\big|\leq3\varsigma(\xi)|v|.
    \label{eq:tangentIncrementStepVI}
\end{equation}
Putting together \eqref{eq:tangentIncrementStepVI} and \eqref{eq:normalIncrementStepVI} we conclude that 
\begin{equation}
    \big|\tilde f(x)-\tilde f(y)-\langle \pi_{V_i}Df_\xi(y),x-y\rangle\big|\leq (n+3)\varsigma(\xi)\lvert x-y\rvert
    \label{eq:tildeAffineIncrementStepVI}
\end{equation}

We now compare $\tilde f$ with functions in $\mathcal F$ on the small balls Let $A\in\mathcal F$. Since $\psi_h=1$ on $B(y,25\xi r_h)$, by
\eqref{eq:defineSecondPerturbation} we have
\begin{equation}
    \hat f_{\xi,\iota}
    =
    (1-\varsigma(\xi))
    (\tilde f+\widetilde\Phi_{h,\iota})
    \qquad\text{on }B(y,25\xi r_h).
    \label{eq:hatLocalStepVI}
\end{equation}
Using \eqref{eq:massComparisonYStepVI},
\eqref{eq:qellBallsInsideYStepVI}, \eqref{eq:tildeAffineIncrementStepVI} and since $\mathcal F$ contains linear functions, we get
\begin{equation}
    \begin{split}
    &\inf_{A\in \mathcal F}\fint_{B(y,25\xi r_h)}
    \frac{|\hat f_{\xi,\iota}-A|}{25\xi r_h}\,d\mu\geq(1-\varsigma(\xi))\inf_{A\in \mathcal F}\fint_{B(y,25\xi r_h)}
    \frac{|\widetilde{\Phi}_{h,\iota}-A|}{25\xi r_h}\,d\mu-2(n+3)\varsigma(\xi).
    \end{split}
    \label{eq:lowerBySmallBallsStepVI}
\end{equation}
By \eqref{eq:PhiTildeValuesActualSecondPerturbation} and by the fact that
$\widetilde\Phi_{h,\iota}$ is $1/4$-Lipschitz, for every
$w\in B(q_{\ell,m_h}^{c_h}(z),\mathfrak r\xi r_h/16)$ we have
\begin{equation}
    \bigg|
    \widetilde\Phi_{h,\iota}(w)
    -
    \frac{\mathfrak r\xi r_h}{16}\beta_\ell^\iota
    \bigg|
    \leq
    \frac{\mathfrak r\xi r_h}{50}.
    \label{eq:PhiAlmostConstantStepVI}
\end{equation}
Using \eqref{eq:PhiAlmostConstantStepVI}, \eqref{eq:qellSeparationOriginalScaleStepII},
\eqref{eq:qellCloseYStepVI}, and \cref{evadeaffinepro}, applied with $B:=B(y,25\xi r_h)$, $\zeta:=\frac{\mathfrak r}{400}$ and $\eta:=\frac{\mathfrak r}{400}$ we obtain, that
\begin{equation}
    \sum_{\iota=1}^{d}\Omega_{1,\mu,\mathcal F}(\hat f_{\xi,\iota};y,25\xi r_h)\geq\frac{d\,\mathfrak r}{20000(d+1)^2D}-2(n+3)d\varsigma(\xi).
    \label{eq:almostClaimStepVI}
\end{equation}
Finally, by the choice of $\xi$ we infer that
\begin{equation}
    \sum_{\iota=1}^{d}
    \Omega_{1,\mu,\mathcal F}(\hat f_{\xi,\iota};y,25\xi r_h)
    \geq
    \frac{d\,\mathfrak r}{40000(d+1)^2D}.
    \label{eq:claimStepVIproved}
\end{equation}
This proves \eqref{eq:claimStepVI}, and therefore, by Step V, the proof of the
lemma is complete.
\end{proof}

\subsection{Proof of the main qualitative result}

\begin{proposizione}\label{coveringtangents}
Suppose $\mu$ is a Radon measure such that
$$
\limsup_{r\to 0}\frac{\mu(B(x,2r))}{\mu(B(x,r))}<\infty
\qquad\text{for $\mu$-almost every }x\in \R^n.
$$
For every $D,A\geq 1$ and every $\mathfrak r>0$, let $E_{D,\mathfrak r,A}$ be the set of those $x\in\supp\mu$ for which
\begin{equation}
    \Tan(\mu,x)\cap \bigcup_{k\in\N}\bigcap_{j\geq k}\mathfrak S_d(V(\mu,x),2^{-j},D,\mathfrak r,A)\neq \emptyset.
    \label{definitionEDrA}
\end{equation}
Then, the set $E_{D,\mathfrak r,A}$ is $\mu$-measurable.

Assume that for $\mu$-almost every $x\in\R^n$ there exist $D,A\geq1$ and $\mathfrak r>0$ such that $x\in E_{D,\mathfrak r,A}$. Then
$$
\mu\Big(\R^n\setminus\bigcup_{m_1,m_2,m_3\in\N}E_{m_1,m_32^{-m_2},m_3}\Big)=0.
$$
\end{proposizione}

\begin{proof}
By \cref{measurability}, the sets $E_{m_1,m_32^{-m_2},m_3}$ are $\mu$-measurable for every $m_1,m_2,m_3\in\N$. Fix $x\in\R^n$ such that $x\in E_{D,\mathfrak r,A}$ for some $D,A\geq1$ and $\mathfrak r>0$. Choose $m_3\in\N$ in such a way that
$A\leq m_3$ and $\mathfrak r\leq m_3$.
Then choose $m_2\in\N$ in such a way that
$$
m_32^{-m_2}\leq \mathfrak r\leq m_32^{-m_2+1}.
$$
By \cref{badtangentsrelationships}, for every $k\in\Z$ we have
$$
\mathfrak S_d(V(\mu,x),2^{-k},D,\mathfrak r,A)
\subseteq
\mathfrak S_d(V(\mu,x),2^{-k},A^2D,m_32^{-m_2},A),
$$
and finally choose $m_1\in\N$ such that $A^2D\leq m_1$.
Using the monotonicity in $D$ and in the doubling constant given by \cref{badtangentsrelationships}, we obtain
$$
\mathfrak S_d(V(\mu,x),2^{-k},A^2D,m_32^{-m_2},A)
\subseteq
\mathfrak S_d(V(\mu,x),2^{-k},m_1,m_32^{-m_2},m_3)
$$
for every $k\in\Z$. Therefore
$$
x\in E_{m_1,m_32^{-m_2},m_3}.
$$
Hence
$$
\{x:\text{there exist }D,A\geq1\text{ and }\mathfrak r>0\text{ such that }x\in E_{D,\mathfrak r,A}\}
\subseteq
\bigcup_{m_1,m_2,m_3\in\N}E_{m_1,m_32^{-m_2},m_3}.
$$
By assumption, the set on the left-hand side has full $\mu$-measure. Since the union on the right-hand side is $\mu$-measurable thanks to the discussion above, we conclude the proof of the proposition.
\end{proof}

We now prove the main qualitative theorem. The argument has two rather
different parts. The first one is local and infinitesimal: the perturbation
lemma turns the presence of the tangents in
$\mathfrak S_d(V(\mu,x),\xi,D,\mathfrak r,A)$ into an obstruction to the
vanishing of the coefficients $\Omega_{1,\mu,\mathcal F}$. More precisely,
after a countable decomposition of the measure, the lemma shows that on each
piece the functions for which the coefficients vanish on a set of positive
measure form a meagre set. This is the Baire category part of the proof.

The second part is global and uses a beautiful result of G. Alberti and
A. Marchese, namely \cite[Proposition 8.4]{AlbertiMarchese}. Once the bad functions have been produced separately on the pieces of
the decomposition, one has to paste them together without destroying the lower
bounds for $\Omega_{1,\mu,\mathcal F}$. The extension result of Alberti and
Marchese makes this possible in a particularly efficient way: each local
function can be extended to a globally Lipschitz function which agrees with it
on the relevant compact set and is smooth outside that set. Consequently, when the extensions are summed,
all the terms coming from the other pieces are differentiable at the points of
the piece under consideration, and hence they are invisible to
$\Omega_{1,\mu,\mathcal F}$ at infinitesimal scales, since here $\mathcal F$
contains the affine functions. This is the only patching mechanism needed in
the proof, and it is what allows the local residual statement to be converted
into a single Lipschitz function with big $\Omega_{1,\mu,\mathcal F}$ coefficients at small scales on the whole $E$.

\begin{teorema}\label{mainv2}
Suppose $\mu$ is a Radon measure on $\R^n$ such that
$$
\limsup_{r\to 0}\frac{\mu(B(x,2r))}{\mu(B(x,r))}<\infty
\qquad\mu\text{-almost everywhere.}
$$
Suppose $\mathcal F$ is a vector subspace of $L^1_{\mathrm{loc}}(\R^n)$ of
dimension $d$ containing affine functions. Suppose there exists a Borel set
$E$ such that for $\mu$-almost every $x\in E$ there are $D,A\geq 1$ and
$\mathfrak r>0$ such that
$$
  \Tan(\mu,x)\cap \bigcup_{k\in\N}\bigcap_{j\geq k}\mathfrak S_d(V(\mu,x),2^{-j},D,\mathfrak r,A)\neq \emptyset.
$$
Then there exists a Lipschitz function $f$ such that
$$
\limsup_{r\to 0}\Omega_{1,\mu,\mathcal F}(f;x,r)>0
\qquad\text{for $\mu$-almost every $x\in E$}.
$$
\end{teorema}

\begin{proof}
We divide the proof in several steps.

\textbf{Step I. The countable splitting.}
First of all, we observe that, thanks to \cref{coveringtangents}, we have
$$
\mu\left(E\setminus
\bigcup_{m_1,m_2,m_3\in\N}E_{m_1,m_3 2^{-m_2},m_3}\right)=0,
$$
where the sets $E_{D,\mathfrak r,A}$ are those introduced in \cref{coveringtangents}. Arguing as in Step 1 of the
proof of \cite[Proposition 4.1]{AlbertiMarchese}, and using the strong
locality principle \cite[Proposition 2.9]{AlbertiMarchese}, the locality of
tangents \cite[\S 2.3(4)]{Preiss1987GeometryDensities},
\cref{lemma:reg:ball:dens}, the Lebesgue differentiation theorem
\cite[Theorems 2.8.17, 2.9.8]{Federer1996GeometricTheory}, and
\cref{prop:doubllimit}, we can find countably many pairwise disjoint compact
sets $K_i\subseteq E$, continuous vector fields
$\mathscr E_i:=\{\mathfrak e_1^i,\ldots,\mathfrak e_n^i\}$, and measures
$$
\mu_i:=\mu\trace K_i
$$
such that
$$
\mu\trace E=\sum_{i\in\N}\mu_i,
$$
and, for every $i\in\N$, there are $A_i,D_i,k'\in\N$ and a rational number
$0<\mathfrak r_i<1/2(d+1)$ such that
\begin{enumerate}
    \item $\mu_i$ has compact support and $\diam\supp\mu_i\leq 1$;

    \item for $\mu_i$-almost every $x\in \R^n$ there holds
    $$
    \limsup_{r\to0}
    \frac{\mu_i(B(x,2r))}{\mu_i(B(x,r))}
    \leq A_i;
    $$

    \item the vector fields $\mathscr E_i$ coincide $\mu_i$-almost everywhere
    with the vector fields constructed in \cref{campi};

    \item for $\mu_i$-almost every $x\in\R^n$, we have
    $$
    \Tan(\mu_i,x)\cap
     \Tan(\mu,x)\cap \bigcap_{j\geq k'}\mathfrak S_d(V(\mu,x),2^{-j},D,\mathfrak r,A)\neq \emptyset.
    $$
\end{enumerate}

\medskip

\textbf{Step II. The absurd assumption.} Throughout this and the next step, $i\in\N$ will be fixed. We prove that the typical
function in $X_{\mu_i,\mathscr E_i}$ satisfies
$$
\limsup_{r\to0}
\Omega_{1,\mu_i,\mathcal F}(f;x,r)>0
\qquad\text{for $\mu_i$-almost every }x\in\R^n,
$$
where $X_{\mu_i,\mathscr E_i}$ is the metric space introduced in \cref{spaceoffunctions}.
For $\rho>0$, set
$$
U_\rho f(x):=
\sup_{0<s\leq\rho}
\Omega_{1,\mu_i,\mathcal F}(f;x,s).
$$
We recall that the maps $U_\rho$ are of Baire class $1$, as maps from
$X_{\mu_i,\mathscr E_i}$ into the space of Borel functions endowed with
convergence in $\mu_i$-measure on compact sets, see \cref{operators}. Therefore the set
$\mathscr R_i^q$ of functions which are continuity points of $U_{2^{-q}}$, is residual for
every $q\in\N$ in $X_{\mu_i,\mathscr E_i}$.
Let $\mathscr R_i:=\cap_{q\in \N} \mathscr R_i^q$, notice that $\mathscr R_i$ is still residual in $X_{\mu_i,\mathscr E_i}$ and pick a function $f\in\mathscr R_i$. Suppose by contradiction that
$$
B_f:=\Big\{x\in\R^n:\limsup_{r\to0}\Omega_{1,\mu_i,\mathcal F}(f;x,r)=0\Big\}$$
has positive $\mu_i$-measure.

\textbf{Step III. The residuality of bad functions for $\mu_i$.}
First of all let 
$$
\mathfrak d_\rho(x):=
\sup_{\ell\in\N}
\frac{\mu_i(B(x,2^{-\ell+1}\rho))}
{\mu_i(B(x,2^{-\ell-2}\rho))}.
$$
By \cref{lemma:reg:ball:dens}, the function $\mathfrak d_\rho$ is Borel. For
every $\Lambda>0$ and $j\in\N$, let us define
$$
E_{\Lambda,j}:=
\{x\in\R^n:\mathfrak d_{2^{-j}}(x)<\Lambda\}.
$$
Thanks to item 2 of step I, it is immediate to see that 
$$\mu_i(\R^n\setminus \bigcup_{\Lambda\in \N}\bigcup_{j_0\in \N} \bigcap_{j\geq j_0}E_{\Lambda,j})=0.$$
Moreover, by item 4 of Step I, for $\mu_i$-almost every $x\in\R^n$ we have
$$
x\in \mathfrak G_{\mathfrak r_i,D_i,2^{-k},A_i}
\qquad\text{for every }k\geq k',
$$
where $\mathfrak G_{\mathfrak r_i,D_i,2^{-k},A_i}$ denotes the set of points
whose tangents satisfy the condition appearing in the statement of
\cref{perturbationlemma} with parameters
$\mathfrak r_i,D_i,2^{-k},A_i$.
Since $B_f$ has positive $\mu_i$-measure, by the preceding observations and
inner regularity we can find $\Lambda\in\N$, $j_0\in\N$ and a
compact set $C\subseteq B_f$ with $\mu_i(C)>0$ such that $C\subseteq E_{\Lambda,j}\cap \bigcap_{k\geq k'}\mathfrak G_{\mathfrak r_i,D_i,2^{-k},A_i}$ for every $j\geq j_0$ and thanks to Severini-Egoroff theorem we can also assume that 
\begin{equation}
    \begin{split}
        \Omega_{1,\mu_i,\mathcal F}(f;x,2^{-j-1})\leq \frac{\mathfrak r_i}{160000\Lambda(d+1)^2D_i}\qquad\text{for every }x\in C\text{ and every }j\geq j_0.
        \label{stimaomegasuC}
    \end{split}
\end{equation}
Fix $k_0\geq k'$ in such a way that 
$$2^{-k_0}\leq \frac{\mathfrak r_i}{40000(d+1)^2D_i}.$$
Applying \cref{perturbationlemma} to the
measure $\mu_i$, to the compact set $C$ and
with
$$
A=A_i,\qquad
D=D_i,\qquad
\mathfrak r=\mathfrak r_i,\qquad
\xi=2^{-k}\text{ with $k\geq k_0$},\qquad
r=2^{-j-1}\text{ with $j\geq j_0$},
$$
we obtain a compact set $C_{j,k}\subseteq C$ and a function
$f_{j,k}\in X_{\mu_i,\mathscr E_i}$ such that
\begin{equation}
    \mu_i(C_{j,k})
    \geq
    \frac{1}{4(d+1)A_i^2}\mu_i(C),
    \label{eq:measureCjkMainV2}
\end{equation}
$$
\|f-f_{j,k}\|_\infty
\leq
2(\|f\|_\infty+d+1)\varsigma(2^{-k}),
$$
and, for every $y\in C_{j,k}$,
\begin{equation}
    \sup_{0<s<2^{-j-1}}
    \sup_{z\in \mathfrak G_{\mathfrak r_i,D_i,2^{-k},A_i}:\,y\in B(z,s/5)}
    \Omega_{1,\mu_i,\mathcal F}(f_{j,k};z,s)
    \geq  \frac{\mathfrak r_i}{40000(d+1)^2D_i}.
    \label{eq:offcenterMainV2}
\end{equation}

We now pass from the off-center estimate to a centered one. If $y\in B(z,s/5)$,
then
$$
B(y,s/2)\subseteq B(z,s)\subseteq B(y,2s).
$$
Therefore, for every Lipschitz function $h$ and every $A\in\mathcal F$,
$$
\fint_{B(y,2s)}
\frac{|h-A|}{2s}\,d\mu_i
\geq
\frac{\mu_i(B(y,s/2))}{2\mu_i(B(y,2s))}
\fint_{B(z,s)}
\frac{|h-A|}{s}\,d\mu_i.
$$
Taking the infimum over $A\in\mathcal F$, we obtain
\begin{equation}
    \Omega_{1,\mu_i,\mathcal F}(h;y,2s)
    \geq
    \frac{\mu_i(B(y,s/2))}{2\mu_i(B(y,2s))}
    \Omega_{1,\mu_i,\mathcal F}(h;z,s).
    \label{eq:centerChangeMainV2}
\end{equation}
Let $y\in C$ and $0<s<2^{-j-1}$. Choose $\ell\in\N$ such that $2^{-j-\ell-1}<s\leq 2^{-j-\ell}$. Then
$$
\frac{\mu_i(B(y,2s))}{\mu_i(B(y,s/2))}
\leq
\frac{\mu_i(B(y,2^{-j-\ell+1}))}
{\mu_i(B(y,2^{-j-\ell-2}))}
\leq
\mathfrak d_{2^{-j}}(y)
\leq
\Lambda,
$$
because $C\subseteq E_{\Lambda,j}$. Hence, by
\eqref{eq:centerChangeMainV2} and \eqref{eq:offcenterMainV2}, we get
\begin{equation}
    U_{2^{-j}}f_{j,k}(y)
    \geq\frac{\mathfrak r_i}{80000\Lambda(d+1)^2D_i}
    \qquad\text{for every }y\in C_{j,k}.
    \label{eq:UperturbedMainV2}
\end{equation}
Let us observe that
$$
\lim_{k\to \infty}\|f-f_{j,k}\|_\infty
\leq\lim_{k\to \infty}
2(\|f\|_\infty+d+1)\varsigma(2^{-k})=0,
$$
and since $f$ is a continuity point of $U_{2^{-j}}$, we have
\begin{equation}
    \lim_{k\to \infty}U_{2^{-j}}f_{j,k}(w)=U_{2^{-j}}f(w)\qquad\text{for $\mu_i$-almost every $w\in C$.}
    \label{convergenzadegliU}
\end{equation}
Thanks to our choice of $C$ and of $j_0$, see \eqref{stimaomegasuC}, this implies in particular that for every $j\geq j_0$ we have 
$$\limsup_{k\to \infty}U_{2^{-j}}f_{j,k}(w)\leq \frac{\mathfrak r_i}{160000\Lambda(d+1)^2D_i}\qquad\text{for every }w\in C.$$
In particular, applying Severini-Egoroff, we infer that there must exist $k\in \N$ and a Borel set $E\subseteq C$ with 
$$\mu_i(C\setminus E)\leq \frac{1}{100(d+1)A_i^2}\mu_i(C)\quad\text{and}\quad U_{2^{-j}}f_{j,k}(w)\leq \frac{\mathfrak r_i}{100000\Lambda(d+1)^2D_i}\quad\text{for every }w\in E.$$
This however is in contradiction with the existence of $C_{j,k}$. This finally implies $\mu(B_f)=0$ for the typical $f\in X_{\mu_i,\mathscr E_i}$ and hence, for the typical Lipschitz function in $X_{\mu_i,\mathscr E_i}$ there holds 
$$\limsup_{r\to0}\Omega_{1,\mu_i,\mathcal F}(f;x,r)>0\qquad\text{for $\mu_i$-almost every }x\in \R^n.$$

\medskip

\textbf{Step IV. Conclusion of the proof.} In this final step we patch together the counterexamples we obtained from the previous step. For every $i\in\N$, we choose $f_i\in X_{\mu_i,\mathscr E_i}$ such that
$$
\limsup_{r\to0}
\Omega_{1,\mu_i,\mathcal F}(f_i;x,r)>0
\qquad\text{for $\mu_i$-almost every }x\in\R^n.
$$
Replacing $f_i$ by a positive multiple of $f_i-f_i(x_i)$, for some
$x_i\in\supp\mu_i$, we may assume that $f_i$ is $1$-Lipschitz and
$\|f_i\|_\infty\leq 2$ on $\supp\mu_i$, as $\diam \supp\mu_i\leq 1$ by item 1 in step I. By means of
\cite[Proposition 8.4]{AlbertiMarchese}, for every $i\in\N$ we can find a
$2$-Lipschitz function $g_i$ such that
\begin{enumerate}
    \item[(a)] $g_i$ agrees with $f_i$ on $\supp\mu_i$ and is smooth on
    $\R^n\setminus K_i$;
    \item[(b)] $\|g_i\|_\infty\leq 3$;
    \item[(c)] $|f_i(x)-g_i(x)|\leq \dist(x,\supp\mu_i)^{10}$ for every $x\in\R^n$.
\end{enumerate}
Let
$$ g:=\sum_{i\in\N}2^{-i}g_i.$$
Then $g$ is Lipschitz. We claim that $g$ satisfies the conclusion.

Fix $i_0\in\N$. Since $g_{i_0}=f_{i_0}$ on $\supp\mu_{i_0}$, we have
\begin{equation}
    \limsup_{r\to0}
    \Omega_{1,\mu_{i_0},\mathcal F}(g_{i_0};x,r)>0
    \qquad\text{for $\mu_{i_0}$-almost every }x\in\R^n.
    \label{eq:positiveForMuiMainV2}
\end{equation}
Moreover, since $\mu_{i_0}=\mu\trace K_{i_0}$, the Lebesgue differentiation
theorem gives
$$
\lim_{r\to0}
\frac{\mu_{i_0}(B(x,r))}{\mu(B(x,r))}
=
1
\qquad\text{for $\mu_{i_0}$-almost every }x\in\R^n.
$$
For every such point $x$, every $r>0$, and every $A\in\mathcal F$, we have
$$
\fint_{B(x,r)}
\frac{|g_{i_0}-A|}{r}\,d\mu
\geq
\frac{\mu_{i_0}(B(x,r))}{\mu(B(x,r))}
\fint_{B(x,r)}
\frac{|g_{i_0}-A|}{r}\,d\mu_{i_0}.
$$
Taking the infimum over $A\in\mathcal F$ and using
\eqref{eq:positiveForMuiMainV2}, we get
\begin{equation}
    L(x):=
    \limsup_{r\to0}
    \Omega_{1,\mu,\mathcal F}(g_{i_0};x,r)>0
    \qquad\text{for $\mu_{i_0}$-almost every }x\in\R^n.
    \label{eq:positiveForMuMainV2}
\end{equation}
Arguing as in \cite[Step 4, proof of Theorem 4.1]{AlbertiMarchese}, for every
$x\in\supp\mu_{i_0}$ the function
$$
\tilde g:=\sum_{i\neq i_0}2^{-i}g_i
$$
is differentiable at $x$. Fix a point $x\in\supp\mu_{i_0}$ where
\eqref{eq:positiveForMuMainV2} holds and where $\tilde g$ is differentiable.
Let
$$
\ell_x(y):=\tilde g(x)+D\tilde g(x)(y-x).
$$
Choose $\rho=\rho(x)>0$ such that
$$
|\tilde g(y)-\ell_x(y)|
\leq
\frac{2^{-i_0}L(x)}{100}|y-x|
\qquad\text{for every }y\in B(x,\rho).
$$
Then, for every $0<s<\rho$ and every $A\in\mathcal F$,
\begin{equation}
    \begin{split}
    \fint_{B(x,s)}
    \frac{|g-A|}{s}\,d\mu
    &\geq
    \fint_{B(x,s)}
    \frac{|2^{-i_0}g_{i_0}+\ell_x-A|}{s}\,d\mu
    -
    \frac{2^{-i_0}L(x)}{100}.
    \end{split}
    \label{eq:subtractDifferentiableMainV2}
\end{equation}
Since $\ell_x$ is affine and $\mathcal F$ contains the affine functions, we
have $A-\ell_x\in\mathcal F$ whenever $A\in\mathcal F$. Taking the infimum in
\eqref{eq:subtractDifferentiableMainV2} over $A\in\mathcal F$, and using that
$\mathcal F$ is a vector space, we obtain
$$
\Omega_{1,\mu,\mathcal F}(g;x,s)
\geq
2^{-i_0}
\Omega_{1,\mu,\mathcal F}(g_{i_0};x,s)
-
\frac{2^{-i_0}L(x)}{100}.
$$
Taking the limsup as $s\to0$, we infer that
$$
\limsup_{s\to0}
\Omega_{1,\mu,\mathcal F}(g;x,s)
\geq
\frac{99}{100}2^{-i_0}L(x)>0.
$$
Hence
$$
\limsup_{r\to0}
\Omega_{1,\mu,\mathcal F}(g;x,r)>0
\qquad\text{for $\mu_{i_0}$-almost every }x\in\R^n,
$$
and since $i_0\in\N$ was arbitrary the conclusion follows by choosing $f:=g$. This concludes the proof.
\end{proof}

\subsection{Consequences of \cref{mainv2}}

In this section we use \cref{mainv2} to prove some rectifiability criteria. First however, we prove that under certain regularity assumptions for the measure, we fall in the assumptions of \cref{mainv2}.

\begin{proposizione}\label{upperdensityunrectifiableimpliestangentcondition}
Let $\mu$ be a Radon measure on $\R^n$ such that
$$
\limsup_{r\to0}\frac{\mu(B(x,2r))}{\mu(B(x,r))}<\infty
\qquad\text{for $\mu$-almost every }x.
$$
Let $E$ be a Borel set. Assume that for $\mu$-almost every $x\in E$ there
exists $\beta(x)>\dim V(\mu,x)$ such that
$$
\liminf_{R\to0}\inf_{0<r\leq R}
\frac{\mu(B(x,R))r^{\beta(x)}}{\mu(B(x,r))R^{\beta(x)}}>0.
$$
Then, for $\mu$-almost every $x\in E$, there are $D,A\geq1$ and
$\mathfrak r>0$ such that
$$
\Tan(\mu,x)\subseteq 
\bigcup_{k\in\Z}\bigcap_{j\geq k}
\mathfrak S_d(V(\mu,x),2^{-j},D,\mathfrak r,A).
$$
\end{proposizione}

\begin{proof}
By Lebesgue differentiation theorem and the strong locality principle of the decomposability bundle, see \cite[Proposition 2.9]{AlbertiMarchese}, we can assume from now that $\mu$ is a compactly supported finite measure, and that
\begin{enumerate}
\item there exists $q=0,\ldots, n$ such that $\dim V(\mu,x)=q$ for $\mu$-almost every $x\in \R^n$;
\item there is $A\in \N$ such that  $\displaystyle{\limsup_{s\to 0}\frac{\mu(B(x,2s))}{\mu(B(x,s))}\leq \tfrac{A}{2}}$ for $\mu$-almost every $x\in \R^n$. 
\item for $\mu$-almost every $x\in E$ there
exists $\beta(x)>q$ such that
$\displaystyle{\liminf_{R\to0}\inf_{0<r\leq R}
\frac{\mu(B(x,R))r^{\beta(x)}}{\mu(B(x,r))R^{\beta(x)}}>0}$.
\end{enumerate}
Notice that thanks to \cref{lemma:reg:ball:dens}, the map 
$$x\mapsto \liminf_{R\to0}\inf_{0<r\leq R}
\frac{\mu(B(x,R))r^{q+\delta}}{\mu(B(x,r))R^{q+\delta}},$$
is Borel for every $\delta>0$. This shows thanks to item 3, and Lebesgue differentiability theorem that we can reduce to assuming that there holds
\begin{itemize}
    \item[3a.] there are $\delta>0$ and $C>0$ such that for $\mu$-almost every $x\in E$ there holds
$$\liminf_{R\to0}\inf_{0<r\leq R}
\frac{\mu(B(x,R))r^{q+\delta}}{\mu(B(x,r))R^{q+\delta}}>\frac{C}{2}.$$
\end{itemize}
By Severini-Egoroff and items 1. and 3. above for every $\varepsilon>0$ there exists a Borel set $E_\varepsilon$ and $r_\varepsilon>0$ such that
$\mu(\R^n\setminus E_\varepsilon)\leq \varepsilon\mu(\R^n)$ and for every $0<R<r_\varepsilon$ we have
$$\frac{\mu(B(x,R))}{\mu(B(x,r))}\geq C\Big(\frac{R}{r}\Big)^{q+\delta}\qquad\text{for every }0<r<R\text{ and }\qquad\frac{\mu(B(x,2R))}{\mu(B(x,R))}\leq A.$$
Let $r_j$ be any infinitesimal sequence, let $x\in E_\varepsilon$ and let 
$$
\nu_j:=\frac{1}{\mu(B(x,r_j))}T_{x,r_j}\mu.
$$
By the doubling assumption, and thanks to \cite[Proposition 1.12]{Preiss1987GeometryDensities} passing to a subsequence, we may assume that there exists a non-zero Radon measure $\nu$ such that 
$\nu_{j,x}\rightharpoonup\nu$. In particular, by definition we have 
$\nu_x\in\Tan(\mu,x)$. By \cref{th:plit} we know that there exists a Radon measure $\eta\in \mathscr{M}(V(\mu,x)^\perp)$ such that $\nu_x=\Haus^q\trace V(\mu,x)\otimes \eta_x$ and the measure $\eta$ is doubling, i.e.:
$$\eta_x(B(z,2r))\leq 2^{-q}A^{\log_2(2n)+1}\eta_x(B(z,r))\qquad \text{ for every $z\in \supp\eta$ and every $r>0$}.$$

We claim that $\nu$ satisfies a reverse doubling estimate. Let $p\in\supp\nu$ and $0<r<R$.  Indeed, let $p\in\supp\nu$ and $s>0$. By locality of tangents, see \cite[\S 2.3(4)]{Preiss1987GeometryDensities}, we know that $\nu_x\in \Tan(\mu\trace E_\varepsilon,x)$ and thus, thanks to \cref{puntivicini} there exist points $y_j\in E_\varepsilon$ such that $(y_j-x)/r_{j,x}\to p$. Thus, for $\Leb^1$-almost every $0<r\leq \frac{1-\varepsilon}{1+\varepsilon}R$, Portmanteau's theorem implies that
\begin{equation}
    \begin{split}
        \frac{\nu(B(p,R))}{\nu(B(p,r))}=&\lim_{j\to \infty}\frac{\nu_j(B(p,R))}{\nu_j(B(p,r))}=\lim_{j\to \infty}\frac{\mu(B(x+r_jp,Rr_j))}{\mu(B(x+r_jp,rr_j))}\\
        \geq& \lim_{j\to \infty}\frac{\mu(B(y_j,(1-\varepsilon)Rr_j))}{\mu(B(y_j+r_jp,(1+\varepsilon)rr_j))}\geq C\Big(\frac{R}{r}\Big)^{q+\delta}\Big(\frac{1-\varepsilon}{1+\varepsilon}\Big)^{q+\delta}.
    \end{split}
\end{equation}
The arbitrariness of $\varepsilon$ shows that for every $p\in \supp\nu$, every $R>0$ and $\Leb^1$-almost every $r\leq R$ there holds 
$$\frac{\nu(B(p,R))}{\nu(B(p,r))}\geq C\Big(\frac{R}{r}\Big)^{q+\delta}.$$
The conclusion for every $r\leq R$ follows by monotonicity of the measure.

We now prove that
$\eta$ satisfies a reverse doubling estimate with positive exponent $\delta$. Let $z\in\supp\eta$ and $0<r<R$. Since
$$B_{V(\mu,x)}(0,r)\times B_{V(\mu,x)^\perp}(z,r)\subseteq B((0,z),2r)\quad\text{and}\quad B((0,z),\tfrac{R}{2})\subseteq B_{V(\mu,x)}(0,\tfrac{R}{2})\times B_{V(\mu,x)^\perp}(z,\tfrac{R}{2}),$$
we get
$$\omega_q r^q\eta(B(z,r))\leq\nu(B((0,z),2r))\qquad\text{and}\qquad\nu(B((0,z),R/2))\leq\omega_q R^q\eta(B(z,R)).$$
In particular this implies that
$$\frac{\eta(B(z,R))}{\eta(B(z,r))}\geq \frac{C}{2^{2(q+\delta)}}\Big(\frac{R}{r}\Big)^\delta.$$
We now choose $\mathfrak r>0$ small enough, depending only on the doubling
constant of $\eta$, $\delta$, $q$ and on $d$, such that 
$$\mathfrak r:=\frac{1}{4}\Big(\frac{C}{2(d+1)}A^{-2(\log(2n)+1)}\Big)^{\max\{1,1/\delta\}},$$
and the following holds. For every $z\in\supp\eta$ and every $\xi>0$, there
exist points
$$
z_1,\ldots,z_{d+1}\in\supp\eta\cap B(z,\xi/8)
$$
such that
$$
|z_\ell-z_{\ell'}|\geq 4\mathfrak r\xi
\qquad\text{whenever }\ell\neq\ell'.
$$
Indeed, if no such points existed, then a maximal $4\mathfrak r\xi$-separated
subset of $\supp\eta\cap B(z,\xi/8)$ would have at most $N\leq d$ points $w_j$. Hence
$\supp\eta\cap B(z,\xi/8)$ would be covered by at most $N$ balls of radius
$4\mathfrak r\xi$. For each of the points $w_j$ we have, by the reverse doubling estimate for
$\eta$ applied with radii $4\mathfrak r\xi$ and $\xi/4$,
$$
\eta(B(w_j,4\mathfrak r\xi))
\leq
\frac{2^{2(q+\delta)}}{C}(16\mathfrak r)^\delta
\eta(B(w_j,\xi/4)).
$$
Since $w_j\in B(z,\xi/8)$, we also have $B(w_j,\xi/4)\subseteq B(z,\xi/2)$, and hence, by doubling twice, we infer that
$$\eta(B(w_j,\xi/4))\leq \eta(B(z,\xi/2))\leq 2^{-2q}A^{2(\log_2(2n)+1)}\eta(B(z,\xi/8)).$$
Therefore
\begin{equation}
\begin{split}
\eta(B(z,\xi/8))\leq&
\sum_{j=1}^{N}\eta(B(w_j,4\mathfrak r\xi))\leq \frac{2^{2(q+\delta)}}{C}(16\mathfrak r)^\delta\sum_{j=1}^N
\eta(B(w_j,\xi/4))\\
\leq& d\frac{2^{2(q+\delta)}}{C}(16\mathfrak r)^\delta 2^{-2q}A^{2(\log_2(2n)+1)}\eta(B(z,\xi/8)),
\end{split}
\end{equation}
which results in a contradiction thanks to the choice of $\mathfrak r$.

We now check the mass comparison required in the definition of
$\mathfrak S_d$. Since $\eta$ is doubling and
$z_\ell\in B(z,\xi/8)$, we have $B(z,50\xi)\subseteq B(z_\ell,51\xi)$. Choose $N_0\in\N$ such that
$2^{N_0}\frac{\mathfrak r}{32}\geq 51$, observe that $N_0$ is independent on $\xi$ and notice that iterating the doubling inequality $N_0$ times, we get
$$
\eta(B(z,50\xi))
\leq
\eta(B(z_\ell,51\xi))
\leq (2^{-qN_0}A^{(\log_2(2n)+1)N_0})\eta(B(z_\ell,\mathfrak r\xi/32)).
$$
This finally implies that
$$
\nu\in
\mathfrak S_d(V(\mu,x),\xi,D,\mathfrak r,A^{(\log_2(2n)+1)N_0}))
\qquad\text{for every }\xi>0.
$$
Since the infinitesimal sequence $r_j$ was arbitrary, the same conclusion holds for
every $\nu\in\Tan(\mu,x)$. Hence
$$
\Tan(\mu,x)\subseteq
\bigcup_{k\in\Z}\bigcap_{j\geq k}
\mathfrak S_d(V(\mu,x),2^{-j},D,\mathfrak r,A^{(\log_2(2n)+1)N_0})),
$$
and this proves the claim for $\mu$-almost every $x\in E_\varepsilon$. Since
$\varepsilon>0$ was arbitrary, the proposition follows.
\end{proof}

From the above result we immediately obtain the following consequence.

\begin{teorema}
    \label{rectthm}
Let $\mu$ be a Radon measure on $\R^n$ such that for $\mu$-almost every $x\in \R^n$ the following hold
\begin{enumerate}
    \item $\displaystyle{\limsup_{r\to0}\frac{\mu(B(x,2r))}{\mu(B(x,r))}<\infty}$;
    \item there
exists $\beta(x)>\dim V(\mu,x)$ such that
$\displaystyle{\liminf_{R\to0}\inf_{0<r\leq R}
\frac{\mu(B(x,R))r^{\beta(x)}}{\mu(B(x,r))R^{\beta(x)}}>0.}$
\end{enumerate}
Then, if $\mathcal F$ is a vector subspace of $L^1_{\mathrm{loc}}(\R^n)$ of dimension $d$ containing affine function, then there exists a $1$-Lipschtiz function $f:\R^n\to \R$ such that 
$$\limsup_{r\to 0}\Omega_{1,\mu,\mathcal F}(f;x,r)>0\qquad\text{ for $\mu$-almost every $x\in K$.}$$
\end{teorema}

\begin{proof}
   Thanks to \cref{upperdensityunrectifiableimpliestangentcondition}, the hypothesis of \cref{mainv2} are satisfied and thus the conclusion immediately falls.  
\end{proof}

We are finally ready to state our rectifiability result and some consequences.

\begin{corollario}
    \label{rectthm2}
Suppose $\mu$ is a Radon measure on $\R^n$ and let $\alpha>0$. Suppose that 
$$
\Theta^{\alpha,*}(\mu,x)>0
\qquad\text{and}\qquad
\limsup_{r\to0}\frac{\mu(B(x,2r))}{\mu(B(x,r))}<\infty
$$
for $\mu$-almost every $x\in \R^n$. Suppose moreover that
$$
\liminf_{R\to0}\inf_{0<r\leq R}
\frac{\mu(B(x,R))r^\alpha}{\mu(B(x,r))R^\alpha}>0
\qquad\text{for $\mu$-almost every }x\in\R^n.
$$
Suppose that $\mathcal F$ is a vector subspace of $L^1_{\mathrm{loc}}(\R^n)$ of dimension $d$ and assume that for every $1$-Lipschitz function $f:\R^n\to \R$ we have 
$$
\limsup_{r\to 0}\Omega_{1,\mu,\mathcal F}(f;x,r)=0
\qquad\text{for $\mu$-almost every }x\in \R^n.
$$
Then $\alpha\in \N$ and $\mu$ is $\alpha$-rectifiable.
\end{corollario}

\begin{proof} Let us denote $\widetilde{\mathcal {F}}:=\mathscr A+\mathcal F$, the sum of affine functions and $\mathcal F$ as vector subspaces of $L^1_{loc}(\R^n)$.
As a first step, let us show that 
$\Theta^{\alpha,*}(\mu,x)<\infty$ for $\mu$-almost every $x\in \R^n$. Indeed, fix $x$ such that
$$
\liminf_{R\to0}\inf_{0<r\leq R}
\frac{\mu(B(x,R))r^\alpha}{\mu(B(x,r))R^\alpha}>0.
$$
Then, there are $c(x)>0$ and $R_x>0$ such that, for every $0<r\leq R<R_x$, there holds
$$\mu(B(x,r))\leq c(x)^{-1}\left(\frac{r}{R}\right)^\alpha\mu(B(x,R)),$$
and in particular $\Theta^{\alpha,*}(\mu,x)<\infty$. This shows that $0<\Theta^{\alpha,*}(\mu,x)<\infty$
for $\mu$-almost every $x\in\R^n$.

Let $E_{j}:=\{x\in E:j^{-1}\leq \Theta^{\alpha,*}(\mu,x)\leq j\}$. Then \cite[\S 2.10.19(1), (3)]{Federer1996GeometricTheory} and the Lebesgue differentiation theorem imply that 
$$\mu\trace E_j\leq 2^\alpha j\Haus^\alpha\trace E_j$$
and that for every $B$ measurable and $V$ open containing $B$ we have 
$$\Haus^\alpha\trace E_j(B)\leq j\mu\trace E_j(V).$$
Thus $\Haus^\alpha\trace E_j$ is a $\sigma$-finite measure and $\mu\trace E_j$ and $\Haus^\alpha\trace E_j$ are mutually absolutely continuous. This implies that since $\mu$ is Radon, $\Haus^\alpha\trace E_j$ is Radon and hence by Radon-Nikodym decomposition theorem, there exists $\rho_j\in L^1(\mu\trace E_j)$ such that 
$$\mu\trace E_j=\rho_j\Haus^\alpha\trace E_j.$$
By \cref{abscontinuous} we infer that $\dim V(\mu,x)\leq \alpha$ for $\mu$-almost every $x\in E_j$. By the arbitrariness of $j\in\N$ this shows in particular that 
$$\dim(V(\mu,x))\leq \alpha\qquad\text{for $\mu$-almost every }x\in E.$$

Assume first that $\alpha\notin\N$, then 
\begin{equation}
    \dim(V(\mu,x))< \alpha\qquad\text{for $\mu$-almost every }x\in E.
    \label{strictinequality}
\end{equation}
holds trivially. Thus the assumptions of \cref{rectthm} are satisfied with $\beta(x):=\alpha$.
Hence there exists a $1$-Lipschitz function $f:\R^n\to\R$ such that
$$
\limsup_{r\to0}\Omega_{1,\mu,\widetilde{\mathcal{F}}}(f;x,r)>0
\qquad\text{for $\mu$-almost every }x\in\R^n,
$$
contradicting the hypothesis. Therefore $\alpha\in\N$.

It remains to prove that $\mu$ is $\alpha$-rectifiable. We can assume without loss of generality that $\mu$ is purely unrectifiable. This is an immediate consequence of the fact that all the quantities above are stable under restricting $\mu$ to subsets thanks to Lebesgue differentiability theorem. Since $\mu$ is $\alpha$-purely unrectifiable, \cref{criteriodirett}
gives
$$
\dim V(\mu,x)<\alpha
\qquad\text{for $\mu$-almost every }x\in \R^n.
$$
Therefore \cref{rectthm}, applied to the measure $\mu$ with
$\beta(x):=\alpha$, gives a $1$-Lipschitz function $f:\R^n\to\R$ such that
$$
\limsup_{r\to0}\Omega_{1,\mu,\mathcal F}(f;x,r)\geq \limsup_{r\to0}\Omega_{1,\mu,\widetilde{\mathcal F}}(f;x,r)>0
\qquad\text{for $\mu$-almost every }x\in \R^n,
$$
thus resulting in a contradiction and in particular $\mu$ is $\alpha$-rectifiable.
\end{proof}

\begin{corollario}
    \label{rectthm3}
Suppose $\mu$ is a Radon measure on $\R^n$ and let $\alpha>0$. Suppose that
$$
0<\Theta^\alpha_*(\mu,x)\leq \Theta^{\alpha,*}(\mu,x)<\infty
\qquad\text{for $\mu$-almost every }x\in\R^n.
$$
Suppose that $\mathcal F$ is a vector subspace of $L^1_{\mathrm{loc}}(\R^n)$ of dimension $d$ and assume that for every $1$-Lipschitz function $f:\R^n\to \R$ we have 
$$
\limsup_{r\to 0}\Omega_{1,\mu,\mathcal F}(f;x,r)=0
\qquad\text{for $\mu$-almost every }x\in \R^n.
$$
Then $\alpha\in \N$ and $\mu$ is $\alpha$-rectifiable.
\end{corollario}

\begin{proof}
We check that the assumptions of \cref{rectthm2} are satisfied. Indeed
$$\limsup_{r\to 0}\frac{\mu(B(x,2r))}{\mu(B(x,r))}\leq\frac{\Theta^{\alpha,*}(\mu,x)}{2^\alpha\Theta^{\alpha}_*(\mu,x)}.$$
Now, notice that for $\mu$-almost every $x\in \R^n$ there exists $r_x>0$ such that for every $0<r<r_x$, we have 
$$\frac{\Theta^\alpha_*(\mu,x)}{2}\leq \frac{\mu(B(x,r))}{r^\alpha}\leq 2\Theta^{\alpha,*}(\mu,x).$$
This shows in particular that for every $0<r\leq R<r_x$ we have
$$
\frac{\mu(B(x,R))r^\alpha}{\mu(B(x,r))R^\alpha}
\geq\frac{1}{4}\frac{\Theta^{\alpha}_*(\mu,x)}{\Theta^{\alpha,*}(\mu,x)}.
$$
Thus, all the assumptions of \cref{rectthm2} hold. Applying \cref{rectthm2}, we
obtain that $\alpha\in\N$ and that $\mu$ is $\alpha$-rectifiable.
\end{proof}



\section{Quantitative results}
\label{sectionquantitativeresults}

The purpose of this section is to prove \cref{teoremaintro}. The main difficulty is to turn the qualitative mechanism developed in \cref{sezionequalitativa} into a statement which is uniform across locations and scales. In this section we develop quantitative arguments that are a quantitative analogue of those that can be found in \cref{sezionequalitativa}.

The central notion is that of a cube which is quantitatively invariant along a given direction. Roughly speaking, invariance means that near every point of the cube one can find a curve directed inside a prescribed cone which remains close to the support of the measure for most of its length. The maximal number of quantitatively independent invariant directions plays the role of a quantitative decomposability bundle. Invariance in $j$ independent directions forces the support to be dense, in the right coarse sense, along the $j$-dimensional affine planes generated by those directions; in particular, $k$ independent invariant directions yield the flatness needed for the bilateral weak geometric lemma.

It remains to show that the cubes with fewer than $k$ invariant directions form Carleson families. Arguing by contradiction, we extract from any non-Carleson family arbitrarily long and well-separated towers of cubes. Failure of invariance on these towers is converted into a quantitative loss of width: we construct nested open sets which retain most of the relevant $\mu$-mass but have small intersection with curves directed inside prescribed cones. A quantitative version of the joining-cones argument of G. Alberti, M. Csörnyei and D. Preiss then upgrades the estimates obtained for finitely many narrow cones to the wider transverse cones needed in the perturbation argument.

Finally, the resulting width functions are used first to cancel the derivatives of a Lipschitz function in the transverse directions and then to insert perturbations which cannot be approximated by the finite-dimensional space $\mathcal F$. Iterating this construction along the tower produces a single Lipschitz function whose $\Omega$ coefficients are large on arbitrarily many separated layers, contradicting the Carleson estimate in the GWALA condition. Consequently, all the strata with fewer than $k$ invariant directions are Carleson, whereas the remaining cubes satisfy the required bilateral flatness estimate. This proves the fixed-threshold OUWGL condition and hence uniform rectifiability.

\subsection{The WALA condition and its consequences}

First of all, let us introduce the condition at the center of this work.

\begin{definizione}\label{defWALA}Let $\mu$ be an $\alpha$-dimensional AD regular measure on $\R^n$ and let us fix a vector space $\mathcal F$ of $L^1_{\mathrm{loc}}(\R^n)$ of dimension $d$ and let $q\in [1,\infty]$.
    We say that $\mu$ satisfies the $\text{GWALA}_q$ condition if for every $1$-Lipschitz function $f$ we have that 
    \begin{equation}
           \mathscr{B}(f,\varepsilon):=\mathscr{B}_{q,\mu,\mathcal F}(f,\varepsilon):=\{(x,r)\in \supp(\mu)\times(0,\infty):\Omega_{q,\mu,\mathcal F}(f;x,r)> \varepsilon\},
           \label{BADsetWALA}
    \end{equation}   
    is a Carleson set, whose Carleson norm depends on $\varepsilon$ but not on $f$. More precisely we will say that $\mu$ satisfies the $\text{GWALA}_q$ condition with constant $\{C(t):t>0\}$ if for every Lipschitz function $f:\R^n\to \R$, every $x\in \supp\mu$ and $R>0$ we have 
    $$\int_0^R\int_{B(x,R)}\mathbb{1}_{\mathscr{B}(f,\varepsilon)}(y,s)d\mu(y)\frac{ds}{s}\leq \Lip(f)C(\varepsilon)\mu(B(x,R)).$$
    When $\mathcal F$ coincides with the affine functions, then the above condition is typically called simply WALA condition, see \cite[Part I, Chapter 2, Definition 2.47]{DavidSemmes}. 
\end{definizione}

\begin{osservazione}\label{osservazioneriduzionewala}
    Let us observe that the $\text{GWALA}_1$ condition is the weakest condition. Indeed, thanks to \cref{omega_monotone}, we know that for every Lipschitz function $f$ we have
    $$\mathscr{B}_{q,\mu,\mathcal F}(f,\varepsilon)\subseteq \mathscr{B}_{1,\mu,\mathcal F}(f,\varepsilon).$$
Thus, if we prove that the $\text{GWALA}_1$ condition implies uniform rectifiability, we would prove that $\text{GWALA}_q$, for every $q\in [1,\infty]$ implies uniform rectifiability. For this reason, from now on we will drop the dependence on $q$ of the condition $\text{GWALA}$ and we will only study the case $q=1$. 
\end{osservazione}

We first collect the first low-hanging fruit ripened in Section \cref{sezionequalitativa}.

\begin{proposizione}\label{prop:rettificabilitya}
    Suppose $\alpha>0$ and let $\mu$ be an $\alpha$-dimensional AD-regular measure with regularity constant $D$ satisfying the GWALA condition. Then $\alpha\in \N$ and $\mu$ is $\alpha$-rectifiable.
\end{proposizione}

\begin{proof}
    First, we notice that by \cref{propcompatibilityOmega}, the coefficients $\Omega_{1,\mu,\mathcal F}$ are compatible in the sense of \cref{defscalecompatible}. This implies by \cref{lemcontinuouspackinglimitzero} that for $\mu$-almost every $x\in \R^n$ and every Lipschitz function $f$ we have 
    $$\lim_{r\to 0}\Omega_{1,\mu,\mathcal F}(f;x,r)=0.$$
    By \cref{rectthm3} this immediately concludes the proof.
\end{proof}

\begin{osservazione}
Just from the quantitative proposition one can prove several consequences. First, with a standard argument one can show abundance of flat cubes, i.e. one can show that for every cube there exists a subcube with comparable diameter with small bilateral $\beta$. Second, one can show that the measures satisfying the GWALA condition satisfy a form of compactness, in the sense that if one takes a sequence of measures with the GWALA with uniformly bounded constants, then all the limits still satisfy the GWALA with slightly bigger constants. These consequences are not proved in this work since they are not needed to reach the proof of uniform rectifiability. 
\end{osservazione}

\subsection{Non-Carleson families}

In this subsection we record a quantitative form of the failure of the Carleson packing condition. Starting from a non-Carleson family of dyadic cubes $\mathscr F$, we show that, for every $M\in\mathbb N$, one can find a cube $P\in\mathscr F$ and a finite tower of $M$ generations below $P$ whose bottom layer covers almost all of $P$ in $\mu$-measure. The construction also allows us to impose separation of scales and a mild geometric separation from the boundaries of the ancestors, properties that will be needed in the later constructions.

\begin{lemma}\label{nonClayers}
 Suppose $\mu$ is a $k$-AD-Radon measure with regularity constant $D$.
Suppose that $\mathscr F \subseteq \Delta_\mu$ is not Carleson. Then, for every $L,M,W\in\mathbb N$ and $\eta>0$, we can find a cube $P\in\mathscr F$ and $M+1$ finite families $\mathfrak{L}_0,\dots,\mathfrak{L}_M\subset\mathscr{F}(P)$ so that
\begin{itemize}
\item[(i)] $\mathfrak{L}_0=\{P\}$.
\item[(ii)] No cube appears in more than one of the families $\mathfrak{L}_0,\dots,\mathfrak{L}_M$.
\item[(iii)] The cubes in each family $\mathfrak{L}_m$, for $m=0,\dots,M$, are pairwise disjoint and for every $Q,Q'\in \mathfrak L_m$ we have 
\begin{equation}
    \ell(Q)=\ell(Q')\qquad\text{or, if $\ell(Q)<\ell(Q')$, then }\qquad\frac{\ell(Q)}{\ell(Q')}\leq \frac{1}{W};
    \label{separazionediscale}
\end{equation}
\item[(iv)] Each cube $Q'\in\mathfrak{L}_{m}$, for every $m=1,\dots,M$, is contained in a unique strictly larger cube $Q\in \mathfrak{L}_{m-1}$ with 
$$\diam(Q)\geq L\,\diam(Q').$$
\item[(v)] $\displaystyle\sum_{Q\in\mathfrak{L}_M}\mu(Q)\ge(1-\eta)\mu(P)\,.$
\end{itemize}
\end{lemma}

\begin{proof}
Let us choose $W_0\in\mathbb N$ in such a way that 
$W_0\geq W$ and $\oldC{cubi}^{-2}2^{W_0}\geq L$.
For every $\iota=0,\ldots,W_0-1$ let
$$
\tilde{\Delta}_\mu^\iota:=\bigcup_{j\in \mathbb Z}\Delta_{\iota+W_0j,\mu}.
$$
Since 
$$
\Delta_{\mu}=\bigcup_{\iota=0}^{W_0-1}\tilde{\Delta}_\mu^\iota,
$$
there exists $\iota_0\in\{0,\ldots,W_0-1\}$ such that $\mathscr F\cap \tilde{\Delta}_\mu^{\iota_0}$ is not Carleson. Indeed, otherwise $\mathscr F$ would be a finite union of Carleson families, and hence Carleson. Replacing $\mathscr F$ by $\mathscr F\cap \tilde{\Delta}_\mu^{\iota_0}$, we may assume without loss of generality that
\begin{equation}
    \mathscr F\subseteq \tilde\Delta_{\mu}^{\iota_0}.
    \label{inclusione}
\end{equation}

By applying \cite[Lemma 7]{NTV} to this reduced family, we find $P\in \mathscr F$ and $M+1$ finite families $\mathfrak L_0,\dots,\mathfrak L_M\subset \mathscr F(P)$ satisfying items (i), (ii), the pairwise disjointness part of item (iii), item (v), and such that each cube $Q'\in\mathfrak L_m$, $m=1,\dots,M$, is contained in a unique strictly larger cube $Q\in\mathfrak L_{m-1}$.

It remains to check the quantitative scale conclusions. Let $Q,Q'\in\mathfrak L_m$ and suppose $\ell(Q)<\ell(Q')$. Since both cubes belong to $\tilde\Delta_\mu^{\iota_0}$, their generations differ by at least $W_0$. Hence
$$
\frac{\ell(Q)}{\ell(Q')}\leq 2^{-W_0}\leq 2^{-W}\leq \frac{1}{W},
$$
which proves \eqref{separazionediscale}. Finally, let $Q'\in\mathfrak L_m$, $m=1,\dots,M$, and let $Q\in\mathfrak L_{m-1}$ be the unique strictly larger cube containing $Q'$. If $Q\in\Delta_{j(Q),\mu}$ and $Q'\in\Delta_{j(Q'),\mu}$, then \eqref{inclusione} gives
$$
j(Q')-j(Q)\geq W_0.
$$
Therefore, by the size estimates for dyadic cubes,
$$
\frac{\diam Q}{\diam Q'}
\geq
\frac{\oldC{cubi}^{-1}2^{-j(Q)}}{\oldC{cubi}2^{-j(Q')}}
=
\oldC{cubi}^{-2}2^{j(Q')-j(Q)}
\geq
\oldC{cubi}^{-2}2^{W_0}
\geq L.
$$
This proves item (iv) and concludes the proof.
\end{proof}

The previous lemma gives a tower of cubes with quantitative scale separation. In the next variant we refine the construction by discarding the terminal cubes whose ancestor chain passes too close to the boundary of a previous ancestor. The loss of mass is controlled by the small-boundary property of the David cubes. After this pruning, the remaining layers still carry almost all the mass of the top cube and satisfy an additional separation property for the enlarged balls.

\begin{lemma}\label{nonClayersv2}
Suppose $\mu$ is a $k$-AD-Radon measure with regularity constant $D$. Suppose that $\mathscr F \subseteq \Delta_\mu$ is not Carleson. Then, for every $L,M,W,\gimel \in\mathbb N$ and $0<\eta<1$, we can find a cube $P\in\mathscr F$ and $M+1$ finite families $\mathfrak{L}_0,\dots,\mathfrak{L}_M\subset\mathscr{F}(P)$ so that
\begin{itemize}
\item[(i)] $\mathfrak{L}_0=\{P\}$.
\item[(ii)] No cube appears in more than one of the families $\mathfrak{L}_0,\dots,\mathfrak{L}_M$.
\item[(iii)] The cubes in each family $\mathfrak{L}_m$, for $m=0,\dots,M$, are pairwise disjoint and for every $Q,Q'\in \mathfrak L_m$ we have 
\begin{equation}
    \ell(Q)=\ell(Q')\qquad\text{or, if $\ell(Q)<\ell(Q')$, then }\qquad\frac{\ell(Q)}{\ell(Q')}\leq \frac{1}{W};
    \label{separazionediscale2asfasdf}
\end{equation}
\item[(iv)] Each cube $Q'\in\mathfrak{L}_{m}$, for every $m=1,\dots,M$, is contained in a unique strictly larger cube $Q\in \mathfrak{L}_{m-1}$ with 
$$\diam(Q)\geq L\,\diam(Q').$$
Moreover, if $Q\in \mathfrak L_{m-1}$ and $Q'\in \mathfrak L_m$ are such that $\gimel B_{Q'}\cap \gimel B_Q\neq \emptyset$, then 
$$\diam(Q)\geq L\diam(Q').$$
\item[(v)] For every $m=1,\ldots,M$ and every $Q\in \mathfrak L_m$, if $Q^-\in \mathfrak L_{m-1}$ is the unique cube containing $Q$, then
$$10\gimel B_Q\cap \supp\mu\subseteq Q^-.$$
\item[(vi)] $\displaystyle\sum_{Q\in\mathfrak{L}_M}\mu(Q)\ge(1-\eta)\mu(P).$
\end{itemize}
\end{lemma}

\begin{proof}
Recall that, by \cref{evev}, for every $\tau>0$ we have
$$\mu(\partial (Q,\tau))\leq \oldC{cubi}^2\tau^{1/\oldC{cubi}}\mu(Q).$$
Let us fix $\tau_0\in (0,1)$ such that
$M\oldC{cubi}^2(2\tau_0)^{1/\oldC{cubi}}\leq \eta/2$.
Choose $\Lambda\in \mathbb N$ so large that
\begin{equation}\label{eq:choice-Lambda}
\Lambda\geq \max\{L,2\tau_0^{-1},20\gimel\tau_0^{-1},8\gimel L\tau_0^{-1}\}.
\end{equation}
and such that the hypotheses of \cref{nonClayers} are satisfied with $\Lambda$ in place of $L$. We apply \cref{nonClayers} to the family $\mathscr F$ with the parameters $\Lambda,M,W$ and $\eta/2$. We obtain a cube $P\in\mathscr F$ and $M+1$ finite families $\mathfrak L_0,\dots,\mathfrak L_M\subset\mathscr F(P)$ satisfying items (i)--(v) of \cref{nonClayers}.

Now we construct sublayers of $\mathfrak L_m$ by pruning away descendants close to the boundaries of their ancestors. For every $Q\in\mathfrak L_M$, let
$$
\mathfrak f_0(Q)\supset \mathfrak f_1(Q)\supset\cdots\supset \mathfrak f_M(Q)=Q
$$
be its ancestor chain, with $\mathfrak f_m(Q)\in\mathfrak L_m$. We say that $Q$ is \emph{good} if for every $m=1,\dots,M$,
\begin{equation}\label{eq:goodness}
\mathfrak f_m(Q)\cap \partial(\mathfrak f_{m-1}(Q),\tau_0)=\emptyset.
\end{equation}
Let $\mathfrak G_M\subset \mathfrak L_M$ be the family of all good cubes. For $m=0,\dots,M-1$, define $\mathfrak G_m$ to be the family of all ancestors in $\mathfrak L_m$ of cubes in $\mathfrak G_M$.

Since $\mathfrak G_m\subseteq \mathfrak L_m$, items (ii) and (iii) are inherited by the layers $\mathfrak G_m$. The first part of item (iv) is also inherited from the layers $\mathfrak L_m$, since $\Lambda\geq L$. Once the lower bound for $\mathfrak G_M$ is proved, we will also have $\mathfrak G_M\neq\emptyset$, and hence $\mathfrak G_0=\{P\}$.

We first prove item (v). Let $m=1,\dots,M$, let $Q\in\mathfrak G_m$, and let $Q^-\in\mathfrak G_{m-1}$ be the unique cube containing $Q$. Since $Q$ is an ancestor of a good cube, we have
$$
Q\cap \partial(Q^-,\tau_0)=\emptyset.
$$
Moreover, by \eqref{eq:choice-Lambda},
$$
10\gimel\diam Q\leq \frac{10\gimel}{\Lambda}\diam Q^-\leq \tau_0\diam Q^-.
$$
Thus, if a point of $10\gimel B_Q\cap\supp\mu$ belonged to $\supp\mu\setminus Q^-$, then $Q$ would meet $\partial(Q^-,\tau_0)$, a contradiction. Hence
$10\gimel B_Q\cap \supp\mu\subseteq Q^-$.

We now estimate the mass lost in the pruning. For each $m=1,\dots,M$, let $\mathcal B_m\subset \mathfrak L_M$ be the family of cubes $Q$ for which
\begin{equation}\label{eq:Bm}
\mathfrak f_m(Q)\cap \partial(\mathfrak f_{m-1}(Q),\tau_0)\neq\emptyset.
\end{equation}
We claim that
\begin{equation}\label{eq:Bm-estimate}
\sum_{S\in\mathcal B_m}\mu(S)\leq \oldC{cubi}^2(2\tau_0)^{1/\oldC{cubi}}\mu(P)\qquad\text{for every }m=1,\dots,M.
\end{equation}
Fix $m$ and $R\in\mathfrak L_{m-1}$. Let $S\in\mathcal B_m$ be such that $\mathfrak f_{m-1}(S)=R$, and set $T:=\mathfrak f_m(S)$. Then $T\in\mathfrak L_m$, $T\subset R$, and by the definition of $\mathcal B_m$, we have $T\cap \partial(R,\tau_0)\neq\emptyset$.
Moreover, by the first part of item (iv) for the layers $\mathfrak L_m$,
$$
\diam(T)\leq \Lambda^{-1}\diam R\leq \frac{\tau_0}{2}\diam R.
$$
Hence $T\subseteq \partial(R,2\tau_0)$. For each fixed such cube $T$, the terminal cubes $S\in\mathfrak L_M$ with $\mathfrak f_m(S)=T$ are pairwise disjoint and contained in $T$. Therefore
\begin{equation}
    \begin{split}
        \sum_{\substack{S\in \mathcal B_m\\ \mathfrak f_{m-1}(S)=R}}\mu(S)
        &\leq
        \sum_{\substack{T\in\mathfrak L_m\\ T\subset R,\ T\cap\partial(R,\tau_0)\neq\emptyset}}
        \sum_{\substack{S\in\mathfrak L_M\\ \mathfrak f_m(S)=T}}\mu(S)\leq
        \sum_{\substack{T\in\mathfrak L_m\\ T\subset R,\ T\cap\partial(R,\tau_0)\neq\emptyset}}\mu(T)
        \leq \mu(\partial(R,2\tau_0))\\
        &\leq \oldC{cubi}^2(2\tau_0)^{1/\oldC{cubi}}\mu(R).
    \end{split}
\end{equation}
Summing over the pairwise disjoint cubes $R\in\mathfrak L_{m-1}$ gives \eqref{eq:Bm-estimate}.
Every cube in $\mathfrak L_M\setminus \mathfrak G_M$ fails \eqref{eq:goodness} at least at one level $m$, and hence belongs to $\bigcup_{m=1}^M\mathcal B_m$. By \eqref{eq:Bm-estimate} and the choice of $\tau_0$,
\begin{equation}\label{eq:loss}
\sum_{Q\in\mathfrak L_M\setminus\mathfrak G_M}\mu(Q)
\leq \sum_{m=1}^M \sum_{Q\in\mathcal B_m}\mu(Q)
\leq M\oldC{cubi}^2(2\tau_0)^{1/\oldC{cubi}}\mu(P)
\leq \frac{\eta}{2}\mu(P).
\end{equation}
Since \cref{nonClayers} gives
$$
\sum_{Q\in\mathfrak L_M}\mu(Q)\geq \bigg(1-\frac{\eta}{2}\bigg)\mu(P),
$$
we obtain
$$
\sum_{Q\in\mathfrak G_M}\mu(Q)\geq  \sum_{Q\in\mathfrak L_M}\mu(Q)-\sum_{Q\in\mathfrak L_M\setminus\mathfrak G_M}\mu(Q)
\geq \bigg(1-\frac{\eta}{2}\bigg)\mu(P)-\frac{\eta}{2}\mu(P)=(1-\eta)\mu(P).
$$
This proves item (vi). In particular $\mathfrak G_M\neq\emptyset$, and hence $\mathfrak G_0=\{P\}$, so item (i) also holds for the pruned layers.

We are left to check the second part of item (iv). Fix $m\in\{1,\dots,M\}$. Let $Q'\in\mathfrak G_m$ and let $R\in\mathfrak G_{m-1}$ be its unique parent. Let $Q\in\mathfrak G_{m-1}$ be such that
\begin{equation}\label{eq:AB-intersect}
\gimel B_Q\cap \gimel B_{Q'}\neq\emptyset.
\end{equation}
We prove that $\diam(Q)\geq L\diam(Q')$. If $Q=R$, this follows from the first part of item (iv). Suppose then that $Q\neq R$. Since $Q$ and $R$ are distinct cubes in the same layer $\mathfrak G_{m-1}$, they are disjoint. Arguing by contradiction, suppose that
$$
\diam(Q)<L\diam(Q').
$$
By \eqref{eq:AB-intersect}, we get
$$
\dist(Q',Q)\leq 2\gimel(\diam(Q')+\diam(Q))\leq 2\gimel(1+L)\diam(Q')\leq 4\gimel L\diam(Q').
$$
Since $Q\subseteq \supp\mu\setminus R$, this gives 
$$\dist(Q',\supp\mu\setminus R)\leq 4\gimel L\diam(Q').$$
Since $R$ is the parent of $Q'$, the first part of item (iv) gives
$\diam(Q')\leq \Lambda^{-1}\diam(R)$ and thus by \eqref{eq:choice-Lambda}, we infer that
$$
\dist(Q',\supp\mu\setminus R)\leq \tau_0\diam(R).
$$
Hence $Q'\cap \partial(R,\tau_0)\neq\emptyset$, contradicting the goodness of $Q'$. This concludes the proof of the proposition.
\end{proof}

\subsection{Invariant cubes}

We introduce a quantitative notion of invariance of a cube along a direction. This notion is meant to localize, at the level of dyadic cubes, the splitting of tangent measures given by \cref{th:plit}. More precisely, it reflects the basic feature of the decomposability bundle: directions belonging to the bundle are directions along which the measure admits local decompositions into curves, and therefore directions along which the support behaves as if it were locally invariant. Here this is encoded by requiring that, near most points of the cube, there is a $C(e,\sigma)$-curve which stays close to $\supp\mu$ for most of its length.

This is one of the central definitions of the paper. It is the mechanism that will allow us to construct the perturbation in the proof of the WALA conjecture.

\begin{definizione}\label{definizionecubieinvarianti}
Suppose $\mu$ is a $k$-AD-regular measure with regularity constant $D$. Fix
$e\in \mathbb S^{n-1}$ and the parameters
$\sigma,\varepsilon,\delta,\zeta,\mathfrak m\in (0,1)$ and $A\geq 1$. We say that a cube $Q\in \Delta_\mu$ is
$(e,\sigma,\varepsilon,\delta,\zeta,\mathfrak m,A)$-invariant if, denoting by
$G_{e,\sigma,\varepsilon,\delta,\zeta,\mathfrak m,A}(Q)$ the set of those
$x\in AB_Q$ for which there is a $C(e,\sigma)$-curve passing through
$B(x,\zeta \diam Q)$ such that
$$
\Haus^1\Big(2AB_Q\cap \gamma\setminus B(\supp\mu,\delta \diam Q)\Big)
\leq \varepsilon\diam Q,
$$
we have that $G_{e,\sigma,\varepsilon,\delta,\zeta,\mathfrak m,A}(Q)$ is
$\mathfrak m\diam Q$-dense in $AB_Q\cap \supp\mu$.
We denote by $\mathfrak F_{e,\sigma,\varepsilon,\delta,\zeta,\mathfrak m,A}$
the family of cubes that are \emph{not}
$(e,\sigma,\varepsilon,\delta,\zeta,\mathfrak m,A)$-invariant.
\end{definizione}

We also need a multi-directional version of the previous definition. A cube is invariant in $j$ independent directions if it is invariant along each of those directions and the directions are quantitatively transverse. The transversality condition prevents the corresponding cones from overlapping and will allow us to move independently along the selected directions.

\begin{definizione}\label{definvariantiasognsaodn}
Suppose $\mu$ is a $k$-AD-regular measure with regularity constant $D$. Fix
$j\in\{1,\ldots,n\}$ and choose parameters $\varepsilon,\delta,\zeta,\mathfrak m\in (0,1)$ and $A\geq 1$. We further choose
$$
\mathfrak{v}=(v_1,\ldots,v_j)\in (\mathbb S^{n-1})^{j},\qquad
\boldsymbol{\vartheta}=(\vartheta_1,\ldots,\vartheta_{j})\in (0,1)^{j},
\qquad
\boldsymbol{\sigma}=(\sigma_1,\ldots,\sigma_j)\in (0,1)^{j},
$$
with
$$\max_{i=1,\ldots,j}\sigma_i\leq \min_{i=1,\ldots,j}\frac{\vartheta_i}{16}.$$
We say that a cube $Q\in \Delta_\mu$ is
$(j,\mathfrak{v},\boldsymbol{\vartheta},\boldsymbol{\sigma},\varepsilon,\delta,\zeta,\mathfrak m,A)$-invariant
if
\begin{equation}
\dist(v_i,\mathrm{span}(\{v_h:h\neq i\})\geq \vartheta_{i},
\qquad\text{for every }i\in \{1,\ldots,j\},
\label{indipendenzavettori}
\end{equation}
and the cube $Q$ is
$(v_i,\sigma_i,\varepsilon,\delta,\zeta,\mathfrak m,A)$-invariant for every
$i=1,\ldots,j$.
When $\vartheta_i\equiv \vartheta$ and $\sigma_i\equiv \sigma$ for every $i=1,\ldots, j$, we say that $Q$ is $(j,\mathfrak{v},\vartheta,\sigma, \varepsilon,\delta,\zeta,\mathfrak m,A)$-invariant.
\end{definizione}

\begin{osservazione}
Notice that if $e_1,\ldots,e_j\in \mathbb S^{n-1}$ are such that 
\begin{equation}
\dist(e_i,\mathrm{span}(\{e_h:h\neq i\}))\geq \vartheta,
\qquad\text{for every }i\in \{1,\ldots,j\},
\nonumber
\end{equation}
and $\sigma\leq \vartheta/16$, then 
$$C(e_{i_1},\sigma)\cap C(e_{i_2},\sigma)=\{0\}\qquad\text{whenever $i_1\neq i_2$.}$$
Indeed, suppose by contradiction that there exists
$x\in C(e_{i_1},\sigma)\cap C(e_{i_2},\sigma)$ with $x\neq 0$, and set $u:=x/|x|$ and
$V:=\mathrm{span}(\{e_h:h\neq i_1\})$. Since $e_{i_2}\in V$ and $u\in C(e_{i_2},\sigma)$, we have
$|P_{V^\perp}u|\leq \sigma$. Moreover, since $u\in C(e_{i_1},\sigma)$, we have
$$
\dist(e_{i_1},V)\leq \dist(e_{i_1},\mathrm{span}(u))+|P_{V^\perp}u|\leq 2\sigma<\vartheta,
$$
contradicting the assumption.
\end{osservazione}

\begin{definizione}
    Let $j=1,\ldots, n$ and suppose $\mathfrak e=(e_1,\ldots, e_j)\in (\mathbb S^{n-1})^j$. Then we let 
    $$\mathrm{span}(\mathfrak e):=\mathrm{span}(\{e_1,\ldots, e_j\}).$$
\end{definizione}

\begin{proposizione}\label{densitalungolecurve}
Let $\gamma$ be a $C(e,\sigma)$-curve such that
$$
\Haus^1(2AB_R\cap \gamma\setminus B(\supp\mu,\delta \diam R))
\leq \varepsilon\diam R.
$$
Then, for every $q\in \gamma\cap 2(A-\varepsilon)B_R$, we have
$B(q,(\delta+\varepsilon)\diam R)\cap \supp\mu\neq \emptyset$.
\end{proposizione}

\begin{proof}
Suppose by contradiction that, for some $q\in\gamma\cap 2(A-\varepsilon)B_R$,
$$
B(q,(\delta+\varepsilon)\diam R)\cap\supp\mu=\emptyset.
$$
Then, we would have
$$
B(q,\varepsilon\diam R)\cap B(\supp\mu,\delta\diam R)=\emptyset.
$$
Since $q\in 2(A-\varepsilon)B_R$, there also holds that $B(q,\varepsilon\diam R)\subseteq 2AB_R$. Therefore
$$
\Haus^1(2AB_R\cap \gamma\setminus B(\supp\mu,\delta\diam R))
\geq
\Haus^1(\gamma\cap B(q,\varepsilon\diam R)).
$$
Since $q\in\gamma$ and $\gamma$ is a full $C(e,\sigma)$-curve, 
$$
\Haus^1(\gamma\cap B(q,\varepsilon\diam R))\geq 2\varepsilon\diam R,
$$
contradicting the assumed bound.
\end{proof}

We now show that invariance in $j$ independent directions forces the support to be quantitatively dense along the $j$-dimensional affine planes spanned by those directions and passing through points of the support. In this sense, invariant cubes carry a local fibration by parallel affine planes.

\begin{proposizione}\label{propo:piano-da-invarianza}
Suppose $\mu$ is a $k$-AD-regular measure with regularity constant $D$. Let
$1\leq j\leq n$, let $\mathfrak e=(e_1,\ldots,e_j)\in(\mathbb S^{n-1})^j$, and suppose that
$Q\in\Delta_\mu$ is
$(j,\mathfrak e,\vartheta,\sigma,\varepsilon,\delta,\zeta,\mathfrak m,A)$-invariant.
Assume that
$$
A\vartheta\geq 16n,\qquad \sigma\leq \frac{\vartheta}{16n},
\qquad \varepsilon+\delta+\zeta+\mathfrak m\leq \frac{1}{16n^2}.
$$
Then, for every $x\in \frac{A}{16n}B_Q\cap\supp\mu$, writing
$V_\mathfrak{e}:=\mathrm{span}(\mathfrak e)$, we have
$$
\sup_{z\in (x+V_\mathfrak{e})\cap B(x,\frac{A\vartheta}{16n}\diam Q)}
\frac{\dist(z,\supp\mu\cap B(x,\tfrac{A\vartheta}{8n}\diam Q))}
{\frac{A\vartheta}{8n}\diam Q}
\leq 8n\Big(\frac{\sigma}{\vartheta}
+n\frac{\varepsilon+\delta+\zeta+\mathfrak m}{A\vartheta}\Big).
$$
\end{proposizione}

\begin{proof}
Fix $x\in \tfrac{A}{16n}B_Q\cap\supp\mu$ and set $\Pi_Q:=x+\mathrm{span}(e_1,\ldots,e_j)$, let $z\in \Pi_Q\cap B(x,\tfrac{A\vartheta}{16n}\diam Q)$ and write
$$
z=x+\sum_{i=1}^j t_i e_i.
$$
We first estimate the coefficients $t_i$. Fix $i\in\{1,\ldots,j\}$ and let
$\pi_i$ be the orthogonal projection onto
$\mathrm{span}(\{e_h:h\neq i\})^\perp$. Then $\pi_i e_h=0$ for every $h\neq i$, and by assumption on $Q$, see \eqref{indipendenzavettori}, we have 
$$
|\pi_i e_i|=\dist(e_i,\mathrm{span}(\{e_h:h\neq i\}))\geq \vartheta .
$$
Therefore $\pi_i(z-x)=t_i\pi_i e_i$, and hence
$\vartheta |t_i|\leq |\pi_i(z-x)|\leq |z-x|$.
Since $z\in B(x,\tfrac{A\vartheta}{16n}\diam Q)$, we have
$$
|t_i|\leq \frac{A}{16n}\diam Q.
$$
For every $\iota=0,\ldots,j$, we define
$z_\iota:=x+\sum_{i=1}^\iota t_i e_i$ and notice that
\begin{equation}
|z_j-x|\leq \sum_{i=1}^j |t_i|\leq j\frac{A}{16n}\diam Q\leq \frac{A}{16}\diam Q.
\end{equation}
Thus all the points $z_0,\ldots,z_j$ lie inside $AB_Q$. We now construct inductively points $x_0,\ldots,x_j\in\supp\mu$ satisfying
$$
|x_i-z_i|\leq 2\sigma\sum_{h=1}^i |t_h|+i(\varepsilon+\delta+\zeta+\mathfrak m)\diam Q,
\qquad i=0,\ldots,j.
$$
Set $x_0=x=z_0$. Suppose $x_{i-1}$ has been constructed. Then $x_{i-1}\in AB_Q$, indeed
\begin{equation}
\begin{split}
|x-x_{i-1}|&\leq |x-z_{i-1}|+|x_{i-1}-z_{i-1}|\leq \frac{A}{8}\diam Q+2\sigma\sum_{h=1}^{i-1}|t_h|+(i-1)(\varepsilon+\delta+\zeta+\mathfrak m)\diam Q\\
&\leq \frac{A}{8}\diam Q+2\sigma(j-1)\frac{A}{8n}\diam Q+(j-1)(\varepsilon+\delta+\zeta+\mathfrak m)\diam Q<\frac{A}{4}\diam Q.
\nonumber
\end{split}
\end{equation}
Since $x\in \tfrac{A}{16n}B_Q$ and $A\geq 16n$, this implies $x_{i-1}\in AB_Q$. Since $Q$ is $(e_i,\sigma,\varepsilon,\delta,\zeta,\mathfrak m,A)$-invariant, there exists
$y_{i-1}\in G_{e_i,\sigma,\varepsilon,\delta,\zeta,\mathfrak m,A}(Q)$ such that
$$
|y_{i-1}-x_{i-1}|\leq \mathfrak m\diam Q .
$$
By the definition of the good set, there is a $C(e_i,\sigma)$-curve $\gamma_i$
passing through $B(y_{i-1},\zeta\diam Q)$ and satisfying
$$
\Haus^1\Big(2AB_Q\cap \gamma_i\setminus B(\supp\mu,\delta\diam Q)\Big)\leq \varepsilon\diam Q.
$$
Choose $p_{i-1}\in\gamma_i$ with
$$
|p_{i-1}-y_{i-1}|\leq \zeta\diam Q.
$$
We use the standard parametrization of a $C(e_i,\sigma)$-curve by its
$e_i$-coordinate $\gamma(t)=te_i+\eta(t)$ with $\eta(t)-\eta(s)\in e_i^\perp$ for every $s,t\in \R$.
Let $p_i$ be the point of $\gamma_i$ whose $e_i$-coordinate is shifted by $t_i$ with respect to $p_{i-1}$. Without loss of generality we can assume that $p_{i-1}=\gamma(0)$ and hence $p_i=\gamma(t_i)$.
Since $\sigma\leq 1/(16n)$, this point satisfies
$$
|p_i-p_{i-1}-t_i e_i|\leq 2\sigma |t_i|.
$$
Moreover, since $|t_i|\leq \tfrac{A}{16n}\diam Q$, we have
\begin{equation}
\begin{split}
|p_i-x|&\leq |p_{i-1}-x|+|t_i|(1+2\sigma)\leq |p_{i-1}-y_{i-1}|+|y_{i-1}-x_{i-1}|+|x-x_{i-1}|+|t_i|(1+2\sigma)\\
&\leq (\zeta+\mathfrak m)\diam Q+\frac{A}{4}\diam Q+(1+2\sigma)\frac{A}{16n}\diam Q<\frac{A}{3}\diam Q.
\nonumber
\end{split}
\end{equation}
Since $x\in \tfrac{A}{16n}B_Q$ and $A\geq 16n$, this implies $p_i\in AB_Q$. Hence, by \cref{densitalungolecurve}, there exists $x_i\in \supp\mu$ with
$$
|x_i-p_i|\leq (\delta+\varepsilon)\diam Q.
$$
Hence
\begin{equation}
\begin{split}
|x_i-z_i|&\leq |p_i-x_i|+|p_i-z_i|\leq (\delta+\varepsilon)\diam Q+|p_i-z_{i-1}-t_ie_i|\\
&\leq (\delta+\varepsilon)\diam Q+|p_i-p_{i-1}-t_ie_i|+|p_{i-1}-y_{i-1}|+|y_{i-1}-x_{i-1}|+|x_{i-1}-z_{i-1}|\\
&\leq 2\sigma |t_i|+(\varepsilon+\delta+\zeta+\mathfrak m)\diam Q+2\sigma\sum_{h=1}^{i-1}|t_h|+(i-1)(\varepsilon+\delta+\zeta+\mathfrak m)\diam Q\\
&=2\sigma\sum_{h=1}^{i}|t_h|+i(\varepsilon+\delta+\zeta+\mathfrak m)\diam Q.
\nonumber
\end{split}
\end{equation}
This concludes the induction. This shows in particular that
\begin{equation}
\begin{split}
|x_j-z|&\leq 2\sigma\sum_{h=1}^j |t_h|+j(\varepsilon+\delta+\zeta+\mathfrak m)\diam Q\leq \Big(\frac{A\sigma}{4}+j(\varepsilon+\delta+\zeta+\mathfrak m)\Big)\diam Q.
\nonumber
\end{split}
\end{equation}
Moreover, since $\sigma\leq \vartheta/(16n)$, $j\leq n$, $\varepsilon+\delta+\zeta+\mathfrak m\leq 1/(16n^2)$ and $A\vartheta\geq 16n$, we have
$$
|x_j-z|\leq \frac{A\vartheta}{16n}\diam Q.
$$
Since $|z-x|\leq \frac{A\vartheta}{16n}\diam Q$, this implies
$$
x_j\in\supp\mu\cap B(x,\tfrac{A\vartheta}{8n}\diam Q).
$$
Therefore
$$
\dist(z,\supp\mu\cap B(x,\tfrac{A\vartheta}{8n}\diam Q))\leq \Big(\frac{A\sigma}{4}+j(\varepsilon+\delta+\zeta+\mathfrak m)\Big)\diam Q.
$$
Dividing by $\tfrac{A\vartheta}{8n}\diam Q$, we get
$$
\frac{\dist(z,\supp\mu\cap B(x,\tfrac{A\vartheta}{8n}\diam Q))}{\tfrac{A\vartheta}{8n}\diam Q}
\leq
2n\frac{\sigma}{\vartheta}
+
8nj\frac{\varepsilon+\delta+\zeta+\mathfrak m}{A\vartheta}.
$$
Since $j\leq n$, taking the supremum over $z\in \Pi_Q\cap B(x,\frac{A\vartheta}{16n}\diam Q)$ concludes the proof.
\end{proof}

The next corollary is a quantified dyadic analogue of \cref{abscontinuous}. It says that a $k$-AD-regular measure cannot be locally too close, at sufficiently small error, to affine planes of dimension strictly larger than $k$. Consequently, no cube can be invariant in more than $k$ independent directions.

\begin{corollario}\label{noinvariantmorethank}
Suppose $\mu$ is a $k$-AD-regular measure with regularity constant $D$ and let $k<j\leq n$. Then there are no cubes $Q\in\Delta_\mu$ that are
$(j,\mathfrak e,\vartheta,\sigma,\varepsilon,\delta,\zeta,\mathfrak m,A)$-invariant, where $\mathfrak e=(e_1,\ldots,e_j)$, provided
$$
A\vartheta\geq 16n,\qquad\text{and}\qquad
16n\Big(\frac{\sigma}{\vartheta}
+n\frac{\varepsilon+\delta+\zeta+\mathfrak m}{A\vartheta}\Big)\leq \frac{1}{2\cdot 3^k 16^n D^2}.
$$
\end{corollario}

\begin{proof}
Set
$$
\eta:=16n\Big(\frac{\sigma}{\vartheta}
+n\frac{\varepsilon+\delta+\zeta+\mathfrak m}{A\vartheta}\Big)\qquad\text{and}\qquad \eta_0:=\frac{1}{2\cdot 3^k 16^n D^2}.
$$
Suppose by contradiction that such a cube $Q$ exists. Choose $x_Q\in Q\cap\supp\mu$, set
$V:=\mathrm{span}(\mathfrak e)$, and let
$r:=\frac{A\vartheta}{16n}\diam Q$. By \cref{propo:piano-da-invarianza}, for every $z\in (x_Q+V)\cap B(x_Q,r)$ there exists
$y_z\in\supp\mu\cap B(x_Q,2r)$ such that
$$
|z-y_z|\leq 2\eta r.
$$
Let $\mathcal Z$ be a maximal $8\eta r$-separated subset of
$(x_Q+V)\cap B(x_Q,r/2)$. Then
$$
\mathrm{Card}(\mathcal Z)\geq (16\eta)^{-j}.
$$
For $\eta\leq \eta_0$, the points $y_z$ relative to the points $z$ of the maximal set $\mathcal Z$, are $4\eta r$-separated. Hence the balls
$B(y_z,\eta r)$ are pairwise disjoint and contained in $B(x_Q,3r)$. By AD-regularity, we see that
$$
D^{-1}\mathrm{Card}(\mathcal Z)(\eta r)^k
\leq
\mu(B(x_Q,3r))
\leq
D(3r)^k.
$$
This implies in particular that
$$
1\leq 3^kD^2 16^j\eta^{j-k},
$$
which by the choice of $\eta\leq \eta_0$ and the fact that $k<j\leq n$, results in a contradiction. 
\end{proof}

The next proposition show that if a cube is invariant in $j$ independent directions, then, at a sufficiently coarse scale, the support is locally organized as a product. The proposition makes this precise by constructing, near each point of the cube, a separated set $E_{Q,x}\subset \mathrm{span}(\mathfrak e)^\perp$ such that $\supp\mu$ is close to $\mathrm{span}(\mathfrak e)\otimes E_{Q,x}$. The AD-regularity of $\mu$ then forces $E_{Q,x}$ to have coarse dimension $k-j$, as expressed by the upper and lower cardinality estimates.

\begin{proposizione}
\label{prop:discrete-product-from-partial-invariance}
Let $\mu$ be a $k$-AD-regular measure in $\mathbb R^n$ with regularity constant $D$. Let
$Q\in\Delta_\mu$ be
$(j,\mathfrak e,\vartheta,\sigma,\varepsilon,\delta,\zeta,\mathfrak m,A)$-invariant, where
$\mathfrak e=(e_1,\ldots,e_j)$ and $1\leq j\leq k$. Assume that
$$
A\vartheta\geq 16n,\qquad16n\Big(\frac{\sigma}{\vartheta}
+n\frac{\varepsilon+\delta+\zeta+\mathfrak m}{A\vartheta}\Big)\leq \frac{\vartheta}{2^{10+12k}D^2},
$$
and let us define 
$$\rho_0:=2^{9+12k}D^2\Big(\frac{\sigma}{\vartheta}+n\frac{\varepsilon+\delta+\zeta+\mathfrak m}{A\vartheta}\Big)A\vartheta \diam Q.$$
Then, for every $x\in \supp\mu\cap \frac{A}{32n}B_Q$, there exists a discrete set
$E_{Q,x}\subset \mathrm{span}(\mathfrak e)^\perp$ that is $\rho_0$-separated such that, writing
$K_x:=x+\mathrm{span}(\mathfrak e)\otimes E_{Q,x}$, one has
\begin{equation}
d_{\Haus,B(x,\frac{A\vartheta}{256n}\diam Q)}
(\supp\mu,K_x)\leq 2\rho_0.
\label{eq:discrete-product-hausdorff}
\end{equation}
Moreover, for every $r$ such that $\rho_0\leq r\leq\frac{A\vartheta}{1024n}\diam Q$, one has
\begin{equation}
2^{-8k}D^{-2}\Big(\frac{r}{\rho_0}\Big)^{k-j}\leq\mathrm{Card}\big(E_{Q,x}\cap B(p,r)\big),
\label{eq:discrete-lower}
\end{equation}
for every $p\in E_{Q,x}\cap B_{\mathfrak e}\big(0,\frac{A\vartheta}{128n}\diam Q\big)$, and
\begin{equation}
\mathrm{Card}\big(E_{Q,x}\cap B(p,r)\big)\leq2^{5k}D^2\Big(\frac{r}{\rho_0}\Big)^{k-j}
\label{eq:discrete-upper}
\end{equation}
for every $p\in E_{Q,x}$.
\end{proposizione}

\begin{proof}
Fix $x\in \supp\mu\cap \frac{A}{32n}B_Q$ and assume without loss of generality that $x=0$. Translating back gives the set
$K_x=x+\mathrm{span}(\mathfrak e)\otimes E_{Q,x}$ appearing in the statement.

By \cref{propo:piano-da-invarianza}, for every
$y\in \supp\mu\cap \frac{A}{16n}B_Q$ and every
$z\in (y+\mathrm{span}(\mathfrak e))\cap B(y,\frac{A\vartheta}{16n}\diam Q)$, one has
\begin{equation}
\dist\Big(z,\supp\mu\cap B(y,A\vartheta\diam Q)\Big)
\leq
\Big(
\frac{\sigma}{\vartheta}
+
n\frac{\varepsilon+\delta+\zeta+\mathfrak m}{A\vartheta}
\Big)A\vartheta\diam Q.
\label{eq:discrete-fibre-estimate}
\end{equation}
Since $x=0$ and $x\in \frac{A}{32n}B_Q$, while $\vartheta<1$, we have
$$
\supp\mu\cap B(0,\tfrac{A\vartheta}{64n}\diam Q)
\subset
\supp\mu\cap \frac{A}{16n}B_Q.
$$
Thus \eqref{eq:discrete-fibre-estimate} applies with base point any
$y\in \supp\mu\cap B(0,\tfrac{A\vartheta}{64n}\diam Q)$. Let
$$
E:=
\pi_{\mathrm{span}(\mathfrak e)^\perp}
\Big(\supp\mu\cap B(0,\tfrac{A\vartheta}{64n}\diam Q)\Big),
$$
and let $E_{Q,x}$ be a maximal $\rho_0$-separated subset of $E$. Then
\begin{equation}
E\subset \bigcup_{p\in E_{Q,x}} B(p,\rho_0).
\label{eq:E-net-cover}
\end{equation}

\medskip

\noindent \textbf{Step I. Local product approximation.}
Let $y\in \supp\mu\cap B(0,\frac{A\vartheta}{256n}\diam Q)$. Then
$\pi_{\mathrm{span}(\mathfrak e)^\perp}(y)\in E$, and by \eqref{eq:E-net-cover} there exists
$p\in E_{Q,x}$ such that
$$
|\pi_{\mathrm{span}(\mathfrak e)^\perp}(y)-p|\leq \rho_0.
$$
Hence
\begin{equation}
\dist\big(y,\mathrm{span}(\mathfrak e)\otimes E_{Q,x}\big)
\leq \rho_0.
\label{eq:first-product-inclusion}
\end{equation}
Conversely, let
$y\in (\mathrm{span}(\mathfrak e)\otimes E_{Q,x})\cap
B_{\mathfrak e}(0,\frac{A\vartheta}{256n}\diam Q)$. Then there exists $p\in E_{Q,x}$ such that
$\pi_{\mathrm{span}(\mathfrak e)^\perp}(y)=p$. Since $E_{Q,x}\subset E$, there exists
$y_p\in \supp\mu\cap B(0,\frac{A\vartheta}{64n}\diam Q)$ such that
$p=\pi_{\mathrm{span}(\mathfrak e)^\perp}(y_p)$. Moreover,
$$
|y-y_p|
\leq
|y|+|y_p|
\leq
\frac{A\vartheta}{256n}\diam Q+
\frac{A\vartheta}{64n}\diam Q
<
\frac{A\vartheta}{16n}\diam Q.
$$
Thus \eqref{eq:discrete-fibre-estimate} applies with base point $y_p$ and point $y$, and gives
\begin{equation}
\dist(y,\supp\mu)
\leq
\Big(
\frac{\sigma}{\vartheta}
+
n\frac{\varepsilon+\delta+\zeta+\mathfrak m}{A\vartheta}
\Big)A\vartheta\diam Q
\leq \rho_0.
\label{eq:second-product-inclusion}
\end{equation}
Combining \eqref{eq:first-product-inclusion} and \eqref{eq:second-product-inclusion}, we obtain
$$
d_{\Haus,B(0,\frac{A\vartheta}{256n}\diam Q)}
(\supp\mu,\mathrm{span}(\mathfrak e)\otimes E_{Q,x})\leq 2\rho_0.
$$

\medskip

\noindent \textbf{Step II. Coarse upper regularity of $E_{Q,x}$.}
Fix $p\in E_{Q,x}$ and let $\rho_0\leq r\leq\frac{A\vartheta}{1024n}\diam Q$. Write
$$
E_{Q,x}\cap B(p,r)=\{p_\alpha:\alpha\in I\}.
$$
For every $\alpha\in I$, choose
$y_\alpha\in \supp\mu\cap B(0,\frac{A\vartheta}{64n}\diam Q)$ such that
$p_\alpha=\pi_{\mathrm{span}(\mathfrak e)^\perp}(y_\alpha)$. Since $p\in E_{Q,x}\subset E$, choose
$y_p\in \supp\mu\cap B(0,\frac{A\vartheta}{64n}\diam Q)$ such that
$p=\pi_{\mathrm{span}(\mathfrak e)^\perp}(y_p)$.

Let $\{g_\beta:\beta\in J\}$ be a maximal $\rho_0$-separated set in
$\mathrm{span}(\mathfrak e)\cap
B(\pi_{\mathrm{span}(\mathfrak e)}(y_p),r)$. Then
\begin{equation}
\mathrm{Card}(J)\geq
\Big(\frac r{\rho_0}\Big)^j.
\label{eq:plane-grid-cardinality}
\end{equation}
For every $\alpha\in I$ and $\beta\in J$, set
$z^0_{\alpha,\beta}:=p_\alpha+g_\beta$. Then $z^0_{\alpha,\beta}\in y_\alpha+\mathrm{span}(\mathfrak e)$ and
$$
|z^0_{\alpha,\beta}-y_\alpha|
=
|g_\beta-\pi_{\mathrm{span}(\mathfrak e)}(y_\alpha)|
\leq
r+|y_p-y_\alpha|
<
\frac{A\vartheta}{16n}\diam Q.
$$
Applying \eqref{eq:discrete-fibre-estimate}, for every pair $(\alpha,\beta)$ there exists
$z_{\alpha,\beta}\in \supp\mu$ such that
\begin{equation}
|z_{\alpha,\beta}-z^0_{\alpha,\beta}|
\leq
2^{-9-12k}D^{-2}\rho_0.
\label{eq:grid-error}
\end{equation}
Since the two factors are orthogonal and both families $\{p_\alpha\}_{\alpha\in I}$ and
$\{g_\beta\}_{\beta\in J}$ are $\rho_0$-separated, the points $z^0_{\alpha,\beta}$ are
$\rho_0$-separated. By \eqref{eq:grid-error}, the points $z_{\alpha,\beta}$ are
$\rho_0/2$-separated. Moreover,
$$
|z_{\alpha,\beta}-y_p|
\leq |z^0_{\alpha,\beta}-y_p|+|z_{\alpha,\beta}-z^0_{\alpha,\beta}|
\leq \sqrt2 r+\rho_0\leq 4r.
$$
Since $\mu$ is $k$-AD-regular, the number of $\rho_0/2$-separated points of
$\supp\mu$ in $B(y_p,4r)$ is at most
$32^kD^2(r\rho_0^{-1})^k$. Consequently,
$$
\mathrm{Card}(I)\mathrm{Card}(J)
\leq
2^{5k}D^2
\Big(\frac r{\rho_0}\Big)^k.
$$
Combining this with \eqref{eq:plane-grid-cardinality}, we get
$$
\mathrm{Card}\big(E_{Q,x}\cap B(p,r)\big)
=
\mathrm{Card}(I)
\leq
2^{5k}D^2
\Big(\frac r{\rho_0}\Big)^{k-j}.
$$
This proves \eqref{eq:discrete-upper}.

\medskip

\noindent \textbf{Step III. Coarse lower regularity of $E_{Q,x}$.}
Fix $p\in E_{Q,x}\cap B(0,\frac{A\vartheta}{128n}\diam Q)$ and let
$\rho_0\leq r\leq\frac{A\vartheta}{1024n}\diam Q$. We prove that
$$
\mathrm{Card}\big(E_{Q,x}\cap B(p,r)\big)
\geq
2^{-8k}D^{-2}
\Big(\frac r{\rho_0}\Big)^{k-j}.
$$
Since $p\in E_{Q,x}\subset E$, choose
$y_p\in \supp\mu\cap B(0,\frac{A\vartheta}{64n}\diam Q)$ such that
$p=\pi_{\mathrm{span}(\mathfrak e)^\perp}(y_p)$. Then
$p\in y_p+\mathrm{span}(\mathfrak e)$ and
$$
|p-y_p|\leq |y_p|<\frac{A\vartheta}{16n}\diam Q.
$$
By \eqref{eq:discrete-fibre-estimate}, there exists
$\widetilde y_p\in\supp\mu$ such that
\begin{equation}
|\widetilde y_p-p|
\leq
2^{-9-12k}D^{-2}\rho_0.
\label{eq:ytilde-close-p}
\end{equation}
If $r<4\rho_0$, then the desired estimate follows from
$p\in E_{Q,x}\cap B(p,r)$ and the fact that
$2^{-8k}D^{-2}
(\frac r{\rho_0})^{k-j}\leq 1$.
Thus we may assume that $r\geq4\rho_0$.

Let $u\in \supp\mu\cap B(\widetilde y_p,r/8)$. By \eqref{eq:ytilde-close-p} and $r\geq4\rho_0$,
$$
|u|
\leq
|p|+|\widetilde y_p-p|+|u-\widetilde y_p|
\leq
\frac{A\vartheta}{128n}\diam Q+\frac r4+\frac r8
<
\frac{A\vartheta}{64n}\diam Q.
$$
Hence $\pi_{\mathrm{span}(\mathfrak e)^\perp}(u)\in E$. Choose now
$p_u\in E_{Q,x}$ such that
$|\pi_{\mathrm{span}(\mathfrak e)^\perp}(u)-p_u|\leq\rho_0$. Since
$p\in\mathrm{span}(\mathfrak e)^\perp$,
$$
|p_u-p|\leq \rho_0+|u-p|
\leq \rho_0+\frac r8+|\widetilde y_p-p|
\leq \frac{5r}{8}<r.
$$
Therefore $p_u\in E_{Q,x}\cap B(p,r)$ and
\begin{equation}
\supp\mu\cap B(\widetilde y_p,r/8)
\subset
\bigcup_{q\in E_{Q,x}\cap B(p,r)}
\Big(
\supp\mu\cap B(\widetilde y_p,r/8)
\cap
\{w:\dist(w,q+\mathrm{span}(\mathfrak e))\leq \rho_0\}
\Big).
\label{eq:lower-tube-cover}
\end{equation}
Fix $q\in E_{Q,x}\cap B(p,r)$. The set
$$
B(\widetilde y_p,r/8)
\cap
\{w:\dist(w,q+\mathrm{span}(\mathfrak e))\leq \rho_0\}
$$
can be covered by at most
$2^{3j}(r/\rho_0)^j$ balls of radius $2\rho_0$. Discarding those which do not meet
$\supp\mu$, the upper $k$-AD-regularity gives
\begin{equation}
\mu\Big(
\supp\mu\cap B(\widetilde y_p,r/8)
\cap
\{w:\dist(w,q+\mathrm{span}(\mathfrak e))\leq \rho_0\}
\Big)
\leq
2^{5k}D r^j\rho_0^{k-j}.
\label{eq:lower-tube-measure}
\end{equation}
On the other hand, since $\widetilde y_p\in\supp\mu$, the lower $k$-AD-regularity gives
\begin{equation}
\mu(B(\widetilde y_p,r/8))
\geq
2^{-3k}D^{-1}r^k.
\label{eq:lower-ball-measure}
\end{equation}
Combining \eqref{eq:lower-tube-cover}, \eqref{eq:lower-tube-measure}, and
\eqref{eq:lower-ball-measure}, we get
$$
2^{-3k}D^{-1}r^k
\leq
\mathrm{Card}\big(E_{Q,x}\cap B(p,r)\big)
2^{5k}D r^j\rho_0^{k-j},
$$
and this proves \eqref{eq:discrete-lower} and concludes the proof.
\end{proof}

\subsection{Quantitative width loss}

The purpose of this subsection is to prove a quantitative mechanism which turns non-invariance into thinness.

\begin{lemma}\label{lemmaraggipallefini}
Suppose $\mu$ is a $k$-dimensional AD-regular measure. Let $x\in \supp\mu$ and let $r>0$. Let $E_1,E_2\subseteq \R^n$ be Borel sets. Let $\lambda,\eta\in (0,1)$ be such that $4\eta<\lambda$ and define
$$
M:=\Big\lceil \log\Big(\frac{1}{1-\eta}\Big)^{-1}
\log\Big(\frac{1-2\eta}{1-\lambda/2}\Big)\Big\rceil.
$$
Then, there exists $(1-\lambda/2)r\leq s\leq (1-2\eta)r$ such that
\begin{equation}
\mu(B(x,s/(1-\eta))\setminus B(x,s))
\leq \frac{3}{M}\mu(B(x,(1-\eta)r)\setminus B(x,(1-\lambda)r)),
\label{stimapallefinitotale}
\end{equation}
\begin{equation}
\mu(E_1\cap (B(x,s/(1-\eta))\setminus B(x,s)))
\leq \frac{3}{M}\mu(E_1\cap (B(x,(1-\eta)r)\setminus B(x,(1-\lambda)r))),
\label{stimapallefiniU}
\end{equation}
\begin{equation}
\mu(E_2\cap (B(x,s/(1-\eta))\setminus B(x,s)))
\leq \frac{3}{M}\mu(E_2\cap (B(x,(1-\eta)r)\setminus B(x,(1-\lambda)r))).
\label{stimapallefiniV}
\end{equation}
 Moreover, we can estimate $M$ with
    \begin{equation}
        \frac{\lambda}{2\eta}-2\leq M\leq \frac{\lambda}{\eta}+1.
        \label{stimaMlemma}
    \end{equation}
\end{lemma}

\begin{proof}
The estimates in \eqref{stimaMlemma} follow from elementary bounds on the logarithm. For $j=0,\ldots,M-1$ set
$$
s_j:=(1-\lambda/2)(1-\eta)^{-j}r,
\qquad
A_j:=B(x,s_j/(1-\eta))\setminus B(x,s_j),
$$
and
$$
\Omega:=B(x,(1-\eta)r)\setminus B(x,(1-\lambda)r).
$$
By the definition of $M$, we have $(1-\lambda/2)r\leq s_j\leq (1-2\eta)r$ for every $j=0,\ldots,M-1$. Moreover, the annuli $A_j$ are pairwise disjoint and contained in $\Omega$. Set
$$
\nu_0:=\mu,\qquad \nu_1:=\mu\llcorner E_1,\qquad \nu_2:=\mu\llcorner E_2.
$$
We claim that there exists $j\in\{0,\ldots,M-1\}$ such that
$$
\nu_i(A_j)\leq \frac{3}{M}\nu_i(\Omega)\qquad\text{for every }i=0,1,2.
$$
Indeed, otherwise, for every $j$ there is $i\in\{0,1,2\}$ such that
$\nu_i(A_j)>3\nu_i(\Omega)/M$. Defining
$$
I_i:=\{j:\nu_i(A_j)>3\nu_i(\Omega)/M\},
$$
the disjointness of the $A_j$ gives
$$
\mathrm{Card}(I_i)\frac{3}{M}\nu_i(\Omega)
<
\sum_{j\in I_i}\nu_i(A_j)
\leq \nu_i(\Omega),
$$
hence $\mathrm{Card}(I_i)<M/3$ for $i=0,1,2$. Since by assumption
$\{0,\ldots,M-1\}\subseteq I_0\cup I_1\cup I_2$, this gives
$$
M\leq \mathrm{Card}(I_0)+\mathrm{Card}(I_1)+\mathrm{Card}(I_2)<M,
$$
a contradiction. Taking $s=s_j$ for the index given by the claim yields
\eqref{stimapallefinitotale}, \eqref{stimapallefiniU}, and \eqref{stimapallefiniV}.
\end{proof}

In the next proposition, we start from many layers of cubes which are non-invariant in the directions of a finite set $\mathscr E$, and we construct an open set $\Omega$ which still contains a very big portion of a given set $U$ in measure. At the same time, $\Omega$ is thin in the directions of $\mathscr E$: for every $e\in\mathscr E$, every $C(e,\sigma)$-curve loses a definite proportion of its length when restricted to $\Omega$.

The proof is a thinning procedure. We first fix one direction $e\in\mathscr E$ and construct an open set $\Omega_e$ by induction on the layers. At each step, we choose a Vitali subfamily of cubes and add to $\Omega_e$ balls centered at the witnesses $z(Q,e)$ of non-invariance. These balls are chosen precisely because every $C(e,\sigma)$-curve which crosses one of them has to spend a definite amount of length in $2AB_Q\setminus B(\supp\mu,\delta\diam Q)$. Thus the length of $\gamma\cap\Omega_e$ can be paid for by portions of $\gamma$ which stay away from $\supp\mu$. The scale separation between the layers and the disjointness in the Vitali families control the overlap of these portions.

The radii of the balls are chosen by \cref{lemmaraggipallefini}. The annuli around the chosen radii are not removed from the construction; they are used to control the errors. More precisely, they control the cubes in the following layers which meet the boundary of the balls, and they also control the loss produced when the old radii are enlarged during the induction. This gives the boundary estimate needed to continue the construction and the measure estimate which shows that $\Omega_e$ captures most of $U$, up to the terminal defect.

After the construction is performed for each $e\in\mathscr E$, we take $\Omega:=\bigcap_{e\in\mathscr E}\Omega_e$. The width estimate in each direction comes from the corresponding construction of $\Omega_e$, while the measure and boundary errors are obtained by summing over the finite set of directions.

\begin{proposizione}\label{stepbasethinning}
Suppose $\mu$ is a $k$-AD-regular measure with regularity constant $D$.
Given $\xi,\xi'\in(0,1)$, fix parameters $\sigma,\varepsilon,\delta,\zeta,\mathfrak m\in (0,1)$, $A\geq 1$ and $N,\mathfrak N,\qof\in\mathbb N$. Assume that
\begin{equation}
\zeta \leq \delta\leq 2^{-24},
\qquad
\mathfrak N\geq 2N+1,
\qquad
\frac{4\zeta}{\varepsilon\sqrt{1-\sigma^2}}\leq 1-\xi',
\label{sceltaparametri1}
\end{equation}
where
\begin{equation}
\qof\geq 4\lceil \log_2(80\delta^{-1}\zeta^{-1}A)\rceil\qquad\text{and}\qquad N\geq \frac{\log(\xi/(72\sqrt{\delta}D^22^k))}
{\log(1-\frac{\zeta^k}{40^{k+1}A^kD^2}(1-\sqrt{\delta}))}.
\label{sceltaparametri}
\end{equation}
Suppose moreover that $\mathscr E\subseteq \mathbb S^{n-1}$ is a finite set of directions and that $\mathfrak F\subseteq \Delta_\mu$ is a family of cubes contained in some open set $\Omega_\mathfrak s$, with $10AB_Q\subseteq \Omega_\mathfrak s$ for every $Q\in\mathfrak F$, and that
$$
\mathfrak F=\bigcup_{m=0}^{\mathfrak N}\mathfrak L_m.
$$
Let $U\subseteq \Omega_{\mathfrak s}$ be a measurable set.  Assume that the families $\mathfrak L_0,\ldots,\mathfrak L_{\mathfrak N}$ and the set $U$ satisfy the following properties:
\begin{itemize}
\item[(i)] No cube appears in more than one of the families $\mathfrak L_0,\ldots,\mathfrak L_{\mathfrak N}$.

\item[(ii)] For every $m=0,\ldots,\mathfrak N$, the cubes in $\mathfrak L_m$ are pairwise disjoint. Moreover, whenever $Q,Q'\in \mathfrak L_m$, either
$\ell(Q)=\ell(Q')$, or, if $\ell(Q)<\ell(Q')$, then
\begin{equation}
\frac{\ell(Q)}{\ell(Q')}\leq \frac{1}{2^\qof}.
\label{separazionediscale4}
\end{equation}

\item[(iii)] For every $m=1,\ldots,\mathfrak N$, each cube $Q'\in\mathfrak L_m$ is contained in a unique strictly larger cube $Q\in\mathfrak L_{m-1}$ such that
$$
\diam(Q)\geq 2^\qof\diam(Q').
$$
Further, whenever $m=1,\ldots,\mathfrak N$, $Q\in\mathfrak L_{m-1}$ and $Q'\in\mathfrak L_m$ satisfy
$10AB_{Q'}\cap 10AB_Q\neq \emptyset$, then
$$
\diam(Q)\geq 2^\qof\diam(Q').
$$

\item[(iv)] For every $m=1,\ldots,\mathfrak N$, every $Q\in \mathfrak L_m$, and every $e\in\mathscr E$, the cube $Q$ is not
$(e,\sigma,\varepsilon,\delta,\zeta,\mathfrak m,A)$-invariant, i.e. there exists a point $z(Q,e)\in AB_Q$ such that every $C(e,\sigma)$-curve $\gamma$ intersecting $B(z(Q,e),\zeta\diam Q)$ satisfies
$$
\Haus^1(2AB_Q\cap \gamma\setminus B(\supp\mu,\delta\diam Q))>\varepsilon\diam Q.
$$
\item[(v)] For every $m=1,\ldots,\mathfrak N$ and every $Q\in \mathfrak L_m$, if $Q'\in \mathfrak L_{m-1}$ is the unique cube containing $Q$, then
$$10A B_Q\cap \supp\mu\subseteq Q'.$$
\item[(vi)] For every $m=1,\ldots,\mathfrak N$, every $e\in\mathscr E$, and every $Q\in\mathfrak L_m$, we have
$$ \mu(B(z(Q,e),\zeta\diam Q/2)\cap U) \geq \tfrac12\mu(B(z(Q,e),\zeta\diam Q/2)). $$
\end{itemize}
Then, there exists an open set $\Omega\subseteq \Omega_{\mathfrak s}$ such that:
\begin{enumerate}
\item $\mu(U\setminus \Omega)\leq \mathrm{Card}(\mathscr E)\xi\mu(U)+2\mathrm{Card}(\mathscr E)\mu(U\setminus\bigcup\{Q:Q\in\mathfrak L_{\mathfrak N}\})$;
\item defining $\partial(\Omega,\mathfrak L_{2N+1})$ as the family of cubes $Q\in\mathfrak L_{2N+1}$ such that $10AB_Q\cap \Omega\neq \emptyset$ and $40AB_Q\not\subseteq \Omega$ we have 
\begin{equation}
\begin{split}
    &\mu(U\cap\bigcup\{10AB_Q:Q\in \partial(\Omega,\mathfrak L_{2N+1})\})\\
    &\qquad\qquad\qquad\qquad\leq\mathrm{Card}(\mathscr E)\xi\mu(U)+\mathrm{Card}(\mathscr E)\mu(U\setminus\bigcup\{Q:Q\in\mathfrak L_{\mathfrak N}\}).
\end{split}    
\end{equation}

\item there holds $\Omega\subseteq \bigcup_{Q\in\mathfrak L_1}\frac32 AB_Q$, and, for every $e\in\mathscr E$, and every $C(e,\sigma)$-curve $\gamma$, we have
$$
\Haus^1(\gamma\cap \Omega\cap 2AB_Q)
\leq
(1-\xi')\Haus^1(2(1+2^{-\qof})AB_Q\cap \gamma),\qquad\text{for every $Q\in \mathfrak L_1$.}
$$
\item For every $e\in\mathscr E$ and every $C(e,\sigma)$-curve $\gamma$, we have
$$
\Haus^1(\Omega\cap\gamma)
\leq
(1-\xi')\Haus^1(\Omega_\mathfrak s\cap\gamma).
$$
\end{enumerate}
\end{proposizione}

\begin{proof} We divide the proof into several steps. 

\textbf{Step I. Reduction to a single $e$.} First of all, notice that thanks to the choice of $\qof$, if $Q'\in\mathfrak L_{m+1}$ is contained in $Q\in\mathfrak L_m$, then 
\begin{equation} 80AB_{Q'}\subseteq B(Q,\delta\zeta\diam Q)\subseteq B(\supp\mu,\delta\zeta\diam Q). 
\label{eq:refined-layer-containment}
\end{equation}
Indeed, by the scale separation between consecutive layers, $\diam Q'\leq 2^{-\qof}\diam Q$ and the choice of $\qof$ gives 
$$ 80A\diam Q' \leq 80A2^{-\qof}\diam Q \leq \delta\zeta\diam Q. $$
This proves \eqref{eq:refined-layer-containment}. We shall also use the following estimate. Suppose $Q\in\mathfrak L_m$, $Q'\in\mathfrak L_{m+1}$ and 
$10AB_Q\cap 10AB_{Q'}\neq\emptyset$. Then, item (iii) gives
\begin{equation} 
\diam Q'\leq 2^{-\qof}\diam Q. 
\label{eq:refined-touching-scale}
\end{equation} 
In what follows, for every $e\in \mathscr E$ we shall construct an open set $\Omega_e\subseteq \Omega_{\mathfrak s}$ such that 
\begin{enumerate} 
\item[(\hypertarget{propomegaei}{i})] 
$\mu(U\setminus \Omega_e)\leq \xi\mu(U)+2\mu(U\setminus\bigcup\{Q:Q\in\mathfrak L_{\mathfrak N}\})$;
\item[(\hypertarget{propomegaeii}{ii})]
$\mu(U\cap\bigcup\{10AB_Q:Q\in \partial(\Omega_e,\mathfrak L_{2N+1})\})\leq \xi\mu(U)+\mu(U\setminus\bigcup\{Q:Q\in\mathfrak L_{\mathfrak N}\})$;
\item[(\hypertarget{propomegaeiii}{iii})] for every $C(e,\sigma)$-curve $\gamma$ we have 
$\Haus^1(\Omega_e\cap \gamma) \leq (1-\xi')\Haus^1(\Omega_\mathfrak s\cap \gamma)$; 
\item[(\hypertarget{propomegaeiv}{iv})] for every $Q_0\in\mathfrak L_1$ and every $C(e,\sigma)$-curve $\gamma$ we have 
$$\Haus^1(\gamma\cap \Omega_e\cap 2AB_{Q_0}) \leq (1-\xi')\Haus^1(2(1+2^{-\qof})AB_{Q_0}\cap \gamma).$$ 
\end{enumerate}
Defining $$ \Omega:=\bigcap_{e\in \mathscr E}\Omega_e, $$ we immediately obtain items 1, 3, and 4 of the proposition. 
Moreover, if $Q\in\partial(\Omega,\mathfrak L_m)$, then $10AB_Q\cap\Omega_e\neq\emptyset$ for every $e\in\mathscr E$, while $40AB_Q\not\subseteq\Omega_e$ for at least one $e\in\mathscr E$.
Hence
$$ \partial(\Omega,\mathfrak L_m)\subseteq \bigcup_{e\in\mathscr E}\partial(\Omega_e,\mathfrak L_m), $$
and item 2 follows as well. Therefore, from now on we fix $e\in\mathscr E$ and construct $\Omega_e$.

\medskip

\textbf{Step II. Construction of open sets with big measure and small width.}
In this step we fix $e\in\mathscr E$ and construct an open set $\Omega_e\subseteq\Omega_{\mathfrak s}$ satisfying \hyperlink{propomegaei}{(i)}, \hyperlink{propomegaeii}{(ii)}, \hyperlink{propomegaeiii}{(iii)}, and \hyperlink{propomegaeiv}{(iv)}. We introduce the following notation
$$
\mathfrak d:=\frac{\delta\zeta}{80A},\qquad
p_0:=\frac{\zeta^k}{40^{k+1}A^kD^2},\qquad
\theta:=(1-p_0)+p_0\sqrt{\delta}.
$$

\smallskip

\noindent\underline{\textsf{Base step of the construction of $\Omega_e$.}}
By Vitali's covering lemma, we choose a subfamily $\mathfrak B_1\subseteq\mathfrak L_2$ such that the balls $2AB_Q$, $Q\in\mathfrak B_1$, are pairwise disjoint, while the balls $10AB_Q$, $Q\in\mathfrak B_1$, cover $\bigcup\{AB_Q:Q\in\mathfrak L_2\}$. In particular, they cover all of $U$ up to the error $\mu(U\setminus\bigcup\{Q:Q\in\mathfrak L_{\mathfrak N}\})$.

We now choose the radius of the ball centered at $z(Q,e)$. The choice is made in two steps. The first application of \cref{lemmaraggipallefini} produces an auxiliary corona $\tilde A(Q,1)$ whose mass is controlled by the larger corona $A^\circ(Q,1)$. This auxiliary corona will be used to pay for the cubes of the next layer which intersect the boundary of the ball that we are constructing. The second application chooses the actual radius $r_{Q,1}$ inside the gap determined by $\tilde r_{Q,1}$ and produces the thinner corona $A(Q,1)$. This latter corona is the one that will be propagated in the induction when the old radii are enlarged. Thus $\tilde A(Q,1)$ controls boundary cubes, while $A(Q,1)$ controls the future loss produced by updating the radius. Let us proceed with the construction.

For each $Q\in\mathfrak B_1$, we first apply \cref{lemmaraggipallefini} to $z(Q,e)$ with $r:=\zeta\diam Q$, $\lambda:=4\sqrt{\delta}$, $\eta:=\delta$, $E_1:=U$ and $E_2:=U_1:=\bigcup\{R:R\in\mathfrak L_3\}$. We obtain $\tilde r_{Q,1}$ such that
$$
(1-2\sqrt{\delta})\zeta\diam Q\leq \tilde r_{Q,1}\leq (1-2\delta)\zeta\diam Q.
$$
We introduce some notation. Set
$$\tilde A(Q,1):=B\Big(z(Q,e),\frac{\tilde r_{Q,1}}{1-\delta}\Big)\setminus B(z(Q,e),\tilde r_{Q,1}),$$
and
$$A^\circ(Q,1):=B(z(Q,e),\zeta\diam Q)\setminus B(z(Q,e),(1-4\sqrt{\delta})\zeta\diam Q).$$
Then, there holds
$$
\mu(\tilde A(Q,1))\leq 3\sqrt{\delta}\mu(A^\circ(Q,1)),\quad
\mu(U\cap\tilde A(Q,1))\leq 3\sqrt{\delta}\mu(U\cap A^\circ(Q,1)),
$$
and
$$
\mu(U_1\cap\tilde A(Q,1))\leq 3\sqrt{\delta}\mu(U_1\cap A^\circ(Q,1)).
$$

We apply \cref{lemmaraggipallefini} a second time to $z(Q,e)$ with $r:=\tilde r_{Q,1}/(1-\delta)$, $\lambda:=\delta/36$, $\eta:=\delta\mathfrak d/144$, $E_1:=U$ and $E_2:=U_1$. This gives $r_{Q,1}$ such that
$$
(1+\frac{\delta}{72})\tilde r_{Q,1}\leq r_{Q,1}\leq (1-\frac{\delta\mathfrak d}{72})\frac{\tilde r_{Q,1}}{1-\delta}.
$$
Hence, defined  
$$
A(Q,1):=B\Big(z(Q,e),\frac{r_{Q,1}}{1-\frac{\delta\mathfrak d}{144}}\Big)\setminus B(z(Q,e),r_{Q,1}),
$$
we infer that
$$
\mu(A(Q,1))\leq 12\sqrt{\delta}\mathfrak d\,\mu(A^\circ(Q,1)),\quad
\mu(U\cap A(Q,1))\leq 12\sqrt{\delta}\mathfrak d\,\mu(U\cap A^\circ(Q,1)),
$$
and
$$
\mu(U_1\cap A(Q,1))\leq 12\sqrt{\delta}\mathfrak d\,\mu(U_1\cap A^\circ(Q,1)).
$$
We are ready to define the first open set. We let $\Omega_0:=\emptyset$, for notation's sake and we set
$$
\Omega_1:=\bigcup_{Q\in\mathfrak B_1}B(z(Q,e),r_{Q,1}).
$$
Since $z(Q,e)\in AB_Q$ and $r_{Q,1}\leq \zeta\diam Q\leq \diam Q$, the balls defining $\Omega_1$ are contained in the pairwise disjoint balls $2AB_Q$.

Let $Q'\in\partial(\Omega_1,\mathfrak L_3)$. Then $10AB_{Q'}$ meets $A(Q,1)$ for some $Q\in\mathfrak B_1$. Since $10AB_{Q'}\cap10AB_Q\neq\emptyset$, \eqref{eq:refined-touching-scale} gives $\diam Q'\leq 2^{-\qof}\diam Q\leq \mathfrak d^4\diam Q$. The choice of $r_{Q,1}$ then gives
$$
r_{Q,1}+80A\diam Q'<\frac{\tilde r_{Q,1}}{1-\delta}
\qquad\text{and}\qquad
r_{Q,1}-80A\diam Q'>\tilde r_{Q,1}.
$$
Thus
$40AB_{Q'}\subseteq \tilde A(Q,1)$. Consequently,
$$
\bigcup\{10AB_{Q'}:Q'\in\partial(\Omega_1,\mathfrak L_3)\}
\subseteq
\bigcup_{Q\in\mathfrak B_1}\tilde A(Q,1).
$$
This however implies that
\begin{equation}
\begin{split}
    &\mu(\bigcup\{10 AB_Q:Q\in \partial (\Omega_1,\mathfrak L_3)\})\leq \sum_{Q\in \mathfrak B_1}\mu(\tilde{A}(Q,1))
    \leq 3\sqrt{\delta}\sum_{Q\in \mathfrak B_1}\mu(A^\circ(Q,1))\\
&\leq3\sqrt{\delta}D\zeta^k\sum_{Q\in \mathfrak B_1}\diam Q^k\leq3\sqrt{\delta}D^22^{k+1} \sum_{Q\in \mathfrak B_1} \frac{1}{2}\mu(B(z(Q,e),\zeta \diam Q/2))\\
&\leq 3\sqrt{\delta}D^22^{k+1} \sum_{Q\in \mathfrak B_1}\mu(U\cap B(z(Q,e),\zeta \diam Q/2))\leq 36\sqrt{\delta}D^22^{k+1}\mu(U).
    \nonumber
\end{split}    
\end{equation}
Similarly we obtain the bounds 
\begin{equation}
\begin{split}
    &\mu(U\cap\bigcup\{10 AB_Q:Q\in \partial (\Omega_1,\mathfrak L_3)\})\leq 3\sqrt{\delta}\sum_{Q\in \mathfrak B_1}\mu(U\cap A^\circ(Q,1))\leq 3\sqrt{\delta}D^22^k\mu(U),\\
    &\mu(U_1\cap\bigcup\{10 AB_Q:Q\in \partial (\Omega_1,\mathfrak L_3)\})\leq 3\sqrt{\delta}\sum_{Q\in \mathfrak B_1}\mu(U_1\cap A^\circ(Q,1))\leq 3\sqrt{\delta}D^22^k\mu(U).
    \nonumber
\end{split}    
\end{equation}

Since $z(Q,e)\in AB_Q$ and $\tfrac{\zeta}{2}\diam Q\leq r_{Q,1}\leq \zeta \diam Q\leq \diam Q$, each ball $B(z(Q,e),r_{Q,1})$ is contained in $2AB_Q$. Hence the balls $\{B(z(Q,e),r_{Q,1})\}_{Q\in \mathfrak B_1}$ are pairwise disjoint. Therefore, by the AD-regularity of $\mu$,
\begin{equation}
    \begin{split}
        &\mu(\Omega_{1}\cap U)
        =\sum_{Q\in \mathfrak B_1}\mu(B(z(Q,e),r_{Q,1})\cap U)\geq 
       \frac{1}{2} \sum_{Q\in \mathfrak B_1}\mu(B(z(Q,e),\tfrac{\zeta}{2}\diam Q))\\
\geq& \frac{\zeta^k}{2^kD}\sum_{Q\in \mathfrak B_1}\diam Q^k
       \geq p_0\sum_{Q\in \mathfrak B_1}\mu(10AB_Q)\geq p_0\mu(U)-p_0\mu(U\setminus\bigcup\{Q:Q\in\mathfrak L_{\mathfrak N}\}).
        \label{primasinmasomega}
    \end{split}
\end{equation}
Hence
$$
\mu(U\setminus\Omega_1)\leq\theta\mu(U)+2\mu(U\setminus\bigcup\{Q:Q\in\mathfrak L_{\mathfrak N}\}).$$

Finally, let $\gamma$ be a $C(e,\sigma)$-curve and set
$$
\mathfrak B_1(\gamma):=\{Q\in\mathfrak B_1:\gamma\cap B(z(Q,e),r_{Q,1})\neq\emptyset\}.
$$
For $Q\in\mathfrak B_1(\gamma)$, define $E_Q:=2AB_Q\cap\gamma\setminus B(\supp\mu,\delta\diam Q)$. Since $Q$ is not $(e,\sigma,\varepsilon,\delta,\zeta,\mathfrak m,A)$-invariant, $\Haus^1(E_Q)>\varepsilon\diam Q$. Thus by \cref{propo:curvainpalla}, we have
$$
\Haus^1(\gamma\cap B(z(Q,e),r_{Q,1}))
\leq \frac{2r_{Q,1}}{\sqrt{1-\sigma^2}}
\leq \frac{2\zeta}{\varepsilon\sqrt{1-\sigma^2}}\Haus^1(E_Q).
$$
The sets $E_Q$, $Q\in\mathfrak B_1(\gamma)$, are disjoint and contained in $\gamma\cap\Omega_{\mathfrak s}$. Therefore
$$
\Haus^1(\gamma\cap\Omega_1)
\leq \frac{2\zeta}{\varepsilon\sqrt{1-\sigma^2}}\Haus^1(\gamma\cap\Omega_{\mathfrak s})
\leq (1-\xi')\Haus^1(\gamma\cap\Omega_{\mathfrak s}).
$$

\noindent\underline{\textsf{Inductive step of the construction of $\Omega_{\ell+1}$.}}
Let us assume that for some $\ell\geq 1$ we have already constructed families of cubes
$\mathfrak B_h\subseteq \mathfrak L_{2h}$ for $h=1,\ldots,\ell$, radii $r_{Q,\ell}$ for every $Q\in \mathfrak B_h$, and an open set
$$
\Omega_{\ell}:=\bigcup_{h=1}^{\ell}\bigcup_{Q\in \mathfrak B_h}B(z(Q,e),r_{Q,\ell}),
$$
in such a way that the following properties hold.

\begin{enumerate}
\item[(1)] For every $m=1,\ldots,\ell$, the balls $2AB_Q$, with $Q\in \mathfrak B_m$, are pairwise disjoint.

\item[(2)] For every $1\leq m\leq \ell$ and every $Q\in \mathfrak B_m$ there are parameters
$$
\eta_{h,m}:=\frac{\delta\mathfrak d}{144}\Big(\frac{\sqrt{\delta}\mathfrak d}{144}\Big)^{h-m},
\qquad h=m,\ldots,\ell,
$$
such that, for every $h=m,\ldots,\ell-1$,
\begin{equation}
r_{Q,h}\leq r_{Q,h+1}\leq \frac{1-\eta_{h+1,m}}{1-\eta_{h,m}}\,r_{Q,h}\leq \zeta \diam Q,
\nonumber
\end{equation}

\item[(3)] For every $1\leq m_1<m_2\leq \ell$, every $Q\in \mathfrak B_{m_1}$ and every $Q'\in \mathfrak B_{m_2}$, one has
\begin{equation}
B(z(Q,e),r_{Q,m_2})\cap 2AB_{Q'}=\emptyset.
\label{eq:interlayer-separation-h}
\end{equation}

\item[(4)] For every $1\leq m\leq \ell$, every $Q\in \mathfrak B_m$, and every $h=m,\ldots,\ell$, define
$$
A(Q,h):=
B\Bigl(z(Q,e),\frac{r_{Q,h}}{1-\eta_{h,m}}\Bigr)\setminus B(z(Q,e),r_{Q,h}).
$$
Then for every $h=m,\ldots,\ell-1$ one has $A(Q,h+1)\subseteq A(Q,h)$, and
\begin{equation}
\begin{split}
    \mu(A(Q,h+1))\leq& 4\mathfrak d\,\mu(A(Q,h)),\\
\mu(U\cap A(Q,h+1))\leq& 4\mathfrak d\,\mu(U\cap A(Q,h)).
\label{eq:shell-decay-h}
\end{split}
\end{equation}

\item[(5)] The open set $\Omega_{\ell}$ satisfies
$$
\mu(U\setminus\Omega_{\ell})\leq \theta^\ell\mu(U)+2\mu\Big(U\setminus\bigcup\{Q:Q\in\mathfrak L_{\mathfrak N}\}\Big),
$$
and for every $C(e,\sigma)$-curve $\gamma$ one has
\begin{equation}
\Haus^1(\gamma\cap \Omega_{\ell})
\leq(1-\xi')
\Haus^1(\gamma\cap \Omega_{\mathfrak s}).
\label{eq:inductive-length-loss-h}
\end{equation}
In addition, for every $Q\in \mathfrak L_1$ we have
$$\Haus^1(\gamma\cap \Omega_{\ell}\cap 2AB_{Q})\leq (1-\xi')\Haus^1(2(1+2^{-\qof})AB_Q\cap \gamma).$$
\item[(6)] There holds that
$$\mu(U\cap \bigcup\{10AB_Q:Q\in \partial(\Omega_\ell,\mathfrak L_{2\ell+1})\})\leq36\sqrt{\delta} D^22^k\theta^\ell\mu(U)+\mu(U\setminus\bigcup\{Q:Q\in\mathfrak L_{\mathfrak N}\}).$$
\end{enumerate}
We now construct the next open set $\Omega_{\ell+1}$, the family $\mathfrak B_{\ell+1}$ and the new radii $r_{Q,\ell+1}$.
As a first step, we update the old radii. 
Fix $1\leq m\leq \ell$ and $Q\in \mathfrak B_m$. We apply \cref{lemmaraggipallefini} to the point $z(Q,e)$ with
$$
r:=\frac{r_{Q,\ell}}{1-\eta_{\ell,m}},
\quad
\lambda:=\frac{\eta_{\ell,m}}{3},
\quad
\eta:=\frac{\sqrt{\delta}\eta_{\ell,m}}{12},\quad E_1=U,\quad E_2:=\bigcup\{Q:Q\in \mathfrak L_{2\ell+3}\}=:U_{\ell+1}.
$$
the lemma yields a radius $\tilde r_{Q,\ell+1}$ such that
$$
\Bigl(1+\frac{\eta_{\ell,m}}{6}\Bigr)r_{Q,\ell}
\leq
\frac{1-\eta_{\ell,m}/6}{1-\eta_{\ell,m}}\,r_{Q,\ell}
\leq
\tilde r_{Q,\ell+1}
\leq
\Bigl(1-\frac{\sqrt{\delta}\,\eta_{\ell,m}}{6}\Bigr)\frac{r_{Q,\ell}}{1-\eta_{\ell,m}},
$$
and denoted
$$\tilde{A}(Q,\ell+1):=B\Bigl(z(Q,e),\tilde r_{Q,\ell+1}/(1-\tfrac{\sqrt{\delta}\eta_{\ell,m}}{12})\Bigr)\setminus B\bigl(z(Q,e),\tilde r_{Q,\ell+1}\bigr),$$
we have
\begin{equation}
\begin{split}
\mu(\tilde{A}(Q,\ell+1))\leq&
3\sqrt{\delta}
\mu\Big(B\Big(z(Q,e),\frac{r_{Q,\ell}}{1-\eta_{\ell,m}}\Big)\setminus B(z(Q,e),r_{Q,\ell})\Big)=3\sqrt{\delta}\mu(A(Q,\ell)),
\label{stmimgfj}
\end{split}
\end{equation}
and similarly 
\begin{equation}
\begin{split}
\mu(U\cap \tilde{A}(Q,\ell+1))\leq&3\sqrt{\delta}\mu(U\cap A(Q,\ell)),\\
\mu(U_{\ell+1}\cap \tilde{A}(Q,\ell+1))\leq&3\sqrt{\delta}\mu(U_{\ell+1}\cap A(Q,\ell)),
\label{stmimgfj2}
\end{split}
\end{equation}
The corona $\tilde A(Q,\ell+1)$ is a controlled subcorona of the old corona $A(Q,\ell)$. It will be used to estimate the mass of the cubes of the next relevant layer which meet the boundary of the ball obtained after the radius of $Q$ has been updated. We now choose the final radius $r_{Q,\ell+1}$ with an additional gap, so that a small neighbourhood of the new boundary is still contained in $\tilde A(Q,\ell+1)$. This produces a second controlled corona, $A(Q,\ell+1)$, which will be used in the next steps of the induction when the radius is enlarged again. Thus we apply again \cref{lemmaraggipallefini} to the point $z(Q,e)$ with
$$
r:=\frac{\tilde{r}_{Q,\ell+1}}{1-\tfrac{\sqrt{\delta}\eta_{\ell,m}}{12}},
\quad
\lambda:=\frac{\eta}{3}=\frac{\sqrt{\delta}\eta_{\ell,m}}{36},
\quad
\eta:=\frac{\mathfrak d}{4}\lambda=\eta_{\ell+1,m},\quad E_1=U,\quad E_2:=U_{\ell+1}.
$$
With these choices we find a radius  $r_{Q,\ell+1}$ satisfying
\begin{equation}
    \begin{split}
    (1+\frac{\sqrt{\delta}\eta_{\ell,m}}{36})\tilde{r}_{Q,\ell+1}\leq (1-\lambda/2)r\leq r_{Q,\ell+1}\leq (1-2\eta)r=\frac{1-\frac{\sqrt{\delta}\mathfrak d}{72}\eta_{\ell,m}}{1-\tfrac{\sqrt{\delta}\eta_{\ell,m}}{12}}\tilde{r}_{Q,\ell+1}.
    \label{stimadaaosaoisgdn}
    \end{split}
\end{equation}
and by \eqref{stmimgfj} and \eqref{stmimgfj2}, defining 
$$A(Q,\ell+1):=B\Bigl(z(Q,e),\frac{ r_{Q,\ell+1}}{1-\eta_{\ell+1,m}}\Bigr)\setminus B\bigl(z(Q,e), r_{Q,\ell+1}\bigr),$$
we have
\begin{align*}
    \mu(A(Q,\ell+1))\leq &36\sqrt{\delta}\mathfrak d \mu(A(Q,\ell)),\\
    \mu(U\cap A(Q,\ell+1))\leq &36\sqrt{\delta}\mathfrak d \mu(U\cap A(Q,\ell)),\\
    \mu(U_{\ell+1}\cap 
    A(Q,\ell+1))\leq& 36\sqrt{\delta}\mathfrak d \mu(U_{\ell+1}\cap A(Q,\ell)).
\end{align*}
Notice that, we still have
$r_{Q,\ell+1}\leq \zeta \diam Q$.

\smallskip

Now we proceed with constructing the new generation of balls of $\Omega_{\ell+1}$.
Let $Q\in \mathfrak B_m$ for some $1\leq m\leq \ell$, and let $Q'\in \mathfrak L_{2(\ell+1)}$ satisfy
$$
10AB_Q\cap 10AB_{Q'}\neq \emptyset.
$$
By iterating the scale separation between the layers, we have
$\diam Q'\leq \bigl(2^{-\qof}\bigr)^{2(\ell+1-m)}\diam Q$ and since  
$2^{-\qof}\leq (\delta\zeta/80A)^4=\mathfrak d^4$,
it follows that
\begin{equation}
\diam Q'\leq \mathfrak d^{8(\ell+1-m)}\diam Q.
\label{eq:diamQprime-h}
\end{equation}
In the notations above, we claim that if
$10AB_{Q'}\cap B(z(Q,e),r_{Q,\ell+1})\neq \emptyset$ and $10AB_{Q'}\not\subseteq  B(z(Q,e),r_{Q,\ell+1})$
then
\begin{equation}
    40A B_{Q'}\subseteq \tilde{A}(Q,\ell+1).
    \label{inclusionebordizx}
\end{equation}
In order to see this, we notice that since 
\begin{equation}
    \begin{split}
        80A\diam Q'\leq&80 A\frac{\delta^8   \zeta^8}{(80A)^8} \mathfrak d^{8(\ell-m)}\diam Q
        \leq\zeta^7\frac{\delta^8 }{(80A)^7}\mathfrak d^{8(\ell-m)}\zeta\diam Q\\
        \leq& \frac{\zeta^5\delta^5}{(80A)^5}\mathfrak d^{6(\ell-m)}\eta_{\ell,m}r_{Q,\ell}<\frac{\zeta^5\delta^5}{(80A)^5}\eta_{\ell,m}r_{Q,\ell},
        \label{stimadiametriextrema}
    \end{split}
\end{equation}
Hence, this implies that 
$$r_{Q,\ell+1}+80 A\diam Q'\leq r_{Q,\ell+1}+\frac{\zeta^5\delta^5}{(80A)^5}\eta_{\ell,m}r_{Q,\ell}<\frac{r_{Q,\ell}}{1-\frac{\sqrt{\delta}\eta_{\ell,m}}{12}}<\frac{\tilde r_{Q,\ell+1}}{1-\frac{\sqrt{\delta}\eta_{\ell,m}}{12}},$$
where the last inequality follows from the fact that $r_{Q,\ell}<r_{Q,\ell+1}$ and
\begin{equation}
    \begin{split}
    r_{Q,\ell+1}+\frac{\zeta^5\delta^5}{(80A)^5}\eta_{\ell,m}r_{Q,\ell}\leq \frac{1-\frac{\sqrt{\delta}\mathfrak d}{72}\eta_{\ell,m}}{1-\tfrac{\sqrt{\delta}\eta_{\ell,m}}{12}}\tilde{r}_{Q,\ell+1}+\frac{\zeta^5\delta^5}{2(80A)^5}\eta_{\ell,m}\tilde r_{Q,\ell+1}<\frac{1}{1-\tfrac{\sqrt{\delta}\eta_{\ell,m}}{12}}\tilde{r}_{Q,\ell+1}
    \end{split}
\end{equation}
On the other hand, \eqref{stimadaaosaoisgdn} and \eqref{stimadiametriextrema} implies also that
\begin{equation}
    \begin{split}
        r_{Q,\ell+1}-80 A\diam Q'\geq &r_{Q,\ell+1}-\frac{\zeta^5\delta^5}{(80A)^5}\eta_{\ell,m}r_{Q,\ell}\\
        \geq& (1+\frac{\sqrt{\delta}\eta_{\ell+1,m}}{36})\tilde{r}_{Q,\ell+1}-\frac{\zeta^5\delta^5}{(80A)^5}\eta_{\ell,m}r_{Q,\ell}>\tilde{r}_{Q,\ell+1}.
        \nonumber
    \end{split}
\end{equation}
This finally proves that \eqref{inclusionebordizx} holds. 
We now define $\mathfrak B'_{\ell+1}$ to be the family of cubes $Q'\in \mathfrak L_{2(\ell+1)}$ such that
$$
10AB_{Q'}\cap B(z(Q,e),r_{Q,\ell+1})=\emptyset
$$
for every $1\leq m\leq \ell$ and every $Q\in \mathfrak B_m$.
Applying Vitali's covering lemma, we find a subfamily
$\mathfrak B_{\ell+1}\subseteq \mathfrak B'_{\ell+1}$ such that the balls $2AB_Q$, $Q\in \mathfrak B_{\ell+1}$, are pairwise disjoint, while the balls $10AB_Q$, $Q\in \mathfrak B_{\ell+1}$, cover
$$
\bigcup\{2AB_Q:Q\in \mathfrak B'_{\ell+1}\}\cap \supp\mu.
$$

We are now ready to introduce the new generation of $\Omega_{\ell+1}$. For every $Q\in\mathfrak B_{\ell+1}$ we repeat the two-step choice obtained by means of \cref{lemmaraggipallefini}, with $U_{\ell+1}$ in place of $U_1$. Thus we obtain radii $\tilde r_{Q,\ell+1}$ and $r_{Q,\ell+1}$ such that
$$
(1-2\sqrt{\delta})\zeta\diam Q\leq \tilde r_{Q,\ell+1}\leq (1-2\delta)\zeta\diam Q
$$
$$
(1+\frac{\delta}{72})\tilde r_{Q,\ell+1}\leq r_{Q,\ell+1}\leq (1-\frac{\delta\mathfrak d}{72})\frac{\tilde r_{Q,\ell+1}}{1-\delta}.
$$
Defined
$$
\tilde A(Q,\ell+1):=B\Big(z(Q,e),\frac{\tilde r_{Q,\ell+1}}{1-\delta}\Big)\setminus B(z(Q,e),\tilde r_{Q,\ell+1}),
$$
$$
A^\circ(Q,\ell+1):=B(z(Q,e),\zeta\diam Q)\setminus B(z(Q,e),(1-4\sqrt{\delta})\zeta\diam Q),
$$
$$
A(Q,\ell+1):=B\Big(z(Q,e),\frac{r_{Q,\ell+1}}{1-\frac{\delta\mathfrak d}{144}}\Big)\setminus B(z(Q,e),r_{Q,\ell+1}),
$$
arguing as in the base step and the same choices of constants, we have 
$$
\mu(\tilde A(Q,\ell+1))\leq 3\sqrt{\delta}\mu(A^\circ(Q,\ell+1)),$$
$$\mu(U\cap\tilde A(Q,\ell+1))\leq 3\sqrt{\delta}\mu(U\cap A^\circ(Q,\ell+1)),$$
$$
\mu(U_{\ell+1}\cap\tilde A(Q,\ell+1))\leq 3\sqrt{\delta}\mu(U_{\ell+1}\cap A^\circ(Q,\ell+1)).
$$
and similarly
$$\mu(A(Q,\ell+1))\leq 12\sqrt{\delta}\mathfrak d\mu(A^\circ(Q,\ell+1)),$$
$$\mu(U\cap A(Q,\ell+1))\leq 12\sqrt{\delta}\mathfrak d\mu(U\cap A^\circ(Q,\ell+1)),$$
$$\mu(U_{\ell+1}\cap A(Q,\ell+1))\leq 12\sqrt{\delta}\mathfrak d\mu(U_{\ell+1}\cap A^\circ(Q,\ell+1)).$$
Finally, arguing exactly as in the base step, if $Q'\in\mathfrak L_{2\ell+3}$ and, for some $Q\in\mathfrak B_{\ell+1}$, we have
$10AB_{Q'}\cap A(Q,\ell+1)\neq\emptyset$, then
$40AB_{Q'}\subseteq \tilde A(Q,\ell+1)$ and we define
$$
\Omega_{\ell+1}:=\bigcup_{m=1}^{\ell+1}\bigcup_{Q\in\mathfrak B_m}B(z(Q,e),r_{Q,\ell+1}).
$$
We now check that the induction hypotheses are satisfied. Items (1), (2), (3), and (4) follow from the construction. We now verify item (5), namely the measure gain and the length loss, and then item (6).

\smallskip

\noindent\textsf{$\bullet$ Measure gain.}
In order to simplify notation, we  introduce the auxiliary open sets
$$
\widetilde{\Omega}_{0,\ell+1}
:=
\bigcup_{m=1}^{\ell}\bigcup_{Q\in \mathfrak B_m}
B\Bigl(z(Q,e),\frac{r_{Q,\ell+1}}{1-\eta_{\ell+1,m}}\Bigr)
\qquad\text{and}\qquad
\Omega^\star_{0,\ell+1}
:=
\bigcup_{m=1}^{\ell}\bigcup_{Q\in \mathfrak B_m}
B(z(Q,e),r_{Q,\ell+1}).
$$
As a first step, we prove that
\begin{equation}
    \mu(U\setminus \Omega_{\ell+1})
    \leq
    (1-p_0)\mu(U\setminus \Omega_{\ell})
+p_0\,\mu(U\cap(\widetilde{\Omega}_{0,\ell+1}\setminus \Omega_{\ell}))
+p_0\mu(U\setminus\bigcup\{Q:Q\in\mathfrak L_{\mathfrak N}\}).
    \label{iterazimigliU-h}
\end{equation}
We claim that
\begin{equation}
\mu(U\setminus\widetilde{\Omega}_{0,\ell+1})\leq \mu(\bigcup\{Q:Q\in \mathfrak B'_{\ell+1}\})+\mu(U\setminus\bigcup\{Q:Q\in\mathfrak L_{\mathfrak N}\}).
\label{eq:asodifaosidfnaosidnfasodighss1241}
\end{equation}
Indeed, it is enough to prove the inclusion outside the terminal defect. Let
$$
x\in U\cap\bigcup\{Q:Q\in\mathfrak L_{\mathfrak N}\}
\quad\text{and}\quad
x\notin\widetilde{\Omega}_{0,\ell+1}.
$$
Let $Q'\in\mathfrak L_{2(\ell+1)}$ be the ancestor of the terminal cube containing $x$. If $Q'\notin\mathfrak B'_{\ell+1}$, then there are $1\leq m\leq\ell$ and $Q\in\mathfrak B_m$ such that
$$
10AB_{Q'}\cap B(z(Q,e),r_{Q,\ell+1})\neq\emptyset.
$$
If $10AB_{Q'}\subseteq B(z(Q,e),r_{Q,\ell+1})$, then $x\in\Omega^\star_{0,\ell+1}\subseteq\widetilde{\Omega}_{0,\ell+1}$, a contradiction. Otherwise, by \eqref{inclusionebordizx}, we have
$$
40AB_{Q'}\subseteq\tilde A(Q,\ell+1)\subseteq\widetilde{\Omega}_{0,\ell+1},
$$
again a contradiction. Hence $Q'\in\mathfrak B'_{\ell+1}$, and this proves the claim. On the other hand, by construction of $\mathfrak B_{\ell+1}$, the family
$\{10AB_Q:Q\in \mathfrak B_{\ell+1}\}$ covers the set
$\bigcup\{2AB_Q:Q\in \mathfrak B'_{\ell+1}\}\cap \supp\mu$.
Hence, defined 
$$\Omega_{0,\ell+1}^{\ell+1}:=\bigcup_{Q\in \mathfrak B_{\ell+1}}B(z(Q,e),r_{Q,\ell+1}),$$
thanks to the AD-regularity of $\mu$, we infer that 
\begin{equation}
    \begin{split}
    \mu(U\cap  \Omega_{0,\ell+1}^{\ell+1})
        \geq  &\sum_{Q\in \mathfrak B_{\ell+1}}\mu(U\cap B(z(Q,e),\zeta \diam Q/2))
        \geq \frac{\zeta^k}{4^{k+1}D}\sum_{Q\in \mathfrak B_{\ell+1}}(\diam Q)^k\\
        \geq& \frac{\zeta^k}{40^{k+1}A^kD^2}\sum_{Q\in \mathfrak B_{\ell+1}}\mu(10AB_Q)\geq
p_0\mu(\bigcup\{Q:Q\in \mathfrak B'_{\ell+1}\})\\
\overset{\eqref{eq:asodifaosidfnaosidnfasodighss1241}}{\geq}& p_0\mu(U\setminus \widetilde{\Omega}_{0,\ell+1})-p_0\mu(U\setminus\bigcup\{Q:Q\in\mathfrak L_{\mathfrak N}\}).
    \end{split}
    \label{eq:gain-new-layer-U-h}
\end{equation}
Since $\Omega_{\ell}\subseteq \widetilde{\Omega}_{0,\ell+1}$, we can write
\begin{equation}
    \mu(U\setminus \widetilde{\Omega}_{0,\ell+1})
    \geq
    \mu(U\setminus \Omega_{\ell})
    -
\mu(U\cap(\widetilde{\Omega}_{0,\ell+1}\setminus \Omega_{\ell})).
    \label{eq:UminuswidetildeU-h}
\end{equation}
Plugging \eqref{eq:UminuswidetildeU-h} into \eqref{eq:gain-new-layer-U-h}, we deduce that
\begin{equation}
\mu(U\cap \Omega_{0,\ell+1}^{\ell+1})\geq p_0\mu(U\setminus \Omega_{\ell})-p_0\mu(U\cap(\widetilde{\Omega}_{0,\ell+1}\setminus\Omega_{\ell}))-p_0\mu(U\setminus\bigcup\{Q:Q\in\mathfrak L_{\mathfrak N}\}).
    \label{eq:gain-new-layer-U-v2-h}
\end{equation}

Now notice that $\Omega_{\ell}\subseteq \Omega^\star_{0,\ell+1}$ and $\Omega^\star_{0,\ell+1}\cap \Omega_{0,\ell+1}^{\ell+1}=\emptyset$.
Therefore
\begin{equation}
    \begin{split}
      &\mu(U\setminus \Omega_{\ell+1})=\mu(U\setminus \Omega^\star_{0,\ell+1})-\mu(U\cap \Omega_{0,\ell+1}^{\ell+1})\leq\mu(U\setminus \Omega_{\ell})-\mu(U\cap\Omega_{0,\ell+1}^{\ell+1})\\
       &\overset{\eqref{eq:gain-new-layer-U-v2-h}}{\leq}
      (1-p_0)\mu(U\setminus \Omega_{\ell})+p_0\mu(U\cap(\widetilde{\Omega}_{0,\ell+1}\setminus \Omega_{\ell}))+p_0\mu(U\setminus\bigcup\{Q:Q\in\mathfrak L_{\mathfrak N}\}),
    \end{split}
\end{equation}
thus concluding the proof of  \eqref{iterazimigliU-h}.
We now estimate the remainder
$r_\ell:=\mu(U\cap(\widetilde{\Omega}_{0,\ell+1}\setminus \Omega_{\ell}))$. Let us notice that 
\begin{equation}
    \begin{split}
\widetilde{\Omega}_{0,\ell+1}\setminus \Omega_{\ell}=
(\widetilde{\Omega}_{0,\ell+1}\setminus \Omega^\star_{0,\ell+1})
\cup(\Omega^\star_{0,\ell+1}\setminus\Omega_{\ell})\subseteq\bigcup_{m=1}^{\ell}\bigcup_{Q\in \mathfrak B_m}
        A(Q,\ell+1)\cup A(Q,\ell),
        \nonumber
    \end{split}
\end{equation}
from which
\begin{equation}
    \begin{split}
        r_\ell
        \leq
        \sum_{m=1}^{\ell}\sum_{Q\in \mathfrak B_m}\mu(U\cap A(Q,\ell+1))
        +
        \sum_{m=1}^{\ell}\sum_{Q\in \mathfrak B_m}\mu(U\cap A(Q,\ell)).
    \end{split}
    \label{eq:rllbound0-h}
\end{equation}
Since the cube $Q$ is first selected at generation $m$, by construction of $r_{Q,m}$ we have
$$
A(Q,m)=
B\Bigl(z(Q,e),\frac{r_{Q,m}}{1-\tfrac{\delta\mathfrak d}{144}}\Bigr)\setminus B(z(Q,e),r_{Q,m}),
$$
and therefore, the choice of $r_{Q,m}$ yields
\begin{equation}
    \mu(U\cap A(Q,m))\leq 12\sqrt{\delta}\mathfrak d\mu(U\cap A^\circ(Q,m)).
    \label{eq:initial-U-annulus-h}
\end{equation}
Notice now that $A^\circ(Q,m)\subseteq 2AB_Q$ for every $Q\in \mathfrak B_m$, since $z(Q,e)\in AB_Q$ and $\zeta\leq 1\leq A$. In addition, if $Q,Q'\in \mathfrak B_m$ are distinct, then by construction
$2AB_Q\cap 2AB_{Q'}=\emptyset$
, and hence
$A^\circ(Q,m)\cap  A^\circ(Q',m)=\emptyset$. Finally, for every $Q\in \mathfrak B_m$, again by construction we have
$2AB_Q\cap \Omega_{m-1}=\emptyset$.
Putting together the three observations above, we infer that
\begin{equation}
    \sum_{Q\in \mathfrak B_m}\mu(U\cap A^\circ(Q,m))
    \leq \mu(U\setminus \Omega_{m-1}).
    \label{eq:sum-initial-U-annuli-h}
\end{equation}
Combining \eqref{eq:initial-U-annulus-h} and \eqref{eq:sum-initial-U-annuli-h}, we obtain
\begin{equation}
    \sum_{Q\in \mathfrak B_m}\mu(U\cap A(Q,m))
    \leq 12\sqrt{\delta}\mathfrak d\mu(U\setminus \Omega_{m-1}).
    \label{eq:initial-shell-sum-h}
\end{equation}
By \eqref{eq:shell-decay-h}, for every $Q\in \mathfrak B_m$ and every $\ell\geq m$ we have
\begin{equation}
    \mu(U\cap A(Q,\ell))\leq (4\mathfrak d)^{\ell-m}\mu(U\cap A(Q,m)).
    \label{eq:U-annulus-decay-step-h}
\end{equation}
Hence, using \eqref{eq:initial-shell-sum-h}, we infer that for every $\ell\geq m$
\begin{equation}
    \sum_{Q\in \mathfrak B_m}\mu(U\cap A(Q,\ell))
    \leq 3\sqrt{\delta} (4\mathfrak d)^{\ell+1-m}\mu(U\setminus \Omega_{m-1}).
    \label{eq:sum-U-annulus-ellplusone-h}
\end{equation}
Therefore, by \eqref{eq:rllbound0-h}, we obtain the sought estimate for $r_\ell$
\begin{equation}
    \begin{split}
        r_\ell
        \leq
        12\sqrt{\delta}(1+4\mathfrak d)
        \sum_{m=1}^{\ell}(4\mathfrak d)^{\ell+1-m}\mu(U\setminus \Omega_{m-1}).
    \end{split}
    \label{eq:rllbound-h}
\end{equation}
Plugging \eqref{eq:rllbound-h} into \eqref{iterazimigliU-h}, we finally obtain
\begin{equation}
\begin{split}
    \mu(U\setminus \Omega_{\ell+1})\leq(1-p_0)\mu(U\setminus \Omega_{\ell})&+12p_0\sqrt{\delta}(1+4\mathfrak d)\sum_{m=1}^{\ell}(4\mathfrak d)^{\ell+1-m}\mu(U\setminus \Omega_{m-1})\\
&+p_0\mu(U\setminus\bigcup\{Q:Q\in\mathfrak L_{\mathfrak N}\}).
    \label{eq:closed-recursion-X-h}
\end{split}
\end{equation}
Let us now prove that 
\begin{equation}
    \mu(U\setminus \Omega_{\ell+1})\leq\theta^{\ell+1}\mu(U)+2\mu(U\setminus\bigcup\{Q:Q\in\mathfrak L_{\mathfrak N}\}).
    \label{exponentialdecaya_n}
\end{equation}
Indeed, by \eqref{primasinmasomega}, we notice that 
\begin{equation}
    \begin{split}
        \mu(U\setminus \Omega_{1})&=
        \mu(U)-\mu(U\cap \Omega_1)\leq
        (1-p_0)\mu(U)+p_0\mu(U\setminus\bigcup\{Q:Q\in\mathfrak L_{\mathfrak N}\})\\
        &\leq\theta\mu(U)+2\mu(U\setminus\bigcup\{Q:Q\in\mathfrak L_{\mathfrak N}\}).
        \nonumber
    \end{split}
\end{equation}
Secondly, let us suppose that \eqref{exponentialdecaya_n} holds for $\mu(U\setminus \Omega_m)$ for every $1\leq m\leq \ell$. Then, by \eqref{eq:closed-recursion-X-h}, we have
\begin{equation}
    \begin{split}
&\mu(U\setminus \Omega_{\ell+1})\leq(1-p_0)\theta^\ell\mu(U)+12p_0\sqrt{\delta}(1+4\mathfrak d)\Big(\sum_{m=1}^\ell(4\mathfrak d)^{\ell+1-m}\theta^{m-1}\Big)\mu(U)\\
&\qquad\qquad\qquad\qquad+2(1-p_0)\mu(U\setminus\bigcup\{Q:Q\in\mathfrak L_{\mathfrak N}\})+p_0\mu(U\setminus\bigcup\{Q:Q\in\mathfrak L_{\mathfrak N}\})\\
&\qquad\qquad\qquad\qquad+24p_0\sqrt{\delta}(1+4\mathfrak d)\sum_{m=1}^\ell(4\mathfrak d)^{\ell+1-m}\mu(U\setminus\bigcup\{Q:Q\in\mathfrak L_{\mathfrak N}\})\\
&\qquad\qquad\quad\leq\theta^{\ell+1}\mu(U)+2\mu(U\setminus\bigcup\{Q:Q\in\mathfrak L_{\mathfrak N}\}).
    \end{split}
\end{equation}
where in the last inequality we used that $\mathfrak d\leq \delta\leq 2^{-24}$ and hence
$$
12\sqrt{\delta}(1+4\mathfrak d)(4\mathfrak d)\frac{1}{\theta-4\mathfrak d}\leq \sqrt{\delta}
\qquad\text{and}\qquad
24\sqrt{\delta}(1+4\mathfrak d)\frac{4\mathfrak d}{1-4\mathfrak d}\leq 1.
$$
This concludes the proof of the measure gain.

\smallskip

\noindent\textsf{$\bullet$ Length loss.}
Let $\gamma$ be a $C(e,\sigma)$-curve. 
We let
$$
\mathfrak B(\gamma):=
\bigcup_{m=1}^{\ell+1}\{Q\in \mathfrak B_m:\gamma\cap B(z(Q,e),r_{Q,\ell+1})\neq\emptyset\},
$$
and for every $Q\in \mathfrak B(\gamma)$ we define
$$E_Q:=\gamma\cap 2AB_Q\setminus B(\supp\mu,\delta \diam Q).$$
Since every cube in every $\mathfrak B_m$ is not
$(e,\sigma,\varepsilon,\delta,\zeta,\mathfrak m,A)$-invariant, we have that $\Haus^1(E_Q)>\varepsilon \diam Q$ for every $Q\in \mathfrak B(\gamma)$. On the other hand, by \cref{propo:curvainpalla} we have
\begin{equation}
    \begin{split}
        \Haus^1(\gamma\cap B(z(Q,e),r_{Q,\ell+1}))\leq \frac{2r_{Q,\ell+1}}{\sqrt{1-\sigma^2}}\leq \frac{2\zeta}{\sqrt{1-\sigma^2}}\diam Q\leq \frac{2\zeta}{\varepsilon\sqrt{1-\sigma^2}}\Haus^1(E_Q),
        \label{eq:length-vs-bad-piece}
    \end{split}
\end{equation}
for every $Q\in \mathfrak B(\gamma)$. Summing \eqref{eq:length-vs-bad-piece} over all $Q\in \mathfrak B(\gamma)$, we infer that 
$$\Haus^1(\gamma\cap \Omega_{\ell+1})\leq \sum_{Q\in \mathfrak B(\gamma)}\Haus^1(B(z(Q,e),r_{Q,\ell+1})\cap \gamma)\leq \frac{2\zeta}{\varepsilon\sqrt{1-\sigma^2}}\sum_{Q\in \mathfrak B(\gamma)}\Haus^1(E_Q).$$
For every $x\in \gamma$ we let 
$\mathfrak Q^\flat(x):=\{Q\in \mathfrak B(\gamma):x\in E_Q\}$,
and we notice that
$$\sum_{Q\in \mathfrak B(\gamma)}\Haus^1(E_Q)=\int \Big(\sum_{Q\in \mathfrak B(\gamma)}\mathbb{1}_{E_Q}(z)\Big)d\Haus^1\trace \gamma(z)=\int \mathrm{Card}(\mathfrak Q^\flat(z))d\Haus^1\trace\Big( \bigcup_{Q\in \mathfrak B(\gamma)}E_Q \Big)(z).$$
From this, we immediately see that 
$$\sum_{Q\in \mathfrak B(\gamma)}\Haus^1(E_Q)\leq M \Haus^1(\bigcup_{Q\in \mathfrak B(\gamma)}E_Q)\leq M\Haus^1(\gamma \cap \Omega_{\mathfrak s}),$$
where $M:=\sup_{z\in \bigcup_{Q\in \mathfrak B(\gamma)}E_Q}\mathrm{Card}(\mathfrak Q^\flat(z))$ and where the last inequality comes from the fact that by construction we have 
$$2AB_Q\subseteq \Omega_{\mathfrak s}\qquad\text{for every $m=1,\ldots,\ell+1$ and every $Q\in \mathfrak B_m$}.$$
In order to give a bound to $M$, we fix $z\in \bigcup_{Q\in \mathfrak B(\gamma)}E_Q$ and we observe that if $Q\in \mathfrak Q^\flat(z)$, then  $z\in 2A B_Q$ and 
$$\frac{\dist(z,\supp\mu)}{2A}\leq \diam Q\leq \frac{\dist(z,\supp\mu)}{\delta}.$$
Recall that for every $\mathfrak B_m$ the cubes satisfy $2AB_S\cap 2AB_{S'}=\emptyset$ for every $S,S'\in \mathfrak B_m$ and every $m=1\ldots,\ell+1$, hence for every $z$ and $Q$ as above, for each generation, there is at most one cube belonging to $\mathfrak Q^\flat(z)$. Thanks to the scale separation between the layers, we thus infer that 
$$M\leq 1+\frac{1}{\qof}\log_2(2A/\delta)\leq 2,$$
thanks to the choice of $\qof$.
From this we finally infer thanks to the choice of the parameters that 
$$\Haus^1(\gamma\cap \Omega_{\ell+1})\leq \frac{4\zeta}{\varepsilon\sqrt{1-\sigma^2}}\Haus^1(\gamma \cap \Omega_{\mathfrak s})\leq (1-\xi')\Haus^1(\gamma \cap \Omega_{\mathfrak s}).$$
Let us now show the small intersection property with the balls relative to $\mathfrak L_1$. 
Let us now pick $Q_0\in \mathfrak L_1$ and a $C(e,\sigma)$-curve $\gamma$ and let us define
$$
\mathfrak B(\gamma;Q_0):=
\bigcup_{m=1}^{\ell+1}\{Q\in \mathfrak B_m:\gamma\cap 2AB_{Q_0}\cap B(z(Q,e),r_{Q,\ell+1})\neq\emptyset\}.$$
Notice that, with the same argument we used above we infer that
$$ \Haus^1(\gamma\cap \Omega_{\ell+1}\cap 2AB_{Q_0})\leq \frac{4\zeta}{\varepsilon \sqrt{1-\sigma^2}} \Haus^1\Big(\bigcup_{Q\in \mathfrak B(\gamma;Q_0)}E_Q\Big).$$
Notice that if $Q\in\mathfrak B(\gamma;Q_0)$, then
$B(z(Q,e),r_{Q,\ell+1})\cap 2AB_{Q_0}\neq\emptyset$ and since  $z(Q,e) \in AB_Q$, we deduce that $2AB_Q\cap 2AB_{Q_0}\neq\emptyset$.
Moreover $Q\in\mathfrak L_{2m}$ for some $m\geq1$, while $Q_0\in\mathfrak L_1$. By the refined scale separation from Step I, we get that $\diam Q\leq 2^{-\qof}\diam Q_0$.
Therefore
$$
2AB_Q\subseteq 2(1+2^{-\qof})AB_{Q_0}.
$$
Consequently, we obtain thanks to the choice of the parameters that
$$
\Haus^1(\gamma\cap \Omega_{\ell+1}\cap 2AB_{Q_0})
\leq
\tfrac{4\zeta}{\varepsilon\sqrt{1-\sigma^2}}
\Haus^1(\gamma\cap 2(1+2^{-\qof})AB_{Q_0})\leq
(1-\xi')\Haus^1(\gamma\cap 2(1+2^{-\qof})AB_{Q_0}).
$$
This proves the desired small intersection estimate for every $Q_0\in\mathfrak L_1$. This concludes the proof of item (5).

\medskip

Let us conclude the induction with the proof of item (6). If
$Q'\in \partial(\Omega_\ell,\mathfrak L_{2\ell+1})$, then by the construction of the radii we have
\begin{equation}
    \begin{split}
        U\cap \bigcup\{10AB_Q:Q\in \partial(\Omega_\ell,\mathfrak L_{2\ell+1})\}
        \subseteq
        \bigcup_{m=1}^{\ell}\bigcup_{Q\in \mathfrak B_m}(U\cap \tilde A(Q,\ell)).
    \end{split}
\end{equation}

Fix now $m\in\{1,\ldots,\ell\}$ and $Q\in \mathfrak B_m$.
Notice, that since for every $m=1,\ldots,\ell$, every $Q\in \mathfrak B_m$ there holds
\begin{equation}
\mu(U\cap \tilde A(Q,\ell))
\leq
3\sqrt{\delta}\,(36\sqrt{\delta}\mathfrak d)^{\ell-m}\mu(U\cap A^\circ(Q,m)).
\label{eq:Ushellclosed}
\end{equation}
Summing \eqref{eq:Ushellclosed} over $Q\in \mathfrak B_m$ and then over $m=1,\ldots,\ell$, we infer that
\begin{equation}
    \begin{split}
        &\mu\Big(U\cap \bigcup\{10AB_Q:Q\in \partial(\Omega_\ell,\mathfrak L_{2\ell+1})\}\Big)\leq
    3\sqrt{\delta}\sum_{m=1}^{\ell}(36\sqrt{\delta}\mathfrak d)^{\ell-m}\sum_{Q\in \mathfrak B_m}\mu(U\cap A^\circ(Q,m)).
    \end{split}
\end{equation}
Now, for every fixed $m$, the sets $A^\circ(Q,m)$, $Q\in \mathfrak B_m$, are pairwise disjoint and disjoint from $\Omega_{m-1}$, hence
$$
\sum_{Q\in \mathfrak B_m}\mu(U\cap A^\circ(Q,m))
\leq
\mu(U\setminus \Omega_{m-1}).
$$
Using the measure gain already proved,
$$\mu(U\setminus \Omega_{m-1})\leq\theta^{m-1}\mu(U)+2\mu(U\setminus\bigcup\{Q:Q\in\mathfrak L_{\mathfrak N}\}),$$
we obtain
\begin{equation}
    \begin{split}
        &\qquad\qquad\qquad\qquad\qquad\mu\Big(U\cap \bigcup\{10AB_Q:Q\in \partial(\Omega_\ell,\mathfrak L_{2\ell+1})\}\Big)\\
        &\leq3\sqrt{\delta}\sum_{m=1}^{\ell}(36\sqrt{\delta}\mathfrak d)^{\ell-m}\theta^{m-1}\mu(U)+6\sqrt{\delta}\sum_{m=1}^{\ell}(36\sqrt{\delta}\mathfrak d)^{\ell-m}
    \mu\Big(U\setminus\bigcup\{Q:Q\in\mathfrak L_{\mathfrak N}\}\Big)\\
&\qquad\qquad\leq
\frac{3\sqrt{\delta}}{\theta-36\sqrt{\delta}\mathfrak d}\,\theta^\ell\,\mu(U)
+
\frac{6\sqrt{\delta}}{1-36\sqrt{\delta}\mathfrak d}
\mu\Big(U\setminus\bigcup\{Q:Q\in\mathfrak L_{\mathfrak N}\}\Big)\\
&\qquad\qquad\qquad\leq
36\sqrt{\delta}D^22^k\theta^\ell\mu(U)
+
\mu\Big(U\setminus\bigcup\{Q:Q\in\mathfrak L_{\mathfrak N}\}\Big).
\nonumber
    \end{split}
\end{equation}

\textbf{Step III. Conclusion of the proof.}
By \eqref{exponentialdecaya_n} and by the choice of the parameters in \eqref{sceltaparametri}, we have
\begin{equation}
    \begin{split}
        \mu(U\setminus \Omega_N)\leq&36\sqrt{\delta} D^22^k\theta^N\mu(U)+2\mu(U\setminus\bigcup\{Q:Q\in\mathfrak L_{\mathfrak N}\})\\
        \leq&\xi\mu(U)+2\mu(U\setminus\bigcup\{Q:Q\in\mathfrak L_{\mathfrak N}\}).
    \end{split}
\end{equation}
Let us set $\Omega_e:=\Omega_N$. Then, item \hyperlink{propomegaei}{(i)} is satisfied by the above computation. Item \hyperlink{propomegaeii}{(ii)} follows thanks to item (6) of the induction hypothesis with $\ell=N$ and the choice of $N$, 
and items \hyperlink{propomegaeiii}{(iii)} and \hyperlink{propomegaeiv}{(iv)} follow from item (5) of the induction hypothesis. 
\end{proof}

We now iterate the previous thinning step. Starting from the original stack of layers, we apply \cref{stepbasethinning} inside $\Omega_{\mathfrak s}$ and obtain a smaller open set $\Omega_1$. We then retain only those cubes in a deeper layer whose enlarged balls are contained in $\Omega_1$, and restart the same construction on their descendants. Repeating this procedure produces nested open sets
$$
\Omega_{\mathfrak s}=\Omega_0\supseteq\Omega_1\supseteq\cdots\supseteq\Omega_K.
$$
At each restart, items \textnormal{(i)}--\textnormal{(vi)} of \cref{stepbasethinning} are preserved for the current family. The ancestor-containment condition in item \textnormal{(v)} is used precisely to propagate the terminal density condition \textnormal{(vi)} to the restarted family. At the end, we splice together the layers that survive at the different stages. The resulting families $\mathfrak G_h$ still satisfy the structural assumptions \textnormal{(i)}--\textnormal{(v)}, while the open sets $\Omega_h$ give an iterated length loss along all $C(e,\sigma)$-curves. The terminal layer is not completely preserved, but the discarded part is quantitatively controlled by the mass lost outside $\Omega_K$ and by the boundary errors produced at each stage.

\begin{proposizione}\label{productionopens}
Suppose $\mu$ is a $k$-AD-regular measure with regularity constant $D$. Fix parameters
$\sigma,\varepsilon,\delta,\zeta,\mathfrak m,\beth,\xi'\in (0,1)$ and $\mathfrak N,K\in \N\setminus\{0\}$ and $A,\qof\geq 1$.
Further, fix a finite subset $\mathscr E \subseteq \mathbb S^{n-1}$ and assume also that
$$
\zeta \leq \delta\leq 1/2^{24},\quad \qof\geq 4\lceil \log_2(80\delta^{-1}\zeta^{-1}A)\rceil, \quad \frac{4\zeta}{\varepsilon \sqrt{1-\sigma^2}}\leq 1-\xi',
$$
where
\begin{equation}
\begin{split}
 \mathfrak N\geq 8(K+1)(2N+1)\log_2(160\delta^{-1}\zeta^{-1}A)  \qquad\text{and}\qquad N:= \Big\lceil\frac{\log(\beth/(144\sqrt{\delta}D^22^k))}
{\log(1-\frac{\zeta^k}{40^{k+1}A^kD^2}(1-\sqrt{\delta}))}\Big\rceil.
\end{split}
    \label{sceltaparametriv2}
\end{equation}
Suppose moreover that $\mathfrak F\subseteq \Delta_\mu$ is a family of cubes contained in some open set $\Omega_\mathfrak s$ with $10AB_Q\subseteq \Omega_\mathfrak s$ for every $Q\in \mathfrak F$, and that
$$ \mathfrak F=\bigcup_{m=0}^{\mathfrak N} \mathfrak L_m. $$ 
Let $U$ be a measurable set such that
$$
\mu(U\setminus\bigcup\{Q:Q\in\mathfrak L_{\mathfrak N}\})\leq \tfrac{\beth}{4}\mu(U).
$$
and assume that the families $\mathfrak L_0,\ldots,\mathfrak L_{\mathfrak N}$ and the set $U$ satisfy the hypothesis (i) to (vi) of \cref{stepbasethinning}.

Then, there are families of dyadic cubes $\mathfrak G_0,\ldots, \mathfrak G_\mathfrak N$ such that $\mathfrak G_h\subseteq \mathfrak L_h$ for every $h=0,\ldots, \mathfrak N$ and nested open sets 
$$\Omega_\mathfrak s=:\Omega_0\supseteq \Omega_1\supseteq \ldots \supseteq\Omega_K,$$
such that 
\begin{enumerate}
    \item for every $h=1,\ldots, \mathfrak N$ every cube $Q\in \mathfrak G_h$ satisfies
$40AB_Q\subseteq \Omega_{\min \{\lfloor h/(2N+1)\rfloor,K\}}$ 
\item the layers $\mathfrak G_h$ still satisfy the hypothesis in items (i) to (v) of \cref{stepbasethinning};
\item for every $h=0,\ldots, K-1$, every $e\in \mathscr E$ and every $\gamma\in C(e,\sigma)$ we have
\begin{equation}
\Haus^1(\gamma\cap \Omega_{h+1})
\leq(1-\xi')\Haus^1(\gamma\cap \Omega_h),
\nonumber
\end{equation}
\item for every $h=0,\ldots, K-1$ there holds
$$\Omega_{h+1}\subseteq \bigcup_{Q\in\mathfrak G_{h(2N+1)+1}}\tfrac{3}{2}AB_Q$$
and for every $Q\in\mathfrak G_{h(2N+1)+1}$, for every $e\in \mathscr E$ and every $C(e,\sigma)$-curve $\gamma$ we have that 
$$\Haus^1(\gamma \cap \Omega_{h+1}\cap 2AB_Q)\leq
(1-\xi')\Haus^1(\gamma \cap \Omega_h \cap 2(1+2^{-\qof})AB_Q).$$
\item $\mu(U\setminus \Omega_K)\leq 2\mathrm{Card}(\mathscr{E})K\beth\mu(U)$.
\end{enumerate}
\end{proposizione}

\begin{proof}
For every $h=0,\ldots,K$ we shall construct an open set $\Omega_h$, a measurable set $U^{(h)}$ and families $ \mathfrak L_m^{(h)}\subseteq \mathfrak L_{h(2N+1)+m}$, 
where $m=0,\ldots,\mathfrak N-h(2N+1)$, such that if we define 
$$ \mathfrak F^{(h)}:=\bigcup_{m=0}^{\mathfrak N-h(2N+1)}\mathfrak L_m^{(h)}, $$
then $\mathfrak F^{(h)}$ and $U^{(h)}$ satisfy the hypotheses of \cref{stepbasethinning} in the open set $\Omega_h$.
We begin by setting 
$$ \Omega_0:=\Omega_{\mathfrak s}, \qquad U^{(0)}:=U, \qquad \mathfrak L_m^{(0)}:=\mathfrak L_m \quad\text{for every }m=0,\ldots,\mathfrak N. $$ 
By the assumptions of the present proposition, $\mathfrak F^{(0)}$ and $U^{(0)}$ satisfy the hypotheses of \cref{stepbasethinning}.

Let us move to the inductive step. Let us suppose $\mathfrak N-h(2N+1)\geq 2N+1$ and let us assume that we have constructed a family of nested sets 
$$\Omega_h\subseteq \ldots\subseteq \Omega_1\subseteq \Omega_0=\Omega_\mathfrak s,$$
and families of cubes 
$$ \mathfrak F^{(h)}=\bigcup_{m=0}^{\mathfrak N-h(2N+1)}\mathfrak L_m^{(h)},$$
such that $\mathfrak F^{(h)}$ and $U^{(h)}$ satisfy the hypothesis of \cref{stepbasethinning} for the open set $\Omega_h$.
We now construct $\Omega_{h+1}$, $U^{(h+1)}$, $\mathfrak F^{(h+1)}$ and $\mathfrak L_m^{(h+1)}$ with $m=0,\ldots, \mathfrak N-(h+1)(2N+1)$.

Let us apply \cref{stepbasethinning} to $\Omega_{h}$, to the family $\mathfrak{F}^{(h)}$ and to the set $U^{(h)}$ with parameters $\xi:=\beth/2$, $\xi'$ and let us note that we obtain an open set $\Omega_{h+1}\subseteq \Omega_h$ satisfying the following properties 
\begin{enumerate}
\item[a.]$\mu(U^{(h)}\setminus \Omega_{h+1})\leq\frac12\mathrm{Card}(\mathscr E)\beth\mu(U^{(h)})+2\mathrm{Card}(\mathscr E)\mu(U^{(h)}\setminus\bigcup\{Q:Q\in\mathfrak L_{\mathfrak N-h(2N+1)}^{(h)}\})$;
\item[b.] defining $\partial (\Omega_{h+1},\mathfrak L_{2N+1}^{(h)})$ as the family of cubes $Q$ in $\mathfrak L_{2N+1}^{(h)}$ such that $10 AB_Q\cap \Omega_{h+1}\neq \emptyset$ and $40 AB_Q\not\subseteq \Omega_{h+1}$, we have
\begin{equation}
    \begin{split}
        &\mu(U^{(h)}\cap\bigcup\{10AB_Q:Q\in \partial(\Omega_{h+1}, \mathfrak L_{2N+1}^{(h)})\})\\
        &\qquad\qquad\qquad\leq\frac12\mathrm{Card}(\mathscr E)\beth\mu(U^{(h)})+\mathrm{Card}(\mathscr E)\mu\Big(U^{(h)}\setminus\bigcup\{Q:Q\in\mathfrak L_{\mathfrak N-h(2N+1)}^{(h)}\}\Big).
        \nonumber
    \end{split}
\end{equation}
\item[c.] There holds
$\Omega_{h+1}\subseteq \bigcup_{Q\in \mathfrak L_1^{(h)}}\tfrac{3}{2}AB_Q$
and for every $e\in \mathscr E$ and every $Q\in \mathfrak L_1^{(h)}$ we have
$$\Haus^1(\gamma\cap \Omega_{h+1}\cap 2AB_{Q})
     \leq (1-\xi')\Haus^1(2(1+2^{-\qof})AB_Q\cap \gamma),$$
     for every $C(e,\sigma)$-curve; 
\item[d.] for every $e\in \mathscr E$ and every $C(e,\sigma)$-curve $\gamma$ we have
$$\Haus^1(\gamma\cap \Omega_{h+1})\leq (1-\xi')\Haus^1(\gamma\cap \Omega_{h}).$$
\end{enumerate}
In what is left we need to construct the families $\mathfrak{L}_m^{(h+1)}$ for $m=0,\ldots,\mathfrak N-(h+1)(2N+1)$ of cubes contained in $\Omega_{h+1}$ and prove they satisfy the hypothesis of \cref{stepbasethinning}.
We define
$$\mathfrak L_0^{(h+1)}:=\{Q\in \mathfrak L_{2N+1}^{(h)}:40AB_Q\subseteq \Omega_{h+1}\},$$
and for every $m=1,\ldots,\mathfrak N-(h+1)(2N+1)$ we let
$$\mathfrak L_m^{(h+1)}:=\Bigl\{Q\in \mathfrak L_{2N+1+m}^{(h)}:\text{the ancestor of }Q\text{ in }\mathfrak L_{2N+1}^{(h)}\text{ belongs to }\mathfrak L_0^{(h+1)}\Bigr\}.$$
We define 
$$U^{(h+1)}:=U^{(h)}\cap\bigcup\{Q:Q\in\mathfrak L_0^{(h+1)}\}.$$
We claim that the restarted family 
$$ \mathfrak F^{(h+1)} := \bigcup_{m=0}^{\mathfrak N-(h+1)(2N+1)}\mathfrak L_m^{(h+1)}, $$
together with $U^{(h+1)}$, still satisfies the hypotheses of the \cref{stepbasethinning} inside $\Omega_{h+1}$ with parameters $\xi:=\beth/2$ and $\xi'$.
The properties (i), (ii), (iii), (iv) and (v) are inherited immediately, since we are only passing to subfamilies.

It remains to check property (vi) for the restarted families $\mathfrak L_m^{(h+1)}$.
Fix $m\in\{1,\ldots,\mathfrak N-(h+1)(2N+1)\}$, $Q\in \mathfrak L_m^{(h+1)}$ and $e\in \mathscr E$, and denote by $Q'\in \mathfrak L_0^{(h+1)}$ the ancestor of $Q$ in $\mathfrak L_{2N+1}^{(h)}$. By property (v), iterated along the chain of ancestors from $Q$ to $Q'$, we have
$$
B(z(Q,e),\zeta \diam Q/2)\cap\supp\mu\subseteq Q'.
$$
Since $U^{(h+1)}=U^{(h)}\cap\bigcup\{S:S\in\mathfrak L_0^{(h+1)}\}$ and $Q'\in\mathfrak L_0^{(h+1)}$, it follows that
$$
B(z(Q,e),\zeta \diam Q/2)\cap U^{(h)}
\subseteq
B(z(Q,e),\zeta \diam Q/2)\cap U^{(h+1)}
$$
up to a $\mu$-null set. Therefore, by property (vi) for $\mathfrak F^{(h)}$ and $U^{(h)}$, we infer that
$$
\mu\Bigl(B(z(Q,e),\zeta \diam Q/2)\cap U^{(h+1)}\Bigr)
\geq
\frac12\,\mu(B(z(Q,e),\zeta \diam Q/2)).
$$
This shows that property (vi) is inherited. Moreover,
$$
U^{(h+1)}\setminus\bigcup\{Q:Q\in \mathfrak L_{\mathfrak N-(h+1)(2N+1)}^{(h+1)}\}
\subseteq
U^{(h)}\setminus\bigcup\{Q:Q\in \mathfrak L_{\mathfrak N-h(2N+1)}^{(h)}\}.
$$
Thus the terminal-defect term does not increase in the restart, and the induction on $h$ is complete.

\medskip
For every $h=0,\ldots,\mathfrak N$, set 
$$
\mathfrak G_h:=
\begin{cases}
\mathfrak L_m^{(r)}, 
&\text{if } h=r(2N+1)+m,\text{ where } 0\leq r\leq K-1,\text{ and }0\leq m\leq 2N,\\
\mathfrak L_m^{(K)}, 
&\text{if } h=K(2N+1)+m,\text{ where }\ 0\leq m\leq \mathfrak N-K(2N+1).
\end{cases}
$$
We now check that $\mathfrak G_h$ satisfies the properties 1 to 5. Let $Q\in\mathfrak G_h$ with $h\geq1$. If $h/(2N+1)\geq1$, there is nothing to prove thanks to the construction of the $\mathfrak G_h$. If instead $1\leq h<2N+1$, let $Q_0\in\mathfrak L_0$ be the ancestor of $Q$. By the scale separation along the chain from $Q$ to $Q_0$, we have $40AB_Q\subseteq10AB_{Q_0}$. Since $10AB_{Q_0}\subseteq\Omega_{\mathfrak s}=\Omega_0$, we get
$$
40AB_Q\subseteq\Omega_0=\Omega_{\min\{\lfloor h/(2N+1)\rfloor,K\}}.
$$
This proves item 1.

The layers $\mathfrak G_h$ satisfy the assumption (i) to (v) of the hypothesis of \cref{stepbasethinning} thanks to the above induction argument. Items 3 and 4 are a direct consequence of item c and d above. 
Indeed, for every $h=0,\ldots,K-1$ we have
$$
\mathfrak G_{h(2N+1)+1}=\mathfrak L_1^{(h)}.
$$
Hence item c gives
$$
\Omega_{h+1}\subseteq\bigcup_{Q\in\mathfrak G_{h(2N+1)+1}}\frac32AB_Q.
$$
Moreover, if $Q\in\mathfrak G_{h(2N+1)+1}$, then $Q\in\mathfrak L_1^{(h)}$, and therefore $10AB_Q\subseteq\Omega_h$. Since $\qof\geq1$,
$$
2(1+2^{-\qof})AB_Q\subseteq3AB_Q\subseteq10AB_Q\subseteq\Omega_h.
$$
Thus
$$
2(1+2^{-\qof})AB_Q\cap\gamma
=
\Omega_h\cap2(1+2^{-\qof})AB_Q\cap\gamma.
$$
Combining this with item c gives item 4.

We now move to the proof of item 5, i.e. we now show that $\Omega_K$ still covers most of $U$.
For every $h=0,\ldots,K-1$, let us observe that 
\begin{equation} 
\begin{split} 
U^{(h)}\setminus U^{(h+1)}
\subseteq&
\Bigl(U^{(h)}\setminus\bigcup\{Q:Q\in \mathfrak L_{\mathfrak N-h(2N+1)}^{(h)}\}\Bigr)\cup
\bigl(U^{(h)}\setminus \Omega_{h+1}\bigr)\\
&\cup
\Big(U^{(h)}\cap\bigcup\{10AB_Q:Q\in\partial(\Omega_{h+1},\mathfrak L_{2N+1}^{(h)})\}\Big). 
\end{split} 
\label{eq:discarded-Uh} 
\end{equation} 
Indeed, let $x\in U^{(h)}\setminus U^{(h+1)}$. If
$$
x\notin\bigcup\{Q:Q\in \mathfrak L_{\mathfrak N-h(2N+1)}^{(h)}\},
$$
then there is nothing to prove. Otherwise, we may choose $R\in\mathfrak L_{\mathfrak N-h(2N+1)}^{(h)}$ such that $x\in R$. Since $x\notin U^{(h+1)}$, the ancestor $Q'\in\mathfrak L_{2N+1}^{(h)}$ of $R$ does not belong to $\mathfrak L_0^{(h+1)}$, and hence $40AB_{Q'}\not\subseteq \Omega_{h+1}$. 
If $10AB_{Q'}\cap\Omega_{h+1}=\emptyset$, then $x\in U^{(h)}\setminus\Omega_{h+1}$. Otherwise $Q'\in\partial(\Omega_{h+1},\mathfrak L_{2N+1}^{(h)})$, and this proves \eqref{eq:discarded-Uh}. 

Since the sets $U^{(h)}$ are decreasing and $U^{(K)}\subseteq\Omega_K$, we have
\begin{equation} 
\begin{split} 
U\setminus\Omega_K \subseteq \bigcup_{h=0}^{K-1}\bigl(U^{(h)}\setminus U^{(h+1)}\bigr). 
\nonumber
\end{split}
\end{equation} 
Using \eqref{eq:discarded-Uh}, item a. and item b., we obtain
\begin{equation} 
\begin{split}
&\mu(U\setminus\Omega_K)\leq\sum_{h=0}^{K-1}
\mu\Bigl(U^{(h)}\setminus\bigcup\{Q:Q\in \mathfrak L_{\mathfrak N-h(2N+1)}^{(h)}\}\Bigr)+
\sum_{h=0}^{K-1}\mu(U^{(h)}\setminus\Omega_{h+1})\\
&\qquad\qquad\qquad\qquad\qquad+
\sum_{h=0}^{K-1}
\mu\Bigl(U^{(h)}\cap\bigcup\{10AB_Q:Q\in\partial(\Omega_{h+1},\mathfrak L_{2N+1}^{(h)})\}\Bigr)\\
&\leq
K\mu\Big(U\setminus\bigcup\{Q:Q\in\mathfrak L_{\mathfrak N}\}\Big)
+
\mathrm{Card}(\mathscr E)K\beth\,\mu(U)+3\mathrm{Card}(\mathscr E)K
\mu\Big(U\setminus\bigcup\{Q:Q\in\mathfrak L_{\mathfrak N}\}\Big)\\
&\qquad\qquad\qquad\qquad\qquad\leq
2\mathrm{Card}(\mathscr E)K\beth\,\mu(U),
\nonumber
\end{split}
\end{equation}
because $U^{(h)}\subseteq U$ for every $h=0,\ldots,K$. This proves item 5 and concludes the proof.
\end{proof}

\subsection{Joining cones}

In this subsection we pass from the estimates produced by the thinning procedure to the estimates that are needed in the perturbation argument. The construction of the previous subsection gives open sets which have small intersection with curves going in a direction of a lot of very thin cones. This is not sufficient for the next step. Indeed, the perturbations are built from width functions, and the Lipschitz constants of these functions depend on the aperture of the cone. If the cone is very narrow, the Lipschitz constant along the invariant directions is too large. Thus, we need smallness not only for a finite family of narrow cones, but for the whole transverse cone covered by such small cones.

The difficulty is that this implication is not automatic. Indeed, a curve going in the large cone might have derivative oscillating in a very fast fashion between one thin cone and another. If the curve zig-zags too rapidly, it is not clear if we can get a big portion of this curve going in the direction of one of the two thin cones. If this was the case we would know that the curve we are studying looses a little bit of length when intersected with the open sets.

This is where the joining cones lemma by Alberti--Csörnyei--Preiss lemma enters. Suppose that a convex cone $C$ is the union of two overlapping convex subcones $C^+$ and $C^-$. The lemma allows one to find, inside every $C$-curve, a compact set with uniformly controlled length on which the curve coincides with a curve whose increments lie either in $C^+$ or in $C^-$. Hence estimates known for the two smaller cones can be transferred to the larger cone, with a quantitative loss. We organize this transfer through a finite binary tree of convex cones. The root is the large cone for which we want the width estimate, while the leaves are contained in the narrow cones where the thinning estimates are already available. Iterating the two-cone argument and going back up from the leaves to the root, gives the desired estimate for the large cone.

We then make this applicable to the families of cubes. First, we decompose $\Delta_\mu$ according to the maximal number of independent invariant directions. If the family $\mathscr I_j$ is not Carleson, then, after pigeonholing the possible invariant directions, we find arbitrarily long layers of cubes with one fixed invariant $j$-tuple. Finally, we show that the existence of a huge amount of cubes non-invariant along a fixed amount of directions will imply that we can generate many open sets that have small width with respect to very wide cones.

\begin{lemma}\label{lem:pigeonhole-direction-small-cone}
Let $0<\sigma'<\sigma$ and $e\in\mathbb S^{n-1}$. There exists a constant $C_n$ depending only on $n$, and a finite set $\mathcal E_{\sigma'}(e,\sigma)\subset C(e,\sigma)\cap\mathbb S^{n-1}$ such that
$$
\mathrm{Card}\bigl(\mathcal E_{\sigma'}(e,\sigma)\bigr)\leq \frac{C_n}{(\sigma')^{n-1}},
$$
and for every $v\in C(e,\sigma)\cap\mathbb S^{n-1}$ there exists $\bar e\in\mathcal E_{\sigma'}(e,\sigma)$ such that $v\in C(\bar e,\sigma')$.
\end{lemma}

\begin{proof}This is an immediate consequence of the fact that $\Haus^{n-1}\trace \mathbb{S}^{n-1}$ is $(n-1)$-AD-regular.
\end{proof}

The following lemma is due to G. Alberti, M. Csörnyei and D. Preiss. It is the curve-selection lemma used in their joining-cones argument; see \cite[Lemma 2.35]{acpnull}. More precisely, their Section 2.7 shows how, given two one-sided subcones of a convex cone separated by a strip of width $\alpha$, every curve directed in the original cone can be replaced, on a subset of uniformly positive one-dimensional measure, by a curve directed by one of the two subcones. This is then used in \cite[Lemma 2.36]{acpnull} to join the corresponding zero-width conditions. The numbering refers to the version of the manuscript in the author's possession, dating back to 2015, as \cite{acpnull} at the time of writing is still unpublished. We will use this proposition to do the exact same thing. The strength of their argument is that they get quantification of the portion of the measure covered by the curves going in the smaller cones. This is of course of capital importance in a quantified argument.

\begin{lemma}\label{lem:acp-two-convex-slope-cones}
Let $C$ be a closed convex cone and let $e\in\mathbb S^{n-1}$ be such that
$\langle x,e\rangle>0$ for every 
$x\in C\setminus\{0\}$.
Let us choose $u\in e^\perp\cap\mathbb S^{n-1}$, $c\in\mathbb R$ and $\alpha>0$ in such a way that $c-\alpha\leq0\leq c+\alpha$ and define
$$
C^+:=\{x\in C:\langle x,u\rangle\geq (c-\alpha)\langle x,e\rangle\},
\qquad
C^-:=\{x\in C:\langle x,u\rangle\leq (c+\alpha)\langle x,e\rangle\}.
$$
Let us define
$$
K_0:=\max\Big\{1,\sup_{x\in C\setminus \{0\}}\frac{\lvert \langle u,x\rangle\rvert}{\langle e,x\rangle}\Big\}.
$$
Let $I$ be a compact interval and let $\gamma:I\to\mathbb R^n$ be a $C$-curve so that $\gamma(t)=te+\eta(t)$, with $\eta(s)-\eta(t)\in e^\perp$ for every $s,t\in I$.
Then, there exists a compact set $S\subseteq I$ and a Lipschitz curve $\tilde \gamma:I\to\mathbb R^n$
such that
$$
\tilde \gamma(t)=\gamma(t)
\qquad\text{for every }t\in S,
$$
and $\tilde \gamma$ is either a $C^+$-curve or a $C^-$-curve, and
$$
\Leb^1(S)\geq \Leb^1(I)\min\left\{1,\frac{\alpha}{2K_0}\right\}.
$$
\end{lemma}

\begin{proof}
Let us assume $I=[a,b]$.
We prove the $C^-$ alternative under the assumption that
$\langle \gamma(b)-\gamma(a),u\rangle\leq c(b-a)$.
The $C^+$ alternative is obtained by reversing the inequalities. Let us first notice that if
$$
\langle \gamma(s)-\gamma(t),u\rangle\leq (c+\alpha)(s-t)
\qquad\text{for every }a\leq t\leq s\leq b,
$$
then there is nothing to prove. We may therefore suppose that there exist $a\leq t<s\leq b$ such that
$$
\langle \gamma(s)-\gamma(t),u\rangle>(c+\alpha)(s-t).
$$
By continuity of the function $(s,t)\mapsto \langle u,\gamma(s)-\gamma(t)\rangle$, there exist $a\leq t<s\leq b$ such that
$$
\langle \gamma(s)-\gamma(t),u\rangle=(c+\alpha)(s-t).
$$
Since $\gamma(s)-\gamma(t)\in C$ and $\langle \gamma(s)-\gamma(t),e\rangle=s-t$, the definition of $K_0$ gives
$c+\alpha\leq K_0$.
Define
$$
S:=
\left\{
t\in[a,b]:
\langle \gamma(s)-\gamma(t),u\rangle\leq (c+\alpha)(s-t)
\text{ for every }s\in[t,b]
\right\}.
$$
Since $\gamma$ is continuous, then $S$ is compact and thanks to the above discussion, we have $S\neq [a,b]$ and $S\neq \emptyset$.
We set
$$
\tilde \gamma(t):=\gamma(t)
\qquad\text{for every }t\in S.
$$
On every connected component $(a_*,b_*)$ of $(a,b)\setminus S$ with $a_*,b_*\in S$, we define
$\tilde{\gamma}$ affinely by
$$
\tilde{\gamma}(t):=
\frac{b_*-t}{b_*-a_*}\gamma(a_*)+
\frac{t-a_*}{b_*-a_*}\gamma(b_*).
$$
If $a\notin S$, let $t_0:=\min S$ and set
$$
v:=\frac{\gamma(b)-\gamma(a)}{b-a}.
$$
Then $v\in C$, $\langle v,e\rangle=1$, and by the assumption
$\langle v,u\rangle\leq c\leq c+\alpha$. Thus $v\in C^-$. We define
$$
\tilde{\gamma}(t):=\gamma(t_0)+(t-t_0)v
\qquad\text{for }a\leq t\leq t_0.
$$
This makes $\tilde{\gamma}$ continuous and Lipschitz on $[a,b]$.
We now prove that $\tilde{\gamma}$ is a $C^-$-curve. Let
$$
f(t):=\langle \gamma(t),u\rangle-(c+\alpha)t.
$$
By the definition of $S$, if we set
$$
F(t):=\max_{s\in[t,b]} f(s),
$$
then $S=\{t\in[a,b]: f(t)=F(t)\}$.
Let $(a_*,b_*)$ be a connected component of $[a,b]\setminus S$ with endpoints in $S$.
We prove that
$$
f(a_*)=f(b_*).
$$

First observe that $F$ is non-increasing. We claim that $F$ is constant on $(a_*,b_*)$.
Indeed, let $a_*<x<y<b_*$. Suppose by contradiction that
$F(x)>F(y)$.
Since
$F(y)=\max_{s\in[y,b]} f(s)$,
the maximum $F(x)$ cannot be attained in $[y,b]$. Hence it is attained at some point
$z\in[x,y)$. Therefore
$f(z)=F(x)$
.
Since $z\in[x,y)$, we also have
$F(z)\leq F(x)$
,
because $F$ is non-increasing, and
$$
F(z)\geq f(z)=F(x).
$$
Thus
$F(z)=f(z)$
,
so $z\in S$ and this is a contradiction.
Thus $F$ is constant on $(a_*,b_*)$.
Since $f$ is continuous, the function $F$ is continuous and hence
$$
F(a_*)=F(b_*).
$$
Finally, since $a_*,b_*\in S$, we have obtained the sought identity
$f(a_*)=f(b_*)$. Equivalently,
$$
\langle \gamma(b_*)-\gamma(a_*),u\rangle
=(c+\alpha)\langle \gamma(b_*)-\gamma(a_*),e\rangle.
$$
Since $\gamma(b_*)-\gamma(a_*)\in C$, this gives
$\gamma(b_*)-\gamma(a_*)\in C^-$.
Therefore every affine piece of $\tilde{\gamma}$ inside $[a,b]\setminus S$ has derivative in $C^-$.
If $s<t$ are two points of $S$, then by the defining property of $S$,
$$
\langle \gamma(t)-\gamma(s),u\rangle\leq (c+\alpha)(t-s).
$$
Since $\gamma(t)-\gamma(s)\in C$ and $\langle \gamma(t)-\gamma(s),e\rangle=t-s$, the same argument used above yields
$\gamma(t)-\gamma(s)\in C^-$.
Finally, take arbitrary $a\leq s<t\leq b$. The interval $[s,t]$ is decomposed into pieces on which
the corresponding increments of $\tilde\gamma$ belong to $C^-$. Therefore
$$
\tilde \gamma(t)-\tilde \gamma(s)\in C^-.
$$
Thus $\tilde \gamma$ is a $C^-$-curve. It remains to estimate $\Leb^1(S)$. Since $f(a_*)=f(b_*)$ on every connected component $(a_*,b_*)$ of $[a,b]\setminus S$ with endpoints in $S$, and since the possible initial component satisfies the corresponding inequality, we have
$$
\int_{[a,b]\setminus S}\langle \gamma'(t),u\rangle\,dt
\geq  (c+\alpha)\Leb^1([a,b]\setminus S).
$$
By assumption we have that 
$\langle \gamma(b)-\gamma(a),u\rangle\leq c(b-a)$, and by construction of $\tilde \gamma$ we obtain
\begin{equation}
\begin{split}
c(b-a)
\geq&
\langle \gamma(b)-\gamma(a),u\rangle=\int_{[a,b]} \langle \gamma'(t),u\rangle dt=\int_S\langle \gamma'(t),u\rangle dt+\int_{[a,b]\setminus S}\langle \gamma'(t),u\rangle dt\\
\geq&
(c+\alpha)\Leb^1([a,b]\setminus S)
+
\int_S \langle \gamma'(t),u\rangle\,dt.
\end{split}
\end{equation}
Equivalently, we can rewrite the above inequality as
$$
\alpha(b-a)
\leq
\int_S \bigl((c+\alpha)-\langle \gamma'(t),u\rangle\bigr)\,dt.
$$
For almost every $t\in S$ we have $\gamma'(t)\in C$ and
$\langle \gamma'(t),e\rangle=1$.
Hence
$|\langle \gamma'(t),u\rangle|\leq K_0$.
Since also $c+\alpha\leq K_0$, we get
$$
(c+\alpha)-\langle \gamma'(t),u\rangle\leq 2K_0
\qquad\text{for almost every }t\in S.
$$
Therefore
$\alpha(b-a)\leq 2K_0\Leb^1(S)$,
and hence
$$
\Leb^1(S)\geq \frac{\alpha}{2K_0}(b-a).
$$
This concludes the proof.
\end{proof}

Here we split the family $\Delta_\mu$ according to the number of quantitatively independent invariant directions carried by each cube. The last condition says that this number is maximal in a quantitative sense, namely no further invariant direction can be added outside a fixed cone around the span of the existing ones.

\begin{proposizione}\label{rioridinocubi}
Suppose $\mu$ is a $k$-AD-regular measure with regularity constant $D$ and choose parameters $\sigma,\varepsilon,\vartheta,\delta,\zeta,\mathfrak m\in (0,1)$ and $A\geq 1$ satisfying the following properties
$$A\vartheta^n\geq 16n,\qquad \text{and} \qquad 16n\Big(\frac{\sigma}{\vartheta^n}
+n\frac{\varepsilon+\delta+\zeta+\mathfrak m}{A\vartheta^n}\Big)\leq \frac{1}{2\cdot 3^k 16^n D^2}.$$
Then, we can write $\Delta_\mu$ as a disjoint union of families of cubes $\mathscr I_j:=\mathscr I(j,\vartheta,\sigma,\varepsilon,\delta,\zeta,\mathfrak m,A)$, for $j=0,\ldots, k$, we have
\begin{enumerate}
\item for every $Q\in \mathscr I_j$ there are
$\mathfrak{v}=(v_1,\ldots,v_j)\in (\mathbb S^{n-1})^j$, such that $Q$ is $(j,\mathfrak{v},\vartheta^j,\sigma,\varepsilon,\delta,\zeta,\mathfrak m,A)$-invariant;
\item for every
$e\in X_c(\mathrm{span}(\mathfrak v),\vartheta) \cap \mathbb{S}^{n-1}$, the cube $Q$ is not
$(e,\sigma,\varepsilon,\delta,\zeta,\mathfrak m,A)$-invariant.
\end{enumerate}
\end{proposizione}

\begin{proof}
Let us define $\widetilde{\mathscr I}_0:=\Delta_\mu$. For every $j=1,\ldots,k+1$ let us define
$$
\widetilde{\mathscr I}_j:=\{Q\in \Delta_\mu:Q\text{ is }(j,\mathfrak e,\vartheta^j,\sigma,\varepsilon,\delta,\zeta,\mathfrak m,A)\text{-invariant for some }\mathfrak e\in (\mathbb S^{n-1})^j\}.
$$
Notice that thanks to \cref{noinvariantmorethank}, and to the choice of the parameters with $\vartheta^n$, we have that $\widetilde{\mathscr I}_{k+1}=\emptyset$. We define
$$
\mathscr I_j:=
\widetilde{\mathscr I}_j\setminus
\bigcup_{\ell=j+1}^{k+1}\widetilde{\mathscr I}_\ell,
\qquad j=0,\ldots,k.
$$
This gives a disjoint decomposition of $\Delta_\mu$.
Let us verify item 1. and item 2. Item 1. is satisfied by construction. Let $Q\in\mathscr I_j$, and let $\mathfrak v=(v_1,\ldots,v_j)$ be such that $Q$ is
$(j,\mathfrak v,\vartheta^j,\sigma,\varepsilon,\delta,\zeta,\mathfrak m,A)$-invariant.

We claim that if $e\in X_c(\mathrm{span}(\mathfrak v),\vartheta)$, then the tuple
$(v_1,\ldots,v_j,e)$ is $\vartheta^{j+1}$-separated. First, by assumption we have $\dist(e,\mathrm{span}(\mathfrak v))\geq \vartheta$. It remains to check the old directions. Fix $i=1,\ldots,j$ and set $V_i:=\mathrm{span}\{v_\ell:\ell\neq i\}$. Let $\pi_i$ be the orthogonal projection onto $V_i^\perp$, and define
$p:=\pi_i v_i$ and $q:=\pi_i e$. Then $|p|=\dist(v_i,V_i)\geq \vartheta^j$.
Moreover, since $\mathrm{span}(V_i,v_i)=V_i+\mathrm{span}(p)$, we have
$$
\vartheta \leq \dist(e,\mathrm{span}(\mathfrak v))=\dist(e,\mathrm{span}(V_i,v_i))
=
\dist(q,\mathrm{span}(p)).
$$
In particular $q\neq 0$. Similarly, $\mathrm{span}(V_i,e)=V_i+\mathrm{span}(q)$, and therefore
$$
\dist(v_i,\mathrm{span}(V_i,e))
=
\dist(p,\mathrm{span}(q)).
$$
We now compare the two distances. Since $p,q\neq 0$, the orthogonal projection formula gives
$$
\dist(p,\mathrm{span}(q))^2
=
|p|^2-\frac{\langle p,q\rangle^2}{|q|^2}, \qquad\text{and}\qquad
\dist(q,\mathrm{span}(p))^2
=
|q|^2-\frac{\langle p,q\rangle^2}{|p|^2}.
$$
Multiplying the second identity by $|p|^2/|q|^2$, we obtain
$$
\frac{|p|^2}{|q|^2}\dist(q,\mathrm{span}(p))^2
=
|p|^2-\frac{\langle p,q\rangle^2}{|q|^2}
=
\dist(p,\mathrm{span}(q))^2.
$$
Since $|q|\leq |e|=1$, it follows that
$\dist(p,\mathrm{span}(q))
\geq
|p|\,\dist(q,\mathrm{span}(p))$. Consequently,
$$
\dist(v_i,\mathrm{span}(V_i,e))
=
\dist(p,\mathrm{span}(q))
\geq
\vartheta^j\vartheta
=
\vartheta^{j+1}.
$$
Since this holds for every $i=1,\ldots,j$, the tuple $(v_1,\ldots,v_j,e)$ is $\vartheta^{j+1}$-separated.

Now suppose by contradiction that $Q$ is $(e,\sigma,\varepsilon,\delta,\zeta,\mathfrak m,A)$-invariant for some
$e\in X_c(\mathrm{span}(\mathfrak v),\vartheta)$. By the claim, the tuple $(\mathfrak v,e)$ is $\vartheta^{j+1}$-separated. Hence $Q$ is
$(j+1,(\mathfrak v,e),\vartheta^{j+1},\sigma,\varepsilon,\delta,\zeta,\mathfrak m,A)$-invariant. Therefore $Q\in\widetilde{\mathscr I}_{j+1}$, contradicting
$Q\in\mathscr I_j\subseteq \widetilde{\mathscr I}_j\setminus \widetilde{\mathscr I}_{j+1}$. This proves item 2. and concludes the proof.
\end{proof}

The next proposition constructs a finite binary tree of convex cones which decomposes the large cone $C(e,\theta)$ into smaller cones contained in a fixed $\sigma$-net of directions. Each splitting has the precise form needed to apply the joining-cones lemma: the two children are obtained by cutting the parent cone with two overlapping half-spaces, and the ratio $\alpha_Q/K_Q$ is bounded from below uniformly along the whole tree. Thus estimates available on the terminal cones $C(e_i,\sigma)$ can be propagated back through the tree, with a controlled quantitative loss, to obtain the corresponding estimate on the original cone $C(e,\theta)$.

\begin{proposizione}\label{alberodiconi}
Let us fix $e\in S^{n-1}$, two parameters $0<\sigma<\theta<1$, with $0<\sigma<1/10$, and set
$V:=e^\perp$. Fix an orthonormal basis $v_1,\dots,v_{n-1}$ of $V$ and let
$\{e_i\}_{i=1}^N\subset S^{n-1}\cap C(e,\theta)$ be a maximal
$\sigma/4$-separated set of $S^{n-1}\cap C(e,\theta)$. 
Assume that
$$
\frac{\sigma\theta\sqrt{n-1}}{1600\sqrt{1-\theta^2}}\leq \frac12.
$$
Then, there is a finite rooted binary tree $\mathcal T$ of depth smaller than
$n\log_2(\frac{6400\theta\sqrt{n-1}}{\sigma\sqrt{1-\theta^2}})$
with the following properties.
Each node $\omega\in\mathcal T$ is associated to a convex set
$\Omega_\omega\subset V\cap B(0,\theta/\sqrt{1-\theta^2})$ and to the convex cone
$$
\Gamma_\omega:=\{\lambda(e+y):\lambda\geq0,\ y\in\Omega_\omega\}.
$$
The root $\emptyset$ satisfies $\Gamma_{\emptyset}=C(e,\theta)$ and every
non-terminal node $\omega$ has exactly two children $\omega+,\omega-$, with
$$
\Omega_\omega=\Omega_{\omega+}\cup\Omega_{\omega-},
\qquad
\Gamma_\omega=\Gamma_{\omega+}\cup\Gamma_{\omega-}.
$$
In particular $\Gamma_\omega$ is convex for every $\omega\in\mathcal T$.
For every non-terminal node $\omega$ there are
$\mathfrak e_\omega\in \Gamma_\omega\cap S^{n-1}$, $u_\omega\in \mathfrak e_\omega^\perp$ and
$\alpha_\omega>0$ such that
$$
\{x\in\Gamma_\omega:\langle x,u_\omega\rangle\geq -\alpha_\omega\langle x,\mathfrak e_\omega\rangle\}\subset\Gamma_{\omega+},
\qquad
\{x\in\Gamma_\omega:\langle x,u_\omega\rangle\leq \alpha_\omega\langle x,\mathfrak e_\omega\rangle\}\subset\Gamma_{\omega-},
$$
and $\Gamma_\omega$ is the union of these two sets. Moreover, if we let
$$
K_\omega:=\max\Big\{1,\sup_{x\in\Gamma_\omega\setminus\{0\}}
\frac{|\langle x,u_\omega\rangle|}{\langle \mathfrak e_\omega,x\rangle}\Big\},
\qquad\text{then}\qquad
\frac{\alpha_\omega}{2K_\omega}\geq
\frac{\sigma(1-\theta^2)^2}{6400\sqrt{n-1}} .
$$
Finally, for every terminal node $\omega$ there exists $i$ such that
$\Gamma_\omega\subset C(e_i,\sigma)$.
\end{proposizione}

\begin{proof}
Let us observe that every non-zero $x\in C(e,\theta)$ is uniquely written as
$x=\lambda(e+y)$, where $\lambda>0$, $y\in V$ and
$|y|\leq\theta/\sqrt{1-\theta^2}$. If $\Omega\subset V$ is convex, then
$\Gamma(\Omega):=\{\lambda(e+y):\lambda\geq0,\ y\in\Omega\}$ is a convex cone,
and notice that
$$
\Gamma(B(0,\theta/\sqrt{1-\theta^2}))=C(e,\theta).
$$
We first define the sets associated to the first $n-1$ non-root generations of
the tree as follows. Set
$$
\eta_0:=\frac{\sigma}{1600\sqrt{n-1}}
\qquad\text{and}\qquad
\Omega_\varnothing:=V\cap B(0,\theta/\sqrt{1-\theta^2}).
$$
Suppose that $0\leq k\leq n-2$ and that $\Omega_\omega$ has been defined for
every word $\omega\in\{+,-\}^k$. Define the two children of $\omega$ by
$$
\Omega_{\omega+}:=\Omega_\omega\cap\{\langle y,v_{k+1}\rangle\geq-\eta_0\},
\qquad
\Omega_{\omega-}:=\Omega_\omega\cap\{\langle y,v_{k+1}\rangle\leq\eta_0\}.
$$
Thus every node has exactly two children, the children are convex, and their
union is the parent. At the end of this first induction with $0\leq k\leq n-1$, the cells are indexed by
words
$\omega=(\omega_1,\dots,\omega_{n-1})\in\{+,-\}^{n-1}$ and each cell is a
fattened orthant inside $V\cap B(0,\theta/\sqrt{1-\theta^2})$.

For each non-terminal node constructed above, we construct
$\mathfrak e_\omega,u_\omega,\alpha_\omega$. Suppose first that $\omega$ is a
node of length $k$, with $0\leq k\leq n-2$. Thus the cone relative to $\omega$
is split in the coordinate $v_{k+1}$. Set
$$
\mathfrak e_\omega:=e,
\qquad
u_\omega:=v_{k+1},
\qquad
\alpha_\omega:=\eta_0(1-\theta^2).
$$
Then $\mathfrak e_\omega\in\Gamma_\omega\cap S^{n-1}$ and
$u_\omega\in\mathfrak e_\omega^\perp$. If
$x=\lambda(e+y)\in\Gamma_\omega$ and
$\langle x,u_\omega\rangle\geq-\alpha_\omega\langle x,\mathfrak e_\omega\rangle$, then
$\langle x,u_\omega\rangle\geq-\alpha_\omega|x|$, and therefore
$$
\langle y,v_{k+1}\rangle
\geq
-\alpha_\omega\textstyle{\sqrt{1+|y|^2}}
\geq
-\eta_0.
$$
Thus $x\in\Gamma_{\omega+}$. The inequality
$\langle x,u_\omega\rangle\leq\alpha_\omega\langle x,\mathfrak e_\omega\rangle$ gives similarly
$x\in\Gamma_{\omega-}$. This concludes the proof of the base step of the
construction.

We now continue the construction by blocks of $n-1$ generations. Suppose that
$h\geq1$, that $\omega$ is a word of length $h(n-1)$, and that
$\Omega_\omega$ has already been defined. If
\begin{equation}
\sup_{y,y'\in\Omega_\omega}|\langle y-y',v_j\rangle|\leq\frac{\sigma}{16\sqrt{n-1}}\qquad\text{for every }j=1,\dots,n-1,
    \label{stopaconsg}
\end{equation}
we declare $\omega$ terminal and do not split it further. Otherwise we
construct one full block below $\omega$. Suppose that $0\leq k\leq n-2$ and
that $\Omega_{\omega\omega'}$ has been defined for every word
$\omega'\in\{+,-\}^k$. For every such $\omega'$, set
$$
t_{\omega\omega'}:=
\frac12\bigg(
\inf_{y\in\Omega_{\omega\omega'}}\langle y,v_{k+1}\rangle+\sup_{y\in\Omega_{\omega\omega'}}\langle y,v_{k+1}\rangle
\bigg).
$$
Notice that $\lvert t_{\omega\omega'}\rvert\leq \frac{\theta}{\sqrt{1-\theta^2}}$ and define
$$
\Omega_{\omega\omega'+}:=\Omega_{\omega\omega'}\cap\{\langle y,v_{k+1}\rangle\geq t_{\omega\omega'}-\eta_0\},
\qquad
\Omega_{\omega\omega'-}:=
\Omega_{\omega\omega'}\cap
\{\langle y,v_{k+1}\rangle\leq t_{\omega\omega'}+\eta_0\}.
$$
After these $n-1$ successive binary splittings, the next block has been
constructed. Every non-terminal node has
exactly two children, the children are convex, and their union is the parent.

Since $\Omega_{\omega\omega'}$ is compact and convex, there exists $z_{\omega\omega'}\in\Omega_{\omega\omega'}$ such that
$\langle z_{\omega\omega'},v_{k+1}\rangle=t_{\omega\omega'}$. Define
$$
\mathfrak e_{\omega\omega'}:=\frac{e+z_{\omega\omega'}}{|e+z_{\omega\omega'}|},
\qquad
u_{\omega\omega'}:=\frac{v_{k+1}-t_{\omega\omega'}e}{\textstyle{\sqrt{1+t_{\omega\omega'}^2}}}\qquad
\alpha_{\omega\omega'}:=\eta_0(1-\theta^2).
$$
Then $\mathfrak e_{\omega\omega'}\in\Gamma_{\omega\omega'}\cap S^{n-1}$ and $u_{\omega\omega'}\in\mathfrak e_{\omega\omega'}^\perp$ since $\langle z_{\omega\omega'},e\rangle=0$.
Writing $x=\lambda(e+y)\in\Gamma_{\omega\omega'}$, we have
$$
\langle x,v_{k+1}-t_{\omega\omega'}e\rangle
=
\lambda(\langle y,v_{k+1}\rangle-t_{\omega\omega'}).
$$
Since
$|y|\leq\theta/\sqrt{1-\theta^2}$, the inequality
$\langle x,u_{\omega\omega'}\rangle\geq-\alpha_{\omega\omega'}\langle x,\mathfrak e_{\omega\omega'}\rangle$ implies, a fortiori,
$\langle x,u_{\omega\omega'}\rangle\geq-\alpha_{\omega\omega'}|x|$, and hence
$$
\langle y,v_{k+1}\rangle
\geq
t_{\omega\omega'}-\eta_0.
$$
Thus $x\in\Gamma_{{\omega\omega'+}}$.

The inequality
$\langle x,u_{\omega\omega'}\rangle\leq\alpha_{\omega\omega'}\langle x,\mathfrak e_{\omega\omega'}\rangle$ gives, 
$\langle x,u_{\omega\omega'}\rangle\leq\alpha_{\omega\omega'}|x|$, and hence $x\in\Gamma_{\omega\omega'-}$.
Therefore, for every non-terminal node ${\omega\omega'}$, we have
$$
\{x\in\Gamma_{\omega\omega'}:\langle x,u_{\omega\omega'}\rangle\geq-\alpha_{\omega\omega'}\langle x,\mathfrak e_{\omega\omega'}\rangle\}\subset\Gamma_{{\omega\omega'+}},
$$
$$
\{x\in\Gamma_{\omega\omega'}:\langle x,u_{\omega\omega'}\rangle\leq\alpha_{\omega\omega'}\langle x,\mathfrak e_{\omega\omega'}\rangle\}\subset\Gamma_{\omega\omega'-}.
$$
Let us prove that the construction stops. For $h\geq0$, let $\Omega^h$ be a
node of generation $h(n-1)$ and set
$$
W_j^h:=\sup_{y,y'\in\Omega^h}|\langle y-y',v_j\rangle|,
\qquad j=1,\dots,n-1.
$$
Assume that $\Omega^h$ is not terminal, and let $\Omega^{h+1}$ be any
descendant of $\Omega^h$ at generation $(h+1)(n-1)$. We claim that
$$
W_j^{h+1}\leq \frac12 W_j^h+\eta_0
\qquad\text{for every }j=1,\dots,n-1.
$$
Indeed, when the coordinate $v_j$ is split, the interval of possible values of
$\langle y,v_j\rangle$ is cut around its midpoint with an overlap of size
$2\eta_0$. Thus each child has $v_j$-width at most $W_j^h/2+\eta_0$, and the
remaining intersections cannot increase widths. The same estimate holds for
the first block. Hence
$$
W_j^h\leq
2^{-h}\frac{2\theta}{\sqrt{1-\theta^2}}+2\eta_0,
\qquad j=1,\dots,n-1.
$$
If
$h\geq\log_2\frac{3200\theta\sqrt{n-1}}{\sigma\sqrt{1-\theta^2}},$
then
$$
2^{-h}\frac{2\theta}{\sqrt{1-\theta^2}}
\leq
\frac{\sigma}{1600\sqrt{n-1}},
$$
and therefore
$$
W_j^h
\leq
\frac{\sigma}{1600\sqrt{n-1}}
+
\frac{\sigma}{800\sqrt{n-1}}
<
\frac{\sigma}{16\sqrt{n-1}}.
$$
Thus the stopping condition is reached after fewer than
$2+
\log_2\frac{1600\theta\sqrt{n-1}}{\sigma\sqrt{1-\theta^2}}$ blocks. Since each block consists of $n-1$ generations, the tree has depth
smaller than
$n
\log_2(\frac{6400\theta\sqrt{n-1}}{\sigma\sqrt{1-\theta^2}})$.

Notice that if $\omega$ is non-terminal, then
$\Omega_{\omega}=\Omega_{\omega+}\cup\Omega_{\omega-}$
by construction, and therefore
$$
\Gamma_{\omega}=
\Gamma_{\omega+}\cup\Gamma_{\omega-}.
$$

It remains to estimate $K_{\omega}$. For the first block, including the root, $\mathfrak e_{\omega}=e$ and
$u_{\omega}\in e^\perp$, so for $x=\lambda(e+y)$ one has
$$
\frac{|\langle x,u_{\omega}\rangle|}{\langle\mathfrak e_{\omega},x\rangle}\leq|y|\leq\frac{\theta}{\sqrt{1-\theta^2}}.
$$
Let now ${\omega\omega'}$ be a non-terminal node, where as usual $|\omega|=(n-1)h$ and $|\omega'|=k$ with $1\leq k\leq n-1$. Since $\Omega_{\omega\omega'}$ is contained in one of the fattened orthants produced in the first $n-1$ generations, a simple algebraic computation shows that any two points $y,z\in\Omega_{\omega\omega'}$ satisfy
$$
\langle y,v_j\rangle\langle z,v_j\rangle
\geq
-\eta_0\frac{\theta}{\sqrt{1-\theta^2}}
\qquad\text{for every }j=1,\dots,n-1.
$$
In particular, using the definition of $\eta_0$, we have
$$
\langle y,z_{\omega\omega'}\rangle
\geq
-(n-1)\eta_0\frac{\theta}{\sqrt{1-\theta^2}}=-\frac{\sigma\theta\sqrt{n-1}}{1600\sqrt{1-\theta^2}}\geq -\frac{1}{2}.
$$
Thus, for $x=\lambda(e+y)\in\Gamma_{\omega\omega'}$, with the choices of $\mathfrak e_{\omega\omega'}$, $u_{\omega\omega'}$ and $\alpha$ we have
$$
\langle\mathfrak e_{\omega\omega'},x\rangle
=
\lambda\frac{1+\langle z_{\omega\omega'},y\rangle}{\sqrt{1+|z_{\omega\omega'}|^2}}
\geq
\frac{\lambda\sqrt{1-\theta^2}}2,
\qquad
|\langle x,u_{\omega\omega'}\rangle|\leq |x|\leq\frac{\lambda}{\sqrt{1-\theta^2}}.
$$
Hence $K_{\omega\omega'}\leq2/(1-\theta^2)$ for every non-terminal node ${\omega\omega'}$. Since $\alpha_{\omega\omega'}=\eta_0(1-\theta^2)$ for every non-terminal node, we obtain
$$
\frac{\alpha_{\omega\omega'}}{2K_{\omega\omega'}}
\geq
\frac{\eta_0(1-\theta^2)^2}{4}
=
\frac{\sigma(1-\theta^2)^2}{6400\sqrt{n-1}}.
$$
We are left showing that each of the constructed terminal cones lies in some cone $C(e_i,\sigma)$. Let ${\omega}$ be terminal. The map $y\mapsto(e+y)/|e+y|$ is $1$-Lipschitz, because $|e+y|\geq1$. Therefore
$$\operatorname{diam}(\Gamma_{\omega}\cap \mathbb{S}^{n-1})
\leq\operatorname{diam}\Omega_{\omega}\leq\frac{\sigma}{16},$$
where the last inequality comes from the stopping condition \eqref{stopaconsg}.
Choose $\xi_{\omega}\in\Gamma_{\omega}\cap \mathbb{S}^{n-1}$. Since $\{e_i\}_{i=1}^N$ is a maximal $\sigma/4$-net of $\mathbb{S}^{n-1}\cap C(e,\theta)$, it is a $\sigma/4$-covering of that set. Hence there exists $i$ such that $|\xi_{\omega}-e_i|\leq\sigma/4$. Therefore every $\xi\in\Gamma_{\omega}\cap \mathbb{S}^{n-1}$ satisfies
$$
|\xi-e_i|\leq |\xi-\xi_{\omega}|+|\xi_{\omega}-e_i|\leq\frac{\sigma}{16}+\frac{\sigma}{4}<\sigma.
$$
Since $|\xi-e_i|<\sigma$ and $|\xi|=|e_i|=1$, we get
$\langle \xi,e_i\rangle>1-\sigma^2/2\geq\sqrt{1-\sigma^2}$. Hence $\xi\in C(e_i,\sigma)$, and so $\Gamma_{\omega}\subset C(e_i,\sigma)$.
This concludes the proof.
\end{proof}

The next proposition is the step where we pass from estimates on the two subfamilies of curves to estimates on all $C$-curves. Recall that, by construction, $C=C^+\cup C^-$. At each level, the set $\Omega_{h+1}$ is covered by balls of the previous generation, and these balls are much smaller than the balls from the previous generation that they intersect. We assume that inside each such ball the portion of $\Omega_{h+1}$ seen by every $C^+$-curve and every $C^-$-curve is small. The point is to prove the corresponding decay for an arbitrary $C$-curve. The proof is local: inside each covering ball we use the two-leaves lemma to replace a definite piece of a $C$-curve by either a $C^+$-curve or a $C^-$-curve, where the assumed estimate applies. This gives a fixed loss of length in each selected ball. A Vitali selection makes the relevant balls disjoint, so the local gains can be summed. Finally, the scale separation between consecutive generations keeps the enlarged balls inside the slightly enlarged ball $\Lambda(1+\mathfrak t)B_0$, which gives the localized one-step estimate and then its iteration.

\begin{proposizione}\label{prop:two-leaves-iterated-substitution-additive}
Let $C,C^+,C^-,e,u,c,\alpha,K_0$ be as in
\cref{lem:acp-two-convex-slope-cones}, and set
$\mathfrak q:=\min\{1,\tfrac{\alpha}{8K_0}\}$.
Suppose further that $C\subseteq C(e,\theta)$ for some $\theta\in (0,1)$, let $\tfrac43\leq \Lambda\leq 9$, $K\in \N$ and let $\mathfrak t\in (0,1)$ satisfy
$$\mathfrak t\leq \frac{\sqrt{1-\theta^2}}{16}\qquad\text{and}\qquad (1+\mathfrak t)^K\leq \frac{10}{9}.$$
Let $\Omega_0\supseteq \Omega_1\supseteq \ldots \supseteq\Omega_K$ be nested open sets, and let $\mathscr B^0,\ldots,\mathscr B^K$ be families of balls such that, for every $h=0,\ldots,K-1$, the following hold
\begin{enumerate}
  \item $\Omega_{h+1}\subseteq \bigcup\{B:B\in \mathscr B^{h}\}$, $20B\subseteq \Omega_h$ for every $B\in \mathscr B^{h+1}$ and if $B\in \mathscr B^h$ and $B'\in\mathscr B^{h+1}$ satisfy
$10 B\cap B'\neq \emptyset$, then
$$
\diam B'\leq \frac{\mathfrak t}{20}\diam B.
$$
  \item for every $B\in \mathscr B^{h}$ and every $C^+$-curve or $C^-$-curve $\gamma$ we have
$$
\Haus^1(\gamma \cap \Omega_{h+1}\cap \tfrac{4}{3}B)
\leq
\frac{\mathfrak q(1-\theta^2)}{32}
\Haus^1(\gamma \cap \tfrac{4}{3}(1+\mathfrak t)B).
$$
\end{enumerate}
Then for every $h_0=0,\ldots, K-1$ and $h_0\leq h\leq K-1$, every $B_0\in\mathscr B^{h_0}$ and every $C$-curve $\gamma$ we have 
$$\Haus^1(\gamma \cap \Omega_{h+1}\cap \Lambda B_0)\leq\Big(1-\frac{\mathfrak q(1-\theta^2)^{2}}{5120}\Big)\Haus^1(\gamma \cap \Omega_{h}\cap \Lambda(1+\mathfrak t)B_0).$$
In particular, for every $0\leq h_0\leq h_2<h_1\leq K$, every $B_0\in \mathscr B^{h_0}$ and every $C$-curve $\gamma$ we have
$$\Haus^1(\gamma \cap \Omega_{h_1}\cap \Lambda B_0)\leq\Big(1-\frac{\mathfrak q(1-\theta^2)^{2}}{5120}\Big)^{h_1-h_2}\Haus^1(\gamma \cap \Omega_{h_2}\cap \Lambda(1+\mathfrak t)^{h_1-h_2}B_0).$$
In addition, for every $0\leq h_2<h_1\leq K$ and every $C$-curve $\gamma$ we have
$$\Haus^1(\gamma \cap \Omega_{h_1})\leq\Big(1-\frac{\mathfrak q(1-\theta^2)^{2}}{5120}\Big)^{h_1-h_2}\Haus^1(\gamma \cap \Omega_{h_2}).$$
\end{proposizione}

\begin{proof}We can assume without loss of generality that $\gamma$ is parametrized as usual by a line parallel to $e$, namely $\gamma(t)=te+\eta(t)$ where $\eta(t)-\eta(s)\in e^\perp$ for every $s,t\in \R$. 
Let us fix $h_0=0,\ldots,K-1$, let us fix a ball $B_0\in \mathscr B^{h_0}$ and let $\gamma$ be a $C$-curve intersecting $\Lambda B_0$. 

For every $\ell=h_0+1,\ldots,K$, let $\mathscr B^{\ell}(\gamma)$ be the family of balls $B\in \mathscr B^{\ell}$ such that $B\cap \gamma \neq \emptyset$, and let $\mathscr B^{\ell}(B_0,\gamma)$ be the family of balls $B\in\mathscr B^\ell(\gamma)$ such that $B\cap \gamma\cap \Lambda B_0\neq \emptyset$. 
By Vitali's covering argument we can extract subfamilies $\tilde {\mathscr B}^{\ell}(\gamma)\subseteq \mathscr B^\ell(\gamma)$ and $\tilde{\mathscr B}^{\ell}(B_0,\gamma)\subseteq \mathscr B^\ell(B_0,\gamma)$ such that the balls $\tfrac{4}{3}B$ are disjoint and such that 
$$
\bigcup_{B\in \mathscr B^{\ell}(\gamma)}\tfrac{4}{3}B\subseteq\bigcup_{B\in \tilde{\mathscr{B}}^{\ell}(\gamma)} \tfrac{20}{3}B \subseteq \Omega_{\ell-1}
\quad\text{and similarly}\quad
\bigcup_{B\in \mathscr B^{\ell}(B_0,\gamma)}\tfrac{4}{3}B\subseteq\bigcup_{B\in \tilde{\mathscr{B}}^{\ell}(B_0,\gamma)} \tfrac{20}{3}B \subseteq \Omega_{\ell-1}.
$$
Notice that if $B\in {\mathscr B}^\ell(B_0,\gamma)$, then by iterating the hypothesis on the families $\mathscr B^\ell$ we have
$$
\diam B\leq \Big(\frac{\mathfrak t}{20}\Big)^{\ell-h_0}\diam B_0\leq \frac{\mathfrak t}{20}\diam B_0.
$$
This implies in particular that 
\begin{equation}
    \bigcup_{B\in \mathscr B^{\ell}(B_0,\gamma)}\tfrac{20}{3}B\subseteq \Omega_{\ell-1}\cap \Lambda(1+\mathfrak t)B_0.
    \label{inclusionecollasdngi}
\end{equation}
Indeed, if $B\cap\Lambda B_0\neq\emptyset$, then
$\dist(c_B,c_{B_0})\leq \Lambda r(B_0)+r(B)$,
and since $r(B)\leq \frac{\mathfrak t}{20}r(B_0)$ and $\Lambda\geq \tfrac43$, we have
$$
\dist(c_B,c_{B_0})+\frac{20}{3}r(B)
\leq
\Lambda r(B_0)+\frac{23}{3}r(B)
\leq
\Lambda(1+\mathfrak t)r(B_0).
$$

From now on we consider $h\geq h_0$ to be fixed.
Let $B\in \mathscr B^{h+1}(\gamma)$ and notice that by \cref{propo:curvainpalla} there exists an interval $I\subseteq\gamma^{-1}(\tfrac76B)$ such that
$$
\Leb^1(I)\geq \sqrt{1-\theta^2}r(B)/6.
$$
Moreover, $\Haus^1(\gamma\cap \tfrac76B)\geq \frac16 r(B)$, where $r(B)$ is the radius of the ball $B$. Thanks to \cref{lem:acp-two-convex-slope-cones}, possibly choosing the $C^-$ alternative instead of the $C^+$ alternative, we can assume without loss of generality that there exists a $C^+$-curve $\tilde \gamma$ such that we can find a compact set $S\subseteq I$ for which $\tilde \gamma(t)=\gamma(t)$ for every $t\in S$ and 
$$
\Leb^1(S)\geq \mathfrak{q}\Leb^1(I)\geq \mathfrak q\sqrt{1-\theta^2}r(B)/6.
$$
In particular, since the inverse of $\gamma$ is $1$-Lipschitz, and by \cref{propo:curvainpalla} we have
$$
\Haus^1(\gamma \cap \tilde \gamma\cap \tfrac{4}{3}B)\geq \frac{\mathfrak q\sqrt{1-\theta^2}r(B)}{6}
\qquad\text{and}\qquad
\Haus^1(\gamma\cap \tfrac{4}{3}B)\leq\frac{8r(B)}{3\sqrt{1-\theta^2}},
$$
and in particular we infer that 
\begin{equation}
    \Haus^1(\gamma \cap \tilde \gamma\cap \tfrac{4}{3}B)\geq \frac{\mathfrak q(1-\theta^2)}{16}\Haus^1(\gamma\cap \tfrac{4}{3}B).
    \label{stimacollasodifna}
\end{equation}
Similarly, since $\tilde \gamma$ intersects $\frac76B$, we similarly infer that 
\begin{equation}
    \Haus^1(\gamma \cap \tilde \gamma\cap \tfrac{4}{3}B)\geq \frac{\mathfrak q(1-\theta^2)}{16(1+\mathfrak t)}\Haus^1(\tilde\gamma\cap \tfrac{4}{3}(1+\mathfrak t)B).
    \label{equationalskdmflaskmdflaskdmf}
\end{equation}
We observe that
\begin{equation}
\begin{split}
&\Haus^1(\gamma\cap \Omega_{h+1}\cap \tfrac{4}{3}B)
\leq \Haus^1(\tilde \gamma \cap \Omega_{h+1}\cap \tfrac{4}{3}B)+\Haus^1((\gamma\setminus \tilde \gamma) \cap \Omega_{h+1}\cap \tfrac{4}{3}B)\\
&\qquad\leq \frac{\mathfrak q(1-\theta^2)}{32}\Haus^1(\tilde\gamma\cap \tfrac{4}{3}(1+\mathfrak t)B)+\Haus^1((\gamma\setminus \tilde \gamma) \cap \Omega_{h+1}\cap \tfrac{4}{3}B)\\
&\qquad=\frac{\mathfrak q(1-\theta^2)}{32}\Haus^1(\tilde\gamma\cap (1+\mathfrak t)\tfrac{4}{3}B)+\Haus^1((\gamma\setminus \tilde \gamma) \cap \Omega_{h+1}\cap \tfrac{4}{3}B)\\
&\qquad\leq \frac{\mathfrak q(1-\theta^2)}{32}\Haus^1(\tilde\gamma\cap (1+\mathfrak t)\tfrac{4}{3}B)-\Haus^1(\gamma\cap \tilde\gamma \cap \tfrac{4}{3}B)+\Haus^1(\gamma\cap \Omega_{h}\cap \tfrac{4}{3}B).
\nonumber
\end{split}
\end{equation}
Thanks to \eqref{equationalskdmflaskmdflaskdmf} this implies that 
\begin{equation}
    \begin{split}
        \Haus^1(\gamma\cap \Omega_{h+1}\cap \tfrac{4}{3}B)\leq -\Big(1-\frac{1+\mathfrak t}{2}\Big)\Haus^1(\gamma \cap \tilde \gamma\cap \tfrac{4}{3}B)+\Haus^1(\gamma\cap \Omega_{h}\cap \tfrac{4}{3}B)
    \end{split}
\end{equation}
Thanks to \eqref{stimacollasodifna}, we finally deduce that 
\begin{equation}
    \begin{split}
\Haus^1(\gamma\cap \Omega_{h+1}\cap \tfrac{4}{3}B)
&\leq\Haus^1(\gamma\cap \tfrac{4}{3}B)-\frac{\mathfrak q(1-\theta^2)(1-\mathfrak t)}{32}\Haus^1(\gamma \cap \tfrac{4}{3}B)\\
&\leq
\Big(1-\frac{\mathfrak q(1-\theta^2)}{64}\Big)
\Haus^1(\gamma\cap \tfrac{4}{3}B),
\label{eq:stimaapsdofasdf}
\end{split}
\end{equation}
where the last inequality comes from the choice of $\mathfrak t$.

Let us now estimate how long is the portion of the curve $\gamma$ covered by $\tilde{\mathscr B}^{h+1}(B_0,\gamma)$. Thanks to the fact that $\gamma$ is a $C(e,\theta)$-curve and the balls $B\in \tilde{\mathscr B}^{h+1}(B_0,\gamma)$ cover $\gamma\cap \Omega_{h+1}\cap \Lambda B_0$, we have 
\begin{equation}
\begin{split}
\Haus^1(\gamma\cap \Lambda B_0\cap \Omega_{h+1})
\leq& \sum_{B\in \tilde{\mathscr {B}}^{h+1}(B_0,\gamma)}\Haus^1(\tfrac{20}{3}B\cap \gamma\cap  \Omega_{h+1})
\leq \frac{1}{\sqrt{1-\theta^2}}\sum_{B\in \tilde{\mathscr B}^{h+1}(B_0,\gamma)}\tfrac{40}{3}r(B)\\
\leq& \frac{40}{\sqrt{1-\theta^2}}  \Haus^1\Big(\gamma\cap \bigcup_{B\in \tilde{\mathscr B}^{h+1}(B_0,\gamma)}\tfrac{4}{3}B\cap \Lambda(1+\mathfrak t)B_0\Big)\\
=&\frac{40}{\sqrt{1-\theta^2}}  \Haus^1\Big(\gamma\cap \Omega_h\cap\Lambda(1+\mathfrak t)B_0\cap \bigcup_{B\in \tilde{\mathscr B}^{h+1}(B_0,\gamma)}\tfrac{4}{3}B\Big),
\label{asodasgnjnsdf}
\end{split}  
\end{equation}
where in the inequality in the second line above we used the fact that every ball $B\in \tilde{\mathscr {B}}^{h+1}(B_0,\gamma)$ satisfies $\tfrac{20}{3}B\subseteq \Lambda(1+\mathfrak t)B_0$, thanks to the discussion at the beginning of the proof. 
Let us observe that
\begin{equation}
    \begin{split}
      &  \Haus^1(\gamma \cap \Omega_{h+1}\cap \Lambda B_0)
=\Haus^1\Big(\gamma \cap \Omega_{h+1}\cap \Lambda B_0\setminus \bigcup_{B\in \tilde{\mathscr B}^{h+1}(B_0,\gamma)}\tfrac{4}{3}B\Big)\\
&\qquad\qquad\qquad\qquad\qquad\qquad\qquad\qquad+\Haus^1\Big(\gamma \cap \Omega_{h+1}\cap \Lambda B_0\cap  \bigcup_{B\in \tilde{\mathscr B}^{h+1}(B_0,\gamma)}\tfrac{4}{3}B\Big)\\
&\overset{\eqref{eq:stimaapsdofasdf}}{\leq} \Haus^1\Big(\gamma \cap \Omega_{h+1}\cap \Lambda B_0\setminus \bigcup_{B\in \tilde{\mathscr B}^{h+1}(B_0,\gamma)}\tfrac{4}{3}B\Big)+\Big(1-\frac{\mathfrak q(1-\theta^2)}{64}\Big)\sum_{B\in \tilde{\mathscr B}^{h+1}(B_0,\gamma)}\Haus^1(\gamma\cap \tfrac{4}{3}B\cap \Omega_{h}).
\nonumber
    \end{split}
\end{equation}
Since the balls $\tfrac43B$, with $B\in \tilde{\mathscr B}^{h+1}(B_0,\gamma)$, are pairwise disjoint, and thanks to \eqref{inclusionecollasdngi}, we have
\begin{equation}
    \begin{split}
     \sum_{B\in \tilde{\mathscr{B}}^{h+1}(B_0,\gamma)}\Haus^1(\gamma\cap \tfrac{4}{3}B\cap \Omega_{h})\leq \Haus^1\Big(\gamma\cap \Omega_h\cap \Lambda(1+\mathfrak t)B_0\cap \bigcup_{B\in \tilde{\mathscr B}^{h+1}(B_0,\gamma)}\tfrac{4}{3}B\Big)
     \label{asdjnasijgnairsngapsringapweriong}
    \end{split}
\end{equation}
By \eqref{asdjnasijgnairsngapsringapweriong}, and reorganizing terms, we have 
\begin{equation}
    \begin{split}
        &\Haus^1(\gamma \cap \Omega_{h+1}\cap \Lambda B_0)\leq \Haus^1(\gamma \cap \Omega_{h}\cap \Lambda(1+\mathfrak t) B_0)\\
        &\qquad\qquad\qquad\qquad\qquad\qquad\qquad\qquad-\frac{\mathfrak q(1-\theta^2)}{64}\Haus^1\Big(\gamma\cap \Omega_h\cap \Lambda(1+\mathfrak t)B_0\cap \bigcup_{B\in \tilde{\mathscr B}^{h+1}(B_0,\gamma)}\tfrac{4}{3}B\Big)\\
        &\overset{\eqref{asodasgnjnsdf}}{\leq} \Haus^1(\gamma \cap \Omega_{h}\cap \Lambda(1+\mathfrak t) B_0)-\frac{\mathfrak q(1-\theta^2)^{3/2}}{2560}
         \Haus^1(\gamma\cap \Lambda B_0\cap \Omega_{h+1})
        \nonumber
    \end{split}
\end{equation}
Rearranging terms we conclude that 
\begin{equation}
\begin{split}
\Haus^1(\gamma \cap \Omega_{h+1}\cap \Lambda B_0)\leq \Big(1-\frac{\mathfrak q(1-\theta^2)^{2}}{5120}\Big)\Haus^1(\gamma\cap \Lambda(1+\mathfrak t)B_0\cap  \Omega_h).
\nonumber
\end{split}
\end{equation}
This concludes the proof of the first part of the proposition. With the same argument we can also show that for every $h=0,\ldots, K-1$ we have
$$
\Haus^1(\gamma \cap \Omega_{h+1})
\leq \Big(1-\frac{\mathfrak q(1-\theta^2)^{2}}{5120}\Big)\Haus^1(\gamma \cap \Omega_{h}).
$$
Indeed, one applies the preceding argument without localizing inside $\Lambda B_0$, using the family $\tilde{\mathscr B}^{h+1}(\gamma)$ instead of $\tilde{\mathscr B}^{h+1}(B_0,\gamma)$.

Notice that the preceding argument remains valid with $\Lambda$ replaced by any $\Lambda'\in[\tfrac43,10]$, since the separation assumption is stated with $10B$ and the proof only uses the bounds $\tfrac43\leq\Lambda'\leq10$. Let
$$
\rho:=1-\frac{\mathfrak q(1-\theta^2)^2}{5120}.
$$
The base step is proved above. Assume that, for some $d\geq2$, the claim has been shown for every
$1\leq \ell\leq d-1$. Namely, for every $0\leq h_0\leq h\leq K-\ell$ and every
$B_0\in \mathscr B^{h_0}$, we have
$$
\Haus^1(\gamma \cap \Omega_{h+\ell}\cap \Lambda B_0)
\leq
\rho^\ell
\Haus^1(\gamma \cap \Omega_h\cap \Lambda(1+\mathfrak t)^\ell B_0).
$$
Then, for $0\leq h_0\leq h\leq K-d$, using first the induction hypothesis with
$\ell=d-1$ and then the one-step estimate, we obtain
\begin{equation}
\begin{split}
\Haus^1(\gamma \cap \Omega_{h+d}\cap \Lambda B_0)
\leq& \rho^{d-1}\Haus^1(\gamma \cap \Omega_{h+1}\cap \Lambda(1+\mathfrak t)^{d-1}B_0)\\
\leq& \rho^{d}\Haus^1(\gamma \cap \Omega_{h}\cap \Lambda(1+\mathfrak t)^dB_0).
\nonumber
\end{split}
\end{equation}
Similarly, applying the global one-step estimate inductively, for every $h_1>h_2$ we get
$$
\Haus^1(\gamma \cap \Omega_{h_1})
\leq \rho^{h_1-h_2}\Haus^1(\gamma \cap \Omega_{h_2}).
$$
This concludes the proof.
\end{proof}

The following lemma is the abstract device which allows us to iterate the loss of length through the sets $\Omega_h$. The assumptions say that, at each level, the new set $\Omega_{h+1}$ is covered by balls of the next generation, and that inside each of these balls every $C(e,\theta)$-curve sees only an $\alpha$-portion of $\Omega_{h+1}$, measured with respect to the previous level. The conclusion is that each passage from $\Omega_h$ to $\Omega_{h+1}$ removes a fixed proportion of length, and hence the loss accumulates geometrically from level $h_2$ to level $h_1$.

\begin{lemma}\label{lem:iterated-local-decay}
Let $K\in\mathbb N$, $e\in \mathbb S^{n-1}$, $\theta\in (0,1)$, $\alpha\in (0,\frac14)$ and $\mathfrak t\in(0,1)$ be such that
$$
\alpha\leq \frac{\sqrt{1-\theta^2}}{16},
\qquad
\mathfrak t\leq \sqrt{1-\theta^2}
\qquad\text{and}\qquad
(1+\mathfrak t)^K\leq 2.
$$
Let
$\Omega_0\supseteq \Omega_1\supseteq \ldots \supseteq\Omega_K$ be nested open sets, and let $\mathscr B^h$, $h=0,\ldots,K$, be
families of balls. Assume that, for every $h=0,\ldots,K-1$, the
following hold.
\begin{enumerate}
    \item $\Omega_{h+1}\subseteq \bigcup\{B:B\in \mathscr B^{h+1}\}$,
    $20B\subseteq \Omega_h$ for every $B\in \mathscr B^{h+1}$, and if
    $B\in \mathscr B^h$, $B'\in\mathscr B^{h+1}$ and
    $4B\cap B'\neq \emptyset$, then
    $$
    \diam B'\leq \frac{\mathfrak t}{20}\diam B.
    $$

    \item For every $B\in \mathscr B^{h+1}$ and every
    $C(e,\theta)$-curve $\gamma$ one has
    $$
    \Haus^1(\gamma \cap \Omega_{h+1}\cap 2 B)
    \leq
    \alpha
    \Haus^1(\gamma \cap \Omega_h\cap 2 (1+\mathfrak t)B).
    $$
\end{enumerate}
Then, for every $h_1\geq h_2$ in $\{0,\ldots,K\}$, every
$B\in \mathscr B^{h_2}$ and every $C(e,\theta)$-curve $\gamma$,
one has
$$
\Haus^1(\gamma \cap \Omega_{h_1}\cap 2B)
\leq
\Big(1-\frac{\sqrt{1-\theta^2}}{80}\Big)^{h_1-h_2}
\Haus^1(\gamma \cap \Omega_{h_2}\cap 2(1+\mathfrak t)^{h_1-h_2}B).
$$
\end{lemma}

\begin{proof}
Without loss of generality we assume that $\gamma(t)=te+\eta(t)$ with $\eta(t)-\eta(s)\in e^\perp$.

We first prove the following one-step estimate. For every $0\leq h_0\leq h\leq K-1$, every
$B_0\in\mathscr B^{h_0}$, every $C(e,\theta)$-curve $\gamma$ and every
$1\leq\Lambda\leq 2$ with $\Lambda(1+\mathfrak t)\leq2$, one has
$$
\Haus^1(\gamma \cap \Omega_{h+1}\cap 2\Lambda B_0)
\leq
\Big(1-\frac{\sqrt{1-\theta^2}}{80}\Big)
\Haus^1(\gamma \cap \Omega_{h}\cap 2\Lambda (1+\mathfrak t)B_0).
$$

Let $\mathscr B^{h+1}(B_0,\gamma)$ be the family of balls $B\in\mathscr B^{h+1}$ such that
$B\cap\gamma\cap\Omega_{h+1}\cap 2\Lambda B_0\neq\emptyset $. 

For every $B\in\mathscr B^{h+1}(B_0,\gamma)$, by iterating the separation assumption along a chain of balls from level $h_0$ to level $h+1$, we have
$$
\diam B\leq \Big(\frac{\mathfrak t}{20}\Big)^{h+1-h_0}\diam B_0
\leq
\frac{\mathfrak t}{20}\diam B_0.
$$
In particular, $10B\subseteq \Omega_h\cap 2\Lambda(1+\mathfrak t)B_0$. Indeed, $20B\subseteq\Omega_h$. Moreover, if $B\cap 2\Lambda B_0\neq\emptyset$, then $\dist(c_B,c_{B_0})\leq 2\Lambda r(B_0)+r(B)$, and therefore
$$\dist(c_B,c_{B_0})+10r(B)\leq2\Lambda r(B_0)+11r(B)\leq2\Lambda(1+\mathfrak t)r(B_0),$$
because $r(B)\leq \frac{\mathfrak t}{20}r(B_0)$ and $\Lambda\geq1$. We further extract from $\mathscr B^{h+1}(B_0,\gamma)$ a subfamily $\tilde{\mathscr B}^{h+1}(B_0,\gamma)$ such that the balls $2B$ are disjoint and
$$
\gamma\cap\Omega_{h+1}\cap 2\Lambda B_0
\subseteq
\bigcup_{B\in \tilde{\mathscr B}^{h+1}(B_0,\gamma)}10B.
$$
In particular, we infer that 
$$
\bigcup_{B\in \tilde{\mathscr B}^{h+1}(B_0,\gamma)}2B
\subseteq
\bigcup_{B\in \tilde{\mathscr B}^{h+1}(B_0,\gamma)}10B
\subseteq
\Omega_h\cap 2\Lambda(1+\mathfrak t)B_0.
$$
Since $\gamma$ is a $C(e,\theta)$-curve and each $B\in \tilde{\mathscr B}^{h+1}(B_0,\gamma)$ intersects $\gamma$, by \cref{propo:curvainpalla} we have
$$
\Haus^1(\gamma\cap 10B)\leq \frac{20r(B)}{\sqrt{1-\theta^2}}
\qquad\text{and}\qquad
\Haus^1(\gamma\cap 2B)\geq r(B).
$$
Hence
\begin{equation}
    \begin{split}
    \Haus^1(\gamma\cap 2\Lambda B_0\cap \Omega_{h+1})&\leq\sum_{B\in \tilde{\mathscr B}^{h+1}(B_0,\gamma)}\Haus^1(\gamma\cap 10B)\leq\frac{20}{\sqrt{1-\theta^2}}\sum_{B\in \tilde{\mathscr B}^{h+1}(B_0,\gamma)}\Haus^1(\gamma\cap 2B)\\
    &\leq\frac{20}{\sqrt{1-\theta^2}}\Haus^1\Big(\gamma\cap \Omega_h\cap 2\Lambda(1+\mathfrak t)B_0\cap \bigcup_{B\in \tilde{\mathscr B}^{h+1}(B_0,\gamma)}2B\Big).
    \label{ea:sodifngsodfni}
    \end{split}
\end{equation}

Thanks to item 2, we also have, for every $B\in\tilde{\mathscr B}^{h+1}(B_0,\gamma)$,
$$\Haus^1(\gamma\cap \Omega_{h+1}\cap 2 B)\leq\alpha\Haus^1(\gamma\cap \Omega_h\cap 2 (1+\mathfrak t)B).$$
Moreover, again by \cref{propo:curvainpalla}, we know
$$\Haus^1(\gamma\cap 2(1+\mathfrak t)B)\leq\frac{4(1+\mathfrak t)r(B)}{\sqrt{1-\theta^2}}\leq\frac{8r(B)}{\sqrt{1-\theta^2}},$$
while, since we know that $\Haus^1(\gamma\cap 2B)\geq r(B)$, we have
$$\Haus^1(\gamma\cap 2(1+\mathfrak t)B)\leq\frac{8}{\sqrt{1-\theta^2}}\Haus^1(\gamma\cap 2B).$$
Since the balls $2B$, with $B\in \tilde{\mathscr B}^{h+1}(B_0,\gamma)$, are pairwise disjoint, we obtain
\begin{equation}
    \begin{split}
        &\sum_{B\in \tilde{\mathscr B}^{h+1}(B_0,\gamma)}\Haus^1(\gamma\cap \Omega_h\cap 2(1+\mathfrak t)B)\leq\frac{8}{\sqrt{1-\theta^2}}\Haus^1\Big(\gamma\cap \Omega_h\cap 2\Lambda(1+\mathfrak t)B_0\cap\bigcup_{B\in \tilde{\mathscr B}^{h+1}(B_0,\gamma)}2B\Big).
        \nonumber
    \end{split}
\end{equation}
This implies in particular that 
\begin{equation}
    \begin{split}
         &\Haus^1(\gamma \cap \Omega_{h+1}\cap 2\Lambda  B_0)=\Haus^1\Big(\gamma \cap \Omega_{h+1}\cap 2\Lambda  B_0\setminus \bigcup_{B\in \tilde{\mathscr B}^{h+1}(B_0,\gamma)}2B\Big)\\
&\qquad\qquad\qquad\qquad\qquad\qquad\qquad\qquad+\Haus^1\Big(\gamma \cap \Omega_{h+1}\cap 2\Lambda B_0\cap\bigcup_{B\in \tilde{\mathscr B}^{h+1}(B_0,\gamma)}2B\Big)\\
         &\leq
         \Haus^1\Big(\gamma \cap \Omega_{h}\cap 2\Lambda(1+\mathfrak t)B_0
         \setminus \bigcup_{B\in \tilde{\mathscr B}^{h+1}(B_0,\gamma)}2B\Big)+
         \alpha
         \sum_{B\in \tilde{\mathscr B}^{h+1}(B_0,\gamma)}
         \Haus^1(\gamma\cap \Omega_h\cap 2(1+\mathfrak t)B)\\
         &\leq
         \Haus^1(\gamma\cap \Omega_h\cap 2\Lambda(1+\mathfrak t)B_0)-
         \Big(1-\frac{8\alpha}{\sqrt{1-\theta^2}}\Big)
         \Haus^1\Big(\gamma\cap \Omega_h\cap 2\Lambda(1+\mathfrak t)B_0
         \cap\bigcup_{B\in \tilde{\mathscr B}^{h+1}(B_0,\gamma)}2B\Big).
         \nonumber
    \end{split}
\end{equation}
By the choice of $\alpha$, and combining this with \eqref{ea:sodifngsodfni}, we obtain
\begin{equation}
    \begin{split}
    \Haus^1(\gamma \cap \Omega_{h+1}\cap 2\Lambda  B_0)&\leq\Haus^1(\gamma\cap \Omega_h\cap 2\Lambda(1+\mathfrak t)B_0)-\frac{\sqrt{1-\theta^2}}{40}\Haus^1(\gamma \cap \Omega_{h+1}\cap 2\Lambda  B_0).
    \end{split}
\end{equation}
Rearranging terms, and using $\frac{1}{1+x}\leq 1-\frac{x}{2}$ for $0\leq x\leq1$, we get
\begin{equation}
\Haus^1(\gamma \cap \Omega_{h+1}\cap 2\Lambda  B_0)
\leq
\Big(1-\frac{\sqrt{1-\theta^2}}{80}\Big)
\Haus^1(\gamma\cap \Omega_h\cap 2\Lambda(1+\mathfrak t)B_0).
\end{equation}
This proves the one-step estimate.

We now prove the iterated estimate. More precisely, we prove by induction on $d$ that, for every
$0\leq h_0\leq h\leq K-d$, every $B_0\in\mathscr B^{h_0}$, every
$C(e,\theta)$-curve $\gamma$ and every $1\leq\Lambda\leq2$ with
$\Lambda(1+\mathfrak t)^d\leq2$, one has
$$
\Haus^1(\gamma\cap \Omega_{h+d}\cap 2\Lambda B_0)
\leq
\Big(1-\frac{\sqrt{1-\theta^2}}{80}\Big)^d
\Haus^1(\gamma\cap \Omega_h\cap 2\Lambda(1+\mathfrak t)^dB_0).
$$
The case $d=0$ is trivial, and the case $d=1$ is the one-step estimate proved above.

Assume now that $d\geq2$ and that the claim has been proved for $d-1$. Then, applying the induction hypothesis first and then the one-step estimate with $\Lambda$ replaced by $\Lambda(1+\mathfrak t)^{d-1}$, we obtain
\begin{equation}
    \begin{split}
\Haus^1(\gamma\cap \Omega_{h+d}\cap 2\Lambda B_0)
&\leq
\Big(1-\frac{\sqrt{1-\theta^2}}{80}\Big)^{d-1}
\Haus^1(\gamma\cap \Omega_{h+1}\cap 2\Lambda(1+\mathfrak t)^{d-1}B_0)\\
&\leq
\Big(1-\frac{\sqrt{1-\theta^2}}{80}\Big)^d
\Haus^1(\gamma\cap \Omega_h\cap 2\Lambda(1+\mathfrak t)^dB_0).
\end{split}
\end{equation}
This concludes the induction. Taking $h_0=h_2$, $h=h_2$, $B_0=B$, $\Lambda=1$ and $d=h_1-h_2$ gives the desired estimate.
\end{proof}

We now turn the failure of the Carleson packing condition for the family $\mathscr I_j$ into a structured collection of cubes with a fixed set of invariant directions. More precisely, after a finite pigeonholing of the admissible directions, we obtain arbitrarily long nested layers of cubes on which the same $j$-tuple $\mathfrak e$ realizes the invariant structure. At the same time, directions quantitatively transverse to $\mathrm{span}(\mathfrak e)$ remain non-invariant on every cube of the construction.

\begin{proposizione}\label{joiningcones}
 Suppose $\mu$ is a $k$-AD-regular measure with regularity constant $D$ and choose parameters $\sigma,\vartheta,\varepsilon,\delta,\zeta,\mathfrak m\in (0,1)$ and $A\geq 1$, $\qof\in \N$ satisfying the following properties 
$$\vartheta\leq \tfrac{1}{4},\qquad A\vartheta^n\geq 16n, \qquad\text{and}\qquad 16n\Big(\frac{\sigma}{\vartheta^n}
+n\frac{\varepsilon+\delta+\zeta+\mathfrak m}{A\vartheta^n}\Big)\leq \frac{1}{2\cdot 3^k 16^n D^2}.$$
Let $0\leq j<k$ and suppose $\mathscr I_j(2\vartheta,\sigma/2,\varepsilon,\delta,\zeta,\mathfrak m,A)$, one of the families of cubes yielded by \cref{rioridinocubi}, is not Carleson. Then, for every $\mathfrak N\in \N$ and $\eta\in (0,1)$ there exist 
$$\mathfrak e\in (\mathbb{S}^{n-1})^j,\,\text{a cube }R_{\mathfrak N}\in \mathscr I_j(2\vartheta,\sigma/2,\varepsilon,\delta,\zeta,\mathfrak m,A),\text{ a Borel set }U\subseteq R_{\mathfrak N},$$
$$\text{and families of cubes }\mathfrak L_0,\ldots, \mathfrak L_\mathfrak N\subseteq \mathscr I_j(2\vartheta,\sigma/2,\varepsilon,\delta,\zeta,\mathfrak m,A),$$ 
such that the following holds. 

\begin{enumerate}

\item for every $\iota=1,\ldots,\mathfrak N$, every cube $Q\in \mathfrak L_\iota$ is $(j,\mathfrak{e},\vartheta^j,\sigma, \varepsilon,\delta,\zeta,\mathfrak m,A)$-invariant.

\item $\mathfrak L_0=\{R_\mathfrak N\}$ and no cube appears in more than one of the families $\mathfrak L_0,\ldots,\mathfrak L_{\mathfrak N}$.

\item For every $m=0,\ldots,\mathfrak N$, the cubes in $\mathfrak L_m$ are pairwise disjoint. Moreover, whenever $Q,Q'\in \mathfrak L_m$, either
\begin{equation}
    \ell(Q)=\ell(Q'),
    \qquad\text{or, if }\ell(Q)<\ell(Q'),
    \text{ then }\qquad
    \frac{\ell(Q)}{\ell(Q')}\leq \frac{1}{2^\qof}.
    \label{separazionediscalev2asdfasdf}
\end{equation}

\item For every $m=1,\ldots,\mathfrak N$, each cube $Q'\in \mathfrak L_m$ is contained in a unique strictly larger cube $Q\in \mathfrak L_{m-1}$ such that
$$
\diam(Q)\geq 2^\qof\,\diam(Q').
$$
In addition, whenever $m=1,\ldots,\mathfrak N$, $Q\in \mathfrak L_{m-1}$ and $Q'\in \mathfrak L_m$ satisfy
$$
10A B_{Q'}\cap 10A B_Q\neq \emptyset\qquad\text{then}\qquad
\diam(Q)\geq 2^\qof\,\diam(Q').
$$
Finally, for every $m=1,\ldots, \mathfrak N$, $Q\in \mathfrak L_m$ if we let $Q'\in \mathfrak L_{m-1}$ be the ancestor of $Q$, we have 
$$10A B_Q\cap \supp\mu\subseteq Q'.$$
\item For every $m=1,\ldots,\mathfrak N$, every $Q\in \mathfrak L_m$, and every
$e\in X_c(\mathrm{span}(\mathfrak e),3\vartheta)$,
the cube $Q$ is not
$(e,\sigma/2,\varepsilon,\delta,\zeta,\mathfrak m,A)$-invariant, i.e., there exists a point $z(Q,e)\in AB_Q\cap \supp\mu$ such that every $C(e,\sigma/2)$-curve $\gamma$ intersecting $B(z(Q,e),\zeta \diam Q)$ satisfies
$$\Haus^1(2AB_Q\cap \gamma\setminus B(\supp\mu,\delta \diam Q))>\varepsilon\diam Q.$$
\item For every $m=1,\ldots,\mathfrak N$, every $e\in X_c(\mathrm{span}(\mathfrak e),3\vartheta)$ and every $Q\in\mathfrak L_m$ we have
$$\mu(B(z(Q,e),\zeta \diam Q/2)\cap U)\geq \frac{1}{2}\mu(B(z(Q,e),\zeta \diam Q/2)).$$
\item There holds $\mu(R_{\mathfrak N}\setminus U)\leq \eta \mu(R_{\mathfrak N})$ and $\mu(R_{\mathfrak N}\setminus \bigcup\{Q:Q\in \mathfrak L_{\mathfrak N}\})\leq \eta \mu(R_{\mathfrak N})$.
\end{enumerate}
\end{proposizione}

\begin{proof} We divide the proof in steps.

\textbf{Step I. Pigeonholing planes.} Let $0<\rho<1$ be such that $\rho\leq \sigma/16$ and such that the following two stability properties hold: whenever $v_1,\ldots,v_j\in\mathbb S^{n-1}$ are $(2\vartheta)^j$-separated and $e_1,\ldots,e_j\in\mathbb S^{n-1}$ satisfy $|e_i-v_i|\leq\rho$ for every $i=1,\ldots,j$, then $(e_1,\ldots,e_j)$ is $\vartheta^j$-separated and
$$
X_c(\mathrm{span}(e_1,\ldots,e_j),3\vartheta)\subseteq X_c(\mathrm{span}(v_1,\ldots,v_j),2\vartheta).
$$
Let $\mathscr E$ be a finite $\rho$-dense set in $\mathbb S^{n-1}$.
We first observe that, if $v,e\in\mathbb S^{n-1}$ and $|v-e|\leq\rho$, then $C(v,\sigma/2)\subseteq C(e,\sigma)$.  Indeed, if $u\in C(v,\sigma/2)\cap\mathbb S^{n-1}$, then
$$
\dist(u,\mathrm{span}(e))
\leq
\dist(u,\mathrm{span}(v))+\dist(v,\mathrm{span}(e))
\leq
\frac{\sigma}{2}+\rho
\leq
\sigma.
$$

Let us fix now $Q\in \mathscr I_j(2\vartheta,\sigma/2,\varepsilon,\delta,\zeta,\mathfrak m,A)$. By \cref{rioridinocubi}, there are $v_1,\ldots,v_j\in\mathbb S^{n-1}$ such that $Q$ is
$$
(j,(v_1,\ldots,v_j),(2\vartheta)^j,\sigma/2,\varepsilon,\delta,\zeta,\mathfrak m,A)
$$
-invariant. Choose $e_i\in\mathscr E$ with $|e_i-v_i|\leq\rho$ for every $i=1,\ldots,j$. By the choice of $\rho$, the tuple $\mathfrak e:=(e_1,\ldots,e_j)$ is $\vartheta^j$-separated. Moreover, since $C(v_i,\sigma/2)\subseteq C(e_i,\sigma)$ for every $i=1,\ldots,j$, the cube $Q$ is
$$
(j,\mathfrak e,\vartheta^j,\sigma,\varepsilon,\delta,\zeta,\mathfrak m,A)
$$
-invariant.

Furthermore, again by the choice of $\rho$, we have
$X_c(\mathrm{span}(\mathfrak e),3\vartheta)
\subseteq
X_c(\mathrm{span}(v_1,\ldots,v_j),2\vartheta)$.
Hence, by the second property in \cref{rioridinocubi}, applied with parameters $(2\vartheta,\sigma/2,\varepsilon,\delta,\zeta,\mathfrak m,A)$, the cube $Q$ is not
$(e,\sigma/2,\varepsilon,\delta,\zeta,\mathfrak m,A)$-invariant for every
$e\in X_c(\mathrm{span}(\mathfrak e),3\vartheta)$.

Since $\mathscr E$ is finite, the pigeonhole principle guarantees the existence of a single $\mathfrak e\in \mathscr E^j$ for which the corresponding family $\mathscr F_{\mathfrak e}$ of cubes satisfying the two conclusions above is not Carleson.

\textbf{Step II. Construction of the first layers.} Let us apply now \cref{nonClayersv2} to the non-Carleson family $\mathscr F_{\mathfrak e}$ with parameters $M:=\mathfrak{N}$, $W:=2^\qof$, $L:=2^\qof$, $\gimel:=10A$, and choose $\eta'$ so that 
$$\eta'\leq \frac{\eta}{(1+2D^2(\frac{12A}{\zeta})^k)^{\mathfrak N}}.$$
We infer that we can find a cube $R_{\mathfrak N}$ and layers of cubes $\mathfrak H_0,\ldots, \mathfrak H_{\mathfrak N}$ such that
\begin{itemize}
\item[(a)] $\mathfrak{H}_0=\{R_\mathfrak N\}$.
\item[(b)] No cube appears in more than one of the families $\mathfrak{H}_0,\dots,\mathfrak{H}_{\mathfrak{N}}$.
\item[(c)] The cubes in each family $\mathfrak{H}_m$, for $m=0,\dots,\mathfrak{N}$, are pairwise disjoint and for every $Q,Q'\in \mathfrak H_m$ we have 
\begin{equation}
    \ell(Q)=\ell(Q')\qquad\text{or, if $\ell(Q)<\ell(Q')$, then }\qquad\frac{\ell(Q)}{\ell(Q')}\leq \frac{1}{2^\qof};
    \label{separazionediscale2}
\end{equation}
\item[(d)] Each cube $Q'\in\mathfrak{H}_{m}$, for every $m=1,\dots,\mathfrak{N}$, is contained in a unique strictly larger cube $Q\in \mathfrak{H}_{m-1}$ with 
$$\diam(Q)\geq 2^\qof\,\diam(Q').$$
Secondly, if $Q\in \mathfrak H_{m-1}$, then for all $Q'\in \mathfrak H_{m}$ for which $10A B_{Q'}\cap 10A B_Q\neq \emptyset$, then 
$$\diam (Q)\geq 2^\qof\diam (Q').$$
\item[(e)] for every $m=1,\ldots,\mathfrak{N}$, $Q\in \mathfrak H_m$ if we let $Q'\in \mathfrak H_{m-1}$ be the ancestor of $Q$, we have 
$$10A B_Q\cap \supp\mu\subseteq Q';$$
\item[(f)] $\displaystyle\sum_{Q\in\mathfrak{H}_\mathfrak{N}}\mu(Q)\ge(1-\eta')\mu(R_{\mathfrak N}).$
\end{itemize}
We finally choose 
$U:=\bigcup\{Q:Q\in \mathfrak H_{\mathfrak N}\}$ and
notice that thanks to item (f) above, we have that 
$$\mu(R_{\mathfrak N}\setminus U)\leq \eta'\mu(R_{\mathfrak N}).$$

\textbf{Step III. Subselection of cubes satisfying item 6. and 7.} 
For every $m=1,\ldots,\mathfrak N$ we denote by 
$\mathfrak B_m$ the family of cubes $Q\in \mathfrak H_m$ for which there exist $e_Q\in X_c(\mathrm{span}(\mathfrak e),3\vartheta)$ and a point $z_Q\in AB_Q\cap\supp\mu$ satisfying the non-invariance conclusion for $Q$ and $e_Q$, and such that
$$\mu(B(z_Q,\zeta \diam Q/2)\cap U)< \frac{1}{2}\mu(B(z_Q,\zeta \diam Q/2)).$$
Notice that by construction we have $B(z_Q,\zeta \diam Q/2)\cap \supp\mu\subseteq 2AB_Q\cap \supp\mu\subseteq R_{\mathfrak N}$, and hence, defined $E:=R_{\mathfrak N}\setminus U$, we have 
$$
\mu(E\cap B(z_Q,\zeta \diam Q/2))\geq \frac{1}{2}\mu(B(z_Q,\zeta \diam Q/2)).
$$
By the Vitali covering lemma, there exists a subfamily
$\mathfrak S_m\subset\mathfrak B_m$ and a partition $\{\mathfrak F_Q:Q\in \mathfrak S_m\}$ of $\mathfrak B_m$ such that the balls $B(z_Q,\zeta \diam Q/2)$ with $Q\in \mathfrak S_m$ are pairwise disjoint and such that, whenever
$Q'\in\mathfrak{F}_Q$, we have
$$
B(z_Q,\zeta \diam Q/2)\cap B(z_{Q'},\zeta\diam Q'/2)\neq\emptyset
\qquad\text{and}\qquad
\diam Q'\leq \diam Q.
$$
Then we have $Q'\subseteq 6AB_Q$. Since the cubes in $\mathfrak H_m$ are pairwise disjoint, it follows that for every $Q\in \mathfrak S_m$ we have
$$
\sum_{Q'\in \mathfrak{F}_Q}\mu(Q')\leq \mu(6AB_Q)\leq D^2\Big(\frac{12A}{\zeta}\Big)^k\mu(B(z_Q,\zeta \diam Q/2)).
$$
In particular, we infer that 
\begin{equation}
    \begin{split}
        &\sum_{Q\in \mathfrak B_m}\mu(Q)
        \leq D^2\Big(\frac{12A}{\zeta}\Big)^k
        \sum_{Q\in \mathfrak S_m}\mu(B(z_Q,\zeta \diam Q/2))\\
        \leq& 2D^2\Big(\frac{12A}{\zeta}\Big)^k
        \sum_{Q\in \mathfrak S_m}\mu(B(z_Q,\zeta \diam Q/2)\cap E)
        \leq 2D^2\Big(\frac{12A}{\zeta}\Big)^k\eta' \mu(R_{\mathfrak N}).
        \label{badcubemassjoiningcones}
    \end{split}
   \end{equation}
Now, we define $\mathfrak L_0:=\{R_{\mathfrak N}\}$, $\mathfrak B_0:=\emptyset$ and for every $m\geq 1$, we inductively let 
$$
\mathfrak L_m:=\{Q\in \mathfrak H_m\setminus \mathfrak B_m:\text{there exists $Q'\in \mathfrak L_{m-1}$ such that $Q\subseteq Q'$}\}.
$$
In order to ease notations we let
$C_0:=2D^2\Big(\frac{12A}{\zeta}\Big)^k$. We claim that, for every $m=0,\ldots,\mathfrak N$,
\begin{equation}
\mu\Big(R_{\mathfrak N}\setminus\bigcup_{Q\in\mathfrak L_m}Q\Big)
\leq
(1+mC_0)\eta'\mu(R_{\mathfrak N}).
\label{formulacopricopri}
\end{equation}
Indeed, let $x\in R_{\mathfrak N}\setminus\bigcup_{Q\in\mathfrak L_m}Q$. If $x\notin U$, then $x\in R_{\mathfrak N}\setminus U$. Otherwise, since $U=\bigcup\{Q:Q\in\mathfrak H_{\mathfrak N}\}$, there is a chain of ancestors
$Q_i(x)\in\mathfrak H_i$ with $i=0,\ldots,m$ with $x\in Q_i(x)$. Since $Q_0(x)=R_{\mathfrak N}\in\mathfrak L_0$ and $Q_m(x)\notin\mathfrak L_m$, there exists a first index $i\in\{1,\ldots,m\}$ such that
$$
Q_{i-1}(x)\in\mathfrak L_{i-1}
\qquad\text{and}\qquad
Q_i(x)\notin\mathfrak L_i.
$$
By the definition of $\mathfrak L_i$, this implies $Q_i(x)\in\mathfrak B_i$. Hence
$$
R_{\mathfrak N}\setminus\bigcup_{Q\in\mathfrak L_m}Q
\subseteq
(R_{\mathfrak N}\setminus U)\cup
\bigcup_{i=1}^m\bigcup_{Q\in\mathfrak B_i}Q.
$$
Therefore, using \eqref{badcubemassjoiningcones},
$$
\mu\Big(R_{\mathfrak N}\setminus\bigcup_{Q\in\mathfrak L_m}Q\Big)
\leq
\eta'\mu(R_{\mathfrak N})
+
\sum_{i=1}^m C_0\eta'\mu(R_{\mathfrak N})
=
(1+mC_0)\eta'\mu(R_{\mathfrak N}).
$$
This proves \eqref{formulacopricopri}.

\textbf{Step IV. Conclusion.} Since each $Q \in \mathfrak{H}_m$ is $(j,\mathfrak{e},\vartheta^j,\sigma, \varepsilon,\delta,\zeta,\mathfrak m,A)$-invariant by construction in Step II, any sublayer $\mathfrak L_m \subseteq\mathfrak H_m$ inherits this invariance. Thus item 1 follows. Items 2, 3 and 4 follow from the construction of the families $\mathfrak H_m$ and from the definition of the subfamilies $\mathfrak L_m$. Item 5 follows from the construction of the cubes $\mathscr I_j(2\vartheta,\sigma/2,\varepsilon,\delta,\zeta,\mathfrak m,A)$, see \cref{rioridinocubi}. Item 6 follows from the definition of the bad cubes $\mathfrak B_m$: if $Q\in\mathfrak L_m$, then $Q\notin\mathfrak B_m$, and therefore for every $e\in X_c(\mathrm{span}(\mathfrak e),3\vartheta)$ we can choose a point $z(Q,e)\in AB_Q\cap\supp\mu$ satisfying both the non-invariance conclusion and
$$
\mu(B(z(Q,e),\zeta \diam Q/2)\cap U)\geq \frac{1}{2}\mu(B(z(Q,e),\zeta \diam Q/2)).
$$
Finally, by \eqref{formulacopricopri} and the choice of $\eta'$ we obtain
$$
\mu(R_{\mathfrak N}\setminus \bigcup\{Q:Q\in\mathfrak L_{\mathfrak N}\})\leq \eta\mu(R_{\mathfrak N}).
$$
Moreover, since $\mu(R_{\mathfrak N}\setminus U)\leq \eta'\mu(R_{\mathfrak N})\leq \eta\mu(R_{\mathfrak N})$, item 7 follows.
\end{proof}

In the next proposition, we assemble the tools developed in the previous subsections to establish the main geometric setup for the proof of the WALA conjecture. We start from the assumption that the family $\mathscr I_j$ of cubes with $j$ independent invariant directions is not Carleson. By applying the thinning procedure, we initially force a geometric length loss only along curves directed by a finite net of very thin cones. We then apply the joining-cones argument through a finite binary tree to glue these thin cones together, showing that the constructed sets are thin with respect to very wide transverse cones. This produces an arbitrarily deep hierarchy of nested open sets $\Omega_h$ and corresponding dyadic layers $\mathfrak L_h$ that shrink geometrically along all transverse directions.

\begin{proposizione}\label{propfinalejoiningcones}
Suppose $\mu$ is a $k$-AD-regular measure with regularity constant $D$. We fix $\mathfrak N\in\N$ and $\xi\in (0,1)$ and choose parameters $\sigma,\vartheta,\varepsilon,\delta,\zeta,\mathfrak m,\varsigma \in (0,1/2)$ and $A\geq 1$, $\qof\in \N$ satisfying the following properties
$$\vartheta\leq \frac14,\qquad
\zeta \leq \delta\leq 1/2^{24},\qquad
A\vartheta^n\geq 16n;$$
$$
\frac{4\zeta}{\varepsilon\sqrt{1-\sigma^2}}\leq \frac{\sigma\vartheta^6}{2^{20}\sqrt{n-1}},\qquad
16n\Big(\frac{\sigma}{\vartheta^n}
+n\frac{\varepsilon+\delta+\zeta+\mathfrak m}{A\vartheta^n}\Big)
\leq \frac{1}{2\cdot 3^k16^nD^2},
$$
$$
\qof\geq 4\log_2\Big(\frac{12800A}{\vartheta\varsigma\delta \zeta}
\Big\lceil n\log_2\Big(\frac{6400\sqrt{n-1}}{\sigma\vartheta}\Big)\Big\rceil
\Big\lceil\frac{20\log(\xi)\log(\frac{\sigma\vartheta^6}{2^{25}\sqrt{n-1}})}{\varsigma\vartheta}\Big\rceil\Big).
$$
Let $0\leq j<k$ and suppose that $\mathscr I(j,2\vartheta,\sigma/2,\varepsilon,\delta,\zeta,\mathfrak m,A)$, one of the families of cubes yielded by \cref{rioridinocubi}, is not Carleson. Then there exist
$$
\mathfrak e\in(\mathbb S^{n-1})^j,
\quad
R\in \mathscr I(j,2\vartheta,\sigma/2,\varepsilon,\delta,\zeta,\mathfrak m,A),\quad\text{families }\mathfrak L_0,\ldots,\mathfrak L_{\mathfrak N}
\subseteq \mathscr I(j,2\vartheta,\sigma/2,\varepsilon,\delta,\zeta,\mathfrak m,A),
$$
and nested open sets $\Omega_{\mathfrak s}=:\Omega_0\supseteq \Omega_1\supseteq \ldots \supseteq \Omega_{\mathfrak N}$, such that, if $\mathscr E$ is an orthonormal basis of $\mathrm{span}(\mathfrak e)^\perp$, the following properties hold.

\begin{enumerate}
\item for $m=1,\ldots,\mathfrak N$, the cubes $Q\in \mathfrak L_m$ are $(j,\mathfrak e,\vartheta^j,\sigma,\varepsilon,\delta,\zeta,\mathfrak m,A)$-invariant and $40AB_Q\subseteq \Omega_m$;

\item $\mathfrak L_0=\{R\}$ and no cube appears in more than one of the families $\mathfrak L_0,\ldots,\mathfrak L_{\mathfrak N}$;

\item for $m=0,\ldots,\mathfrak N$, the cubes in $\mathfrak L_m$ are pairwise disjoint. Moreover, whenever $Q,Q'\in \mathfrak L_m$, either
\begin{equation}
\ell(Q)=\ell(Q'),
\qquad\text{or, if }\ell(Q)<\ell(Q'),
\text{ then }\qquad
\frac{\ell(Q)}{\ell(Q')}\leq \frac{1}{2^\qof};
\label{separazionediscalev2asdfawrhs}
\end{equation}

\item For every $m=1,\ldots,\mathfrak N$, each cube $Q'\in \mathfrak L_m$ is contained in a unique strictly larger cube $Q\in \mathfrak L_{m-1}$ such that
$$
\diam(Q)\geq 2^\qof\diam(Q').
$$
In addition, if $Q\in \mathfrak L_{m-1}$ and $Q'\in \mathfrak L_m$ satisfy $10A B_{Q'}\cap 10A B_Q\neq \emptyset$, then
$$
\diam(Q)\geq 2^\qof\diam(Q').
$$
Finally, for every $m=1,\ldots,\mathfrak N$ and every $Q\in \mathfrak L_m$, if $Q'\in \mathfrak L_{m-1}$ is the ancestor of $Q$, then
$$
10A B_Q\cap \supp\mu\subseteq Q'.
$$

\item For every $m=1,\ldots,\mathfrak N$, every $Q\in \mathfrak L_m$, and every
$e\in X_c(\mathrm{span}(\mathfrak e),3\vartheta)$, the cube $Q$ is not
$(e,\sigma/2,\varepsilon,\delta,\zeta,\mathfrak m,A)$-invariant.

\item For every $h=0,\ldots,\mathfrak N-1$, every $e\in \mathscr E$ and every $C(e,\sigma)$-curve $\gamma$ we have
$$
\Haus^1(\gamma\cap \Omega_{h+1})
\leq \xi\Haus^1(\gamma\cap \Omega_h).
$$
In particular, for every $e\in \mathscr E$ and every $C(e,\sigma)$-curve $\gamma$ there holds
$\Haus^1(\gamma\cap \Omega_{\mathfrak N})
\leq \xi^{\mathfrak N}\Haus^1(\gamma\cap \Omega_{\mathfrak s})$;

\item for every $h=0,\ldots,\mathfrak N-1$ there are subfamilies $\mathfrak W_{h+1}\subseteq \Delta_\mu$ of $(j,\mathfrak e,\vartheta^j,\sigma,\varepsilon,\delta,\zeta,\mathfrak m,A)$-invariant cubes such that
$$
\Omega_{h+1}\subseteq \bigcup_{Q\in \mathfrak W_{h+1}}\tfrac{3}{2}AB_Q
\qquad\text{and}\qquad
20AB_Q\subseteq \Omega_h
\quad\text{for every }Q\in \mathfrak W_{h+1}.
$$
If $h=1,\ldots,\mathfrak N-1$, $Q'\in \mathfrak W_{h+1}$ and $Q\in \mathfrak W_h$ satisfy $10AB_Q\cap 10AB_{Q'}\neq \emptyset$, then
$$
\diam Q'\leq 2^{-\qof}\diam Q.
$$
Moreover, for every $Q\in \mathfrak W_{h+1}$, every $e\in \mathscr E$ and every $C(e,\sqrt{1-9\vartheta^2})$-curve $\gamma$ we have
$$
\Haus^1(\gamma \cap \Omega_{h+1}\cap 2AB_Q)
\leq
\xi\Haus^1(\gamma \cap \Omega_h \cap 2(1+\varsigma)AB_Q).
$$
In addition, the families $\mathfrak W_h$ are nested, in the sense that if $Q\in \mathfrak W_h$, $Q'\in \mathfrak W_{h+1}$ and $Q\cap Q'\neq \emptyset$, then $Q'\subseteq Q$. Further, every cube $Q\in \mathfrak W_h$ is contained in a cube of $\mathfrak L_h$ and every cube of $\mathfrak L_{h+1}$ is contained in a cube of $\mathfrak W_h$. Finally, for every $h=1,\ldots,\mathfrak N-1$, if $Q\in \mathfrak W_h$ and $Q'\in \mathfrak W_{h-1}$ is the ancestor of $Q$, then $10A B_Q\cap \supp\mu\subseteq Q'$.
\item
$\mu(R\setminus \Omega_{\mathfrak N})\leq 4\xi^{\mathfrak N}\mu(R)$.
\end{enumerate}
\end{proposizione}

\begin{proof}
We divide the proof in several steps.

\medskip

\textbf{Step I. The setup.} In this step, we construct the objects we need to prove our main result by directly applying the results we have obtained so far. 
First, thanks to \cref{joiningcones}, for every $\mathfrak M\in \N$ and $\eta\in (0,1)$ such that 
$$\eta\leq \frac{\xi^\mathfrak{N}}{\bigg\lceil-\frac{2^{25}\sqrt{n-1}}{\sigma\vartheta^{10}}\log\Big(\frac{\sigma\vartheta^6}{2^{20}\sqrt{n-1}}\Big)\bigg\rceil\bigg\lceil \frac{20}{\varsigma\vartheta}\log(\xi)\log\Big(\frac{\sigma\vartheta^6}{2^{25}\sqrt{n-1}}\Big)\bigg\rceil \bigg\lceil n\log_2\Big(\frac{6400\sqrt{n-1}}{\sigma\vartheta}\Big)\bigg\rceil \mathfrak N},$$
there exist $\mathfrak e\in (\mathbb{S}^{n-1})^j$, a cube $R_{\mathfrak M}\in \mathscr I(j,2\vartheta,\sigma/2,\varepsilon,\delta,\zeta,\mathfrak m,A)$, a Borel set $U\subseteq R_{\mathfrak M}$ and families
$$\widetilde{\mathfrak L}_0,\ldots, \widetilde{\mathfrak L}_\mathfrak M\subseteq \mathscr I(j,2\vartheta,\sigma/2,\varepsilon,\delta,\zeta,\mathfrak m,A),$$
such that items 1. to 7. of the thesis of \cref{joiningcones} are satisfied for such $\widetilde{\mathfrak L}_0,\ldots, \widetilde{\mathfrak L}_\mathfrak M$.

In what follows, we assume $\mathscr E$ is a fixed orthonormal basis of $\mathrm{span}(\mathfrak e)^\perp$. Let us observe that for every $e \in \mathscr E$ then $C(e,\sqrt{1-9\vartheta^2})\subseteq X_c(\mathrm{span}(\mathfrak e),3\vartheta)$. For every $e\in\mathscr E$, let us introduce a $\sigma/16$-dense set $\mathscr U_e$ in $C(e,\sqrt{1-9\vartheta^2})\cap\mathbb S^{n-1}$, and set $\mathscr U:=\bigcup_{e\in\mathscr E}\mathscr U_e$. Thanks to \cref{lem:pigeonhole-direction-small-cone}, we have 
$$\mathrm{Card}(\mathscr U)\leq\frac{nC_n}{\sigma^{n-1}},$$
where $C_n$ is a universal constant depending only on $n$.
Let us fix 
$$\mathfrak q_0:=\frac{\sigma\vartheta^4}{25600\sqrt{n-1}},$$
and let us notice that by \cref{alberodiconi} and our assumption on $\vartheta$ and $\sigma$, we have that for every $e\in \mathscr E$, there is a finite rooted binary tree $\mathcal T_{e}$ of depth $\mathfrak P(\mathcal T_{e})$ smaller than $\mathfrak P:=\lceil n\log_2(\frac{6400\sqrt{n-1}}{\sigma\vartheta})\rceil$ with the following properties.
Each node $\mathfrak n\in\mathcal T_{e}$ is associated to a convex cone $\Gamma_\mathfrak n$ satisfying the following properties
\begin{itemize}
    \item[($\alpha$)] the root $\emptyset$ satisfies $\Gamma_{\emptyset}=C(e,\sqrt{1-9\vartheta^2})$ and every non-terminal node $\mathfrak n$ has exactly two children $\mathfrak n^+$, $\mathfrak n^-$, with
$\Gamma_\mathfrak n=\Gamma_{\mathfrak n^+}\cup\Gamma_{\mathfrak n^-}$;
\item[($\beta$)] for every non-terminal node $\mathfrak n$ there are
$\mathfrak e_\mathfrak n\in \Gamma_\mathfrak n\cap S^{n-1}$, $u_\mathfrak n\in \mathfrak e_\mathfrak n^\perp$ and $\alpha_\mathfrak n>0$ such that
$$
\{x\in\Gamma_\mathfrak n:\langle x,u_\mathfrak n\rangle\geq -\alpha_\mathfrak n\langle x, \mathfrak e_\mathfrak n\rangle\}\subset\Gamma_{\mathfrak n^+},
\qquad
\{x\in\Gamma_\mathfrak n:\langle x,u_\mathfrak n\rangle\leq \alpha_\mathfrak n\langle x,\mathfrak e_\mathfrak n\rangle\}\subset\Gamma_{\mathfrak n^-},
$$
and $\Gamma_\mathfrak n$ is the union of these two sets;
\item[($\gamma$)] if we let
$$
K_\mathfrak n:=\max\left\{1,\sup_{x\in\Gamma_\mathfrak n\setminus\{0\}}\frac{|\langle x,u_\mathfrak n\rangle|}{\langle \mathfrak e_\mathfrak n,x\rangle}\right\},\qquad\text{then}\qquad
\frac{\alpha_\mathfrak n}{8K_\mathfrak n}\geq
\frac{\sigma\vartheta^4}{25600\sqrt{n-1}}=\mathfrak q_0;
$$
\item[($\delta$)] for every terminal node $\mathfrak n$ there exists $w\in \mathscr U$ such that
$\Gamma_\mathfrak n\subset C(w,\sigma/2)$.
\end{itemize}
Let us denote by $\mathfrak d$ the depth of the tree $\mathcal T_e$ and denote by $\mathcal T_e^i$ the family of nodes at depth $\mathfrak d-i$. Recall that $\mathfrak d<\mathfrak P$. 
Before proceeding with the proof we need to introduce some notations. We let
\begin{equation}
    \begin{split}
    N:=\bigg\lceil \frac{\log\big(\frac{\sigma^{n-1}\beth}{288D^22^kC_n\sqrt{\delta}}\big)}{\log\big(1-\frac{\zeta^k(1-\sqrt{\delta})}{40^{k+1}A^kD^2}\big)}\bigg\rceil,\,\,
K_1:=\bigg\lceil\frac{-1280\log\Big(\frac{\mathfrak q_0\vartheta^2}{32}\Big)}{\mathfrak q_0\vartheta^6}\bigg\rceil\mathfrak P,\,\,\text{and}\,\, K_2:=\bigg\lceil\frac{20\log(\xi)\log(\frac{\sigma\vartheta^6}{2^{25}\sqrt{n-1}})}{\varsigma\vartheta}\bigg\rceil,
\nonumber
    \end{split}
\end{equation}
where we choose
$$\beth \leq \frac{1}{4C_nK_1K_2\mathfrak N}\xi^{\mathfrak N}\sigma^{n-1},$$
and where $C_n$ is the constant yielded by \cref{lem:pigeonhole-direction-small-cone}. 
Let $K:=K_1K_2\mathfrak N$ and choose $\mathfrak M\in \N$ in such a way that
$$
\mathfrak M\geq 8(K+1)(2N+1)\log_2(160\delta^{-1}\zeta^{-1}A).
$$
We now check that \cref{productionopens} can be applied to the families of cubes $\widetilde{\mathfrak L}_0,\ldots, \widetilde{\mathfrak L}_\mathfrak M$ produced above, with parameters $\sigma/2$, $\beth_0:=\frac{\sigma^{n-1}\beth}{2nC_n}$ and $\xi':=1-\frac{4\zeta}{\varepsilon\sqrt{1-\sigma^2/4}}$. Then, by item 7. of \cref{joiningcones} and the choice of $\eta$, we have
$$ \mu(R_{\mathfrak M}\setminus U)\leq \eta\mu(R_{\mathfrak M}),\qquad
\mu(U\setminus\bigcup\{Q:Q\in\widetilde{\mathfrak L}_{\mathfrak M}\})\leq \eta\mu(R_{\mathfrak M}).
$$

Moreover, the family
$\bigcup_{m=0}^{\mathfrak M}\widetilde{\mathfrak L}_m$ is contained in
$\widetilde{\Omega}_{\mathfrak s}:=\mathrm{int}(20AB_{R_{\mathfrak M}})$ and satisfies
$$
10AB_Q\subseteq \widetilde{\Omega}_{\mathfrak s}
\qquad\text{for every }Q\in\bigcup_{m=0}^{\mathfrak M}\widetilde{\mathfrak L}_m.
$$
Indeed, this is immediate for $Q=R_{\mathfrak M}$. If $Q\neq R_{\mathfrak M}$, then
$Q\subseteq R_{\mathfrak M}$ and the scale separation along the ancestor chain gives
$\diam Q\leq 2^{-\qof}\diam R_{\mathfrak M}$. Since the choice of $\qof$ gives
$10A2^{-\qof}<1$, the inclusion follows. In addition, since $\mu(R_{\mathfrak M}\setminus U)\leq\eta\mu(R_{\mathfrak M})$ and $\eta\leq1/2$, we have $\mu(U)\geq\frac12\mu(R_{\mathfrak M})$. Hence
$$
\mu\Big(U\setminus\bigcup\{Q:Q\in\widetilde{\mathfrak L}_{\mathfrak M}\}\Big)
\leq
\frac{\beth_0}{4}\mu(U).
$$

Indeed, items 1. to 6. of the thesis of \cref{joiningcones} give the hypotheses (i)-(vi) in the assumptions of \cref{stepbasethinning}, because $\mathscr U\subseteq X_c(\mathrm{span}(\mathfrak e),3\vartheta)$. Therefore, thanks to the above observation and the choice of the parameters, we can apply \cref{productionopens} to  $\widetilde{\mathfrak L}_0,\ldots, \widetilde{\mathfrak L}_\mathfrak M$, and this shows that there are families of cubes 
$\mathfrak G_0,\ldots, \mathfrak G_\mathfrak{M}$ such that $\mathfrak G_h\subseteq \widetilde{\mathfrak L}_h$ for every $h=0,\ldots, \mathfrak M$ and nested open sets 
$$\widetilde{\Omega}_\mathfrak s=\widetilde{\Omega}_0\supseteq \widetilde{\Omega}_1\supseteq \ldots \supseteq\widetilde{\Omega}_K,$$
such that 
\begin{enumerate}
    \item[(a)] for every $h=1,\ldots, \mathfrak M$ every cube $Q\in \mathfrak G_h$ is $(j,\mathfrak{e},\vartheta^j,\sigma, \varepsilon,\delta,\zeta,\mathfrak m,A)$-invariant and satisfies
$40AB_Q\subseteq \widetilde{\Omega}_{\min\{\lfloor h/(2N+1)\rfloor,K\}}$.
\item[(b)] the layers $\mathfrak G_h$ still satisfy items (i) to (v) of the assumptions of \cref{stepbasethinning}.

\item[(c)] for every $h=0,\ldots, K-1$, every $w\in \mathscr U$ and every $C(w,\sigma/2)$-curve $\gamma$ we have
\begin{equation}
\Haus^1(\gamma\cap \widetilde{\Omega}_{h+1})
\leq\Big(\frac{4\zeta}{\varepsilon\sqrt{1-\sigma^2}}\Big)\Haus^1(\gamma\cap \widetilde{\Omega}_h),
\nonumber
\end{equation}
\item[(d)] for every $h=0,\ldots, K-1$ we have
$\widetilde{\Omega}_{h+1}\subseteq \bigcup_{Q\in \mathfrak{G}_{h(2N+1)+1}}\tfrac{3}{2}AB_Q$,
and for every $Q\in \mathfrak{G}_{h(2N+1)+1}$, for every $w\in \mathscr U$ and every $C(w,\sigma/2)$-curve $\gamma$ we have that 
$$\Haus^1(\gamma \cap \widetilde{\Omega}_{h+1}\cap 2AB_Q)\leq
\Big(\frac{4\zeta}{\varepsilon\sqrt{1-\sigma^2}}\Big)\Haus^1(\gamma \cap \widetilde{\Omega}_h \cap 2(1+2^{-\qof})AB_Q).$$
\item[(e)]
$\mu(U\setminus \widetilde{\Omega}_K)\leq 2\mathrm{Card}(\mathscr{U})K\beth_0\mu(U)$;
\item[(f)] for every $w\in \mathscr U$ and every $C(w,\sigma/2)$-curve we have
$$\Haus^1(\gamma \cap \widetilde{\Omega}_K)\leq \Big(\frac{4\zeta}{\varepsilon\sqrt{1-\sigma^2}}\Big)^K\Haus^1(\gamma\cap \widetilde{\Omega}_{\mathfrak s});$$
\end{enumerate}
Notice that (f) is a direct consequence of an iteration of item (c). Let us record the following consequence of (a) and (d). If
$Q\in\mathfrak G_{h(2N+1)+1}$, then
\begin{equation}
40AB_Q
\subseteq
\widetilde{\Omega}_{\min\{\lfloor(h(2N+1)+1)/(2N+1)\rfloor,K\}}
\subseteq
\widetilde{\Omega}_h.
\label{eq:forty-ball-contained-in-omega-h}
\end{equation}
Consequently, if
$B=\tfrac32AB_Q$, then
\begin{equation}
20B=30AB_Q\subseteq 40AB_Q\subseteq\widetilde{\Omega}_h.
\label{eq:twenty-B-contained-in-omega-h}
\end{equation}
In addition, we also have that if $Q\in \mathfrak G_{h(2N+1)+1}$ and $Q'\in \mathfrak G_{(h+1)(2N+1)+1}$ are such that $10AB_Q\cap 10AB_{Q'}\neq \emptyset$, there holds 
\begin{equation}
    \diam(Q')\leq 2^{-\qof}\diam(Q).
    \label{bounddiametriperW}
\end{equation}

\medskip

\textbf{Step II. Putting cones together.} Thanks to item ($\delta$), we know that 
for every $\mathfrak n\in \mathcal T_e^0$, there exists $w\in \mathscr U$ such that $\Gamma_\mathfrak{n}\subseteq C(w,\sigma/2)$ and in particular 
\begin{itemize}
    \item[($\mathrm{c}^0$)] for every $h=0,\ldots, K-1$, every $\mathfrak n\in \mathcal{T}_e^0$ and every $\Gamma_{\mathfrak n}$-curve $\gamma$ we have
\begin{equation}
\Haus^1(\gamma\cap \widetilde{\Omega}_{h+1})
\leq\Big(\frac{4\zeta}{\varepsilon\sqrt{1-\sigma^2}}\Big)\Haus^1(\gamma\cap \widetilde{\Omega}_h),
\nonumber
\end{equation}
\item[($\mathrm{d}^0$)] for every $h=0,\ldots, K-1$ we have
$$\widetilde{\Omega}_{h+1}\subseteq \bigcup_{Q\in \mathfrak G_{h(2N+1)+1}}\tfrac{3}{2}AB_Q$$
and for every $Q\in \mathfrak G_{h(2N+1)+1}$, every $\mathfrak n\in \mathcal T_e^0$ and every $\Gamma_\mathfrak{n}$-curve $\gamma$ we have that 
$$\Haus^1(\gamma \cap \widetilde{\Omega}_{h+1}\cap 2AB_Q)\leq
\Big(\frac{4\zeta}{\varepsilon\sqrt{1-\sigma^2}}\Big)\Haus^1(\gamma \cap \widetilde{\Omega}_h \cap 2(1+2^{-\qof})AB_Q).$$
\end{itemize}

We now inductively prove that for all the cones associated to the nodes of
$\mathcal T_e$, we can prove a similar statement to
($\mathrm{c}^0$) and ($\mathrm{d}^0$). Set
$$
\mathfrak t:=20\cdot 2^{-\qof}
\qquad\text{and}\qquad
\shin:=\Big\lceil
-\frac{1280}{\mathfrak q_0\vartheta^6}
\log\Big(\frac{\mathfrak q_0\vartheta^2}{32}\Big)
\Big\rceil.
$$
We prove by induction on $i=1,\ldots,\mathfrak d$ the following claim. If
$\mathfrak n\in\mathcal T_e^i$, $0\leq h_1<h_2\leq K$, and
$h_2-h_1\geq i\shin$, then, for every $Q\in\mathfrak G_{h_1(2N+1)+1}$ and every
$\Gamma_{\mathfrak n}$-curve $\gamma$, we have
\begin{equation}
\Haus^1(\gamma\cap \widetilde{\Omega}_{h_2}\cap 2AB_Q)
\leq
\Big(\frac{\mathfrak q_0\vartheta^2}{32}\Big)^i
\Haus^1(
\gamma\cap\widetilde{\Omega}_{h_1}
\cap 2(1+\mathfrak t)^{h_2-h_1}AB_Q
),
\label{eq:induction-local-cones}
\end{equation}
and
\begin{equation}
\Haus^1(\gamma\cap\widetilde{\Omega}_{h_2})
\leq
\Big(\frac{\mathfrak q_0\vartheta^2}{32}\Big)^i
\Haus^1(\gamma\cap\widetilde{\Omega}_{h_1}).
\label{eq:induction-global-cones}
\end{equation}

\emph{Base step.} Let $\mathfrak n\in\mathcal T_e^1$. Then its children
$\mathfrak n^+$ and $\mathfrak n^-$ belong to $\mathcal T_e^0$. We apply
\cref{prop:two-leaves-iterated-substitution-additive} with
$C:=\Gamma_{\mathfrak n}$, $C^+:=\Gamma_{\mathfrak n^+}$, $C^-:=\Gamma_{\mathfrak n^-}$ and $\Lambda:=\frac43$ and
$$ \mathscr B^h:=\Big\{\tfrac32AB_Q:Q\in\mathfrak G_{h(2N+1)+1}\Big\}.$$
We verify the hypotheses with $\mathfrak q:=\mathfrak q_0$ and $\theta:=\sqrt{1-9\vartheta^2}$. The splitting hypotheses for $C,C^+,C^-$ follow from items $(\alpha)$,
$(\beta)$ and $(\gamma)$ in the construction of the tree. Moreover, by
construction, $\Gamma_{\mathfrak n}\subseteq C(e,\sqrt{1-9\vartheta^2})=C(e,\theta)$. 
By the choice of $\qof$ and the definition of $\mathfrak t$, we have
$$
\mathfrak t=20\cdot2^{-\qof}
\leq
\frac{\sqrt{1-\theta^2}}{16}
\qquad\text{and}\qquad
(1+\mathfrak t)^{K_1K_2}\leq \frac{10}{9}.
$$
The covering hypothesis follows from item (d) of Step I, i.e. 
$$\widetilde{\Omega}_{h+1}\subseteq
\bigcup_{Q\in\mathfrak G_{h(2N+1)+1}}\frac32AB_Q=\bigcup_{B\in\mathscr B^h}B.$$
If $B=\frac32AB_Q\in\mathscr B^h$, then
\eqref{eq:twenty-B-contained-in-omega-h} gives $20B\subseteq\widetilde{\Omega}_h$ and notice that if  $B=\frac32AB_Q\in\mathscr B^h$, $B'=\frac32AB_{Q'}\in\mathscr B^{h+1}$ and $4B\cap B'\neq\emptyset$ then  $10AB_Q\cap10AB_{Q'}\neq\emptyset$. Therefore
\eqref{bounddiametriperW} gives
$$
\diam B'
\leq
2^{-\qof}\diam B
=
\frac{\mathfrak t}{20}\diam B.
$$

It remains to check the estimates for the two children. Since
$\mathfrak n^+,\mathfrak n^-\in\mathcal T_e^0$, item $(\delta)$ gives, for
each choice of the sign, a vector $w^\pm\in\mathscr U$ such that
$\Gamma_{\mathfrak n^\pm}\subseteq C(w^\pm,\sigma/2)$. Thus every $\Gamma_{\mathfrak n^\pm}$-curve is a $C(w^\pm,\sigma/2)$-curve.
Hence item (d) of Step I implies that, for every
$B=\frac32AB_Q\in\mathscr B^h$, we have
$$
\Haus^1(\gamma\cap\widetilde{\Omega}_{h+1}\cap \tfrac43B)
\leq
\frac{4\zeta}{\varepsilon\sqrt{1-\sigma^2}}
\Haus^1(\gamma\cap\widetilde{\Omega}_h\cap \tfrac43(1+\mathfrak t)B),
$$
whenever $\gamma$ is a $\Gamma_{\mathfrak n^+}$-curve or a
$\Gamma_{\mathfrak n^-}$-curve. Moreover,
$$
\frac{4\zeta}{\varepsilon\sqrt{1-\sigma^2}}
\leq
\frac{\sigma\vartheta^6}{2^{20}\sqrt{n-1}}
\leq
\frac{\mathfrak q_0(1-\theta^2)}{32}.
$$
The corresponding global estimate for the two children follows from item (c)
of Step I. Thus all the hypotheses of
\cref{prop:two-leaves-iterated-substitution-additive} are satisfied.
Therefore, for every $h_2>h_1$ with $h_2-h_1\leq K_1K_2$, every
$Q\in\mathfrak G_{h_1(2N+1)+1}$ and every $\Gamma_{\mathfrak n}$-curve
$\gamma$, we have
$$\Haus^1(\gamma\cap\widetilde{\Omega}_{h_2}\cap 2AB_Q)\leq\left(1-\frac{\mathfrak q_0(1-\theta^2)^2}{5120}\right)^{h_2-h_1}\Haus^1(\gamma\cap\widetilde{\Omega}_{h_1}\cap 2(1+\mathfrak t)^{h_2-h_1}AB_Q),$$
and
$$\Haus^1(\gamma\cap\widetilde{\Omega}_{h_2})\leq\left(1-\frac{\mathfrak q_0(1-\theta^2)^2}{5120}\right)^{h_2-h_1}\Haus^1(\gamma\cap\widetilde{\Omega}_{h_1}).$$
By the choice of $\shin$, if $h_2-h_1\geq\shin$, then
$$
\left(1-\frac{\mathfrak q_0(1-\theta^2)^2}{5120}\right)^{h_2-h_1}
\leq
\frac{\mathfrak q_0\vartheta^2}{32}.
$$
This proves the claim for $i=1$.

\emph{Inductive step.} Assume that the claim has been proved for some
$1\leq i<\mathfrak d$. Let $\mathfrak n\in\mathcal T_e^{i+1}$ and let
$\mathfrak n^+$ and $\mathfrak n^-$ be its children. Then $\mathfrak n^+,\mathfrak n^-\in\mathcal T_e^i$.  Let $h_2-h_1\geq (i+1)\shin$. Since the sets $\widetilde{\Omega}_h$ are nested, it is enough to prove the claim with
$h_2=h_1+(i+1)\shin$.
The general case follows by replacing $h_2$ by $h_1+(i+1)\shin$ and enlarging
the ball on the right hand side. If we set $h_*:=h_1+\shin$,
then $h_2-h_*=i\shin$. Applying the inductive hypothesis to the children
$\mathfrak n^+$ and $\mathfrak n^-$ on the interval $h_2$ and $h_*$ gives that,
for every $\Gamma_{\mathfrak n^\pm}$-curve $\gamma$, every
$S\in\mathfrak G_{h_*(2N+1)+1}$, and each choice of the sign,
$$\Haus^1(\gamma\cap\widetilde{\Omega}_{h_2}\cap 2AB_S)\leq\Big(\frac{\mathfrak q_0\vartheta^2}{32}\Big)^i\Haus^1(\gamma\cap\widetilde{\Omega}_{h_*}\cap 2(1+\mathfrak t)^{i\shin}AB_S),$$
and
$$\Haus^1(\gamma\cap\widetilde{\Omega}_{h_2})\leq\Big(\frac{\mathfrak q_0\vartheta^2}{32}\Big)^i\Haus^1(\gamma\cap\widetilde{\Omega}_{h_*}).$$
We now apply \cref{prop:two-leaves-iterated-substitution-additive} to the sets $\widetilde{\Omega}_{h_1}\supseteq\widetilde{\Omega}_{h_1+1}\supseteq\ldots\supseteq\widetilde{\Omega}_{h_*}$, with
$C:=\Gamma_{\mathfrak n}$, $C^+:=\Gamma_{\mathfrak n^+}$, $C^-:=\Gamma_{\mathfrak n^-}$ with $\Lambda:=\frac43(1+\mathfrak t)^{i\shin}$, $\mathfrak q:=\mathfrak q_0$ and $\theta:=\sqrt{1-9\vartheta^2}$
and with
$$\mathscr B^h:=\left\{\tfrac32AB_Q:Q\in\mathfrak G_{h(2N+1)+1}\right\}.$$
The splitting hypotheses follow again from items $(\alpha)$, $(\beta)$ and
$(\gamma)$ in the construction of the tree. The covering, containment and scale
hypotheses are the same as in the base step, by item (d) of Step I,
\eqref{eq:twenty-B-contained-in-omega-h} and \eqref{bounddiametriperW}.
The value of $\Lambda$ is admissible, since
$$(1+\mathfrak t)^{i\shin}\leq(1+\mathfrak t)^{K_1K_2}\leq1+\varsigma \qquad \text{and thus}\qquad \Lambda\leq\frac43(1+\varsigma)<2.$$
The estimates just obtained for the two children give exactly the child
estimates required in the application of
\cref{prop:two-leaves-iterated-substitution-additive}, with this value of
$\Lambda$. Hence, for every
$Q\in\mathfrak G_{h_1(2N+1)+1}$ and every $\Gamma_{\mathfrak n}$-curve
$\gamma$, we obtain
$$
\Haus^1\Big(
\gamma\cap\widetilde{\Omega}_{h_*}
\cap 2(1+\mathfrak t)^{i\shin}AB_Q
\Big)
\leq
\frac{\mathfrak q_0\vartheta^2}{32}
\Haus^1\Big(
\gamma\cap\widetilde{\Omega}_{h_1}
\cap 2(1+\mathfrak t)^{(i+1)\shin}AB_Q
\Big).
$$
Combining this estimate with the estimate for the children gives
$$
\Haus^1(\gamma\cap\widetilde{\Omega}_{h_2}\cap 2AB_Q)
\leq
\Big(\frac{\mathfrak q_0\vartheta^2}{32}\Big)^{i+1}
\Haus^1\Big(
\gamma\cap\widetilde{\Omega}_{h_1}
\cap 2(1+\mathfrak t)^{(i+1)\shin}AB_Q
\Big).
$$
Since here $h_2-h_1=(i+1)\shin$, this is exactly
\eqref{eq:induction-local-cones}. The global estimate follows in the same way:
first use the global inductive estimate for the two children from $h_*$ to
$h_2$, and then the global conclusion of
\cref{prop:two-leaves-iterated-substitution-additive} from $h_1$ to $h_*$.
This proves the inductive step.

We now apply the claim to the root cone
$$
\Gamma_\emptyset=C(e,\sqrt{1-9\vartheta^2}).
$$
Since the depth of $\mathcal T_e$ is at most $\mathfrak P$, and
$K_1=\mathfrak P\shin$, the claim gives the desired estimates on every interval
of length $K_1$. Iterating these estimates $K_2$ times, and using the choices of
$K_2$ and $\qof$, we obtain that if $h_2=h_1+K_1K_2$, then for every
$e\in\mathscr E$ and every $C(e,\sqrt{1-9\vartheta^2})$-curve $\gamma$,
\begin{equation}
\Haus^1(\gamma\cap \widetilde{\Omega}_{h_2})
\leq
\xi\Haus^1(\gamma\cap \widetilde{\Omega}_{h_1}).
\label{eq:stimawidthapertix}
\end{equation}
Similarly, for every $Q\in\mathfrak G_{h_1(2N+1)+1}$ we get
\begin{equation}
\begin{split}
\Haus^1(\gamma\cap\widetilde{\Omega}_{h_2}\cap 2AB_Q)
&\leq
\xi\,
\Haus^1\Big(
\gamma\cap\widetilde{\Omega}_{h_1}
\cap 2(1+\mathfrak t)^{K_1K_2}AB_Q
\Big)\\
&\leq
\xi\,
\Haus^1\Big(
\gamma\cap\widetilde{\Omega}_{h_1}
\cap 2(1+\varsigma)AB_Q
\Big).
\end{split}
\label{stimaitem8}
\end{equation}
With this we are finally ready to move to the final step of the proof.

\textbf{Step III. Conclusion.}
For every $h=0,\ldots,\mathfrak N$ we define
$$
\Omega_h:=\widetilde{\Omega}_{K_1K_2h},
\qquad
\mathfrak L_h:=\mathfrak G_{K_1K_2h(2N+1)},
$$
and for $h=0,\ldots,\mathfrak N-1$, let
$$
\mathfrak W_{h+1}:=\mathfrak G_{K_1K_2h(2N+1)+1}.
$$
Finally, set $R:=R_{\mathfrak M}$.

Item 1 follows from item (a) of Step I. Items 2, 3, 4 and 5 follow directly from the corresponding properties of the families $\mathfrak G_h$, iterating the ancestor relation between consecutive layers. Item 6 follows from \eqref{eq:stimawidthapertix}, because
$\sigma<1/2\leq\sqrt{1-9\vartheta^2}$ and hence
$C(e,\sigma)\subseteq C(e,\sqrt{1-9\vartheta^2})$.

Item 7 gathers the properties of the subfamilies $\mathfrak W_h$. By item (d) of Step I,
$$
\Omega_{h+1}
=
\widetilde{\Omega}_{K_1K_2(h+1)}
\subseteq
\widetilde{\Omega}_{K_1K_2h+1}
\subseteq
\bigcup_{Q\in\mathfrak W_h}\frac32AB_Q.
$$
Moreover, by \eqref{eq:forty-ball-contained-in-omega-h}, for every $Q\in\mathfrak W_{h+1}$ (which is $\mathfrak G_{K_1K_2h(2N+1)+1}$), we have
$$
20AB_Q\subseteq 40AB_Q\subseteq \widetilde{\Omega}_{K_1K_2h} = \Omega_h.
$$
Similarly, \eqref{bounddiametriperW} implies that if $Q\in \mathfrak W_{h}$ and $Q'\in \mathfrak W_{h+1}$ satisfy $10AB_Q\cap 10AB_{Q'}\neq \emptyset$, then iterating the scale separation properly gives
\begin{equation}
    \diam(Q')\leq 2^{-\qof}\diam(Q).
    \nonumber
\end{equation}
The local estimate in item 7 follows immediately from \eqref{stimaitem8}. Let us pass to the hierarchy properties. It is immediate to see that if $Q\in\mathfrak W_h$, $Q'\in\mathfrak W_{h+1}$ and $Q\cap Q'\neq\emptyset$, then, by the scale separation and dyadic nesting, $Q'\subseteq Q$. On the other hand, 
since $\mathfrak W_{h+1} = \mathfrak G_{K_1K_2h(2N+1)+1}$, every cube $Q\in \mathfrak W_{h}$ must be contained in its ancestor in $\mathfrak G_{K_1K_2h(2N+1)}=\mathfrak L_h$
by the hierarchy properties of the nested cubes. Similarly, one sees that every cube of $\mathfrak L_{h+1}$ must be contained in an ancestor cube in $\mathfrak W_{h+1}$. Finally, the fact that for every $h=1,\ldots, \mathfrak N-1$, if $Q\in \mathfrak W_h$ and $Q'\in \mathfrak W_{h-1}$ is the ancestor of $Q$, we have $10A B_Q\cap \supp\mu\subseteq Q'$, follows directly from the fact that the layers $\mathfrak G_h$ have this property and the above inclusions.

It remains to prove item 8. Since $\Omega_{\mathfrak N}=\widetilde{\Omega}_K$ (where $K = K_1K_2\mathfrak N$), we have
\begin{equation}
    \begin{split}
\mu(R\setminus\Omega_{\mathfrak N})
&=
\mu(R_{\mathfrak M}\setminus\widetilde{\Omega}_K)
\leq
\mu\Big(R_{\mathfrak M}\setminus \bigcup\{Q:Q\in\widetilde{\mathfrak L}_{\mathfrak M}\}\Big)
+
\mu\Big(\bigcup\{Q:Q\in\widetilde{\mathfrak L}_{\mathfrak M}\} \setminus \widetilde{\Omega}_K\Big)\\
&\leq
\eta\mu(R_{\mathfrak M})
+
2\mathrm{Card}(\mathscr U)K\beth_0 \mu(R_{\mathfrak M})\leq
\eta\mu(R) + K\beth\mu(R)\\
&\leq
\eta\mu(R)+\frac12\xi^{\mathfrak N}\mu(R)
\leq
4\xi^{\mathfrak N}\mu(R).
\nonumber
\end{split}
\end{equation}
This concludes the proof.
\end{proof}

\subsection{GWALA implies UR}

The goal of this subsection is to complete the proof that GWALA implies uniform
rectifiability. By the qualitative part of the argument, GWALA already implies
that the AD-regular dimension is an integer \(k\) and that the measure is
\(k\)-rectifiable. It remains to prove the quantitative conclusion.

We prove the fixed-threshold OUWGL condition. The argument is by contradiction. If the OUWGL bad cubes were not Carleson, the non-Carleson selection lemmas
would produce long trees of dyadic cubes with big scale separation.

Cubes which are invariant in $k$ quantitatively independent directions are good for OUWGL: by \cref{propo:piano-da-invarianza}, through each relevant point there is a $k$-plane whose points remain close to $\supp\mu$ at the
corresponding scale. Hence these cubes cannot belong to the fixed-threshold
OUWGL bad family, provided the parameters are chosen sufficiently small.

The remaining cubes fail invariance in some direction. For those cubes, the
width-loss construction from the previous section produces open sets which are
large in $\mu$-measure but have small width with respect to curves going very wide cone. These open sets are then used to build Lipschitz width functions used to cancel the gradient along the bad directions and allow us to perturb Lipschitz functions away from the the functions in $\mathcal F$ functions. This contradicts GWALA.

We first record two elementary technical tools. The first lemma is a boundary-preserving mollification lemma: it smooths a Lipschitz function inside a ball, keeps the boundary values fixed, and preserves the relevant directional derivative bounds up to a small error. The second lemma is a controlled Vitali selection lemma for two generations of balls.

\begin{lemma}\label{mollificazioneapproxalbordo}
Let $r>0$, $V\in \mathrm{Gr}(j,n)$, and let $e_1,\ldots,e_{n-j}$ be an orthonormal basis of $V^\perp$. Suppose $f$ is a Lipschitz function on the ball $B(0,r)$ such that 
$$|d_V f|\le 1\qquad\text{and}\qquad|\partial_{e_\iota}f|\le 1\qquad\text{for $\Leb^n$-almost every $x\in B(0,r)$ and every $\iota=1,\ldots,n-j$.}$$
Then, for every $\alpha_1,\alpha_2\in(0,1)$, there exist a smooth function $u$ in the interior of $B(0,r)$ and continuous up to the boundary such that $$\lVert f-u\rVert_\infty\leq \frac{\alpha_1\alpha_2}{18}r\qquad\text{and}\qquad u=f\qquad\text{on }\partial B(0,r).$$
Furthermore, $|d_Vu|\le 1+\alpha_1$, $|\partial_{e_\iota}u|\leq 1+\alpha_1$ in $B(0,r)$ for every $\iota=1,\ldots,n-j$, and 
$$u=f*\rho_\eta\text{ on $B(0,(1-\tfrac{\alpha_2}{2})r)$, where }\eta:=\frac{\alpha_1\alpha_2}{18\sqrt{n-j+1}}r.$$
\end{lemma}

\begin{proof}
Let $B:=B(0,r)$ and 
$$\rho(z)=c_n\exp\Big(-\frac{1}{1-|z|^2}\Big)\mathbb{1}_{\{|z|<1\}},$$ where $c_n^{-1}=\int_{B_1(0)}e^{-\frac{1}{1-|z|^2}}\,dz$, and put $\rho_\eta(z)=\eta^{-n}\rho(z/\eta)$.

We now introduce a function that will yield the mollifying radius inside the ball $B$. This radius will need to go to $0$ as we approach the boundary. 
Let $\psi(t)=e^{-1/t}$ for $t>0$ and $\psi(t)=0$ for $t\le 0$, set $\Theta(t)=\psi(t)/(\psi(t)+\psi(1-t))$, and define
$$
h(t)=\begin{cases}
t, & t\le 1/2,\\
t+(1-t)\Theta(2t-1), & 1/2<t<1,\\
1, & t\ge 1.
\end{cases}
$$
Then $h\in C^\infty(\mathbb R)$ and with few elementary algebraic computations it is possible to show that $\|h'\|_{L^\infty}\le 9$. Finally, for $x\in B$, define 
$$r(x):=\eta h\Big(\frac{2r-2|x|}{\alpha_2r}\Big).$$
Let us observe that $x-r(x)z\in B$ for every $z\in B(0,1)$. Indeed since $\Theta(t)\leq 1$ we have $h(t)\leq 2t$ and thus 
$$r(x)\leq 4\frac{\eta}{\alpha_2 r}(r-|x|)<r-|x|,$$
thanks to the choice of $\eta$, and this proves that $x+r(x)B(0,1)\subseteq B$.
We are ready to define $u$. For every $x\in B$, let 
 $$u(x):=\int_{\mathbb R^n}f(x-r(x)z)\rho(z)\,dz.$$
Notice that for every $x\in B$ we have 
$$   |f(x)-u(x)|\leq \sqrt{n-j+1} r(x),$$
and thus $u=f$ on $\partial B$. Furthermore, if $|x|\le (1-\alpha_2/2)r$, then $2(r-|x|)/\alpha_2\ge r$, hence $r(x)=\eta$ and therefore $u(x)=\int_{\mathbb R^n}f(x-\eta z)\rho(z)\,dz=f*\rho_\eta(x)$. Thus $u=f*\rho_\eta$ on $B(0,(1-\alpha_2/2)r)$.

Let us observe that 
$$\lVert f-u\rVert_\infty\leq \sqrt{n-j+1}\lVert r(\cdot)\rVert_\infty\leq \frac{\alpha_1\alpha_2}{18}r,$$
where the last inequality comes from the choice of $\eta$. 

It remains to estimate the derivatives. Differentiating $u$ along a direction $a\in \mathbb S^n$ that is either a unit vector in $V$ or one of the vectors $e_\iota$, gives by the dominated convergence theorem 
$$\partial_a u(x)=\int_{\mathbb R^n}Df(x-r(x)z)[a-(\partial_a r(x))z]\rho(z)\,dz.$$
In particular, this implies that 
$$\lvert \partial_a u(x)\rvert\leq 1+\sqrt{n-j+1}\lVert \nabla r\rVert_\infty\int |z|\rho(z)dz\leq  1+\sqrt{n-j+1}\lVert \nabla r\rVert_\infty,$$
where in the first inequality above we used the fact that $\lVert \nabla f\rVert_\infty\leq\sqrt{n-j+1}$. 
However
$$\lVert \nabla r\rVert_\infty\leq\frac{2\eta}{\alpha_2r}\lVert h'\rVert_\infty\leq \frac{18\eta}{\alpha_2r}\leq \frac{\alpha_1}{\sqrt{n-j+1}}.$$
Thus $|\partial_a u|\le 1+\alpha_1$ for every unit $a\in V$ and for each $a=e_\iota$. Taking the supremum over unit vectors $a\in V$ gives $|d_Vu|\le 1+\alpha_1$. This concludes the proof. 
\end{proof}

\begin{lemma}\label{lemmacostruzionefammigliepallecontrollate}
Let $\theta\in (0,1/2)$ and $\rho\in [0,1)$ be such that $\rho+2\theta<1$.
Suppose that we are given two finite families of balls
$$
\mathscr B_i:=\{B_j^i:j=1,\ldots,N_i\},\qquad i=1,2,
$$
satisfying the following property: if $B_{j_1}^1\cap B_{j_2}^2\neq \emptyset$, then
$$
\diam B_{j_2}^2\leq \theta\diam B_{j_1}^1.
$$
Then there exists a subfamily $\mathscr C_2\subseteq \mathscr B_2$ satisfying the following properties:
\begin{itemize}
    \item the balls in $\mathscr C_2$ are pairwise disjoint and, for every $B\in \mathscr C_2$, there exists $B'\in \mathscr B_1$ such that
    $$
    B\subseteq (1-\rho)B';
    $$
    \item for every $B'\in \mathscr B_1$ we have
    $$
    (1-\rho-2\theta)B'\cap \bigcup_{B\in \mathscr B_2}B
    \subseteq
    \bigcup_{B\in \mathscr C_2}5B.
    $$
\end{itemize}
\end{lemma}

\begin{proof}
Let
$$
\mathscr F
:=
\{B\in\mathscr B_2:\text{ there exists }B'\in\mathscr B_1
\text{ such that }B\subseteq (1-\rho)B'\}.
$$
If $\mathscr F$ is non-empy, the finite Vitali covering lemma, there exists a pairwise disjoint subfamily
$\mathscr C_2\subseteq\mathscr F$ such that
$$
\bigcup_{B\in\mathscr F}B
\subseteq
\bigcup_{B\in\mathscr C_2}5B.
$$
The first property follows immediately from the definition of $\mathscr F$.

We prove the second property. Fix $B'\in\mathscr B_1$ and take
$$
y\in (1-\rho-2\theta)B'\cap \bigcup_{B\in\mathscr B_2}B.
$$
Then there exists $B\in\mathscr B_2$ such that $y\in B$. Since
$y\in B'\cap B$, the assumption gives
$$
\diam B\leq \theta\diam B'.
$$
Write $B'=B(x',R)$ and $B=B(x,r)$. Then $r\leq \theta R$. Since
$y\in (1-\rho-2\theta)B'$, for every $z\in B$ we have
$$
|z-x'|\leq |z-x|+|x-y|+|y-x'|
\leq r+r+(1-\rho-2\theta)R
\leq (1-\rho)R.
$$
Thus $B\subseteq (1-\rho)B'$, and hence $B\in\mathscr F$. Therefore
$$
y\in B\subseteq \bigcup_{C\in\mathscr C_2}5C.
$$
Notice that  $\mathscr F=\emptyset$, we take $\mathscr C_2=\emptyset$.
Then the first property is void. Moreover, if
$$
y\in (1-\rho-2\theta)B'\cap \bigcup_{B\in\mathscr B_2}B
$$
for some $B'\in\mathscr B_1$, then the argument below shows that the ball
$B\in\mathscr B_2$ containing $y$ belongs to $\mathscr F$, a contradiction.
Thus the second property is also void.
This proves the desired inclusion.
\end{proof}

The proof of the next theorem is the quantitative counterpart of the perturbation argument used in the qualitative part. In the qualitative proof, the failure of rectifiability gives tangent measures which split along the decomposability bundle and which hava a lot of mass spread in the transverse direction. This allows one first to cancel the derivatives in the transverse directions by using width functions, and then to add a perturbation, constant along the tangent directions, whose values cannot be approximated by functions in the finite-dimensional space $\mathcal F$. 

Here the same mechanism is implemented at the level of cubes. The role of the tangent splitting is replaced by the discrete product structure obtained from invariant cubes. If one of the families
$\mathscr I(j,2\vartheta,\sigma/2,\varepsilon,\delta,\zeta,\mathfrak m,A)$ 
were not Carleson, the joining-cones construction would give a long tower of cube families and nested open sets having small width with respect to a very wide cone. This is the quantitative substitute for the transverse thinness used in the qualitative proof.

The proof then follows the same two perturbative steps. First, on a selected layer, the function is slightly smoothed and its transverse derivatives are canceled by subtracting suitable width functions. This produces a new admissible function which is almost affine along the invariant directions and has very small oscillation in the transverse directions at the next layer. Second, using the discrete product structure inside the smaller cubes, one adds a perturbation which is constant along the invariant directions and whose prescribed values evade the finite-dimensional space $\mathcal F$. This forces a definite lower bound for the coefficients $\Omega_{1,\mu,\mathcal F}$ on a positive portion of the next scale.

The construction is iterated along the tower of layers. At each stage one chooses one of finitely many perturbations and keeps the lower bounds already produced at the previous scales, losing only a summable error because of the large separation between consecutive layers. After the last step, one obtains a Lipschitz function for which the coefficients $\Omega_{1,\mu,\mathcal F}$ are large on many pairwise separated families of cubes inside a fixed cube $R$. These cubes therefore form a packing of size comparable to the length of the tower. Choosing the tower sufficiently long contradicts the Carleson estimate coming from the GWALA condition. Hence the original non-Carleson assumption is impossible.

\begin{teorema}\label{teoremasemifinale}
Suppose $\mu$ is a $k$-AD-regular measure with regularity constant $D$.
Let $\mathcal F$ be a vector subspace of $L^1_{\mathrm{loc}}(\R^n)$ of
dimension $d$ containing the affine functions, and let $q\in [1,\infty]$.
Suppose that $\mu$ satisfies the $\mathrm{GWALA}_q$ condition with respect to
$\mathcal F$, with constants $\{C(t):t>0\}$.

Then there exists a constant
$\vartheta_0=\vartheta_0(n,\oldC{cubi},D,k,d,\{C(t):t>0\})$ such that the following holds. Let
$\sigma,\vartheta,\varepsilon,\delta,\zeta,\mathfrak m \in (0,1/2)$ and
$A\geq 1$ be such that $\vartheta\leq \vartheta_0$ and
$$
\vartheta\leq
\frac{1}{16nd^2}
\Big(\frac{1}{2^{100(k+1)}D^5(d+1)^2(d+2)^2}\Big)^{2k+2},
\qquad
\zeta \leq \delta\leq 1/2^{24},
\qquad
A\vartheta^n\geq 16n,
$$
$$
\frac{4\zeta}{\varepsilon\sqrt{1-\sigma^2}}
\leq
\frac{\sigma\vartheta^6}{2^{20}\sqrt{n-1}},
\qquad
16n\Big(\frac{\sigma}{\vartheta^{2n}}
+n\frac{\varepsilon+\delta+\zeta+\mathfrak m}{A\vartheta^{2n}}\Big)
\leq
\frac{1}{2\cdot 3^k16^nD^2}.
$$
Then, for every $0\leq j\leq k-1$, the family $\mathscr I(j,2\vartheta,\sigma/2,\varepsilon,\delta,\zeta,\mathfrak m,A)$ yielded by \cref{rioridinocubi}, is a Carleson family.
\end{teorema}

\begin{proof}
    We argue by contradiction. The proof is split in several steps.

\textbf{Step I. The setup.}
Suppose this is not the case. Then, there exists $0\leq j\leq k-1$ such that
$\mathscr I(j,2\vartheta,\sigma/2,\varepsilon,\delta,\zeta,\mathfrak m,A)$ is not a Carleson family. Then, choose $\xi,\varsigma\in (0,1/2)$ and $\mathfrak N,\qof\in \mathbb{N}$ such that 
    \begin{equation}
        \begin{split}
              &\qquad\qof\geq  4\log_2\Big(\frac{12800A}{\vartheta\varsigma\delta \zeta}
\Big\lceil n\log_2\Big(\frac{6400\sqrt{n-1}}{\sigma\vartheta}\Big)\Big\rceil
\Big\lceil\frac{20\log(\xi)\log(\frac{\sigma\vartheta^6}{2^{25}\sqrt{n-1}})}{\varsigma\vartheta}\Big\rceil\Big)+8\log_2(2^{10}A\oldC{cubi}^2)\\
&+4(k+2)\log_2(2^{100(k+1)}D^{12}d(d+1)^2(d+2)^2)+\log_2(2^{20}d^2n^4)-4\log_2(\sigma+\tfrac{n}{A}(\varepsilon+\delta +\zeta+\mathfrak m)).
\nonumber
        \end{split}
    \end{equation}
    and
    $$\xi\leq \frac{\vartheta^8}{128n^2D^8 5^{4k}}\qquad\text{and}\qquad\mathfrak N:=\Big\lfloor \frac{1}{64(n+1)\vartheta}\Big\rfloor.$$
By \cref{propfinalejoiningcones}, applied with depth $16\mathfrak N$, there
exist
$\mathfrak e\in(\mathbb S^{n-1})^j$, $R\in \mathscr I(j,2\vartheta,\sigma/2,\varepsilon,\delta,\zeta,\mathfrak m,A)$ and families
$\mathfrak L_0,\ldots,\mathfrak L_{16\mathfrak N}\subseteq \mathscr I(j,2\vartheta,\sigma/2,\varepsilon,\delta,\zeta,\mathfrak m,A)$, subfamilies
$$\mathfrak W_0,\ldots,\mathfrak W_{16\mathfrak N-1}\subseteq \mathscr I(j,2\vartheta,\sigma/2,\varepsilon,\delta,\zeta,\mathfrak m,A),$$
and nested open sets
$\Omega_{\mathfrak s}=:\Omega_0\supseteq \Omega_1\supseteq\ldots
\supseteq\Omega_{16\mathfrak N}$,
satisfying the conclusions of \cref{propfinalejoiningcones}. 
In what follows, we let $e_1,\ldots, e_{n-j}$ be an orthonormal basis of $\mathrm{span}(\mathfrak e)^\perp$ and we let $X$ be the family of Lipschitz functions
$f:\R^n\to\R$ such that
\begin{itemize}
    \item[(i)] denoting $\mathscr E=\{e_{1},\ldots, e_{n-j}\}$, we have
    $$|\partial_{e_\iota}f(y)|\leq 1\qquad\text{for every $\iota=1,\ldots,n-j$ and for $\Leb^n$-almost every }y\in\R^n;$$
    \item[(ii)] for $\Leb^n$-almost every $y\in\R^n$ there holds $|d_{\mathrm{span}(\mathfrak e)}f(y)|\leq 1$.    
\end{itemize}
In the following we let $V:=\mathrm{span}(\mathfrak e)$.

\medskip

\textbf{Step II. The first step of the perturbation: smoothing.} Let us fix $\ell=0,\ldots, \mathfrak N-1$ and let $f$ be some function in $X$. Thanks to Vitali's covering theorem, we can extract from $\mathfrak W_{16\ell+1}$ a subfamily of cubes $\mathfrak C_{16\ell+1}$ such that the balls $\tfrac{3}{2}AB_Q$, with $Q\in \mathfrak C_{16\ell+1}$, are disjoint and $\{\tfrac{15}{2}AB_Q:Q\in \mathfrak C_{16\ell+1 }\}$ cover $\Omega_{16\ell+1}$. Let us notice that 
\begin{equation}
    \begin{split}
\sum_{Q\in \mathfrak C_{16\ell+1}}\mu(\tfrac{3}{2}AB_Q)\geq& \frac{1}{D}\sum_{Q\in \mathfrak C_{16\ell+1}}\Big(\tfrac{3}{2}A\diam Q\Big)^k\geq \frac{1}{D5^k}\sum_{Q\in \mathfrak C_{16\ell+1}}\Big(\tfrac{15}{2}A\diam Q\Big)^k\\
\geq &\frac{1}{D^25^k}\sum_{Q\in \mathfrak C_{16\ell+1}}\mu(\tfrac{15}{2}AB_Q)\geq \frac{\mu(\Omega_{16\ell+1})}{D^25^k}\geq \frac{\mu(R)}{2D^25^k}.
\label{perditadimassa1}
    \end{split}
\end{equation}
For each $Q\in \mathfrak C_{16\ell+1}$ we let $\alpha_1,\alpha_2:=2^{-\qof/4}$
in $\tfrac{3}{2}AB_Q$, \cref{mollificazioneapproxalbordo} gives us for every $Q\in \mathfrak C_{16\ell+1}$ a function $u_Q$ satisfying the following properties 
\begin{itemize}
    \item[($\alpha$)] $u_Q=\tfrac{1}{1+2^{-\qof/4}}f$ on $\partial(\tfrac{3}{2}AB_Q)$; 
    \item[($\beta$)] $|d_Vu_Q|\le 1$ and     $|\partial_{e_\iota}u_Q|\leq 1$ in $\tfrac{3}{2}AB_Q$ for every $\iota=1,\ldots,n-j$;
    \item[($\gamma$)] $u_Q=\frac{1}{1+2^{-\qof/4}}\rho_{\frac{2^{-\qof/2}A\diam Q}{12\sqrt{n-j+1}}}*f$ on $B(\mathfrak c(Q),\tfrac{3}{2}(1-2^{-\qof/4-1} )A \diam Q)$;
    \item[($\delta$)]  $\lVert (1+2^{-\qof/4})^{-1}f-u_Q\rVert_\infty\leq \frac{2^{-\qof/2}A\diam Q}{12}$.
\end{itemize}
Thus, in the following, we let 
\begin{equation}
    g(z):=\begin{cases}
        \tfrac{1}{1+2^{-\qof/4}}f(z)&\text{if }z\not\in \bigcup_{Q\in \mathfrak C_{16\ell+1}}\tfrac{3}{2}AB_Q;\\
        u_Q(z) &\text{if }z\in \tfrac{3}{2}AB_Q\text{ and }Q\in \mathfrak C_{16\ell+1}.
    \end{cases}
\end{equation}
The function $g$ belongs to $X$. Indeed, it is continuous by construction and by construction items (i) and (ii) are satisfied. 

\medskip

\textbf{Step III. The first step of the perturbation. Canceling the gradient.}

We have smoothed the function $f$, and we now move to smaller scales in order to cancel the gradient. To do so, we will restrict to scales where the modified function is linear. Let us observe that thanks to item 7 of the thesis of \cref{propfinalejoiningcones}, we know that if $\tfrac{3}{2}AB_Q\cap \tfrac{3}{2}AB_{Q'}\neq \emptyset$ with $Q\in \mathfrak W_{16\ell+2}$ and $Q'\in \mathfrak W_{16\ell+1}$, then 
$$\diam Q\leq 2^{-\qof}\diam Q'.$$
Applying \cref{lemmacostruzionefammigliepallecontrollate} to $\mathscr B_1:=\{\tfrac{3}{2}AB_{Q'}:Q'\in \mathfrak C_{16\ell+1}\}$ and $\mathscr B_2:=\{\tfrac{3}{2}AB_Q:Q\in \mathfrak W_{16\ell+2}\}$ and with $\rho=2^{-\qof/4-1}$ and $\theta=2^{-\qof}$ we find a family of cubes  $\mathfrak C_{16\ell+2}$ in $\mathfrak B_{16\ell+2}$ such that the balls $\{\tfrac{3}{2}AB_Q:Q\in \mathfrak C_{16\ell+2}\}$ are pairwise disjoint and, for every $Q\in \mathfrak C_{16\ell+2}$, there exists $Q'\in \mathfrak C_{16\ell+1}$ such that
    $$
    \frac{3}{2}AB_Q\subseteq \frac{3}{2}(1-2^{-\qof/4-1} )AB_{Q'};
    $$
and for every $Q'\in \mathfrak C_{16\ell+1}$ we have
    $$
    \frac{3}{2}(1-2^{-\qof/4})AB_{Q'}\cap \bigcup_{Q\in \mathfrak W_{16\ell+2}}\frac{3}{2}AB_Q
    \subseteq
    \bigcup_{Q\in \mathfrak C_{16\ell+2}}\frac{15}{2}AB_Q.
    $$
    Recall further that thanks to item 7 of the thesis of \cref{propfinalejoiningcones}, since the balls $\{\frac{3}{2}AB_Q:Q\in \mathfrak W_{h+1}\}$ cover $\Omega_{h+1}$, then 
   \begin{equation}
   \begin{split}
          \tfrac{3}{2}(1-2^{-\qof/4})AB_{Q'}\cap \Omega_{16\ell+2}\subseteq \tfrac{3}{2}(1-2^{-\qof/4})AB_{Q'}\cap \bigcup_{Q\in \mathfrak W_{16\ell+2}}\tfrac{3}{2}AB_Q\cap \Omega_{16\ell+2}
    \subseteq
    \bigcup_{B\in \mathfrak C_{16\ell+2}}\tfrac{15}{2}AB_Q.
    \label{inclusioneottima}
   \end{split}
   \end{equation}
Let us observe that 
\begin{equation}
    \begin{split}
        \mu\Big(\bigcup_{Q\in \mathfrak C_{16\ell+2}} \tfrac{3}{2}AB_Q\Big)=&\sum_{Q\in \mathfrak C_{16\ell+2}}\mu(\tfrac{3}{2}AB_Q)\geq \frac{1}{D}\sum_{Q\in \mathfrak C_{16\ell+2}}\Big(\tfrac{3}{2}A\diam Q\Big)^k\\
        \geq& \frac{1}{5^kD}\sum_{Q\in \mathfrak C_{16\ell+2}}\Big(\tfrac{15}{2}A\diam Q\Big)^k\geq \frac{1}{5^kD^2}\sum_{Q\in \mathfrak C_{16\ell+2}}\mu(\tfrac{15}{2}AB_Q)\\
       \overset{\eqref{inclusioneottima}}{\geq}& \frac{1}{5^kD^2}\sum_{Q'\in \mathfrak C_{16\ell+1}}\mu(\tfrac{3}{2}(1-2^{-\qof/4})AB_{Q'}\cap \Omega_{16\ell+2}).
        \label{perditadimassav2}
    \end{split}
\end{equation}
Thanks to item 7 of the thesis of \cref{propfinalejoiningcones}, recall that for every  $Q'\in \mathfrak W_{16\ell+1}$, we have
$$\tfrac{3}{2}(1-2^{-\qof/4})AB_{Q'}\cap \supp\mu\subseteq R,$$
hence thanks to the disjointness of the balls $\{\tfrac{3}{2}AB_Q:Q\in \mathfrak C_{16\ell+1}\}$, and item 8 of the thesis of \cref{propfinalejoiningcones}, we have
\begin{equation}
    \begin{split}
&\sum_{Q'\in \mathfrak C_{16\ell+1}}\mu(\tfrac{3}{2}(1-2^{-\qof/4})AB_{Q'}\cap \Omega_{16\ell+2})\geq \sum_{Q'\in \mathfrak C_{16\ell+1}}\mu(\tfrac{3}{2}(1-2^{-\qof/4})AB_{Q'})-\mu(R\setminus \Omega_{16\ell+2})\\
&\quad\geq\frac{(1-2^{-\qof/4})^k}{D^25^k} \sum_{Q'\in \mathfrak C_{16\ell+1}}\mu(\tfrac{15}{2}AB_{Q'})-\mu(R\setminus \Omega_{16\ell+2})\overset{\eqref{perditadimassa1}}{\geq}\Big(\frac{1}{4 D^4 5^{2k}}-4\xi^{16\mathfrak N}\Big)\mu(R). 
    \end{split}
\end{equation}
Finally, from what above we have
\begin{equation}
    \begin{split}
        \mu\Big(\bigcup_{Q\in \mathfrak C_{16\ell+2}} \tfrac{3}{2}AB_Q\Big)\geq &\frac{1}{5^kD^2}\Big(\frac{1}{4 D^4 5^{2k}}-4\xi^{16\mathfrak N}\Big)\mu(R)\geq \frac{1}{8D^65^{3k}}\mu(R),
        \label{stimaperditadimassav3}
    \end{split}
\end{equation}
where the last inequality comes from the choice of $\xi$ and $\mathfrak N$.
   
Note that thanks to item 7 of the thesis of \cref{propfinalejoiningcones}, we know that $\diam Q\leq 2^{-\qof}\diam Q'$, to item ($\gamma$) and \cref{proprietamollificata} we infer that 
$$\lvert  u_{Q'}(y)- u_{Q'}(x)-D u_{Q'}(x)[y-x]\rvert\leq \frac{480n\sqrt{n-j+1}}{2^{-\qof/2}A\diam Q'}\lvert y-x\rvert^2,$$
for every $y,x\in \tfrac{3}{2}(1-2^{-\qof/4-1})AB_{Q'}$. In particular, if $y,x\in \tfrac{3}{2}AB_Q$, where $Q\in \mathfrak C_{16\ell+2}$, then by item $4$ of \cref{proprietamollificata} we have
\begin{equation}
    \begin{split}
        &\qquad\qquad\qquad\qquad\lvert  u_{Q'}(y)- u_{Q'}(x)-D u_{Q'}(\mathfrak c(Q))[y-x]\rvert\\
        &\leq \frac{480n\sqrt{n-j+1}}{2^{-\qof/2}A\diam Q'}\lvert y-x\rvert^2
        +\frac{120n\sqrt{n-j+1}}{2^{-\qof/2}A\diam Q'}\lvert x-\mathfrak c(Q)\rvert\lvert y-x\rvert\leq 1000n^2\cdot 2^{-\qof/2}\lvert y-x\rvert.
        \label{stimaappoxbuonas}
    \end{split}
\end{equation}
By the scale separation for the families $\mathfrak W_h$ given by item 7 of the thesis of \cref{propfinalejoiningcones}, if $Q\in \mathfrak W_{16\ell+2}$ and $S\in \mathfrak W_{16\ell+3}$ are such that
$10AB_Q\cap 10AB_S\neq \emptyset$, then
$$\diam S\leq 2^{-\qof}\diam Q.$$
We apply \cref{lemmacostruzionefammigliepallecontrollate} to the families 
$\mathscr B_1:=\{\tfrac{3}{2}AB_Q:Q\in \mathfrak C_{16\ell+2}\}$ and $\mathscr B_2:=\{\tfrac{3}{2}AB_S:S\in \mathfrak W_{16\ell+3}\}$,
with $\rho=\sqrt{\xi}$ and $\theta=2^{-\qof}$. We obtain a subfamily
$\mathfrak C_{16\ell+3}\subseteq \mathfrak W_{16\ell+3}$ such that the balls
$\{\tfrac{3}{2}AB_S:S\in \mathfrak C_{16\ell+3}\}$ are pairwise disjoint and, for every
$S\in \mathfrak C_{16\ell+3}$, there exists $Q\in \mathfrak C_{16\ell+2}$ such that
$$
\tfrac{3}{2}AB_S\subseteq \tfrac{3}{2}(1-\sqrt{\xi})AB_Q.
$$
Moreover, for every $Q\in \mathfrak C_{16\ell+2}$ we have
$$
\tfrac{3}{2}(1-\sqrt{\xi}-2^{1-\qof})AB_Q
\cap
\bigcup_{S\in \mathfrak W_{16\ell+3}}\tfrac{3}{2}AB_S
\subseteq
\bigcup_{S\in \mathfrak C_{16\ell+3}}\tfrac{15}{2}AB_S.
$$
Since $\Omega_{16\ell+3}\subseteq \bigcup_{S\in \mathfrak W_{16\ell+3}}\tfrac{3}{2}AB_S$, arguing as in \eqref{perditadimassav2} and \eqref{perditadimassa1}, we get
\begin{equation}
    \begin{split}
        &\qquad\mu\Big(\bigcup_{S\in \mathfrak C_{16\ell+3}}\tfrac{3}{2}AB_S\Big)\geq\frac{1}{5^kD^2}\sum_{Q\in \mathfrak C_{16\ell+2}}\mu\Big(\tfrac{3}{2}(1-\sqrt{\xi}-2^{1-\qof})AB_Q\cap \Omega_{16\ell+3}\Big)\\
        &\geq\frac{1}{5^kD^2}\Big(\frac{(1-\sqrt{\xi}-2^{1-\qof})^k}{D^25^k}\sum_{Q\in \mathfrak C_{16\ell+2}}\mu(\tfrac{15}{2}AB_Q)-\mu(R\setminus \Omega_{16\ell+3})\Big)\overset{\eqref{stimaperditadimassav3}}{\geq}\frac{1}{32D^{10}5^{5k}}\mu(R),
        \label{stimapalleC16l+2}
    \end{split}
\end{equation}
where the last inequality follows from the choice of $\xi$ and $\qof$. We now define 
\begin{equation}
    \begin{split}
        g_0(x):=\tfrac{\sqrt{1-9\vartheta^2}}{\sqrt{1-9\vartheta^2}+6n\vartheta}\Big(g(x)-\sum_{Q\in \mathfrak C_{16\ell+2}}\sum_{\iota=1}^{n-j}Du_{Q'}(\mathfrak c(Q))[e_\iota]\omega_{\iota,\ell,Q}(x)\Psi_Q(x)\Big),
        \nonumber
    \end{split}
\end{equation}
where, for each $Q\in\mathfrak C_{16\ell+2}$, $Q'\in\mathfrak C_{16\ell+1}$ denotes the cube fixed above, and where  $\omega_{\iota,\ell,Q}:=\omega_{e_\iota,\sqrt{1-9\vartheta^2}}[\Omega_{16\ell+3}\cap \tfrac{3}{2}AB_Q]$ is the width function constructed in \cref{p:width-open} and associated to $\Omega_{16\ell+3}\cap \tfrac{3}{2}AB_Q$ and to the cone $C(e_\iota,\sqrt{1-9\vartheta^2})$, and where $\Psi_Q$ are smooth functions such that $\Psi_Q=0$ on $(\tfrac{3}{2}AB_Q)^c$, $\Psi_Q=1$ on $ \tfrac{3}{2}(1-\sqrt{\xi})AB_{Q}$ and such that $\lVert D\Psi_Q\rVert_\infty\leq 10(\sqrt{\xi} A\diam Q)^{-1}$.

Now we prove that $g_0\in X$. In order to do so, we compute the derivatives along $V$ and along the vectors $e_\iota$ with $\iota=1,\ldots, n-j$ of the function $g_0$ inside the balls $\tfrac{3}{2}AB_Q$ with $Q\in \mathfrak C_{16\ell+2}$.
First of all, let us notice that 
$$g=\tfrac{\sqrt{1-9\vartheta^2}+6n\vartheta}{\sqrt{1-9\vartheta^2}}g_0\qquad\text{on }\R^n\setminus \bigcup_{Q\in \mathfrak C_{16\ell+2}}\tfrac{3}{2}AB_Q,$$
and this shows there is nothing to prove in this case, as $g\in X$. 
Let us assume $y\in \tfrac{3}{2}AB_Q$ for some $Q\in \mathfrak C_{16\ell+2}$. Then, for every direction $w\in \mathbb S^{n-1}$ we have 
\begin{equation}
\begin{split}
     \tfrac{\sqrt{1-9\vartheta^2}+6n\vartheta}{\sqrt{1-9\vartheta^2}}\partial_wg_0(y)=\partial_w g(y)-&\sum_{\iota=1}^{n-j}Du_{Q'}(\mathfrak c(Q))[e_\iota]\partial_w\omega_{\iota,\ell,Q}(y)\Psi_Q(y)\\
    -&\sum_{\iota=1}^{n-j}Du_{Q'}(\mathfrak c(Q))[e_\iota]\omega_{\iota,\ell,Q}(y)\partial_w\Psi_Q(y).
     \nonumber
\end{split}
\end{equation}
Let us first estimate the term
$$II:=\sum_{\iota=1}^{n-j}Du_{Q'}(\mathfrak c(Q))[e_\iota]\omega_{\iota,\ell,Q}\partial_w\Psi_Q.$$
Thanks to item 7 of the thesis of \cref{propfinalejoiningcones} and by \cref{propo:curvainpalla}, we have the following estimate
\begin{equation}
    \begin{split}
        \lVert\omega_{\iota,\ell,Q}\rVert_\infty\leq \xi\sup\{\Haus^1(\gamma\cap 2(1+\zeta)AB_Q):\gamma\text{ is a $C(e_\iota,\sqrt{1-9\vartheta^2})$-curve}\}       \leq \frac{8\xi}{\vartheta} A\diam Q.
        \label{stimasupnorma}
    \end{split}
\end{equation}
Thus, thanks to the choice of $\Psi_Q$ and of $\vartheta$ and $\xi$, we have 
\begin{equation}
    \begin{split}
\lVert\omega_{\iota,\ell,Q}\partial_w\Psi_Q\rVert_\infty\leq&\frac{8\xi}{\vartheta}A\diam Q \cdot \frac{10}{\sqrt{\xi}A\diam Q}\leq 80\frac{\sqrt{\xi}}{\vartheta}\leq \xi^{1/4},
\nonumber
    \end{split}
\end{equation}
And this shows in particular that $\lvert II\rvert\leq n\xi^{1/4}$. Now, if $w\in V$, by item 3 of \cref{p:width-open} and few algebraic computations, we have that 
\begin{equation}
    \begin{split}
     & \tfrac{\sqrt{1-9\vartheta^2}}{\sqrt{1-9\vartheta^2}+6n\vartheta}\Big(\Big\lvert  \partial_w g-\sum_{\iota=1}^{n-j}Du_{Q'}(\mathfrak c(Q))[e_\iota]\partial_w\omega_{\iota,\ell,Q}\Psi_Q\Big\rvert+|II|\Big)\\
     &\qquad\qquad\leq  \tfrac{\sqrt{1-9\vartheta^2}}{\sqrt{1-9\vartheta^2}+6n\vartheta}\Big(1+n\frac{3\vartheta}{\sqrt{1-9\vartheta^2}}+n\xi^{1/4}\Big)\leq 1.
     \end{split}
\end{equation}
This shows that $\lVert d_Vg_0\rVert_\infty\leq 1$. On the other hand, let us assume that $w\in \{e_1,\ldots, e_{n-j}\}$.
Without loss of generality, let us assume that $w=e_1$ and set
$$I:=\partial_{e_1}g-\sum_{\iota=1}^{n-j}Du_{Q'}(\mathfrak c(Q))[e_\iota]\partial_{e_1}\omega_{\iota,\ell,Q}\Psi_Q.$$
By item ($\gamma$) and by item 4 of \cref{proprietamollificata}, since
$\tfrac{3}{2}AB_Q\subseteq \tfrac{3}{2}(1-2^{-\qof/4-1})AB_{Q'}$, we have
\begin{equation}
    \begin{split}
    \big|
    \partial_{e_1}g(y)-Du_{Q'}(\mathfrak c(Q))[e_1]
    \big|
    &\leq
    \frac{120n\sqrt{n-j+1}}{2^{-\qof/2}A\diam Q'}
    |y-\mathfrak c(Q)|\leq180n^2 2^{-\qof/2}.
    \end{split}
    \nonumber
\end{equation}
Therefore, for $\Leb^n$-almost every $y\in \tfrac{3}{2}AB_Q$, we have
\begin{equation}
    \begin{split}
      |I|\leq180n^2 2^{-\frac{\qof}{2}}&+\Big\lvert Du_{Q'}(\mathfrak c(Q))[e_1]\big(1-\partial_{e_1}\omega_{1,\ell,Q}(y)\Psi_Q(y)\big)\Big\rvert\\
      &+\Big\lvert\sum_{\iota=2}^{n-j} Du_{Q'}(\mathfrak c(Q))[e_\iota]\partial_{e_1}\omega_{\iota,\ell,Q}(y)\Psi_Q(y)\Big\rvert\\
      \leq&1+n\frac{3\vartheta}{\sqrt{1-9\vartheta^2}}+180n^2 2^{-\qof/2}.
        \nonumber
    \end{split}
\end{equation}
Hence, in this case we have
\begin{equation}
    \begin{split}
    \lvert \partial_{e_1}g_0(y)\rvert
    &\leq
    \tfrac{\sqrt{1-9\vartheta^2}}{\sqrt{1-9\vartheta^2}+6n\vartheta}
    (|I|+|II|)\\
    &\leq
    \tfrac{\sqrt{1-9\vartheta^2}}{\sqrt{1-9\vartheta^2}+6n\vartheta}
    \Big(
    1+n\frac{3\vartheta}{\sqrt{1-9\vartheta^2}}
    +180n^2 2^{-\qof/2}
    +n\xi^{1/4}
    \Big)
    \leq 1,
    \end{split}
\end{equation}
where by the choices of $\qof$ and $\xi$, we have
$180n^2 2^{-\qof/2}+n\xi^{1/4}
\leq
n\frac{3\vartheta}{\sqrt{1-9\vartheta^2}}$.
Verbatim, one obtains
$$
\lVert\partial_{e_\iota} g_0\rVert_\infty\leq 1
\qquad\text{for every }\iota=1,\ldots,n-j.
$$

We now further prove a more precise result, that tells us that $g_0$ is almost constant along the directions $e_1,\ldots, e_{n-j}$, on the balls $\tfrac{3}{2}AB_S$, with $S\in \mathfrak C_{16\ell+3}$.
Let $S\in\mathfrak C_{16\ell+3}$ and let
$Q\in\mathfrak C_{16\ell+2}$, $Q'\in\mathfrak C_{16\ell+1}$ be such that
$S\subseteq Q\subseteq Q'$. Since $\tfrac{3}{2}AB_S\subseteq \Omega_{16\ell+2}\cap (1-\sqrt{\xi})\tfrac{3}{2}AB_Q$ and, by construction,
$\Psi_Q\equiv1$ on $\frac{3}{2}AB_S$, for $\Leb^n$-almost every $y\in AB_S$ we have
\begin{equation}
    \begin{split}
      \frac{\sqrt{1-9\vartheta^2}+6n\vartheta}{\sqrt{1-9\vartheta^2}}\partial_{e_\iota}g_0(y)=\partial_{e_\iota}g(y)&-Du_{Q'}(\mathfrak c(Q))[e_\iota]-    \sum_{\substack{i=1\\ i\neq \iota}}^{n-j} Du_{Q'}(\mathfrak c(Q))[e_i]\partial_{e_\iota}\omega_{i,\ell,Q}.
      \nonumber
    \end{split}
\end{equation}
As seen above, by item ($\gamma$) and by item 4 of \cref{proprietamollificata}, since
$\tfrac{3}{2}AB_Q\subseteq \tfrac{3}{2}(1-2^{-\qof/4-1})AB_{Q'}$, arguing as above we have
\begin{equation}
    \begin{split}
    \big|\partial_{e_\iota}g(y)-Du_{Q'}(\mathfrak c(Q))[e_\iota]\big|\leq180n^2 2^{-\qof/2}.
\end{split}
    \nonumber
\end{equation}
Moreover, by \cref{p:width-open}, if $i\neq \iota$ then
$$|\partial_{e_\iota}\omega_{i,\ell,Q}(y)|\leq\frac{3\vartheta}{\sqrt{1-9\vartheta^2}}.$$
Thus, for every $\iota=1,\ldots,n-j$ and for $\Leb^n$-almost every $y\in AB_S$,
\begin{equation}
    \begin{split}
        |\partial_{e_\iota} g_0(y)|&\leq\frac{\sqrt{1-9\vartheta^2}}{\sqrt{1-9\vartheta^2}+6n\vartheta}\Big(180n^2 2^{-\qof/2}+n\frac{3\vartheta}{\sqrt{1-9\vartheta^2}}\Big)\leq\frac{6n\vartheta}{\sqrt{1-9\vartheta^2}},
        \label{stimalocalfgrason}
    \end{split}
\end{equation}
where the last inequality follows from the choice of $\qof$. 
Now notice that, for every
$S\in \mathfrak C_{16\ell+3}$ with $S\subseteq Q\subseteq Q'$, where
$Q\in\mathfrak C_{16\ell+2}$ and $Q'\in\mathfrak C_{16\ell+1}$, and every
$y,x\in \frac{3}{2}AB_S$, there holds
\begin{equation}
    \begin{split}
    \Big|g_0(y)-g_0(x)-\frac{\sqrt{1-9\vartheta^2}}{\sqrt{1-9\vartheta^2}+6n\vartheta}    &Du_{Q'}(\mathfrak c(Q))[\pi_V(y-x)]\Big|\\   
    &\qquad\qquad\leq\Big(720n^2\cdot 2^{-\qof/2}+   2n\frac{3\vartheta}{\sqrt{1-9\vartheta^2}}\Big)\lvert y-x\rvert.
    \end{split}
    \label{linearitag_0sdfsf}
\end{equation}
Indeed, for this choice of $y,x$ we have $\Psi_Q\equiv1$, and hence
\begin{equation}
    \begin{split}
        &\Big|g_0(y)-g_0(x)-\frac{\sqrt{1-9\vartheta^2}}{\sqrt{1-9\vartheta^2}+6n\vartheta}        \pi_{\mathrm{span}(\mathfrak e)}(Du_{Q'}(\mathfrak c(Q)))[y-x]\Big|\\
        &\leq\frac{\sqrt{1-9\vartheta^2}}{\sqrt{1-9\vartheta^2}+6n\vartheta}\underbrace{\Big|g(y)-g(x)-Du_{Q'}(\mathfrak c(Q))[y-x]\Big|}_{=:(A)}\\
        &\quad+\frac{\sqrt{1-9\vartheta^2}}{\sqrt{1-9\vartheta^2}+6n\vartheta}\underbrace{
        \Big|Du_{Q'}(\mathfrak c(Q))[\pi_{V^\perp}(y-x)]-\sum_{\iota=1}^{n-j}Du_{Q'}(\mathfrak c(Q))[e_\iota]\big(\omega_{\iota,\ell,Q}(y)-\omega_{\iota,\ell,Q}(x)\big)\Big|}_{=:(B)} .
        \nonumber
    \end{split}
\end{equation}
By \eqref{stimaappoxbuonas}, we have
$(A)\leq 720n^2\cdot 2^{-\qof/2}|y-x|$. Moreover,
\begin{equation}
    \begin{split}
    (B)\leq\sum_{\iota=1}^{n-j}\Big|Du_{Q'}(\mathfrak c(Q))[e_\iota]\Big(\langle e_\iota,y-x\rangle-\big(\omega_{\iota,\ell,Q}(y)-\omega_{\iota,\ell,Q}(x)\big)\Big)\Big|\leq2n\frac{3\vartheta}{\sqrt{1-9\vartheta^2}}|y-x|.
    \end{split}
    \nonumber
\end{equation}
The last inequality follows from \cref{p:width-open}. Indeed, on
$\Omega_{16\ell+3}\cap \tfrac{3}{2}AB_Q$ the width function
$\omega_{\iota,Q}$ has derivative $1$ in the direction $e_\iota$ and transverse Lipschitz constant
$\frac{3\vartheta}{\sqrt{1-9\vartheta^2}}$.

\medskip

\textbf{Step IV. The final perturbation away from $\mathcal F$.} Throughout this step, the cube $Q\in \mathfrak C_{16\ell+3}$ will be considered fixed and so, we will drop the dependence on $Q$ of all the objects. We are going to go to a scale where \cref{prop:discrete-product-from-partial-invariance} applies and at those scales construct the perturbation. 
Before doing so, we need to extract a covering of the ball $\frac{A}{32n}B_Q$ with balls of diameter smaller than $\frac{A\vartheta^j}{256n}\diam Q$.

Applying Vitali's covering theorem, it is immediate to see that there exists a finite family of disjoint balls $B(x_i,\frac{A\vartheta^j}{1280n}\diam Q)$ with $i=1,\ldots, H_Q$ and $x_i\in \supp\mu\cap \frac{A}{32n}B_Q$ such that the balls $B(x_i,\frac{A\vartheta^j}{256n}\diam Q)$ cover $\supp\mu\cap \frac{A}{32n}B_Q$. Of course these balls are all contained in $\tfrac{A}{2}B_Q$. Recall that $Q\in \mathfrak C_{16\ell+3}$ is a $(j,\mathfrak{e},\vartheta^j,\sigma, \varepsilon,\delta,\zeta,\mathfrak m,A)$-invariant cube. Hence, defined 
$$\rho:=2^{9+12k}D^2
\Big(
\frac{\sigma}{\vartheta
^j}
+
n\frac{\varepsilon+\delta+\zeta+\mathfrak m}{A\vartheta^j}
\Big)
A\vartheta^j \diam Q,$$
we have thanks to \cref{prop:discrete-product-from-partial-invariance}, that for every $i=1,\ldots, H_Q$, there exists a discrete set
$E_{Q,x_i}\subset \mathrm{span}(\mathfrak e)^\perp$ that is $\rho$ separated such that, writing
$K_{x_i}:=x_i+\mathrm{span}(\mathfrak e)\otimes E_{Q,x_i}$, one has
\begin{equation}
d_{\Haus,B(x_i,\frac{A\vartheta^j}{256n}\diam Q)}
(\supp\mu,K_{x_i})\leq 2\rho.
\label{eq:discrete-product-hausdorffv2}
\end{equation}
Moreover, for every $r$ such that $\rho \leq r \leq \frac{A\vartheta^j}{1024n}\diam Q$, one has
\begin{equation}
2^{-8k}D^{-2}\Big(\frac{r}{\rho}\Big)^{k-j}
\leq
\mathrm{Card}\big(E_{Q,x_i}\cap B(p,r)\big)\leq 2^{5k}D^2\Big(\frac{r}{\rho}\Big)^{k-j},
\nonumber
\end{equation}
for every $p\in E_{Q,x_i}\cap B_{\mathfrak e}\big(0,\frac{A\vartheta^j}{128n}\diam Q\big)$. We need to extract a more coarse set $F_{Q,x_i}$ from $E_{Q,x_i}$ so that the scale of closeness of $\supp\mu$ to $K_{x_i}$ is smaller than the scale of separation of the points of $F_{Q,x_i}$. In order to do this, we let $F_{Q,x_i}$ be maximal $10\rho$-separated in $E_{Q,x_i}$. Then, by definition, we have 
$E_{Q,x_i}\subset \bigcup_{q\in F_{Q,x_i}} B(q,10\rho)$.

Moreover, for every $r$ such that
$128\rho\leq r\leq \frac{A\vartheta^j}{2048n}\diam Q$, and every
$p\in F_{Q,x_i}\cap B_{\mathfrak e}(0,\frac{A\vartheta^j}{256n}\diam Q)$, one has
\begin{equation}
2^{-13k}D^{-4}\Big(\frac{r}{128\rho}\Big)^{k-j}\leq\mathrm{Card}(F_{Q,x_i}\cap B(p,r))\leq2^{21k}D^4
\Big(\frac{r}{128\rho}\Big)^{k-j}.
\label{eq:coarse-discrete-lower}
\end{equation}
Indeed, the upper bound follows from the $10\rho$-separation of
$F_{Q,x_i}$. Since $r\geq 128\rho$, each ball $B(q,5\rho)$, with
$q\in F_{Q,x_i}\cap B(p,r)$, is contained in $B(p,2r)$, and these balls
are pairwise disjoint. Applying the original lower bound at scale $5\rho$
and the original upper bound at scale $2r$, we get
    \begin{equation}
        \begin{split}
            \mathrm{Card}(F_{Q,x_i}\cap B(p,r))
2^{-8k}D^{-2}5^{k-j}
\leq&
\sum_{q\in F_{Q,x_i}\cap B(p,r)}
\mathrm{Card}(E_{Q,x_i}\cap B(q,5\rho))\\
\leq& \mathrm{Card}(E_{Q,x_i}\cap B(p,2r))\leq 
2^{5k}D^2
\Big(\frac{2r}{\rho}\Big)^{k-j}.    
        \end{split}
    \end{equation}
In particular, we have 
$$
\mathrm{Card}(F_{Q,x_i}\cap B(p,r))
\leq
2^{21k}D^4
\Big(\frac{r}{128\rho}\Big)^{k-j}.
$$
For the lower bound, every point of $E_{Q,x_i}\cap B(p,r/2)$ is within
distance $10\rho$ from some point of $F_{Q,x_i}\cap B(p,r)$. Also, for every $q\in F_{Q,x_i}\cap B(p,r)$, the original
upper bound gives
$$
\mathrm{Card}(E_{Q,x_i}\cap B(q,10\rho))
\leq
2^{5k}D^2 10^{k-j}
\leq
2^{9k-4j}D^2.
$$
Therefore
$$
\mathrm{Card}(F_{Q,x_i}\cap B(p,r))
\geq
\frac{\mathrm{Card}(E_{Q,x_i}\cap B(p,r/2))}
{2^{9k-4j}D^2}.
$$
Using the original lower bound at scale $r/2$, we obtain
$$
\mathrm{Card}(F_{Q,x_i}\cap B(p,r))
\geq
2^{-18k+5j}D^{-4}
\Big(\frac{r}{\rho}\Big)^{k-j}
\geq 
2^{-13k}D^{-4}
\Big(\frac{r}{128\rho}\Big)^{k-j}.
$$
From now on, we define
$$\ayin:=\frac{1}{2^{40k+20}D^{10}d(d+1)^2(d+2)^2}$$
and we let $r_Q:=16d(d+1)^2D^2\ayin^{-1}\rho$.
Notice that $10r_Q\leq \frac{A\vartheta^j}{2560 n}\diam Q$ thanks to the choice of the parameters.

\medskip

We construct the perturbations in each of the balls $B(x_i,\frac{A\vartheta^j}{2560n}\diam Q)$. 
We let $\widehat{F}_{Q,x_i}$ be the family of those $p\in F_{Q,x_i}$ such that $B_{\mathfrak e}(p,r_Q)\cap B_{\mathfrak e}(0,\frac{A\vartheta^j}{5120n}\diam Q)\neq \emptyset$. We now apply Vitali's covering lemma to the family of balls $\mathscr{B}:=\{B_{\mathfrak e}(p,8r_Q):p\in \widehat{F}_{Q,x_i}\cap B_{\mathfrak e}(0,\frac{A\vartheta^j}{5120n}\diam Q)\}$ obtaining a disjoint subfamily $\mathscr F_{Q,x_i}\subseteq \widehat{F}_{Q,x_i}\cap B_{\mathfrak e}(0,\frac{A\vartheta^j}{5120n}\diam Q)$ such that 
\begin{equation}
            B_{\mathfrak e}(0,\tfrac{A\vartheta^j}{5120n}\diam Q)\cap F_{Q,x_i}\subseteq \bigcup_{p\in \widehat{F}_{Q,x_i}\cap B_{\mathfrak e}(0,\frac{A\vartheta^j}{5120n}\diam Q)}B_{\mathfrak e}(p,8r_Q)\subseteq \bigcup_{p\in \mathscr{F}_{Q,x_i}}B_{\mathfrak e}(p,40r_Q).
    \label{includiosasjgnaidjfgnajsnd}
\end{equation}
Since $\mathscr F_{Q,x_i}\subseteq B_{\mathfrak e}(0,\frac{A\vartheta^j}{5120n}\diam Q)$, for every $p\in\mathscr F_{Q,x_i}$ one has
$F_{Q,x_i}\cap B_{\mathfrak e}(p,r_Q)\subseteq \widehat F_{Q,x_i}\cap B_{\mathfrak e}(p,r_Q)$. Thus, thanks to the choice of $r_Q$, we have
$$\mathrm{Card}(\widehat{F}_{Q,x_i}\cap B_{\mathfrak e}(p,r_Q))\geq d+2\qquad\text{for every }p\in \mathscr{F}_{Q,x_i}.$$
Thus, thanks to \eqref{eq:coarse-discrete-lower} and to the choice of $\ayin$, for every $p\in \mathscr F_{Q,x_i}$ we can choose
\begin{equation}
    y_1(p),\ldots,y_{d+1}(p)\in \widehat{F}_{Q,x_i}\cap B_{\mathfrak e}(p,r_Q),
    \label{puntiyintroduzione}
\end{equation}
such that
$$
\lvert y_{j_1}(p)-y_{j_2}(p)\rvert\geq 10\ayin r_Q
\qquad\text{for every $j_1,j_2\in \{1,\ldots,d+1\}$ with $j_1\neq j_2$}.
$$
Indeed, suppose that this is not possible. Then a maximal $10\ayin r_Q$-separated subset of
$F_{Q,x_i}\cap B_{\mathfrak e}(p,r_Q)$ has cardinality at most $d+1$, and therefore
$F_{Q,x_i}\cap B_{\mathfrak e}(p,r_Q)$ is covered by at most $d+1$ balls of radius
$10\ayin r_Q$ centered in $F_{Q,x_i}\cap B_{\mathfrak e}(p,r_Q)$. Since
$$
10\ayin r_Q\geq 128\rho
\qquad\text{and}\qquad
10\ayin r_Q\leq \frac{A\vartheta^j}{2048n}\diam Q,
$$
we may apply \eqref{eq:coarse-discrete-lower} at the scales $r_Q$ and $10\ayin r_Q$. Thus
$$
2^{-13k}D^{-4}
\Big(\frac{r_Q}{128\rho}\Big)^{k-j}
\leq
(d+1)2^{21k}D^4
\Big(\frac{10\ayin r_Q}{128\rho}\Big)^{k-j},
$$
which contradicts the choice of $\ayin$. 

Let $\beta^1,\ldots,\beta^d\in\mathbb R^{d+1}$ be the vectors given by \cref{evadeaffinepro}, with $|\beta^\iota|\leq 1$ and
$\beta^\iota_1=0$ for every $\iota=1,\ldots,d$. For every $p\in \mathscr F_{Q,x_i}$ and every $j=1,\ldots,d+1$, by \eqref{eq:discrete-product-hausdorffv2} we choose a point $z_j(p)\in\supp\mu$ such that $|z_j(p)-x_i-y_j(p)|\leq 2\rho$. If we let
$\lambda:=\frac{5\ayin}{8d(d+1)^2D^2}$, 
applying \cref{evadeaffinepro} to the points $z_1(p),\ldots,z_{d+1}(p)$ gives
$$
\sum_{\iota=1}^{d} \inf_{A\in\mathcal F}\sup_{j=1,\ldots,d+1}\fint_{B(z_j(p),\lambda r_Q)}|\ayin r_Q\beta^\iota_j-A(z)|d\mu(z)\geq\frac{\ayin}{2(d+1)}r_Q.
$$
For every $Q\in \mathfrak C_{16\ell+3}$ and every $i=1,\ldots, H_Q$ we introduce $d$ $\tfrac{1}{5}$-Lipschitz functions $\phi_{Q,i}^\iota:\mathrm{span}(\mathfrak e)^\perp\to \R$, with $\iota=1,\ldots, d$, such that 
$$
\phi_{Q,i}^\iota(y_j(p)):= \ayin r_Q\beta^\iota_j\text{ for every }p\in\mathscr F_{Q,x_i}\text{ and every }j=1,\ldots, d+1,
$$
and there holds $\lVert \phi_{Q,i}^\iota\rVert_\infty\leq \ayin r_Q$. 
In order to prove that such functions exist, we just need to prove that, when restricted to the points 
$$\mathfrak T_{Q,i}:=\bigcup_{p\in \mathscr F_{Q,x_i}}\{y_1(p),\ldots, y_{d+1}(p)\},$$ the function $\phi_{Q,i}^\iota$ is $\tfrac{1}{5}$-Lipschitz and then obtain a good extension via McShane extension theorem.

Let us check the Lipschitz constant of $\phi_{Q,i}^\iota$ on $\mathfrak T_{Q,i}$ and in order to do so pick $y_1,y_2\in \mathfrak T_{Q,i}$. If $y_1=y_{j_1}(p_1)$ and $y_2=y_{j_2}(p_2)$ for $p_1\neq p_2$, then 
$\lvert y_1-y_2\rvert\geq 14 r_Q$, and hence 
$$\lvert \phi_{Q,i}^\iota(y_1)-\phi_{Q,i}^\iota(y_2)\rvert\leq 2\ayin r_Q\leq \frac{\ayin}{7}\lvert y_1-y_2\rvert.$$
On the other hand, if $y_1=y_{j_1}(p)$ and $y_2=y_{j_2}(p)$ then $\lvert y_1-y_2\rvert\geq 10 \ayin r_Q$ and arguing as above we have
$$\lvert \phi_{Q,i}^\iota(y_1)-\phi_{Q,i}^\iota(y_2)\rvert\leq 2\ayin r_Q\leq \frac{1}{5}\lvert y_1-y_2\rvert.$$
We now define the functions $\Phi_{Q,i}^\iota:\R^n\to \R$ as
$$\Phi_{Q,i}^\iota(z):=\phi_{Q,i}^\iota(\pi_{\mathrm{span}(\mathfrak e)^\perp}(z-x_i)),$$
and we notice that the functions $\Phi_{Q,i}^\iota$ are $\tfrac{1}{5}$-Lipschitz and
finally for every $\iota=1,\ldots, d$ define
$$\mathfrak{F}_\ell^\iota(f) :=(1-\vartheta)g_0+\sum_{Q\in \mathfrak C_{16\ell+3}}\sum_{i=1}^{H_Q}\Phi_{Q,i}^\iota\widetilde{\Psi}_{Q,i},$$
where $\widetilde{\Psi}_{Q,i}$ are cutoff functions supported on $B(x_i,\frac{A\vartheta^j}{1280n}\diam Q)$ and $\widetilde{\Psi}_{Q,i}\equiv 1$ on $B(x_i,\frac{A\vartheta^j}{2560n}\diam Q)$ and $\lVert  \nabla \widetilde{\Psi}_{Q,i}\rVert_\infty \leq 25600n(A\vartheta^j \diam Q)^{-1}$.

\medskip

\textbf{Step V. The functions $\mathfrak{F}_\ell^\iota(f)$ are in $X$. } As we did for $g_0$, we let $w\in \mathbb{S}^{n-1}$ be any direction and notice that 
$$\partial_w\mathfrak{F}_\ell^\iota(f)=(1-\vartheta)\partial_w g_0+\sum_{Q\in \mathfrak C_{16\ell+3}}\sum_{i=1}^{H_Q}\partial_w\Phi_{Q,i}^\iota\widetilde{\Psi}_{Q,i}+\sum_{Q\in \mathfrak C_{16\ell+3}}\sum_{i=1}^{H_Q}\Phi_{Q,i}^\iota\partial_w\widetilde{\Psi}_{Q,i}.$$
If $w\in\mathrm{span}(\mathfrak e)$, then $\partial_w\Phi_{Q,i}^\iota=0$ and
\begin{equation}
    \begin{split}
        \lvert\partial_w\mathfrak{F}_\ell^\iota(f)\rvert
        \leq &(1-\vartheta)\lvert\partial_w g_0\rvert+\Big\lvert \sum_{Q\in \mathfrak C_{16\ell+3}}\sum_{i=1}^{H_Q}\Phi_{Q,i}^\iota\partial_w\widetilde{\Psi}_{Q,i}\Big\rvert
       \leq  (1-\vartheta)+\frac{25600n}{A\vartheta^j \diam Q} \ayin r_Q,
    \end{split}
\end{equation}
where in the last inequality we used the fact that the balls $B(x_i,\frac{A\vartheta}{1280n}\diam Q)$ are disjoint. 
However, we immediately infer thanks to the choice of $r_Q$ and $\ayin$ that 
\begin{equation}
    \begin{split}
        \lvert\partial_w\mathfrak{F}_\ell^\iota(f)\rvert \leq (1-\vartheta)+ \frac{\vartheta}{2^{13(k+10)}D^4(d+2)}\leq 1,
    \end{split}
\end{equation}
We are left to discuss the case in which $w\in \{e_1,\ldots, e_{n-j}\}$. Let us assume without loss of generality that $w=e_1$. In this case we have 
\begin{equation}
    \begin{split}
         \lvert\partial_{e_1}\mathfrak{F}_\ell^\iota(f)\rvert \overset{\eqref{stimalocalfgrason}}{\leq}&(1-\vartheta)2n\frac{3\vartheta}{\sqrt{1-9\vartheta^2}}+\Big\lvert \sum_{Q\in \mathfrak C_{16\ell+3}}\sum_{i=1}^{H_Q}\partial_1\Phi_{Q,i}^\iota\widetilde{\Psi}_{Q,i}\Big\rvert+\Big\lvert\sum_{Q\in \mathfrak C_{16\ell+3}}\sum_{i=1}^{H_Q}\Phi_{Q,i}^\iota\partial_1\widetilde{\Psi}_{Q,i}\Big\rvert\\
 \leq &(1-\vartheta)2n\frac{3\vartheta}{\sqrt{1-9\vartheta^2}}+\frac{1}{5}+ \vartheta\leq 1,
    \end{split}
\end{equation}
where the last inequality follows from the choice of the parameters. This proves $\mathfrak{F}_\ell^\iota(f)\in X$ for every $\iota=1,\ldots, d$ and every $\ell=0,\ldots, \mathfrak N-1$.

\medskip

\textbf{Step VI. Computation of the $\Omega$s at the perturbation scale.} In this step we show that for every $Q\in \mathfrak C_{16\ell+3}$, every $i=1,\ldots, H_Q$ and every $z\in \supp\mu\cap B(x_i,\frac{A\vartheta}{10240n}\diam Q)$ we have 
$$
\sum_{\iota=1}^d\Omega(\mathfrak{F}_\ell^\iota(f);z,100r_Q)
\geq
\frac{1}{4D^22^k}\frac{\lambda^k\ayin}{4\cdot 50^{k+1}D^2(d+1)}
=:c(k,D,d).
$$
Let us prove that for every $z\in \supp\mu\cap B(x_i,\frac{A\vartheta^j}{10240n}\diam Q)$ there exists $p(z)\in \mathscr F_{Q,x_i}$ such that
$$|\pi_{\mathrm{span}(\mathfrak e)^\perp}(z-x_i)-p(z)|\leq 50r_Q.$$
Recall that by \eqref{eq:discrete-product-hausdorffv2} there is $q\in K_{x_i}$ such that 
$\lvert z-q\rvert\leq 2\rho$. Further, by definition of $K_{x_i}$ there exists $t\in  \mathrm{span}(\mathfrak e)$ and $\hat p\in E_{Q,x_i}$ such that $q=x_i+t+\hat p$.

Observe that by definition of $F_{Q,x_i}$ there exists $\hat{\hat p}\in F_{Q,x_i}$ such that $\lvert \hat p-\hat {\hat p}\rvert\leq 10\rho$. Thanks to the fact that $z\in B(x_i,\frac{A\vartheta^j}{10240n}\diam Q)$, we know that
$|\hat p|\leq 2\rho+\frac{A\vartheta^j}{10240n}\diam Q$, and as a consequence 
$$
\lvert \hat {\hat p}\rvert\leq \frac{A\vartheta^j}{10240n}\diam Q+12\rho\leq \frac{A\vartheta^j}{5120n}\diam Q,
$$
which in turn implies that $\hat {\hat{p}}\in \widehat F_{Q,x_i}\cap B_{\mathfrak e}(0,\frac{A\vartheta^j}{5120n}\diam Q)$. Thanks to \eqref{includiosasjgnaidjfgnajsnd}, there exists $p(z)\in \mathscr F_{Q,x_i}$ such that 
$$
\lvert p(z)-\hat{\hat p}\rvert\leq 40r_Q.
$$
In particular,
$$
|\pi_{\mathrm{span}(\mathfrak e)^\perp}(z-x_i)-p(z)|\leq 12\rho+40r_Q\leq 50r_Q.
$$
Hence let $w:=x_i+t+p(z)$ and notice that 
$B(w,50r_Q)\subseteq B(z,100r_Q)\subseteq B(x_i,\frac{A\vartheta^j}{5120n}\diam Q)$. The last inclusion follows from some easy algebraic computations. Finally, let us observe that 
\begin{equation}
    \begin{split}
      & \qquad\qquad\qquad100r_Q\mu(B(z,100r_Q)) \sum_{\iota=1}^d\Omega(\mathfrak{F}_\ell^\iota(f);z,100r_Q)\geq\sum_{\iota=1}^d\inf_{A\in \mathcal F}\int_{B(w,50r_Q)}\lvert \mathfrak{F}_\ell^\iota(f)-A\rvert d\mu\\
        \geq &\sum_{\iota=1}^d\inf_{A\in \mathcal F}\Big(\int_{B(w,50r_Q)}\lvert \Phi_{Q,i}^\iota-A\rvert d\mu-\int_{B(w,50r_Q)}\Big\lvert g_0(y)-g_0(w)-\mathfrak c_\vartheta Du_{Q'}(\mathfrak c(Q))[\pi_V(y-w)]\Big\rvert d\mu(y)\Big)\\\\
&\qquad\overset{\eqref{linearitag_0sdfsf}}{\geq} \Big(\sum_{\iota=1}^d\inf_{A\in \mathcal F}\int_{B(w,50r_Q)}\lvert \Phi_{Q,i}^\iota-A\rvert d\mu\Big)-d(720n^2\cdot 2^{-\qof/2}+\tfrac{6n\vartheta}{\sqrt{1-9\vartheta^2}})\mu(B(w,50r_Q))100r_Q,
        \nonumber
    \end{split}
\end{equation}
where for convenience we defined $\mathfrak c_\vartheta:=\tfrac{\sqrt{1-9\vartheta^2}}{\sqrt{1-9\vartheta^2}+6n\vartheta}$. For each $\iota=1,\ldots, j$ we let 
$\tilde{y}_j(w):=x_i+t+y_j(p(z))$, where the points $y(p)$ have been introduced in \cref{puntiyintroduzione}. Notice further that for every $j=1,\ldots,n$ there are points
\begin{equation}
    q_j(w)\in \supp\mu\qquad\text{such that}\qquad\lvert q_j(w)-\tilde{y}_j(w)\rvert\leq 2\rho.
    \label{puntoq_j(w)}
\end{equation}
 Further, observe that thanks to few omitted algebraic computations, we have
$$B(q_j(w),r_Q)\subseteq B(x_i,\frac{A\vartheta^j}{5120n}\diam Q).$$
Thus, we infer that by \cref{evadeaffinepro} that
\begin{equation}
    \begin{split}
      \sum_{\iota=1}^d \inf_{A\in \mathcal F}&\int_{B(w,50r_Q)}\lvert \Phi_{Q,i}^\iota-A\rvert d\mu\\
       \geq & \sum_{\iota=1}^d \inf_{A\in \mathcal F}\Big(\sum_{j=1}^{d+1}\int_{B(q_j(w),\lambda r_Q)}\lvert \ayin r_Q\beta^\iota_j-A\rvert d\mu
-\sum_{j=1}^d\int_{B(q_j(w),\lambda r_Q)}\lvert \ayin r_Q \beta^\iota_j-\Phi_{Q,i}^\iota\rvert d\mu\Big)\\
    \geq& \sum_{\iota=1}^d \inf_{A\in \mathcal F}\Big(\sum_{j=1}^{d+1}\int_{B(q_j(w),\lambda r_Q)}\lvert \ayin r_Q\beta^\iota_j-A\rvert d\mu-\frac{2(d+1)\lambda r_Q}{5} \mu(B(q_j(w),\lambda r_Q))\Big),
\nonumber
    \end{split}
\end{equation}
where we used the fact that $\phi^\iota_{Q,i}(\tilde{y}_j(w))=\ayin r_Q \beta^\iota$ and \eqref{puntoq_j(w)}. In particular, this also implies 
    \begin{equation}
    \begin{split}
\sum_{\iota=1}^d \inf_{A\in \mathcal F}&\int_{B(w,50r_Q)}\lvert \Phi_{Q,i}^\iota-A\rvert d\mu\\
\geq& \sum_{\iota=1}^d \frac{(\lambda r_Q)^k}{D}\Big(\inf_{A\in \mathcal F}\sum_{j=1}^{d+1}\fint_{B(q_j(w),\lambda r_Q)}\lvert \ayin r_Q\beta^\iota_j-A\rvert d\mu\Big)-\frac{2d^2D (\lambda r_Q)^{k+1}}{5}\\
\geq& \frac{\ayin\lambda^k r_Q^{k+1}}{2D(d+1)}-\frac{2d(d+1)D(\lambda r_Q)^{k+1}}{5} \geq \frac{\lambda^k\ayin}{4D(d+1)} r_Q^{k+1}\\
\geq& \frac{\lambda^k\ayin}{4\cdot 50^{k+1}D^2(d+1)}(50r_Q)\mu(B(w,50r_Q)),
\nonumber
    \end{split}
\end{equation}
where the last inequality comes from the choice of $\ayin$ and $\lambda$. From this, we are able to infer that 
\begin{equation}
    \begin{split}
&\sum_{\iota=1}^d\Omega(\mathfrak{F}_\ell^\iota(f);z,100r_Q)\geq\frac{\mu(w,50r_Q)}{\mu(B(z,100r_Q))}\Big(\frac{\lambda^k\ayin}{8\cdot 50^kD^2(d+1)}-d(720n^2\cdot 2^{-\qof/2}+\tfrac{6n\vartheta}{\sqrt{1-9\vartheta^2}})\Big)
\\
        &\qquad\qquad \geq \frac{1}{D^22^k}\Big(\frac{\lambda^k\ayin}{8\cdot 50^kD^2(d+1)}-d(720n^2\cdot 2^{-\qof/2}+\tfrac{6n\vartheta}{\sqrt{1-9\vartheta^2}})\Big)\geq c(k,D,s),
        \nonumber
    \end{split}
\end{equation}
where the last inequality comes from the choice of $\qof$ and $\vartheta$. Thus, we have reached the proof of the sought claim. The last task of this step is to estimate from below the mass of the set  
$$\mathcal B_\ell:=\bigcup_{Q\in \mathfrak C_{16\ell+3}}\bigcup_{i=1}^{H_Q}B\Big(x_i,\frac{A\vartheta^j\diam Q}{10240n}\Big).$$
Let us start by recalling that the balls $B(x_i,\frac{A\vartheta^j\diam Q}{10240n})$ are disjoint by construction and the balls $B(x_i,\frac{A\vartheta^j\diam Q}{256n})$ cover $\frac{A}{32n}B_Q$. Thus 
\begin{equation}
    \begin{split}
        \mu(\tfrac{A}{32n}B_Q)\leq& \sum_{i=1}^{H_Q}\mu(B(x_i,\tfrac{A\vartheta^j\diam Q}{256n}))\leq DH_Q\Big(\frac{A\vartheta^j\diam Q}{256n}\Big)^k\\
        \leq& 40^kDH_Q\Big(\frac{A\vartheta^j\diam Q}{10240n}\Big)^k\leq 40^kD^2\sum_{i=1}^{H_Q}\mu\Big(B\Big(x_i,\frac{A\vartheta^j\diam Q}{10240n}\Big)\Big)
        \label{PD2sordigf}
    \end{split}
\end{equation}
On the other hand 
\begin{equation}
    \begin{split}
        \sum_{Q\in \mathfrak C_{16\ell+3}}\mu(\tfrac{A}{32n}B_Q)\geq& \frac{2^k}{(96n)^kD}\sum_{Q\in \mathfrak C_{16\ell+3}}\Big(\frac{3A}{2}\diam Q\Big)^k\\
        \geq &  \frac{2^k}{(96n)^kD^2}\sum_{Q\in \mathfrak C_{16\ell+3}}\mu\Big(\frac{3A}{2}B_Q\Big)\overset{\eqref{stimapalleC16l+2}}{\geq} \frac{1}{32D^{12}48^k5^{5k}n^k}\mu(R).
        \label{PD1asodf}
    \end{split}
\end{equation}
Putting together \eqref{PD2sordigf} and \eqref{PD1asodf} we have 
$$\frac{1}{32D^{14}1920^k5^{5k}n^k}\mu(R)\leq \sum_{Q\in \mathfrak C_{16\ell+3}}\sum_{i=1}^{H_Q}\mu\Big(B\Big(x_i,\frac{A\vartheta^j\diam Q}{10240n}\Big)\Big)=\mu(\mathcal B_\ell).$$

\textbf{Step VII. Computation of the $\Omega$s at bigger scales.} Let $0\leq a<\ell$, let $Q\in \mathfrak C_{16a+3}$, let $i=1,\ldots,H_Q$, and let
$$z\in \supp\mu\cap B\Big(x_i,\frac{A\vartheta^j\diam Q}{10240n}\Big).$$
In what follows, $r_Q$ denotes the scale associated to $Q$ in Step IV. We claim that, for every $\iota=1,\ldots,d$,
\begin{equation}
\Omega(\mathfrak{F}_\ell^\iota(f);z,100r_Q)\geq
   \mathfrak c(\vartheta,\qof)\Omega(f;z,100r_Q)
   -2^{-13\qof(\ell-a)}(1+16n\xi+\ayin),
   \label{stimafinalestepVII}
\end{equation}
where
\begin{equation}
    \mathfrak c(\vartheta,\qof):=\tfrac{(1-\vartheta)\sqrt{1-9\vartheta^2}}{(1+2^{-\qof/4})(\sqrt{1-9\vartheta^2}+6n\vartheta)}.
    \label{costantecvartheta}
\end{equation}
Indeed, by the definition of $\mathfrak F_\ell^\iota(f)$ and of $g_0$, we have
\begin{equation}
    \begin{split}
\mathfrak{F}_\ell^\iota(f)(x)
=&\tfrac{(1-\vartheta)\sqrt{1-9\vartheta^2}}{\sqrt{1-9\vartheta^2}+6n\vartheta}g(x)-(1-\vartheta)
\sum_{S\in \mathfrak C_{16\ell+2}}\sum_{\alpha=1}^{n-j}
Du_{S'}(\mathfrak c(S))[e_\alpha]\omega_{\alpha,\ell,S}(x)\Psi_S(x)\\
&\qquad\qquad\qquad\qquad+\sum_{S\in \mathfrak C_{16\ell+3}}\sum_{h=1}^{H_S}\Phi_{S,h}^\iota(x)\widetilde{\Psi}_{S,h}(x)\\
=&\tfrac{(1-\vartheta)\sqrt{1-9\vartheta^2}}{\sqrt{1-9\vartheta^2}+6n\vartheta}g(x)+(\Delta)(x),
\nonumber
    \end{split}
\end{equation}
where, for each $S\in\mathfrak C_{16\ell+2}$, $S'\in\mathfrak C_{16\ell+1}$ denotes the cube fixed in Step III. Therefore
\begin{equation}
    \begin{split}
        &\Omega(\mathfrak{F}_\ell^\iota(f);z,100r_Q)
        \geq
        \frac{(1-\vartheta)\sqrt{1-9\vartheta^2}}{\sqrt{1-9\vartheta^2}+6n\vartheta}
        \inf_{A\in\mathcal F}\fint_{B(z,100r_Q)}
        \frac{|g-A|}{100r_Q}\,d\mu-\sup_{y\in B(z,100r_Q)}\frac{|(\Delta)(y)|}{100r_Q}\\
        &\qquad\qquad\qquad=\mathfrak c(\vartheta,\qof)
        \inf_{A\in\mathcal F}\fint_{B(z,100r_Q)}
        \frac{|(1+2^{-\qof/4})g-A|}{100r_Q}\,d\mu-\sup_{y\in B(z,100r_Q)}\frac{|(\Delta)(y)|}{100r_Q}\\
        &\geq
        \mathfrak c(\vartheta,\qof)\Omega(f;z,100r_Q)-\mathfrak c(\vartheta,\qof)\fint_{B(z,100r_Q)}
        \frac{|f-(1+2^{-\qof/4})g|}{100r_Q}\,d\mu
        -\sup_{y\in B(z,100r_Q)}\frac{|(\Delta)(y)|}{100r_Q}.
        \label{stimafinale}
    \end{split}
\end{equation}
We now estimate the two error terms
$$
(I):=\fint_{B(z,100r_Q)}
\frac{|f-(1+2^{-\qof/4})g|}{100r_Q}\,d\mu
\qquad\text{and}\qquad
(II):=\sup_{y\in B(z,100r_Q)}\frac{|(\Delta)(y)|}{100r_Q}.
$$
Let $S\in \mathfrak W_{16\ell+b}$, with $b\in\{1,2,3\}$, and assume that
$2AB_Q\cap 2AB_S\neq\emptyset$. We claim that
\begin{equation}
    \diam S\leq 2^{-14\qof(\ell-a)}\diam Q.
    \label{separazionestepVII}
\end{equation}
Indeed, by the nesting of the families $\mathfrak W_m$, there exists a chain
$$
S=S^{16\ell+b}\subseteq S^{16\ell+b-1}\subseteq\cdots\subseteq S^{16a+4},
\qquad S^m\in\mathfrak W_m.
$$
Since $S\subseteq S^{16a+4}$ and $2AB_Q\cap 2AB_S\neq\emptyset$, we have
$10AB_Q\cap 10AB_{S^{16a+4}}\neq\emptyset$. Since
$Q\in\mathfrak C_{16a+3}\subseteq\mathfrak W_{16a+3}$ and $S^{16a+4}\in\mathfrak W_{16a+4}$, the scale separation in item 7 of \cref{propfinalejoiningcones} gives
$$
\diam S^{16a+4}\leq 2^{-\qof}\diam Q.
$$
Iterating the same scale separation along the chain from level $16a+4$ to level $16\ell+b$, we get
$$
\diam S\leq 2^{-\qof(16\ell+b-16a-3)}\diam Q.
$$
Since $b\geq1$, this gives
$$
\diam S\leq 2^{-\qof(16(\ell-a)-2)}\diam Q
\leq 2^{-14\qof(\ell-a)}\diam Q.
$$

We first estimate $(I)$. By the definition of $g$ and by \cref{mollificazioneapproxalbordo}, if
$B(z,100r_Q)$ meets a ball associated with some $S\in\mathfrak C_{16\ell+1}$ in $y$, then by \eqref{separazionestepVII}
\begin{equation}
    \begin{split}
        |f(y)-(1+2^{-\qof/4})g(y)|
        &\leq \frac{A\diam S}{12}\leq 2^{-14\qof(\ell-a)}\frac{A\diam Q}{12}
    \end{split}
\end{equation}
for every $y\in B(z,100r_Q)$. Hence, thanks to the choice of $r_Q$ and of $\qof$,
\begin{equation}
    (I)\leq
    2^{-14\qof(\ell-a)}\frac{A\diam Q}{1200r_Q}
    \leq
    2^{-13\qof(\ell-a)}.
    \label{stimafinale2}
\end{equation}

We now estimate $(II)$. We split $(\Delta)=(\Delta)_A+(\Delta)_B$, where
\begin{equation}
    \begin{split}
(\Delta)_A(x):=&
-(1-\vartheta)\sum_{S\in \mathfrak C_{16\ell+2}}\sum_{\alpha=1}^{n-j}
Du_{S'}(\mathfrak c(S))[e_\alpha]\omega_{\alpha,\ell,S}(x)\Psi_S(x),
    \end{split}
\end{equation}
and
$$
(\Delta)_B(x):=
\sum_{S\in \mathfrak C_{16\ell+3}}\sum_{h=1}^{H_S}
\Phi_{S,h}^\iota(x)\widetilde{\Psi}_{S,h}(x).
$$
By \eqref{stimasupnorma}, by the disjointness of the balls associated with $\mathfrak C_{16\ell+2}$, and by \eqref{separazionestepVII}, we get
\begin{equation}
    \begin{split}
        \|(\Delta)_A\|_{\infty,B(z,100r_Q)}
        &\leq
        2^{-14\qof(\ell-a)}
        \frac{8n\xi}{3\vartheta}A\diam Q.
        \nonumber
    \end{split}
\end{equation}
Moreover, by the disjointness of the balls supporting the final perturbations and by \eqref{separazionestepVII},
\begin{equation}
    \begin{split}
        \|(\Delta)_B\|_{\infty,B(z,100r_Q)}
        \leq
        \sup_{S\in \mathfrak W_{16\ell+3}:2AB_S\cap2AB_Q\neq\emptyset}\ayin r_S
        \leq
        \ayin 2^{-14\qof(\ell-a)}r_Q.
        \nonumber
    \end{split}
\end{equation}
Therefore, again by the choice of $\qof$, we have
\begin{equation}
    \begin{split}
        (II)
        &\leq
        2^{-14\qof(\ell-a)}
        \frac{8n\xi}{3\vartheta}\frac{A\diam Q}{100r_Q}
        +\ayin 2^{-14\qof(\ell-a)}\leq
        2^{-13\qof(\ell-a)}(16n\xi+\ayin).
        \label{stimafinale3}
    \end{split}
\end{equation}
Putting together \eqref{stimafinale}, \eqref{stimafinale2} and \eqref{stimafinale3}, we obtain
\begin{equation}
\Omega(\mathfrak{F}_\ell^\iota(f);z,100r_Q)\geq
   \mathfrak c(\vartheta,\qof)\Omega(f;z,100r_Q)
   -2^{-13\qof(\ell-a)}(1+16n\xi+\ayin).
\end{equation}

\textbf{Step VIII. The final inductive construction.}
We perform an inductive construction that will yield functions
$f_1,\ldots,f_{\mathfrak N}$ and measurable sets
$E_0,\ldots,E_{\mathfrak N-1}$ with certain properties. We begin by setting
$f_0:=0$ and
$$
\gamma_0:=\frac{1}{32D^{14}1920^k5^{5k}n^k}.
$$
In the following, if $z\in \mathcal B_a$, then $Q(z)$ denotes the unique cube
$Q\in\mathfrak C_{16a+3}$ for which
$$
z\in B\Big(x_i,\frac{A\vartheta^j\diam Q}{10240n}\Big)
$$
for some $i=1,\ldots,H_Q$, and $r_{Q(z)}$ denotes the scale associated with this cube in Step IV.

Suppose that, for some $\ell=0,\ldots,\mathfrak N-1$, we have constructed
$f_0,\ldots,f_\ell$ and measurable sets $E_0,\ldots,E_{\ell-1}$ such that
\begin{itemize}
    \item[(a)] $f_b\in X$ for every $b=0,\ldots,\ell$;

    \item[(b)] for every $0\leq a\leq \ell-1$ and every $z\in E_a$,
    \begin{equation}
    \begin{split}
        \Omega(f_\ell;z,100r_{Q(z)})\geq&
        \mathfrak c(\vartheta,\qof)^{\ell-a-1}\frac{c(k,D,d)}{d}-\sum_{m=a+1}^{\ell-1}
        \mathfrak c(\vartheta,\qof)^{\ell-1-m}
        2^{-13\qof(m-a)}(1+16n\xi+\ayin),
    \end{split}
    \label{inductiveOmega}
    \end{equation}
    where the sum is empty if $\ell=a+1$, and where $\mathfrak c(\vartheta,\qof)$ is the constant introduced in \eqref{costantecvartheta}.

    \item[(c)] for every $a=0,\ldots,\ell-1$ one has
    $E_a\subseteq \mathcal B_a$ and
    $$
    \mu(E_a)\geq \frac{\gamma_0}{d}\mu(R).
    $$
\end{itemize}
For $\ell=0$ these conditions are empty except for $f_0\in X$. We now construct $f_{\ell+1}$ and $E_\ell$. By Step VI, applied to $f_\ell$, for every
$z\in \mathcal B_\ell$ we have
$$
\sum_{\iota=1}^d
\Omega(\mathfrak F_\ell^\iota(f_\ell);z,100r_{Q(z)})
\geq c(k,D,d).
$$
For $\iota=1,\ldots,d$ set
$$
\mathcal B_\ell^\iota:=
\Big\{z\in\mathcal B_\ell:
\Omega(\mathfrak F_\ell^\iota(f_\ell);z,100r_{Q(z)})
\geq \frac{c(k,D,d)}{d}
\Big\}.
$$
Then
$$
\mathcal B_\ell\subseteq\bigcup_{\iota=1}^d\mathcal B_\ell^\iota,
$$
and hence there exists $\iota(\ell+1)\in\{1,\ldots,d\}$ such that
$$
\mu(\mathcal B_\ell^{\iota(\ell+1)})
\geq
\frac{1}{d}\mu(\mathcal B_\ell).
$$
We define $E_\ell:=\mathcal B_\ell^{\iota(\ell+1)}$ and $f_{\ell+1}:=\mathfrak F_\ell^{\iota(\ell+1)}(f_\ell)$. By Step V, $f_{\ell+1}\in X$. Moreover, by the mass estimate in Step VI, we have
$$
\mu(E_\ell)\geq \frac{1}{d}\mu(\mathcal B_\ell)
\geq
\frac{\gamma_0}{d}\mu(R).
$$
Thus items (a) and (c) hold at level $\ell+1$. It remains to verify item (b). Let $0\leq a\leq \ell-1$ and let $z\in E_a$. By Step VII,
\begin{equation}
    \begin{split}
    \Omega(f_{\ell+1};z,100r_{Q(z)})
    &=
    \Omega(\mathfrak F_\ell^{\iota(\ell+1)}(f_\ell);z,100r_{Q(z)})\\
    &\geq
    \mathfrak c(\vartheta,\qof)\Omega(f_\ell;z,100r_{Q(z)})
    -2^{-13\qof(\ell-a)}(1+16n\xi+\ayin).
    \end{split}
\end{equation}
Using the inductive hypothesis \eqref{inductiveOmega}, we infer
\begin{equation}
\begin{split}
\Omega(f_{\ell+1};z,100r_{Q(z)})
\geq&\mathfrak c(\vartheta,\qof)^{\ell-a}\frac{c(k,D,d)}{d}-\sum_{m=a+1}^{\ell}\mathfrak c(\vartheta,\qof)^{\ell-m}2^{-13\qof(m-a)}(1+16n\xi+\ayin).
\end{split}
\end{equation}
This is exactly item (b) at level $\ell+1$ for the old sets $E_a$, $a=0,\ldots,\ell-1$.

For the new set $E_\ell$, by the definition of $E_\ell$ we have, for every $z\in E_\ell$,
$$
\Omega(f_{\ell+1};z,100r_{Q(z)})\geq \frac{c(k,D,d)}{d},
$$
which is item (b) with $a=\ell$, since the corresponding sum is empty. This completes the induction.

\textbf{Step IX. The packing estimate.}
Our final choice is $f:=f_{\mathfrak N}$. Recall that, by the choice of the parameters,
$$
\mathfrak c(\vartheta,\qof)^a\geq \frac12
\qquad\text{for every }0\leq a\leq \mathfrak N,
$$
where $\mathfrak c(\vartheta,\qof)$ is the constant introduced in
\eqref{costantecvartheta}. Therefore, for every
$a=0,\ldots,\mathfrak N-1$ and every $z\in E_a$, the inductive estimate of
Step VIII gives
\begin{equation}
    \begin{split}
        \Omega(f;z,100r_{Q(z)})\geq&\mathfrak c(\vartheta,\qof)^{\mathfrak N-a-1}\frac{c(k,D,d)}{d}-\sum_{m=a+1}^{\mathfrak N-1}\mathfrak c(\vartheta,\qof)^{\mathfrak N-1-m}2^{-13\qof(m-a)}(1+16n\xi+\ayin).
    \end{split}
\end{equation}
By the choice of $\qof$, we have
$$
(1+16n\xi+\ayin)
\sum_{h=1}^{\infty}2^{-13\qof h}
\leq
\frac{c(k,D,d)}{4d}.
$$
Hence, for every $a=0,\ldots,\mathfrak N-1$ and every $z\in E_a$, we have
\begin{equation}
    \Omega(f;z,100r_{Q(z)})
    \geq
    \frac{c(k,D,d)}{4d}.
    \label{OmegaFinalLower}
\end{equation}
We now prove the following claim. For every $a=0,\ldots,\mathfrak N-1$,
every $Q\in\mathfrak C_{16a+3}$, every $i=1,\ldots,H_Q$ and every
$$
z\in E_a\cap B\Big(x_i,\frac{A\vartheta^j\diam Q}{10240n}\Big),
$$
where the points $x_i$ are those introduced in the previous steps, there exists a cube $\bar Q(Q,i,z)\subseteq R$ such that
$$200r_Q\leq \diam \bar Q(Q,i,z)\leq 200\oldC{cubi}^2r_Q\qquad\text{and}\qquad\Omega(f,\bar Q(Q,i,z))\geq\frac{c(k,D,d)}{4dD^2(2\oldC{cubi}^2)^{k+1}}.$$
Indeed, by the choice of $r_Q$ and by the choice of the parameters, we have
$$
B(z,100r_Q)
\subseteq
B\Big(x_i,\frac{A\vartheta^j\diam Q}{1280n}\Big)
\subseteq
\frac{A}{16n}B_Q.
$$
Thus, by item 4. of the thesis of \cref{propfinalejoiningcones}, this gives $\supp\mu\cap B(z,100r_Q)\subseteq R$.
Hence, by \cref{evev}, there exists a subcube $\bar Q_0\subseteq R$ with
$z\in \bar Q_0\subseteq B(z,100r_Q)$ and
$$
\frac{100}{\oldC{cubi}^2}r_Q\leq \diam \bar Q_0\leq 100r_Q.
$$
We let $\bar Q(Q,i,z)$ be a cube of $\Delta_\mu$ containing $\bar Q_0$ and such that
$$
200r_Q\leq \diam \bar Q(Q,i,z)\leq 200\oldC{cubi}^2r_Q.
$$
The cube $\bar Q(Q,i,z)$ is still contained in $R$, by the choice of $\qof$ and by
\cref{evev}. Moreover $B_{\bar Q(Q,i,z)}$ contains $B(z,100r_Q)$. Therefore, using
\eqref{OmegaFinalLower} and AD-regularity, we get
\begin{equation}
    \begin{split}
        &\qquad\qquad\Omega(f,\bar Q(Q,i,z))=     \inf_{L\in\mathcal F}\fint_{B_{\bar Q(Q,i,z)}}\frac{|f-L|}{\diam \bar Q(Q,i,z)}\,d\mu\\
        &\geq\frac{100r_Q\,\mu(B(z,100r_Q))}{\diam \bar Q(Q,i,z)\,\mu(B_{\bar Q(Q,i,z)})}        \inf_{L\in\mathcal F}\fint_{B(z,100r_Q)}\frac{|f-L|}{100r_Q}\,d\mu\geq\frac{c(k,D,d)}{4dD^2(2\oldC{cubi}^2)^{k+1}}.
    \end{split}
\end{equation}
This proves the claim. Define 
$$c_2(\oldC{cubi},D,k,d):=\frac{c(k,D,d)}{4dD^2(2\oldC{cubi}^2)^{k+1}}\qquad\text{and let}\qquad\mathcal S:=\{S\in\Delta_\mu(R):\Omega(f,S)\geq c_2(\oldC{cubi},D,k,d)\}.$$
For $a=0,\ldots,\mathfrak N-1$, define
$$
\mathfrak B_a:=
\Big\{
\bar Q(Q,i,z):
Q\in\mathfrak C_{16a+3},\ i=1,\ldots,H_Q,\ 
z\in E_a\cap B\Big(x_i,\frac{A\vartheta^j\diam Q}{10240n}\Big)
\Big\}.
$$
Then $\mathfrak B_a\subseteq\mathcal S$ and $E_a\subseteq \bigcup\{S:S\in\mathfrak B_a\}$.
We now show that the families $\mathfrak B_a$ are separated in scale. Let
$S\in\mathfrak B_a$ and $T\in\mathfrak B_b$, with
$0\leq a<b\leq\mathfrak N-1$. Write
$S=\bar Q(Q,i,z)$ and $T=\bar Q(Q',i',z')$ with
$Q\in\mathfrak C_{16a+3}$ and $Q'\in\mathfrak C_{16b+3}$. Assume that $S\cap T\neq\emptyset$. By the localization of the cubes
$\bar Q(Q,i,z)$ and $\bar Q(Q',i',z')$, and by the smallness of the parameters
entering the definition of $r_Q$, we have
$10AB_Q\cap 10AB_{Q'}\neq\emptyset$. By the scale separation of the families $\mathfrak W_m$ along the generations
between $16a+3$ and $16b+3$, we get
$$
\diam Q\geq 2^{16\qof(b-a)}\diam Q'.
$$
Since $r_Q$ is a fixed multiple of $\diam Q$, with the same multiplicative
constant for every cube $Q$, this implies
$r_Q\geq 2^{16\qof(b-a)}r_{Q'}$. On the other hand,
since $\diam S\geq 200r_Q$ and $\diam T\leq 200\oldC{cubi}^2r_{Q'}$. Therefore
$$\diam S\geq \oldC{cubi}^{-2}2^{16\qof(b-a)}\diam T.$$
In particular, by the choice of $\qof$, the same cube cannot belong to two
distinct families $\mathfrak B_a$ and $\mathfrak B_b$. Hence
\begin{equation}
    \begin{split}
    \sum_{Q\in \mathcal S}\mu(Q)\geq\sum_{a=0}^{\mathfrak N-1}\sum_{S\in\mathfrak B_a}\mu(S)\geq\sum_{a=0}^{\mathfrak N-1}\mu(E_a)\geq\frac{\mathfrak N}{32dD^{14}1920^k5^{5k}n^k}
    \mu(R).
    \end{split}
    \label{packinglowerboundfinal}
\end{equation}
However, since we assumed that $\mu$ satisfies the GWALA condition, and since the coefficients $\Omega$ are compatible in the sense of \cref{defscalecompatible}, we can apply \cref{lemcontinuousdyadiccompatible}. Therefore the dyadic family $\mathcal S$ is Carleson with constant $C_{\mathrm{dyad}}(c_2(\oldC{cubi},D,k,d))$, where $C_{\mathrm{dyad}}(c_2(\oldC{cubi},D,k,d))$ is the constant yielded by \cref{lemcontinuousdyadiccompatible} from the GWALA constant at the threshold $c_2(\oldC{cubi},D,k,d)$. Hence
$$\sum_{Q\in\mathcal S}\mu(Q)\leq C_{\mathrm{dyad}}(c_2(\oldC{cubi},D,k,d))\mu(R).$$
Combining this estimate with \eqref{packinglowerboundfinal}, we obtain
$$\frac{\mathfrak N}{32dD^{14}1920^k5^{5k}n^k}\leq C_{\mathrm{dyad}}(c_2(\oldC{cubi},D,k,d)).$$
Since however 
$$\mathfrak N=\Big\lfloor \frac{1}{64(n+1)\vartheta}\Big\rfloor,$$
this is impossible if $\vartheta$ is chosen sufficiently small, with the smallness depending only on $\oldC{cubi},D,k,d$ and on the dyadic Carleson constant $C_{\mathrm{dyad}}(c_2(\oldC{cubi},D,k,d))$ given by \cref{lemcontinuousdyadiccompatible}. This contradiction proves that no family
$\mathscr I(j,2\vartheta,\sigma/2,\varepsilon,\delta,\zeta,\mathfrak m,A)$,
with $0\leq j\leq k-1$, can fail to be Carleson. This concludes the proof.
\end{proof}

Finally, we are ready to prove the main theorem of the paper.

\begin{teorema}\label{teoremaGWALAimpliesUR}
Suppose $\mu$ is an $\alpha$-AD-regular measure on $\R^n$ with regularity
constant $D$. Let $\mathcal F$ be a vector subspace of
$L^1_{\mathrm{loc}}(\R^n)$ of dimension $d$ containing the affine functions,
and let $q\in[1,\infty]$. Suppose that $\mu$ satisfies the
$\mathrm{GWALA}_q$ condition with respect to $\mathcal F$. Then
$\alpha\in\N$ and $\mu$ is uniformly rectifiable.
\end{teorema}

\begin{proof}
By \cref{omega_monotone} and \cref{osservazioneriduzionewala}, the condition $\mathrm{GWALA}_q$ implies the
corresponding $\mathrm{GWALA}_1$ condition, possibly with the same constants.
Thus, in what follows, we work only with $q=1$. By \cref{prop:rettificabilitya}, we have $\alpha\in \N$ and that $\mu$ is $\alpha$-rectifiable. 

Let $\tilde\varepsilon_0=\tilde\varepsilon_0(n,k,D)>0$ be the constant given
by \cref{propsmallOUWGL}. We shall prove that the fixed-threshold OUWGL bad
set
$$
\mathscr B^2_{\tilde\varepsilon_0}
:=
\{(x,r)\in\supp\mu\times(0,\infty):
u\beta_\mu^k(x,r)>\tilde\varepsilon_0\}
$$
is Carleson. Then \cref{propsmallOUWGL} gives uniform rectifiability. By \cref{lemcontinuousdyadiccompatible} and \cref{remOUWGLdiscretization}, it is enough to prove that, for some $0<\eta<1$ and some $\Lambda\geq1$ satisfying $4\oldC{cubi}^2(\eta+\Lambda^{-1})\leq\tilde\varepsilon_0$, the dyadic family
$$
\mathscr U_{\eta,\Lambda}
:=
\{Q\in\Delta_\mu:
u\beta_\mu^k(\mathfrak c(Q),\Lambda\diam Q)>\eta\}
$$
is Carleson. We first choose $\eta>0$ so small that
$$8\oldC{cubi}^2\eta\leq\tilde\varepsilon_0.$$
We now choose the parameters
$\sigma,\vartheta,\varepsilon,\delta,\zeta,\mathfrak m\in(0,1/2)$ and
$A\geq1$ in such a way that the hypotheses of \cref{teoremasemifinale} are
satisfied and, moreover, $\vartheta$ is smaller than the constant
$\vartheta_0$ appearing there and futher we also suppose that 
\begin{equation}
8\oldC{cubi}^2\frac{16n}{A\vartheta^k}\leq\tilde\varepsilon_0\qquad\text{and}\qquad16n\Big(
\frac{\sigma/2}{\vartheta^k}
+
n\frac{\varepsilon+\delta+\zeta+\mathfrak m}{A\vartheta^k}
\Big)
\leq \eta.
\label{eq:smallness-for-OUWGL}
\end{equation}
Chosen $\Lambda:=\frac{A\vartheta^k}{16n}$, we have that 
$4\oldC{cubi}^2(\eta+\Lambda^{-1})\leq\tilde\varepsilon_0$. By \cref{rioridinocubi}, the dyadic cubes decompose into the families
$$
\mathscr I(j,2\vartheta,\sigma/2,\varepsilon,\delta,\zeta,\mathfrak m,A),
\qquad j=0,\ldots,k.
$$
The families with $j<k$ are Carleson by \cref{teoremasemifinale}. We claim that
no cube in the remaining family
$\mathscr I(k,2\vartheta,\sigma/2,\varepsilon,\delta,\zeta,\mathfrak m,A)$ belongs to $\mathscr U_{\eta,\Lambda}$. Indeed, let $Q\in \mathscr I(k,2\vartheta,\sigma/2,\varepsilon,\delta,\zeta,\mathfrak m,A)$. By the definition of the family yielded by \cref{rioridinocubi}, there exists
$\mathfrak e=(e_1,\ldots,e_k)\in(\mathbb S^{n-1})^k$ such that $Q$ is
$(k,\mathfrak e,(2\vartheta)^k,\sigma/2,
\varepsilon,\delta,\zeta,\mathfrak m,A)$-invariant. Apply \cref{propo:piano-da-invarianza} with
$x=\mathfrak c(Q)$ and
$V:=\mathrm{span}(\mathfrak e)$. Thanks to the choice of $\Lambda$, we have 
$$\sup_{z\in(\mathfrak c(Q)+V)\cap B(\mathfrak c(Q),\Lambda\diam Q)}\frac{\dist(z,\supp\mu)}{\Lambda\diam Q}\leq16n\Big(\frac{\sigma/2}{(2\vartheta)^k}+n\frac{\varepsilon+\delta+\zeta+\mathfrak m}{A(2\vartheta)^k}\Big)\leq \eta$$
by \eqref{eq:smallness-for-OUWGL}. Hence $u\beta_\mu^k(\mathfrak c(Q),\Lambda\diam Q)\leq\eta$, and therefore $Q\notin\mathscr U_{\eta,\Lambda}$. Thus
$$
\mathscr U_{\eta,\Lambda}
\subseteq
\bigcup_{j=0}^{k-1}
\mathscr I(j,2\vartheta,\sigma/2,\varepsilon,\delta,\zeta,\mathfrak m,A).
$$
The right-hand side is a finite union of Carleson families, again by
\cref{teoremasemifinale}. Hence $\mathscr U_{\eta,\Lambda}$ is Carleson.

By \cref{remOUWGLdiscretization}, assumption (a) of
\cref{lemcontinuousdyadiccompatible} holds with
$\varepsilon=\tilde\varepsilon_0$, because
$$4\oldC{cubi}^2(\eta+\Lambda^{-1})\leq\tilde\varepsilon_0.$$
Therefore, by item (i) of \cref{lemcontinuousdyadiccompatible}, the continuous fixed-threshold bad set
$$
\mathscr B^2_{\tilde\varepsilon_0}
=
\{(x,r)\in\supp\mu\times(0,\infty):u\beta_\mu^k(x,r)>\tilde\varepsilon_0\}
$$
is Carleson. Therefore \cref{propsmallOUWGL} applies, and $\mu$ is uniformly rectifiable.
\end{proof}

We finally obtain the following easy consequence.

\begin{corollario}\label{corollarioGWALAgenerale}
Let $\mu$ be an $\alpha$-AD-regular measure on $\R^n$, let
$\mathcal F\subset L^1_{\mathrm{loc}}(\R^n)$ be a finite-dimensional vector
space, and let $q\in[1,\infty]$. If $\mu$ satisfies $\mathrm{GWALA}_q$ with
respect to $\mathcal F$, then $\alpha\in\N$ and $\mu$ is uniformly
rectifiable.
\end{corollario}

\begin{proof}
Set $\widetilde{\mathcal F}:=\mathcal F+\mathscr A$, where as usual $\mathscr A$ denotes the set of affine functions. Since
$\mathcal F\subseteq\widetilde{\mathcal F}$, for every Lipschitz function
$f$, every $x\in\supp\mu$ and every $r>0$, there holds
$$
\Omega_{q,\mu,\widetilde{\mathcal F}}(f;x,r)
\leq
\Omega_{q,\mu,\mathcal F}(f;x,r).
$$
Hence, for every $\varepsilon>0$, we have the inclusion
$$\{(x,r):\Omega_{q,\mu,\widetilde{\mathcal F}}(f;x,r)>\varepsilon\}\subseteq\{(x,r):\Omega_{q,\mu,\mathcal F}(f;x,r)>\varepsilon\}.$$
Thus $\mu$ satisfies $\mathrm{GWALA}_q$ with respect to $\widetilde{\mathcal F}$ and conclusion follows from \cref{teoremaGWALAimpliesUR}.
\end{proof}

\appendix

\section{Measurability results}\label{appendix}

In this section we collect results and the proof of the measurability results. We took out of the main body of the proofs these propositions since in this case measurability issues are of secondary importance and their presence in the main body might worsen the readability of the paper.

\begin{lemma}\label{lemma:reg:ball:dens}
    Let $r>0$ and denote by $\mathsf E_r:\mathscr{M}(\R^n)\times \R^n\to \R^+$ the function 
    $$\mathsf E_r(\mu,x):=\mu(B(x,r)).$$
    Then, $\mathsf E_r(\mu,x)$ is upper semicontinuous when $\mathscr{M}$ is endowed with the topology of the weak* convergence and as a consequence, for every $\Lambda>0$ the map $\mathsf R_r:\mathscr M(\R^n)\times \R^n\to \R^+$ defined as
    $$\mathsf R_r(\mu,r):=\frac{\mathsf E_{\Lambda r}(\mu,x)}{E_{r}(\mu,x)},\qquad\text{ is Borel.}$$
\end{lemma}

\begin{proof}
    This is a direct consequence of the argument of the proof of \cite[Lemma A.5]{ReversePansu}.
\end{proof}

\begin{lemma}\label{measu:tan}
    Suppose $\mu$ is an asymptotically  Radon measure on $\R^n$. Then, the multimap $x\mapsto \Tan(\mu,x)$ coincides with a strongly universally measurable multimap on a Borel set of full $\mu$-measure.
\end{lemma}

\begin{proof}
Let us first recall that $d$ denotes the metric inducing the weak* convergence on $\mathscr M(\R^n)$. Let $E$ be the Borel set of full $\mu$ measure such that $\Tan(\mu,x)\neq \emptyset$ and let
$$\Gamma:=\left\{(x,\nu)\in E\times \mathscr M(\R^n):\nu\in\Tan(\mu,x)\right\}.$$
It is immediate to check that the map $\Psi:(x,r,c)\mapsto cT_{x,r}\mu$ is continuous with respect to the weak* convergence. We claim that
$$\Gamma=\bigcap_{h=1}^\infty\bigcap_{m=1}^\infty\bigcup_{\substack{q\in\mathbb{Q}\cap(0,1/m)\\c\in\mathbb{Q}\cap(0,\infty)}}\Big\{(x,\nu)\in E\times \mathscr M(\R^n):d(cT_{x,q}\mu,\nu)<\frac{1}{h}\Big\}.
$$
Indeed, suppose first that $\nu\in\Tan(\mu,x)$. Then there are $r_i\to0$ and $c_i>0$ such that $c_iT_{x,r_i}\mu\rightharpoonup\nu$.

By the continuity of the map $\Psi$, for every $i$ we can find $q_i\in\mathbb{Q}\cap(0,1/i)$ and $a_i\in\mathbb{Q}\cap(0,\infty)$ such that
$$d(a_iT_{x,q_i}\mu,c_iT_{x,r_i}\mu)<\frac{1}{i}.$$
It follows that $a_iT_{x,q_i}\mu\rightharpoonup\nu$, and hence $(x,\nu)$ belongs to the right-hand side.
The converse inclusion is immediate. This in particular that the set $\Gamma$ is Borel. Thus, for every closed set $C$ we have that $C\cap \Gamma$ is Borel and hence its projection on the first component $\R^n$ must be a Suslin set. Hence, by definition of measurable multimap, that can be found at the beggining of  \cite[\S 5.1]{Srivastava}, we have $x\mapsto \Tan(\mu,x)$ is strongly universally measurable. 
\end{proof}

\begin{proposizione}\label{prop:borelDelte}
    Suppose $\mu$ is an asymptotically doubling Radon measure on $\mathbb R^n$. Let $R>0$ be fixed. For $r>0$ define $\Delta_R(x,r):=d_R(\mu(B(x,r))^{-1}T_{x,r}\mu,\Tan(\mu,x))$, whenever $\mu(B(x,r))>0$, and extend it as $\Delta_R(x,r)=0$ if $\mu(B(x,r))=0$, where we recall that $d_R$ is introduced in \cref{d_r}. Then, for every fixed $r>0$, we can find a Borel representative for the map $x\mapsto \Delta_R(x,r)$. Moreover, for $\mu$-almost every $x\in\mathbb R^n$,
    $\lim_{r\to 0}\Delta_R(x,r)=0$.
\end{proposizione}

\begin{proof}
Let $E$ be the Borel set of full $\mu$ measure such that $\Tan(\mu,x)\neq \emptyset$, and 
$$\limsup_{r\rightarrow 0} \frac{\mu(B(\xi,2r))}{\mu(B(\xi,r))}<\infty,$$
for every $x\in E$. Let
$$
\Gamma:=\{(x,\nu)\in E\times\mathscr M(\R^n):\nu\in\Tan(\mu,x)\}.
$$
By the the proof of \cref{measu:tan}, $\Gamma$ is Borel. Fix $r>0$ and notice that by \cref{lemma:reg:ball:dens} and the continuity of $(x,r)\mapsto T_{x,r}\mu$, the map $x\mapsto\mu(B(x,r))^{-1}T_{x,r}\mu$ is Borel on $E$. Moreover, its $F_R$-mass is positive there. Since $F_R$ is continuous, the map
$$
x\mapsto
\frac{\mu(B(x,r))^{-1}T_{x,r}\mu}
{F_R\big(\mu(B(x,r))^{-1}T_{x,r}\mu\big)}
$$
is Borel on $E$. For the definition of $F_R$, we refer to \cref{defF_r}. For every $a>0$, the set $P_a:=\{x\in E:\Delta_R(x,r)<a\}$ is the projection onto $\R^n$ of the Borel set of all $(x,\nu)\in\Gamma$ such that $F_R(\nu)=1$ and
$$F_R\Big(\frac{\mu(B(x,r))^{-1}T_{x,r}\mu}{F_R(\mu(B(x,r))^{-1}T_{x,r}\mu)},\nu\Big)<a.$$
Since $F_R(\cdot,\cdot)$ is weak* continuous, $P_a$ is therefore a Suslin set and hence universally measurable. Thus $x\mapsto\Delta_R(x,r)$ is universally measurable on $E$. Extending it by zero outside $E$, we obtain a universally measurable function which, since $\mu$ is Radon, agrees $\mu$-almost everywhere with a Borel function.

We are left to prove the limit. Fix $x\in E$ and suppose by contradiction that there are $\varepsilon>0$ and let $r_i$ be an infinitesimal decreasing sequence such that
$\Delta_R(x,r_i)\geq\varepsilon$ for every $i$. By \cref{precompactnessmeasures}, after passing to a subsequence, we have $\mu(B(x,r_i))^{-1}T_{x,r_i}\mu\rightharpoonup\nu$ for some $\nu\in\Tan(\mu,x)$. By \cref{t:compactdouble}, $0\in\supp\nu$, and hence $F_R(\nu)>0$. Therefore, by \cref{preiss},
$$\lim_{i\to\infty}\Delta_R(x,r_i)=\lim_{i\to\infty}d_R\Big(\frac{T_{x,r_i}\mu}{\mu(B(x,r_i))},\Tan(\mu,x)\Big)=d_R\big(\nu,\Tan(\mu,x)\big)=0,$$
which is a contradiction. Hence $\Delta_R(x,r)\to0$ for every $x\in E$, concluding the proof.
\end{proof}

\begin{proposizione}\label{measurability}
Let $\mu$ be a Radon measure on $\R^n$. For every $\mathfrak r,D>0$ and $A\geq1$, let $E_{D,\mathfrak r,A}$ be the set of those $x\in\supp\mu$ for which
\begin{equation}
\Tan(\mu,x)\cap
\bigcup_{k\in\N}\bigcap_{j\geq k}
\mathfrak S_d(V(\mu,x),2^{-j},D,\mathfrak r,A)
\neq\emptyset,
\end{equation}
where the measures $\mathfrak S_d(V(\mu,x),2^{-j},D,\mathfrak r,A)$ were introduced in \cref{defSV}. Then $E_{D,\mathfrak r,A}$ is $\mu$-measurable.
\end{proposizione}

\begin{proof}
For fixed $\xi,D,\mathfrak r>0$ and $A\geq1$, we claim that
$$
\mathscr Y:=\{(V,\nu)\in\Gr(\R^n)\times\mathscr M(\R^n):
\nu\in\mathfrak S_d(V,\xi,D,\mathfrak r,A)\}
$$
is Suslin. Let us first show that the claim concludes the proof. If $\mathscr Y$ is Suslin, then the set
$$\tilde{\mathscr Y}:=\Big\{(V,\nu):\nu\in\bigcup_{k\in\N}\bigcap_{j\geq k}
\mathfrak S_d(V,2^{-j},D,\mathfrak r,A)\Big\}
$$
is Suslin. Let $E$ be a Borel set of full $\mu$ measure such that $\Tan(\mu,x)\neq \emptyset$, for every $x\in E$. Let
$\Gamma:=\{(x,\nu)\in E\times\mathscr M(\R^n):\nu\in\Tan(\mu,x)\}$ and recall that thanks to the proof of \cref{measu:tan}, we have that $\Gamma$ is a Borel set. Since the map $\Psi:=(x,\nu)\mapsto(V(\mu,x),\nu)$ is Borel, the preimage of $\tilde{\mathscr Y}$ under $\Psi$ is Suslin. Therefore,
$$\Gamma\cap \Psi^{-1}(\tilde{\mathscr Y})=\Gamma\cap\Big\{(x,\nu)\in E\times\mathscr M(\R^n):\nu\in\bigcup_{k\in\N}\bigcap_{j\geq k}\mathfrak S_d(V(\mu,x),2^{-j},D,\mathfrak r,A)\Big\}$$
is Suslin and its projection onto the first component is exactly $E_{D,\mathfrak r,A}\cap E$. It follows that $E_{D,\mathfrak r,A}\cap E$ is analytic and hence universally measurable. Since $\mu(\R^n\setminus E)=0$, the set $E_{D,\mathfrak r,A}$ is $\mu$-measurable.

It remains to prove that $\mathscr Y$ is Suslin and in order to do this, we collect a series of measurability results that we will put together at the end of the argument. Since $\Gr(\R^n)$ is the finite disjoint union of the spaces $\Gr(n,m)$, it is enough to fix $m$. 

We identify measures on $V^\perp$ with measures $\eta\in\mathscr M(\R^n)$ concentrated on $V^\perp$. The set
$$\{(\eta,u)\in\mathscr M(\R^n)\times\R^n:u\in\supp\eta\}=\bigcap_{m\in\N}\{(\eta,u):\eta(U(u,m^{-1}))>0\}$$
is Borel since, for every $r>0$, the map $(\eta,u)\mapsto\eta(U(u,r))$ is lower semicontinuous.

Since $V\mapsto\pi_V$ is continuous by \cref{prop:continuitax}, the map
$(V,u)\mapsto\pi_V(u)$ from $\Gr(n,m)\times\R^n$ to $\R^n$ is continuous. Consequently,
$$\{(V,u)\in\Gr(n,m)\times\R^n:u\in V^\perp\}=\{(V,u):\pi_V(u)=0\}$$
is closed. Let us introduce
\begin{equation}
\begin{split}
\mathscr T:=&\{(V,\eta,u)\in\Gr(n,m)\times\mathscr M(\R^n)\times\R^n:u\in\supp\eta\cap V^\perp\cap B(0,10)\}\\
&\quad=\bigcap_{m\in\N}\{(V,\eta,u):\eta(U(u,m^{-1}))>0\}\cap\{(V,\eta,u):\pi_V(u)=0\}\cap\{(V,\eta,u):|u|\leq10\}.
\nonumber
\end{split}
\end{equation}
Observe that the first sets in the second line are Borel by the lower semicontinuity, while the last two sets are closed. Hence the set on the left-hand side is Borel. Its section corresponding to each fixed $(V,\eta)$ is the compact set $\supp\eta\cap V^\perp\cap B(0,10)$. In addition, the set
$$\{(V,\eta)\in\Gr(n,m)\times\mathscr M(\R^n):\eta(\R^n\setminus V^\perp)=0\}$$
is closed. Indeed, suppose that $V_i\to V$, $\eta_i\rightharpoonup\eta$, and $\eta_i(\R^n\setminus V_i^\perp)=0$. If $\varphi\in C_c(\R^n\setminus V^\perp)$, then $\supp\varphi\cap V_i^\perp=\emptyset$ for all sufficiently large $i$, and therefore
$\int\varphi\,d\eta=\lim_i\int\varphi\,d\eta_i=0$. Applying this to nonnegative functions which equal one on a given compact subset of $\R^n\setminus V^\perp$, and then using inner regularity, gives $\eta(\R^n\setminus V^\perp)=0$.

Moreover, the set of all $(V,\eta,\nu)$ satisfying
$$\eta(\R^n\setminus V^\perp)=0\qquad\text{and}\qquad\nu=\Haus^k\trace V\otimes\eta$$
is closed. Indeed, suppose that $V_i\to V$, $\eta_i\rightharpoonup\eta$, and $\nu_i=\Haus^k\trace V_i\otimes\eta_i\rightharpoonup\nu$. Choose orthogonal maps $R_i$ that converge in the operator norm to $\mathrm{id}$ and such that $R_iV=V_i$. Then $(R_i^{-1})_\#\eta_i\rightharpoonup\eta$, and, for every $\varphi\in C_c(\R^n)$,
$$\int\varphi\,d\big(\Haus^k\trace V\otimes(R_i^{-1})_\#\eta_i\big)=\int_{V^\perp}\int_V\varphi(v+z)\,d\Haus^k(v)\,d((R_i^{-1})_\#\eta_i)(z).$$
The inner integral is a continuous compactly supported function of $z\in V^\perp$. Consequently,
$\Haus^k\trace V\otimes(R_i^{-1})_\#\eta_i\rightharpoonup\Haus^k\trace V\otimes\eta$. Since
$\nu_i=(R_i)_\#(\Haus^k\trace V\otimes(R_i^{-1})_\#\eta_i)$ and $R_i$ converges to $\mathrm{id}$, it follows that $\nu=\Haus^k\trace V\otimes\eta$.

Let us observe that for every open set $U$, the set
$$\{(V,\eta,u):u\in\supp\eta\cap V^\perp\cap B(0,10)\cap U\}$$
is Borel, since it is the intersection of $\mathscr T$ and $\Gr(n,m)\times\mathscr M(\R^n)\times U$, that are Borel sets. Its section corresponding to $(V,\eta)$ is an open subset of the compact set $\supp\eta\cap V^\perp\cap B(0,10)$ and is therefore $\sigma$-compact. Hence, by \cite[Theorem 5.12.1]{Srivastava}, its projection
$$\{(V,\eta):\supp\eta\cap V^\perp\cap B(0,10)\cap U\neq\emptyset\}\qquad\text{is Borel.}$$
Let $(U_h)_{h\in\N}$ be a countable basis of open subsets of $\R^n$. We endow $(\R^n)^\N$ with the product topology, which is Polish. The set of all $(V,\eta,\{w_\ell\}_{\ell\in\N})$ such that every $w_\ell$ belongs to $\supp\eta\cap V^\perp\cap B(0,10)$ and every tail of $\{w_\ell\}$ is dense in this set is
\begin{equation}
    \begin{split}
&\mathscr G:=\bigcap_{\ell\in\N}\big\{(V,\eta,(w_j)):(V,\eta,w_\ell)\in\mathscr T\big\}\\
&\quad\cap\bigcap_{h,N\in\N}\Big(\big\{(V,\eta,(w_j)):\supp\eta\cap V^\perp\cap B(0,10)\cap U_h=\emptyset\big\}\cup\bigcup_{\ell\geq N}\big\{(V,\eta,(w_j)):w_\ell\in U_h\big\}\Big).
\nonumber
\end{split}
\end{equation}
For every $\ell$, the map
$(V,\eta,(w_j))\mapsto(V,\eta,w_\ell)$ is continuous, so the sets in the first intersection are Borel because $\mathscr T$ is Borel. The first set inside the parentheses is Borel by the preceding paragraph, while the sets $\{(V,\eta,(w_j)):w_\ell\in U_h\}$ are open because the coordinate projections $(w_j)\mapsto w_\ell$ are continuous. Thus $\mathscr G$ is Borel. The condition in the second line says exactly that every open set in the fixed basis for the topology meeting $\supp\eta\cap V^\perp\cap B(0,10)$ contains a point of every tail of $\{w_\ell\}$, which is equivalent to the density of every tail.

Since the map $(\nu,u)\mapsto \nu(B(u,s))$ is upper semicontinuous, for every $\iota=1,\ldots,d+1$, the set
$$\{(\nu,z,y)\in\mathscr M(\R^n)\times\R^n\times\R^n:\nu(B(z,50\xi))\leq D\,\nu(B(y,\mathfrak r\xi/32))\}$$
is Borel. Indeed, the map
$$
(\nu,z,y)\mapsto
\big(\nu(B(z,50\xi)),\nu(B(y,\mathfrak r\xi/32))\big)
$$
is Borel, and the set $\{(a,b)\in[0,\infty)^2:a\leq Db\}$ is closed.

Consider now the Polish space
$$
\Gr(n,m)\times\mathscr M(\R^n)\times\mathscr M(\R^n)
\times\big((\R^n)^{d+3}\big)^\N
$$
with coordinates
$(V,\eta,\nu,(w_\ell,z_\ell,y_{\ell,1},\ldots,y_{\ell,d+1})_{\ell\in\N})$. Consider the subset consisting of the tuples for which
$\eta(\R^n\setminus V^\perp)=0$, $\nu=\Haus^k\trace V\otimes\eta$, $\nu$ is $A$-doubling, and one of the following two alternatives holds:

\begin{enumerate}
\item $\supp\eta\cap V^\perp\cap B(0,10)=\emptyset$;
\item every tail of $\{w_\ell\}$ is dense in $\supp\eta\cap V^\perp\cap B(0,10)$ and, for every $\ell\in\N$,
$$
\begin{gathered}
z_\ell\in\supp\eta\cap V^\perp\cap B(0,10),
\qquad |w_\ell-z_\ell|<\xi+2^{-\ell},\\
y_{\ell,\iota}\in\supp\eta\cap V^\perp\cap B(z_\ell,\xi/8)
\qquad\text{for }\iota=1,\ldots,d+1,\\
|y_{\ell,\iota_1}-y_{\ell,\iota_2}|\geq4\mathfrak r\xi
\qquad\text{whenever }\iota_1\neq\iota_2,\\
\nu(B(z_\ell,50\xi))
\leq D\,\nu(B(y_{\ell,\iota},\mathfrak r\xi/32))
\qquad\text{for }\iota=1,\ldots,d+1.
\end{gathered}
$$
\end{enumerate}

The preceding results, together with the observation that the space of $A$-doubling measures is closed, show that this subset is Borel. We claim that its projection onto the coordinates $(V,\nu)$ is exactly
$$
\mathscr Y\cap\big(\Gr(n,m)\times\mathscr M(\R^n)\big).
$$
Indeed, let $\nu\in\mathfrak S_d(V,\xi,D,\mathfrak r,A)$ and write
$\nu=\Haus^k\trace V\otimes\eta$. If
$\supp\eta\cap V^\perp\cap B(0,10)=\emptyset$, the first alternative applies. Otherwise, choose $\{w_\ell\}$ so that every tail is dense in this compact set. If $\Sigma_\xi$ is the set given by \cref{defSV}, gratterethe corresponding points
$y_{\ell,1},\ldots,y_{\ell,d+1}$. Thus $(V,\nu)$ belongs $\mathscr Y\cap\big(\Gr(n,m)\times\mathscr M(\R^n)\big)$.

Conversely, suppose that a tuple satisfying the second alternative is given and set $\Sigma_\xi:=\{z_\ell:\ell\in\N\}$. For every
$w\in\supp\eta\cap V^\perp\cap B(0,10)$, the density of every tail of $\{w_\ell\}$ gives indices $\ell_i\to\infty$ such that $w_{\ell_i}\to w$. Hence
$$
\dist(w,\Sigma_\xi)
\leq\liminf_{i\to\infty}
\big(|w-w_{\ell_i}|+|w_{\ell_i}-z_{\ell_i}|\big)
\leq\xi.
$$
Therefore $\Sigma_\xi$ is $\xi$-dense, and the remaining conditions in \cref{defSV} hold by construction. The first alternative corresponds to taking $\Sigma_\xi=\emptyset$. This proves the claimed projection identity.

The projection of a Borel subset of a Polish space is Suslin. Hence
$\mathscr Y\cap(\Gr(n,m)\times\mathscr M(\R^n))$ is Suslin. Taking the finite union over $k$ proves that $\mathscr Y$ is Suslin and concludes the proof.
\end{proof}

\printbibliography[title=References, heading=bibintoc]

\end{document}